\documentclass[11pt]{article}

\usepackage{color, amsmath, amssymb, amsthm, graphicx, bbm, footnote, mathtools}

\usepackage[top=1in, bottom=1in, left=1in, right=1in]{geometry}

\usepackage[toc,page,header]{appendix}
\usepackage{titletoc}
\usepackage{tikz}
\usepackage{titling}
\usepackage{enumitem}
\allowdisplaybreaks

\usepackage{setspace}
\usepackage{adjustbox}
\usepackage[colorlinks,allcolors=blue]{hyperref}
\usepackage[round]{natbib}  
  \hypersetup{
	breaklinks=true,   
	colorlinks=true,   
	pdfusetitle=true,  
}

\usepackage{apptools}
\AtAppendix{\counterwithin{theorem}{section}}
\AtAppendix{\counterwithin{remark}{section}}

\newtheorem{theorem}{Theorem}
\newtheorem{definition}{Definition}
\newtheorem{proposition}[theorem]{Proposition}
\newtheorem{corollary}[theorem]{Corollary}
\newtheorem{lemma}[theorem]{Lemma}
\newtheorem{remark}{Remark}
\newtheorem{assumption}{Assumption}
\newtheorem{algorithm}{Algorithm}

\newcommand{\possessivecite}[1]{\citeauthor{#1}'s (\citeyear{#1})}

\newcommand{\possessivecites}[2]{\citeauthor{#1}'s (\citeyear{#2})} 

\DeclareMathOperator*{\argmin}{\arg\min} 

\title{ \textsc{Gaussian and Bootstrap Approximations for Suprema of Empirical Processes }}
\author{Alexander Giessing\thanks{Department of Statistics, University of Washington, Seattle, WA. E-mail: giessing@uw.edu.}}
\date{\today}
\begin{document}
	\maketitle

\begin{abstract}
	In this paper we develop non-asymptotic Gaussian approximation results for the sampling distribution of suprema of empirical processes when the indexing function class $\mathcal{F}_n$ varies with the sample size $n$ and may not be Donsker. Prior approximations of this type required upper bounds on the metric entropy of $\mathcal{F}_n$ and uniform lower bounds on the variance of $f \in \mathcal{F}_n$ which, both, limited their applicability to high-dimensional inference problems.
	In contrast, the results in this paper hold under simpler conditions on boundedness, continuity, and the strong variance of the approximating Gaussian process. The results are broadly applicable and yield a novel procedure for bootstrapping the distribution of empirical process suprema based on the truncated Karhunen-Lo{\`e}ve decomposition of the approximating Gaussian process. We demonstrate the flexibility of this new bootstrap procedure by applying it to three fundamental problems in high-dimensional statistics: simultaneous inference on parameter vectors, inference on the spectral norm of covariance matrices, and construction of simultaneous confidence bands for functions in reproducing kernel Hilbert spaces.	
\end{abstract}

\section{Introduction}\label{sec:Intro}

\subsection{Approximating the sampling distribution of suprema of empirical processes}\label{subsec:Intro}
This paper is concerned with non-asymptotic bounds on the Kolmogorov distance between the sampling distribution of suprema of empirical processes and the distribution of suprema of Gaussian proxy processes.
Consider a simple random sample of independent and identically distributed random variables $X_1, \ldots, X_n$ with common law $P$ taking values in the measurable space $(S, \mathcal{S})$. Let $\mathcal{F}_n$ be a class of measurable functions $f : S \rightarrow \mathbb{R}$ and define the empirical process 
\begin{align*} 
	\mathbb{G}_n(f) = \frac{1}{\sqrt{n}}\sum_{i=1}^n \big( f(X_i) - Pf\big), \quad \quad f \in \mathcal{F}_n.
\end{align*}
Also, denote by $\{G(f) : f \in \mathcal{F}_n\}$ a centered Gaussian process with some positive semi-definite covariance function.
Within this minimal setup (and a few additional technical assumptions), the paper addresses the problem of deriving non-asymptotic bounds on 
\begin{align}\label{eq:sec:Intro-2}
 	\varrho_n :=	\sup_{s \geq 0} \left|\mathbb{P}\Big\{ \sup_{f \in \mathcal{F}_n} |\mathbb{G}_n(f)| \leq s \Big\} -\mathbb{P}\Big\{\sup_{f \in \mathcal{F}_n} |G(f)|\leq s \Big\} \right|,
\end{align}
and, consequently, conditions under which an empirical process indexed by $\mathcal{F}_n$ 
is \emph{Gaussian approximable}, i.e. $\varrho_n \rightarrow 0$ as $n \rightarrow \infty$. Evidently, if $\mathcal{F}_n$ is a Donsker class then the empirical process is trivially Gaussian approximable by the centered Gaussian process with covariance function $(f, g) \mapsto P(fg) - (Pf)(Pg)$. Thus, in this paper we are concerned with conditions that are strictly weaker than those that guarantee a central limit theorem.

Gaussian approximation results of this type have gained significant interest recently,
as they prove to be effective in tackling high- and infinite-dimensional inference problems. Applications include simultaneous inference on high-dimensional parameters~\citep{dezeure2017high, zhang2017SimultaneousInf}, inference on nonparametric models~\citep{chernozhukov2014AntiConfidenceBands, chen2015asymptotic, chen2016nonparametric}, testing for spurious correlations~\citep{fan2018discoveries}, testing for shape restrictions~\citep{chetverikov2019testing}, inference on large covariance matrices~\citep{chen2018GaussianApproxUStat, han2018on,lopes2019BootstrappingSpectral}, goodness-of-fit tests for high-dimensional linear models~\citep{jankova2020goodness}, simultaneous confidence bands in functional data analysis~\citep{lopes2020bootstrapping, singh2023kernel}, and error quantification in randomized algorithms~\citep{lopes2023bootstrapping}.

The theoretical groundwork on Gaussian approximation was laid in three seminal papers by~\cite{chernozhukov2013GaussianApproxVec, chernozhukov2014GaussianApproxEmp, chernozhukov2015ComparisonAnti}. Since then, subsequent theoretical works have developed numerous refinements and extensions~\citep[e.g.][]{chernozhukov2016empirical, chernozhukov2019ImprovedCLT, chernozhukov2021NearlyOptimalCLT, deng2020beyond, fang2021HighDimCLT, kuchibohotla2021high, cattaneo2022yurinskiis, lopes2022central, lopes2022improved, bong2023dual}. In this paper, we resolve two limitations of~\possessivecite{chernozhukov2014GaussianApproxEmp} original results that have not been previously addressed:

\begin{itemize}
	\item \emph{Entropy conditions.} The original bounds on the Kolmogorov distance $\varrho_n$ depend on the metric entropy of the function class $\mathcal{F}_n$. These upper bounds are non-trivial only if the metric entropy grows at a much slower rate than the sample size $n$. As a result, most applications only address sparse inference problems involving function classes with either a finite number of functions, discretized functions, or a small VC-index. The new upper bounds in this paper no longer depend on the metric entropy of the function class and therefore open the possibility to tackle non-sparse, high-dimensional inference problems. We call these bounds \emph{dimension-/entropy-free}.
	
	\item \emph{Lower bounds on weak variances.} The original bounds on $\varrho_n$ require a strictly positive lower bound on the \emph{weak variance} of the Gaussian proxy process, i.e. $\inf_{f \in \mathcal{F}_n} \mathrm{Var}\big(G(f)\big) > 0$. This condition limits the scope of the original results to problems with standardized function classes (studentized statistics) and excludes situations with variance decay, typically observed in non-sparse, high-dimensional problems.  The new Gaussian approximation results in this paper only depend on the \emph{strong variance} $\mathrm{Var}\big(\sup_{f \in \mathcal{F}_n} |G(f)|\big)$ and are therefore applicable to a broad range of problems including those with degenerate distributions. We say that these bounds are \emph{weak variance-free}.
\end{itemize}
Even though these limitations (and our solution) are quite technical, the resulting new Gaussian and bootstrap approximations have immediate practical consequences. We present three substantial applications to problems in high-dimensional inference in Section~\ref{sec:Applications} and two toy examples that cast light on why and when the original approximation results by~\cite{chernozhukov2014GaussianApproxEmp} fail and the new results succeed in Appendix~\ref{sec:DerivationIntro}.

Notably, in the special case of inference on spectral norms of covariance matrices, entropy-free bounds on $\varrho_n$ already exist. Namely, for $\mathcal{F}_n = \{ x \mapsto f(x) = (x'u)(x'v): u, v \in S^{d-1}\}$ with $S^{d-1} = \{u \in \mathbb{R}^d : \|u\|_2 =1\}$~\citeauthor{lopes2019BootstrappingSpectral} (\citeyear{lopes2019BootstrappingSpectral, lopes2020bootstrapping, lopes2022central, lopes2022improved}) have devised an approach to bounding $\varrho_n$ that combines a specific variance decay assumption with a truncation argument.
For a carefully chosen  truncation level, the resulting bound on $\varrho_n$ depends only on the sample size $n$ and the parameter characterizing the variance decay. Since their approach is intimately related to the bilinearity of the functions in $\mathcal{F}_n$ and the specific variance decay assumption, it does not easily extend to arbitrary function classes. We therefore develop a different strategy in this paper.

Furthermore, for finite function classes $|\mathcal{F}_n| < \infty$,~\cite{chernozhukov2021NearlyOptimalCLT} and~\cite{deng2020beyond} have been able to slightly relax the requirements on the weak variance of the Gaussian proxy process.  However, their results do not generalize to arbitrary function classes with $|\mathcal{F}_n| = \infty$. Our strategy for replacing the weak variance with the strong variance in the upper bounds on $\varrho_n$ is therefore conceptually completely different from theirs. For details, we refer to our companion paper on Gaussian anti-concentration inequalities~\cite{giessing2022anticoncentration}.

\subsection{Contributions and overview of the results}
This paper consists of three parts which contribute to probability theory (Section~\ref{sec:ResultsMaximaVec}), mathematical statistics and bootstrap methodology (Section~\ref{sec:ResultsSupremaEP}), and high-dimensional inference (Section~\ref{sec:Applications}). The Appendices~\ref{sec:DerivationIntro}--\ref{sec:Proofs-AuxiliaryResults} contain additional supporting results and all proofs.

Section~\ref{sec:ResultsMaximaVec} contains the main mathematical innovations of this paper. We establish dimension- and weak variance-free Gaussian and bootstrap approximations for maxima of sums of independent and identically distributed high-dimensional random vectors (Section~\ref{subsec:FourBasicLemmas}). Specifically, we derive a Gaussian approximation inequality (Proposition~\ref{lemma:CLT-Max-Norm}), a Gaussian comparison inequality (Proposition~\ref{lemma:GaussianComparison}), and two bootstrap approximation inequalities (Propositions~\ref{lemma:Bootstrap-Max-Norm} and~\ref{lemma:Bootstrap-Max-Norm-Asymptotic-Size}). At the core of these four theoretical results is a new proof of a multivariate Berry-Esseen-type bound which leverages two new technical auxiliary results: an anti-concentration inequality for suprema of separable Gaussian processes (Lemma~\ref{lemma:AntiConcentration-SeparableProcess} in Appendix~\ref{subsec:AuxiliaryResults-AntiConcentration}) and a smoothing inequality for partial derivatives of the Ornstein-Uhlenbeck semigroup associated to a multivariate normal measure (Lemma~\ref{lemma:DifferentiatingUnderIntegral} in Appendix~\ref{subsec:AuxiliaryResults-Smoothing}). We conclude this section with a comparison to results from the literature (Section~\ref{subsec:RelationPreviousWork-Math}). 

Section~\ref{sec:ResultsSupremaEP} contains our contributions to mathematical statistics and bootstrap methodology. The results in this section generalize the four basic results of Section~\ref{sec:ResultsMaximaVec} from finite dimensional vectors to empirical processes indexed by totally bounded (and hence separable) function classes. As one would expect, dimension- (a.k.a entropy-) and weak variance-freeness of the finite dimensional results carry over to empirical processes. We establish Gaussian approximation inequalities (Section~\ref{subsec:GaussianApproximation}), Gaussian comparison inequalities (Section~\ref{subsec:GaussianComparison}), and an abstract bootstrap approximation inequality (Section~\ref{subsec:GaussianProcessBootstrap}). The latter result motivates a new procedure for bootstrapping the sampling distribution of suprema of empirical processes. We call this procedure the \emph{Gaussian process bootstrap} (Algorithm~\ref{algorithm:GaussianProcessBootstrap}) and discuss practical aspects of its implementations via the truncated Karhunen-Lo{\`e}ve decomposition of a Gaussian proxy process (Section~\ref{subsec:GaussianProcessBootstrap-Implementation-Consistency}). We include a selective comparison with results from the literature (Section~\ref{subsec:RelationPreviousWork-Stat}).

In Section~\ref{sec:Applications} we showcase the flexibility of the Gaussian process bootstrap by applying it to three fundamental problems in high-dimensional statistics: simultaneous inference on parameter vectors (Section~\ref{subsec:Application-Vector}), inference on the spectral norm of covariance matrices (Section~\ref{subsec:Application-Matrix}), and construction of simultaneous confidence bands for functions in reproducing kernel Hilbert spaces (Section~\ref{subsec:Application-FDA}). For each of these three examples we include a brief comparison with alternative bootstrap methods from the literature to explain in what sense our method improves over (or matches with) existing results. 


\section{Dimension- and weak variance-free results for the maximum norm of sums of i.i.d. random vectors}\label{sec:ResultsMaximaVec}

\subsection{Four basic results}\label{subsec:FourBasicLemmas}
We begin with four propositions on non-asymptotic Gaussian and bootstrap approximations for finite-dimensional random vectors. These propositions are our main mathematical contribution and essential for the statistical results to come. In Section~\ref{sec:ResultsSupremaEP} we lift these propositions to general empirical processes 
which are widely applicable to problems in mathematical statistics. Throughout this section, $\|x\|_\infty = \max_{1 \leq k \leq d}|x_k|$ denotes the maximum norm of a vector $x \in \mathbb{R}^d$.

The first result is a Gaussian approximation inequality. This inequality provides a non-asymptotic bound on the Kolmogorov distance between the laws of the maximum norm of sums of independent and identically distributed random vectors and a Gaussian proxy statistic. 

\begin{proposition}[Gaussian approximation]\label{lemma:CLT-Max-Norm}
	Let $X, X_1, \ldots X_n \in \mathbb{R}^d$ be i.i.d. random vectors with mean zero and positive semi-definite covariance matrix $\Sigma \neq \mathbf{0} \in \mathbb{R}^{d \times d}$. Set $S_n = n^{-1/2} \sum_{i=1}^n X_i$ and $Z \sim N(0, \Sigma)$. Then, for $M \geq 0$, $n \geq 1$,
	\begin{align*}
		&\sup_{s \geq 0} \Big|\mathbb{P}\left\{\|S_n\|_\infty \leq s \right\} -\mathbb{P}\left\{\|Z\|_\infty \leq s \right\} \Big|\\
		&\quad{}\quad{}\quad{} \lesssim \frac{(\mathrm{E}[\|X \|_\infty^3])^{1/3}}{\sqrt{n^{1/3}\mathrm{Var}(\|Z\|_\infty)}} + \frac{\mathrm{E} \left[ \|X\|_\infty^3 \mathbf{1}\{\|X\|_\infty > M\}\right]}{\mathrm{E} \left[ \|X\|_\infty^3\right]}  + \frac{ \mathrm{E}[\|Z\|_\infty] + M}{\sqrt{n\mathrm{Var}(\|Z\|_\infty)}},
	\end{align*}
	where $\lesssim$ hides an absolute constant independent of $n, d, M$, and the distribution of the $X_i$'s.
\end{proposition}

\begin{remark}[Extension to independent, non-identically distributed random vectors]\label{remark:lemma:CLT-Max-Norm-1}
	The proof of Proposition~\ref{lemma:CLT-Max-Norm} involves an inductive argument inspired by Theorem 3.7.1 in~\cite{nourdin2012normal} who attribute it to~\cite{bolthausen1984estimate}. This argument relies crucially on independent and identically distributed data. Generalizing it to independent but non-identically distributed data would require additional assumptions on the variances similar to the uniform asymptotic negligibility condition in the classical Lindeberg-Feller CLT for triangular arrays. For empirical processes such a generalization is less relevant. We therefore leave this to future research.
\end{remark}


\begin{remark}[Extension to non-centered random vectors]\label{remark:lemma:CLT-Max-Norm-2}
	If the covariance matrix $\Sigma$ is strictly positive definite, then Proposition~\ref{lemma:CLT-Max-Norm} also holds for i.i.d. random vectors $X, X_1, \ldots X_n \in \mathbb{R}^d$ with non-zero mean $\mu \in \mathbb{R}^d$  and $Z \sim N(\mu, \Sigma)$. Nevertheless, in the broader context of empirical processes a strictly positive definite covariance (function) is a very strong assumption. Therefore, in this paper, we do not pursue this refinement. 
	Instead, the interested reader may consult the companion paper~\cite{giessing2022anticoncentration}.
\end{remark}

\begin{remark}[Gaussian approximation versus CLT]\label{remark:lemma:CLT-Max-Norm-3}
	Proposition~\ref{lemma:CLT-Max-Norm} is strictly weaker than a multivariate CLT because the Kolmogorov distance between the maximum norm of two sequences of random vectors induces a topology on the set of probability measures on $\mathbb{R}^d$ which is coarser than the topology of convergence in distribution. Indeed, consider $U_n = (X,Y)' \in \mathbb{R}^2$ and $V_n = (Y,X)' \in \mathbb{R}^2$ for $n \geq 1$ and $X, Y$ arbitrary random variables. Then, the Kolmogorov distance between $\|U_n\|_\infty$ and $\|V_n\|_\infty$ is zero for all $n \geq 1$, but $U_n \overset{d}{=} V_n$ if and only if $X$ and $Y$ are exchangeable. For another perspective on the same issue, see Section~\ref{subsec:RelationPreviousWork-Math}.
\end{remark}

Above Gaussian approximation result differs in several ways from related results in the literature. The two most striking difference are the following: First, the non-asymptotic upper bound in Proposition~\ref{lemma:CLT-Max-Norm} does not explicitly depend on the dimension $d$, but only on the (truncated) third moments of $\|X\|_\infty$ and $\|Z\|_\infty$ and the variance of $\|Z\|_\infty$. Under suitable conditions on the marginal and/ or joint distribution of the coordinates of $X$ these three quantities grow substantially slower than the dimension $d$. Among other things, this opens up the possibility of applying this bound in the context of high-dimensional random vectors when $d \gg n$. 
Second, Proposition~\ref{lemma:CLT-Max-Norm} holds without stringent assumptions on the distribution of $X$ and applies even to degenerate distributions that do not have a strictly positive definite covariance matrix $\Sigma$. In particular, the lemma does not require a lower bound on the variances of the coordinates of $X$ or the minimum eigenvalue of $\Sigma$. This is fundamental for an effortless extension to general empirical processes. For lower bounds on the variance of $\|Z\|_\infty$ the reader may consult~\cite{giessing2022anticoncentration}.  For a comprehensive comparison of Proposition~\ref{lemma:CLT-Max-Norm} with previous work, we refer to Section~\ref{subsec:RelationPreviousWork-Math}.

If we are willing to impose a lower bound on the variances of the coordinates of $X$, we can deduce the following useful inequality:

\begin{corollary}\label{corollary:lemma:CLT-Max-Norm}
	Recall the setup of Proposition~\ref{lemma:CLT-Max-Norm}. In addition, suppose that $\sigma_{(1)}^2 = \min_{1 \leq j \leq d} \mathrm{Var}(X^{(j)}) > 0$ and that $X$ has coordinate-wise finite $3 + \delta$ moments, $\delta > 0$. Define $\widetilde{X} = (X^{(j)}/\sigma_{(1)})_{j=1}^d$, $\widetilde{Z} = (Z^{(j)}/\sigma_{(1)})_{j=1}^d$. Then, for all $ n \geq 1$,
	\begin{align*}
		&\sup_{s \geq 0} \Big|\mathbb{P}\left\{\|S_n\|_\infty \leq s \right\} -\mathbb{P}\left\{\|Z\|_\infty \leq s \right\} \Big| \\
		&\quad\quad\quad\lesssim  \frac{1}{n^{1/6}} \left( \mathrm{E}[\|\widetilde{X} \|_\infty^{3 + \delta}] + \mathrm{E}[\|\widetilde{Z}\|_\infty^{3 + \delta}]\right)^{\frac{1}{3 + \delta}}  \mathrm{E}[\|\widetilde{Z}\|_\infty] + \frac{1}{n^{\delta/3}} \left(\frac{\mathrm{E}[\|\widetilde{X}\|_\infty^{3 + \delta}]^{\frac{1}{3 + \delta}}}{\mathrm{E}[\|\widetilde{X} \|_\infty^3]^{1/3}}\right)^3,
	\end{align*}
	where  $\lesssim$ hides an absolute constant independent of $n, d, \sigma_{(1)}^2$, and the distribution of the $X_i$'s.
	If the coordinates in $X$ are equicorrelated with correlation coefficient $\rho \in (0, 1]$, then above inequality holds with $\mathrm{E}[\|\widetilde{Z}\|_\infty]$ replaced by $1/\rho$.
\end{corollary}
Since $\mathrm{E}[\|\widetilde{X} \|_\infty^3]^{1/3} \geq \max_{1 \leq j \leq d} \mathrm{E}[|\widetilde{X}^{(j)}|^2]^{1/2} = \frac{\sigma_{(n)}}{\sigma_{(1)}}$ for $\sigma_{(n)}^2 = \max_{1 \leq j \leq d} \mathrm{Var}(X^{(j)})$, the right hand side of the inequality in Corollary~\ref{corollary:lemma:CLT-Max-Norm} can be easily upper bounded under a variety of moment assumptions on the marginal distributions of the coordinates of $X$. For a few concrete examples relevant in high-dimensional statistics we refer to Lemmas 1--3 in~\cite{giessing2023bootstrap}.

The second basic result is a Gaussian comparison inequality. The novelty of this result is (again) that it holds even for degenerate Gaussian laws with singular covariance matrices and does not explicitly depend on the dimension of the random vectors. 

\begin{proposition}[Gaussian comparison]\label{lemma:GaussianComparison}
	Let $Y, Z \in \mathbb{R}^d$ be Gaussian random vectors with mean zero and positive semi-definite covariance matrices $\Sigma$ and $\Omega\neq \mathbf{0} \in \mathbb{R}^{d \times d}$, respectively. Then,
	\begin{align*}
		&\sup_{s \geq 0} \Big|\mathbb{P}\left\{\|Y\|_\infty \leq s \right\} -\mathbb{P}\left\{\|Z\|_\infty \leq s \right\} \Big| \lesssim \left( \frac{\max_{j,k} |\Omega_{jk} - \Sigma_{jk}|}{ \mathrm{Var}(\|Y\|_\infty) \vee \mathrm{Var}(\|Z\|_\infty)} \right)^{1/3},
	\end{align*}
	where $\lesssim$ hides an absolute constant independent of $d, \Sigma, \Omega$.
\end{proposition}

Since the Gaussian distribution is fully characterized by its first two moments, Proposition~\ref{lemma:CLT-Max-Norm} instantly suggests that it should be possible to approximate the sampling distribution of $\|S_n\|_\infty$ with the sampling distribution of $\|Z_n\|_\infty$, where $Z_n \mid X_1, \ldots, X_n \sim N(0, \widehat{\Sigma}_n)$ and $\widehat{\Sigma}_n$ is a positive semi-definite estimate of $\Sigma$. The third basic result formalizes this idea; it is a simple consequence of the triangle inequality combined with Propositions~\ref{lemma:CLT-Max-Norm} and~\ref{lemma:GaussianComparison}:

\begin{proposition}[Bootstrap approximation of the sampling distribution]\label{lemma:Bootstrap-Max-Norm}
	Let $X, X_1, \ldots X_n \in \mathbb{R}^d$ be i.i.d. random vectors with mean zero and positive semi-definite covariance matrix $\Sigma \neq \mathbf{0} \in \mathbb{R}^{d \times d}$. Let $\widehat{\Sigma}_n \equiv \widehat{\Sigma}_n(X_1, \ldots, X_n)$ be any positive semi-definite estimate of $\Sigma$. Set $S_n = n^{-1/2} \sum_{i=1}^n X_i$, $Z_n \mid X_1, \ldots, X_n \sim N(0, \widehat{\Sigma}_n)$, and $Z \sim N(0, \Sigma)$. Then, for $M \geq 0$, $n \geq 1$,
	\begin{align*}
		&\sup_{s \geq 0} \Big|\mathbb{P}\left\{\|S_n\|_\infty \leq s \right\} -\mathbb{P}\left\{\|Z_n\|_\infty \leq s \mid X_1, \ldots, X_n \right\} \Big|\\
		&\quad{}\quad{} \lesssim \frac{(\mathrm{E}[\|X \|_\infty^3])^{1/3}}{\sqrt{n^{1/3}\mathrm{Var}(\|Z\|_\infty)}} + \frac{\mathrm{E} \left[ \|X\|_\infty^3 \mathbf{1}\{\|X\|_\infty > M\}\right]}{\mathrm{E} \left[ \|X\|_\infty^3\right]} \\
		&\quad{}\quad{} \quad{}\quad{}+ \frac{ \mathrm{E}[\|Z\|_\infty] + M}{\sqrt{n\mathrm{Var}(\|Z\|_\infty)}} + \left( \frac{\max_{j,k} |\widehat{\Sigma}_{n,jk} - \Sigma_{jk}|}{ \mathrm{Var}(\|Z\|_\infty)} \right)^{1/3},
	\end{align*}
	where $\lesssim$ hides an absolute constant independent of $n, d, M$, and the distribution of the $X_i$'s.
\end{proposition}

Typically, statistical applications require estimates of the quantiles of the sampling distribution. 
Since the covariance matrix $\Sigma$ is unknown, the quantiles of $\|Z\|_p$ with $Z \sim N(0, \Sigma)$ are infeasible. Hence, for $\alpha \in (0,1)$ arbitrary, we define the feasible Gaussian bootstrap quantiles as
\begin{align*}
	c_n(\alpha; \widehat{\Sigma}_n) &:= \inf \left\{s \geq 0: \mathbb{P}\left\{\|Z_n\|_\infty \leq s \mid X_1, \ldots, X_n \right\}  \geq \alpha \right\}, \quad \text{where} \quad Z_n \mid X_1, \ldots, X_n \sim N(0, \widehat{\Sigma}_n).
\end{align*}
Since this quantity is random, it is not immediately obvious that it is a valid approximation of the $\alpha$-quantile of the sampling distribution of $\|S_n\|_p$. However, combing Proposition~\ref{lemma:Bootstrap-Max-Norm} with standard arguments~\citep[e.g.][]{chernozhukov2013GaussianApproxVec} we obtain the fourth basic result:

\begin{proposition}[Bootstrap approximation of quantiles]\label{lemma:Bootstrap-Max-Norm-Asymptotic-Size}
	Consider the setup of Lemma~\ref{lemma:Bootstrap-Max-Norm}. Let $(\Theta_n)_{n \geq 1} \in \mathbb{R}$ be a sequence of arbitrary random variables, not necessarily independent of $X, X_1, \ldots, X_n$. Then, for $M \geq 0$, $n \geq 1$,
	\begin{align*}
		&\sup_{\alpha \in (0,1)} \Big|\mathbb{P}\left\{\|S_n\|_\infty + \Theta_n \leq c_n(\alpha; \widehat{\Sigma}_n) \right\}  - \alpha \Big| \\
		&\quad{}\quad{} \lesssim \frac{(\mathrm{E}[\|X \|_\infty^3])^{1/3}}{\sqrt{n^{1/3}\mathrm{Var}(\|Z\|_\infty)}} + \frac{\mathrm{E} \left[ \|X\|_\infty^3 \mathbf{1}\{\|X\|_\infty > M\}\right]}{\mathrm{E} \left[ \|X\|_\infty^3\right]} + \frac{ \mathrm{E}[\|Z\|_\infty] + M}{\sqrt{n\mathrm{Var}(\|Z\|_\infty)}}\\
		& \quad{}\quad{}\quad{}\quad+ \inf_{\delta > 0}\left\{ \left(\frac{\delta }{ \mathrm{Var}( \|Z\|_\infty)} \right)^{1/3}  + \mathrm{P}\left(\max_{j,k} |\widehat{\Sigma}_{n,jk} - \Sigma_{jk}|  > \delta\right)\right\}\\
		&\quad{}\quad{}\quad{} \quad{}\quad\quad+ \inf_{\eta > 0} \left\{ \frac{\eta }{ \sqrt{\mathrm{Var}( \|Z\|_\infty)} } + \mathrm{P}\left(|\Theta_n| > \eta \right)\right\},
	\end{align*}
	where $\lesssim$ hides an absolute constant independent of $n, d, M$, and the distribution of the $X_i$'s.
\end{proposition}
\begin{remark}[On the purpose of the random variables $(\Theta_n)_{n \geq 1}$]\label{remark:lemma:Bootstrap-Max-Norm-Asymptotic-Size}
	In applications $(\Theta_n)_{n\geq 1}$ may be a higher-order approximation error such as the remainder of a first-order Taylor approximation. For concrete examples we refer to Corollary 2 and Theorem 10 in~\cite{giessing2023bootstrap}. The random variables $(\Theta_n)_{n\geq 1}$ may also be taken identical to zero. In this case, the expression in the last line in Lemma~\ref{lemma:Bootstrap-Max-Norm-Asymptotic-Size} vanishes.
\end{remark}

Propositions~\ref{lemma:Bootstrap-Max-Norm} and~\ref{lemma:Bootstrap-Max-Norm-Asymptotic-Size} inherit the dimension- and weak variance-freeness from Propositions~\ref{lemma:CLT-Max-Norm} and~\ref{lemma:GaussianComparison}. Since these propositions do not require lower bounds the minimum eigenvalue of the estimate of $\Sigma$, we can always use the naive sample covariance matrix $\widehat{\Sigma}_n = n^{-1} \sum_{i=1}^n (X_i - \bar{X}_n)(X_i - \bar{X}_n)'$ with $\bar{X}_n = n^{-1} \sum_{i=1}^n X_i$ to estimate $\Sigma$ even if $d \gg n$. If additional information about $\Sigma$ is available (viz. low-rank, bandedness, or approximate sparsity), we can of course use more sophisticated estimators $\widehat{\Sigma}_n$ to improve the non-asymptotic bounds. This will prove particularly effective when we lift Propositions~\ref{lemma:Bootstrap-Max-Norm} and~\ref{lemma:Bootstrap-Max-Norm-Asymptotic-Size} to general empirical processes (see Section~\ref{subsec:GaussianProcessBootstrap}). The reader can find several more examples in Section 4.2 of the companion paper~\cite{giessing2023bootstrap}.

\subsection{Relation to previous work}\label{subsec:RelationPreviousWork-Math}
Here, we compare Propositions~\ref{lemma:CLT-Max-Norm} and~\ref{lemma:GaussianComparison} with the relevant results in the literature. From a mathematical point of view Propositions~\ref{lemma:Bootstrap-Max-Norm} and~\ref{lemma:Bootstrap-Max-Norm-Asymptotic-Size} are just an afterthought. For a comprehensive review of the entire literature on Gaussian and bootstrap approximations of maxima of sums of random vectors see~\cite{chernozhukov2023high-dimensional}.

\begin{itemize}
	\item \emph{On the dependence on the dimension.} As emphasized above, the upper bounds in Propositions~\ref{lemma:CLT-Max-Norm} and~\ref{lemma:GaussianComparison} are dimension-free. This is a significant improvement over Theorems 2.2, 2.1, 3.2, 2.1 in~\cite{chernozhukov2013GaussianApproxVec, chernozhukov2017CLTHighDim, chernozhukov2019ImprovedCLT, chernozhukov2021NearlyOptimalCLT}, Theorem 1.1 and Corollary 1.3 in~\cite{fang2021HighDimCLT}, and Theorems 2.1 and 2.2 in~\cite{lopes2022central}, which all feature logarithmic factors of the dimension. Such bounds generalize poorly to empirical processes since the resulting upper bounds necessarily depend on the $\varepsilon$-entropy of the function class. This precludes (or, at the very least, substantially complicates) applications to objects as basic as the operator norm of a high-dimensional covariance matrix.
	
	\item \emph{On the moment assumptions.} The above mentioned results by~\cite{chernozhukov2013GaussianApproxVec, chernozhukov2017CLTHighDim, chernozhukov2019ImprovedCLT, chernozhukov2021NearlyOptimalCLT},~\cite{fang2021HighDimCLT}, and~\cite{lopes2022central} all require strictly positive lower bounds on the variances of the components of the random vector $X$ and/ or on the minimum eigenvalue of the covariance matrix $\Sigma$. 
	The strictly positive lower bounds are especially awkward if we try to extend their finite dimensional results to general empirical processes: While is often sensible to impose an upper bound on the variances of the increments of an empirical process, lower bounds on the variances are much harder to justify and in general only achievable via a discretization (or thinning) of the function class. The weak variance-free upper bounds in Propositions~\ref{lemma:CLT-Max-Norm} and~\ref{lemma:GaussianComparison} give us greater leeway and increases the scope of our approximation results considerably (see also Appendix~\ref{sec:DerivationIntro}).
	
	\item \emph{On the sharpness of the upper bounds.} The results in~\cite{chernozhukov2021NearlyOptimalCLT},~\cite{fang2021HighDimCLT}, and~\cite{lopes2022central} show that the upper bound in Proposition~\ref{lemma:CLT-Max-Norm} is sub-optimal and that its dependence on sample size $n$ can be improved to $n^{-1/2}$. The proof techniques in~\cite{chernozhukov2021NearlyOptimalCLT} and~\cite{fang2021HighDimCLT} are both based on delicate estimates of Hermite polynomials and thus inherently dimension dependent. Extending their approaches to the coordinate-free Wiener chaos decomposition (which would yield dimension-free results) is a formidable research task. The proof strategy in~\cite{lopes2022central} is very different and requires sub-Gaussian data. The results in the aforementioned papers also show that under additional distributional assumptions the exponent on the upper bound in Proposition~\ref{lemma:GaussianComparison} can be improved to $1/2$. 
	
	\item \emph{Proper generalizations of classical CLTs.} In a certain regard, Propositions~\ref{lemma:CLT-Max-Norm} and~\ref{lemma:GaussianComparison} are strictly weaker than the results in~\cite{chernozhukov2017CLTHighDim, chernozhukov2021NearlyOptimalCLT},~\cite{fang2021HighDimCLT}, and~\cite{lopes2022central}. Proposition~\ref{lemma:CLT-Max-Norm} (and similarly Proposition~\ref{lemma:GaussianComparison}) provides a bound on
	\begin{align*}
		\sup_{s \geq 0} \Big|\mathbb{P}\left\{\|S_n\|_\infty \leq s \right\} -\mathbb{P}\left\{\|Z\|_\infty \leq s \right\} \Big| = \sup_{A \in \mathcal{Q}_d} \Big|\mathbb{P}\left\{ S_n \in A \right\} -\mathbb{P}\left\{ Z \in A \right\} \Big| ,
	\end{align*}
	where $\mathcal{Q}_d$ be the collection of all hypercubes in $\mathbb{R}^d$ with center at the origin.~\cite{chernozhukov2017CLTHighDim},~\cite{fang2021HighDimCLT}, and~\cite{lopes2022central} provide bounds on above quantity not for the supremum over $\mathcal{Q}_d$ but the supremum over the much larger collection of all hyper-rectangles in $\mathbb{R}^d$. When the dimension $d$ is fixed, their results imply convergence in distribution of $S_n$ to $Z$ as $n \rightarrow \infty$ and are therefore stronger than our results (see also Remark~\ref{remark:lemma:CLT-Max-Norm-3}). In particular, their results can be considered proper generalizations of classical CLTs to high dimensions.
	
	The results in this paper depend on a dimension-free anti-concentration inequality for $\|Z\|_\infty$ (or, in other words, for the supremum over $\mathcal{Q}_d$) and a smoothing inequality specific to the map $x \mapsto \|x\|_\infty$ (i.e. Lemmas~\ref{lemma:DifferentiatingUnderIntegral} and~\ref{lemma:AntiConcentration-SeparableProcess}). Since these inequalities do not apply to the class of all hyper-rectangles in $\mathbb{R}^d$, the arguments in this paper cannot be easily modified to yield dimension-free generalizations of the classical CLTs to high dimensions.
\end{itemize}

\section{Approximation results for suprema of empirical processes}\label{sec:ResultsSupremaEP}

\subsection{Empirical process notation and definitions}\label{subsec:NotationEP}
In the previous section we got away with intuitive standard notation. We now introduce precise empirical process notation for the remainder of the paper: We denote by $X, X_1, X_2, \ldots$ a sequence of i.i.d. random variables taking values in a measurable space $(S, \mathcal{S})$ with common distribution $P$, i.e. $X_i : S^\infty \rightarrow S$, $i \geq 1$, are the coordinate projections of the infinite product probability space $(\Omega, \mathcal{A}, \mathbb{P}) = \left(S^\infty, \mathcal{S}^\infty, P^\infty\right)$ with law $\mathbb{P}_{X_i} = P$. If auxiliary variables independent of the $X_i$'s are involved, the underlying probability space is assumed to be of the form $(\Omega, \mathcal{A}, \mathbb{P}) = \left(S^\infty, \mathcal{S}^\infty, P^\infty\right) \times (Z, \mathcal{Z}, Q)$. We define the empirical measures $P_n$ associated with observations $X_1, \ldots, X_n$ as random measures on $\left(S^\infty,\mathcal{S}^\infty\right)$ given by $P_n(\omega) := n^{-1}\sum_{i=1}^n \delta_{X_i(\omega_i)}$ for all $\omega \in S^\infty$, where $\delta_{x}$ is the Dirac measure at $x$. 

For a class $\mathcal{F}$ of measurable functions from $S$ onto the real line $\mathbb{R}$ we define the \emph{empirical process indexed by $\mathcal{F}$} as
\begin{align*}
	\mathbb{G}_n(f) : = \frac{1}{\sqrt{n}}\sum_{i=1}^n (f(X_i) - P f\big), \quad{} f \in \mathcal{F}.
\end{align*}
Further, we denote by $\{G_P(f) : f \in \mathcal{F}\}$ the \emph{Gaussian $P$-bridge process} with  mean zero and the same covariance function $\mathcal{C}_P : \mathcal{F} \times \mathcal{F} \rightarrow \mathbb{R}$ as the process $\{f(X): f \in \mathcal{F}\}$, i.e.
\begin{align}\label{eq:subsec:NotationEP-1}
	(f, g) \mapsto \mathcal{C}_P(f, g) := \mathrm{E}[G_P(f) G_P(g)] = Pfg - (Pf)(Pg).
\end{align}
Moreover, we denote by $\{Z_Q(f) : f \in \mathcal{F}\}$ the \emph{Gaussian $Q$-motion} with mean zero and covariance function covariance function $\mathcal{E}_Q : \mathcal{F} \times \mathcal{F} \rightarrow \mathbb{R}$ given by
\begin{align}\label{eq:subsec:NotationEP-2}
	(f, g) \mapsto \mathcal{E}_Q(f, g) := \mathrm{E}[Z_Q(f) Z_Q(g)] = Qfg.
\end{align}

For probability measures $Q$ on $(S, \mathcal{S})$ we define the $L_q(Q)$-norm, $q \geq 1$, for function $f \in \mathcal{F}$ by $\|f\|_{Q,q} = (Q|f|^q)^{1/q}$, the $L_2(Q)$-semimetric by $e_Q(f, g) := \|f- g\|_{Q,2}$ and the intrinsic standard deviation metric by $d_P(f,g) := e_P(f - Pf, g - Pg) = \sqrt{P(f-g)^2 - (P(f-g))^2}$, $f, g \in \mathcal{F}$. We denote by $L_q(S, \mathcal{S}, Q)$, $q\geq 1$, the space of all real-valued measurable functions $f$ on $(S, \mathcal{S})$ with finite $L_q(Q)$-norm. A function $F: S\rightarrow \mathbb{R}$ is an \emph{envelop} for the class $\mathcal{F}$ satisfies $|f(x)| \leq F(x)$ for all $f \in \mathcal{F}$ and $x \in S$.
For a semimetric space $(T, d)$ and any $\varepsilon > 0$ we write $N(T, d, \varepsilon)$ to denote the $\varepsilon$\emph{-covering numbers} of $T$ with respect to $d$.
For two deterministic sequences $(a_n)_{n\geq 1}$ and $(b_n)_{n\geq 1}$, we write $a_n \lesssim b_n$ if $a_n = o(b_n)$ and $a_n \asymp b_n$ if there exist absolute constants $C_1, C_2 > 0$ such that $C_1 b_n \leq a_n \leq C_2 b_n$ for all $n \geq 1$. For $a, b \in \mathbb{R}$ we write $a \wedge b$ for $\min\{a, b\}$ and $a \vee b$ for $\max\{a, b\}$.

\subsection{Gaussian approximation inequalities}\label{subsec:GaussianApproximation}
We present the first main theoretical result of this paper: a Gaussian approximation inequality for empirical processes indexed by function classes that need not be Donsker. 
This result generalizes Proposition~\ref{lemma:CLT-Max-Norm} from finite dimensional index sets to general function classes. 

\begin{theorem}[Gaussian approximation]\label{theorem:CLT-Max-Norm-Simultaneous} 
	Let $(\mathcal{F}_n, \rho)$ be a totally bounded pseudo-metric space. Further, let $\mathcal{F}_n \subset L_2(S, \mathcal{S}, P)$ have envelope function  $F_n \in  L_3(S, \mathcal{S}, P)$. Suppose that there exist functions $\psi_n, \phi_n$ such that for $\delta > 0$,
	\begin{align}\label{eq:theorem:CLT-Max-Norm-Simultaneous-0}
		\mathrm{E} \|G_P\|_{\mathcal{F}_{n,\delta}'} \lesssim 	\psi_n(\delta) \sqrt{\mathrm{Var}(\|G_P\|_{\mathcal{F}_n})} \quad{} \quad{} \text{and} \quad{} \quad{} \big\|\|\mathbb{G}_n\|_{\mathcal{F}_{n,\delta}'} \big\|_{P,1} \lesssim 	\phi_n(\delta)\sqrt{\mathrm{Var}(\|G_P\|_{\mathcal{F}_n})} ,
	\end{align}
	where $\mathcal{F}_{n, \delta}' = \{ f - g: f, g \in \mathcal{F}_n, \: \rho(f, g) < \delta\|F_n\|_{P,2}\}$. Let $r_n = \inf\big\{  \psi_n(\delta) \vee \phi_n(\delta): \delta > 0\big\}$. Then, for each $M \geq 0$,
	\begin{align*}
		&\sup_{s \geq 0} \Big|\mathbb{P}\left\{\|\mathbb{G}_n\|_{\mathcal{F}_n} \leq s \right\} -\mathbb{P}\left\{\|G_P\|_{\mathcal{F}_n} \leq s \right\} \Big|\\
		&\quad{}\quad{}\quad\lesssim \frac{\|F_n\|_{P, 3}}{\sqrt{n^{1/3}\mathrm{Var}(\|G_P\|_{\mathcal{F}_n})}} + \frac{\|F_n\mathbf{1}\{F_n > M\}\|_{P,3}^3}{\|F_n\|_{P,3}^3} + \frac{ \mathrm{E} \|G_P\|_{\mathcal{F}_n} + M }{\sqrt{n\mathrm{Var}(\|G_P\|_{\mathcal{F}_n})}} + \sqrt{r_n},
	\end{align*}
	where $\lesssim$ hides an absolute constant independent of $n, r_n, M, \mathcal{F}_n, F_n, \psi_n, \phi_n$, and $P$.
\end{theorem}
\begin{remark}[Lower and upper bounds on $\mathrm{Var}(\|G_P\|_{\mathcal{F}_n}$]
	Since $\mathcal{F}_n$ may change with the sample size, so may $\mathrm{Var}(\|G_P\|_{\mathcal{F}_n}$. In a few special cases it is possible to exactly compute the variance~\citep{giessing2023bootstrap}, but in most cases one can only derive lower and upper bounds. In our companion work~\cite{giessing2022anticoncentration} we derive such bounds under mild conditions on the Gaussian $P$-bridge process. For the reader's convenience we include the relevant results from this paper in Appendix~\ref{subsec:AuxiliaryResults-AntiConcentration}.
\end{remark}

The total boundedness of $\mathcal{F}_n$ in above theorem is a standard assumption which allows us to reduce the proof of Theorem~\ref{theorem:CLT-Max-Norm-Simultaneous} to an application of Proposition~\ref{lemma:CLT-Max-Norm} combined with a discretization argument. Since $\mathcal{F}_n$ is totally bounded whenever $\mathrm{E}\|G_P\|_{\mathcal{F}_n} < \infty$ (i.e. Lemma~\ref{lemma:Sudakov} in Appendix~\ref{subsec:GaussianProcesses}), this is a rather mild technical assumption. For a detailed discussion of this result and a comparison with the literature we refer to Section~\ref{subsec:RelationPreviousWork-Stat}. 

Under more specific smoothness assumptions on the Gaussian $P$-bridge process (beyond the abstract control of the moduli of continuity in eq.~\eqref{eq:theorem:CLT-Max-Norm-Simultaneous-0}) Theorem~\ref{theorem:CLT-Max-Norm-Simultaneous} simplifies as follows: 

\begin{corollary}\label{corollary:theorem:CLT-Max-Norm-Simultaneous}
	Let $\mathcal{F}_n \subset L_2(S, \mathcal{S}, P)$ and have envelope $F_n \in  L_3(S, \mathcal{S}, P)$. If the Gaussian $P$-Bridge process $\{G_P(f) : f \in \mathcal{F}_n\}$ has almost surely uniformly $d_P$-continuous sample paths and $\mathrm{E}\|G_P\|_{\mathcal{F}_n} < \infty$, then for each $M \geq 0$,
	\begin{align*}
		&\sup_{s \geq 0} \Big|\mathbb{P}\left\{\|\mathbb{G}_n\|_{\mathcal{F}_n} \leq s \right\} -\mathbb{P}\left\{\|G_P\|_{\mathcal{F}_n} \leq s \right\} \Big|\\
		&\quad{}\quad{}\quad{} \lesssim \frac{\|F_n\|_{P, 3}}{\sqrt{n^{1/3}\mathrm{Var}(\|G_P\|_{\mathcal{F}_n})}} + \frac{\|F_n\mathbf{1}\{F_n > M\}\|_{P,3}^3}{\|F_n\|_{P,3}^3} + \frac{ \mathrm{E} \|G_P\|_{\mathcal{F}_n} + M}{\sqrt{n\mathrm{Var}(\|G_P\|_{\mathcal{F}_n})}},
	\end{align*}
	where $\lesssim$ hides an absolute constant independent of $n, M, \mathcal{F}_n, F_n$, and $P$.
\end{corollary}

\begin{remark}[Further simplification with lower bound on second moment]
	Recall the setup of Corollary~\ref{corollary:theorem:CLT-Max-Norm-Simultaneous}.  If in addition there exists $\kappa_n > 0$ such that $P f^2 \geq \kappa_n^2 > 0$ for all $f \in \mathcal{F}_n$ and $F_n \in  L_{3 + \delta}(S, \mathcal{S}, P)$, $\delta > 0$, then
	\begin{align*}
		&\sup_{s \geq 0} \Big|\mathbb{P}\left\{\|\mathbb{G}_n\|_{\mathcal{F}_n} \leq s \right\} -\mathbb{P}\left\{\|G_P\|_{\mathcal{F}_n} \leq s \right\} \Big|\\
		&\quad\quad\quad\lesssim  \frac{1}{n^{1/6} \kappa^2_n}  \left(\|F_n\|_{P, 3 + \delta}^{3 + \delta}+ \mathrm{E}[\|G_P\|_{\mathcal{F}_n}^{3 + \delta}]\right)^{\frac{1}{3 + \delta}}  \mathrm{E}[\|G_P\|_{\mathcal{F}_n}]  + \frac{1}{n^{\delta/3}} \left(\frac{\|F_n\|_{P, 3 + \delta}}{\|F_n\|_{P, 3}}\right)^3,
	\end{align*}
	where  $\lesssim$ hides an absolute constant independent of $n, d, \kappa_n^2$, and the distribution of the $X_i$'s. This is the empirical process analogue to Corollary~\ref{corollary:lemma:CLT-Max-Norm}. Since the proof of this inequality is is identical to the one of Corollary~\ref{corollary:lemma:CLT-Max-Norm} and only differs in notation, we omit the details.
\end{remark}

Since centered Gaussian processes are either almost surely uniformly continuous w.r.t. their intrinsic standard deviation metric or almost surely discontinuous, the smoothness condition in Corollary~\ref{corollary:theorem:CLT-Max-Norm-Simultaneous} is natural. (For the proofs to go through uniform $d_P$-continuity in probability would be sufficient.) In Lemmas~\ref{lemma:MetricEntropCondition} and~\ref{lemma:ModulusContinuity} in Appendix~\ref{subsec:GaussianProcesses} we provide simple sufficient and necessary conditions for almost sure uniform $d_P$-continuity of $G_P$ on $\mathcal{F}_n$. Importantly, $n \geq 1$ is fixed; the uniform $d_P$-continuity is a point-wise requirement and does not need to hold when taking the limit $n \rightarrow \infty$. This is especially relevant in high-dimensional statistics where Gaussian processes are typically unbounded and have diverging sample paths as $n \rightarrow \infty$, and, more generally, whenever the function classes are non-Donsker. Of course, just as with Proposition~\ref{lemma:CLT-Max-Norm} the conclusion of Corollary~\ref{corollary:theorem:CLT-Max-Norm-Simultaneous} is weaker than a weaker convergence (see Section~\ref{subsec:RelationPreviousWork-Math}).

\subsection{Gaussian comparison inequalities}\label{subsec:GaussianComparison}
Our second main result is a comparison inequality which bounds the Kolmogorov distance between a Gaussian $P$-bridge process and a Gaussian $Q$-motion both indexed by (possibly) different function classes $\mathcal{F}_n$ and $\mathcal{G}_n$. This generalizes Proposition~\ref{lemma:GaussianComparison} to Gaussian processes.

To state the theorem, we introduce the following notation: Given the pseudo-metric spaces $(\mathcal{F}_n, \rho_i)$ and $(\mathcal{G}_n, \rho_i)$, $i \in \{1,2\}$, we write $\pi^1_{\mathcal{G}_n} : \mathcal{F}_n \cup \mathcal{G}_n \rightarrow \mathcal{G}_n$ to denote the projection from $\mathcal{F}_n \cup \mathcal{G}_n$ onto $\mathcal{G}_n$ defined by $\rho_1(h, \pi^1_{\mathcal{G}_n}h) = \inf_{g \in \mathcal{G}_n}\rho_1(h, g)$ for all $h \in \mathcal{F}_n \cup \mathcal{G}_n$. In the case that the image $\pi^1_{\mathcal{G}_n}h$ is not a singleton for some $h \in \mathcal{F}_n \cup \mathcal{G}_n$, we choose any one of the equivalent points in $\mathcal{G}_n$. Similarly, we write $\pi^2_{\mathcal{F}_n}: \mathcal{F}_n \cup \mathcal{G}_n \rightarrow \mathcal{F}_n$ for the projection from $\mathcal{F}_n \cup \mathcal{G}_n$ onto $\mathcal{F}_n$ defined via $\rho_2$.

\begin{theorem}[Gaussian comparison]\label{theorem:GaussianComparison-PQ-Simultaneous}
	Let $(\mathcal{F}_n, \rho_1)$ and $(\mathcal{G}_n, \rho_2)$ be totally bounded pseudo-metric spaces. Further, let $\mathcal{F}_n \subset L_2(S, \mathcal{S}, P)$ and $\mathcal{G}_n \subset L_2(S, \mathcal{S}, Q)$ have envelope functions $F_n \in  L_2(S, \mathcal{S}, P)$ and $G_n \in L_2(S, \mathcal{S}, Q)$, respectively. Suppose that there exist functions $\psi_n, \phi_n$ such that 
	\begin{align*}
		\mathrm{E} \|G_P\|_{\mathcal{F}_{n,\delta}'} \lesssim 	\psi_n(\delta)\sqrt{\mathrm{Var}(\|G_P\|_{\mathcal{F}_n})} \quad{} \quad{} \text{and} \quad{} \quad{} \mathrm{E} \|Z_Q\|_{\mathcal{G}_{n,\delta}'} \lesssim 	\phi_n(\delta)\sqrt{\mathrm{Var}(\|Z_Q\|_{\mathcal{G}_n})},
	\end{align*}
	where $\mathcal{F}_{n, \delta}' = \{ f - g: f, g \in \mathcal{F}_n, \: \rho_1(f, g) < \delta\|F_n\|_{P,2}\}$ and $\mathcal{G}_{n,\delta}' = \{ f - g: f, g \in \mathcal{G}_n, \: \rho_2(f, g) < \delta\|G_n\|_{Q,2}\}$. Let $r_n = \inf\big\{  \psi_n(\delta) \vee \phi_n(\delta): \delta > 0\big\}$. Then,
	\begin{align*}
		&\sup_{s \geq 0} \Big|\mathbb{P}\left\{\|G_P\|_{\mathcal{F}_n} \leq s \right\} -\mathbb{P}\left\{\|Z_Q\|_{\mathcal{G}_n} \leq s \right\}\Big| \lesssim \left( \frac{\Delta_{P,Q}(\mathcal{F}_n, \mathcal{G}_n)}{\mathrm{Var}(\|G_P\|_{\mathcal{F}_n} ) \vee \mathrm{Var}(\|Z_Q\|_{\mathcal{G}_n} )} \right)^{1/3} + \sqrt{r_n},
	\end{align*}
	where 
	\begin{align}\label{eq:theorem:GaussianComparison-PQ-Simultaneous-0}
		\begin{split}
		\Delta_{P,Q}(\mathcal{F}_n, \mathcal{G}_n) &:= \sup_{f_1, f_2 \in \mathcal{F}_n} \big| \mathcal{C}_P(f_1, f_2) - \mathcal{E}_Q(\pi^1_{\mathcal{G}_n} f_1), \pi^1_{\mathcal{G}_n}f_2) \big|\\
		&\quad \quad \quad \bigvee \sup_{g_1, g_2 \in \mathcal{G}_n }  \big|\mathcal{E}_Q(g_1, g_2) - \mathcal{C}_P(\pi^2_{\mathcal{F}_n} g_1, \pi^2_{\mathcal{F}_n} g_2) \big|,
		\end{split}	
	\end{align}
	where $\lesssim$ hides an absolute constant independent of $n, \mathcal{F}_n, \mathcal{G}_n, F_n, G_n, \psi_n, \phi_n, P$ and $Q$.
\end{theorem}
\begin{remark}[Comparison inequality for Gaussian $Q$-bridge processes]\label{remark:theorem:GaussianComparison-PQ-Simultaneous}
	Theorem~\ref{theorem:GaussianComparison-PQ-Simultaneous} holds without changes (and identical proof) for Gaussian $Q$-bridge processes $\{G_Q(f) : f \in \mathcal{F}_n\}$. This easily follows from the almost sure representation $G_Q(f)  = Z_Q(f) + (Pf)Z$ for all $f \in \mathcal{F}_n$, where $Z \sim N(0,1)$ is independent of $Z_Q(f)$ for all $f \in \mathcal{F}_n$. We state the comparison inequality for $Q$-motions because this is the version that we use in the context of the Gaussian process bootstrap in Section~\ref{subsec:GaussianProcessBootstrap}.
\end{remark}
Informally speaking, Theorem~\ref{theorem:GaussianComparison-PQ-Simultaneous} states that the distributions of the suprema of any two Gaussian processes are close whenever their covariance functions are not too different. Note that we compare two Gaussian processes with different covariance functions $\mathcal{C}_P$ and $\mathcal{E}_Q$ that differ not only in their measures ($P$ and $Q$) but also their functional form (viz. eq.~\eqref{eq:subsec:NotationEP-1} and~\eqref{eq:subsec:NotationEP-2} in Section~\ref{subsec:NotationEP}). This turns out to be essential for developing alternatives to the classical Gaussian multiplier bootstrap procedure for suprema of empirical processes proposed by~\cite{chernozhukov2014AntiConfidenceBands} and~\cite{chernozhukov2014GaussianApproxEmp}. We discuss this in detail in Section~\ref{subsec:GaussianProcessBootstrap}.

Under additional smoothness conditions on the Gaussian processes and assumptions on the function classes Theorem~\ref{theorem:GaussianComparison-PQ-Simultaneous} simplifies. The special cases $\mathcal{G}_n = \mathcal{F}_n$ and $\mathcal{G}_n \subseteq \mathcal{F}_n$ are of particular interest in the context of the bootstrap. They are the content of the next two corollaries.

\begin{corollary}\label{corollary:theorem:GaussianComparison-PQ-Simultaneous}
	Let $\mathcal{F}_n \subset L_2(S, \mathcal{S}, P) \cap  L_2(S, \mathcal{S}, Q)$. If the Gaussian $P$-Bridge process $\{G_P(f) : f \in \mathcal{F}_n\}$ and the Gaussian $Q$-motion $\{Z_Q(f) : f \in \mathcal{F}_n\}$ have almost surely uniformly continuous sample paths w.r.t. their respective standard deviation metrics $d_P$ and $d_Q$, and $\mathrm{E}\|G_P\|_{\mathcal{F}_n} \vee \mathrm{E}\|Z_Q\|_{\mathcal{F}_n}< \infty$, then
	\begin{align*}
		&\sup_{s \geq 0} \Big|\mathbb{P}\left\{\|G_P\|_{\mathcal{F}_n} \leq s \right\} -\mathbb{P}\left\{\|Z_Q\|_{\mathcal{F}_n} \leq s \right\}\Big| \lesssim \left( \frac{\sup_{f, g \in \mathcal{F}_n} \big| \mathcal{C}_P(f,g) - \mathcal{E}_Q(f, g) \big|}{ \mathrm{Var}(\|G_P\|_{\mathcal{F}_n} ) \vee \mathrm{Var}(\|Z_Q\|_{\mathcal{F}_n} )} \right)^{1/3},
	\end{align*}
	where 
	 $\lesssim$ hides an absolute constant independent of $n, \mathcal{F}_n, P$, and $Q$.
\end{corollary}

\begin{corollary}\label{corollary:theorem:GaussianComparison-PQ-Simultaneous-2}
	Let $\mathcal{F}_n \subset L_2(S, \mathcal{S}, P)$ and $\mathcal{G}_n \subseteq \mathcal{F}_n$ be a $\delta \|F\|_{P,2}$-net of $\mathcal{F}_n$ with respect to $d_P$, $\delta > 0$. If the Gaussian $P$-Bridge process $\{G_P(f) : f \in \mathcal{F}_n\}$ has almost surely uniformly continuous sample paths w.r.t. metric $d_P$ and $\mathrm{E}\|G_P\|_{\mathcal{F}_n} < \infty$, then
	\begin{align*}
		&\sup_{s \geq 0} \Big|\mathbb{P}\left\{\|G_P\|_{\mathcal{F}_n} \leq s \right\} -\mathbb{P}\left\{\|G_P\|_{\mathcal{G}_n} \leq s \right\}\Big| \lesssim \left( \frac{ \delta\|F\|_{P,2} \sup_{f \in \mathcal{F}_n} \sqrt{Pf^2}}{ \mathrm{Var}(\|G_P\|_{\mathcal{F}_n} ) \vee \mathrm{Var}(\|G_P\|_{\mathcal{G}_n} )} \right)^{1/3},
	\end{align*}
	where $\lesssim$ hides an absolute constant independent of $\delta, n, \mathcal{F}_n$, and $P$.
\end{corollary}

\subsection{Gaussian process bootstrap}\label{subsec:GaussianProcessBootstrap}
In this section we develop a general framework for bootstrapping the distribution of suprema of empirical processes. We recover the classical Gaussian multiplier bootstrap as a special case of this more general framework. For concrete applications to statistical problems we refer to Section~\ref{sec:Applications}.

The following abstract approximation result generalizes Proposition~\ref{lemma:Bootstrap-Max-Norm} to empirical processes.

\begin{theorem}[Abstract bootstrap approximation]\label{theorem:Abstract-Bootstrap-Sup-Empirical-Process}
	Let $(\mathcal{F}_n, \rho_1)$ and $(\mathcal{G}_n, \rho_2)$ be totally bounded pseudo-metric spaces. Further, let $\mathcal{F}_n \subset L_2(S, \mathcal{S}, P)$ and $\mathcal{G}_n \subset L_2(S, \mathcal{S}, Q)$ have envelope functions $F_n \in  L_3(S, \mathcal{S}, P)$ and $G_n \in L_2(S, \mathcal{S}, Q)$, respectively. Suppose that there exist functions $\psi_n, \phi_n$ such that for $\delta > 0$,
	\begin{align}\label{eq:theorem:Abstract-Bootstrap-Sup-Empirical-Process-0}
		\mathrm{E} \|G_P\|_{\mathcal{F}_{n,\delta}'}  \vee  \big\|\|\mathbb{G}_n\|_{\mathcal{F}_{n,\delta}'} \big\|_{P,1} \lesssim 	\psi_n(\delta)\sqrt{\mathrm{Var}(\|G_P\|_{\mathcal{F}_n})}, \quad \mathrm{E} \|Z_Q\|_{\mathcal{G}_{n,\delta}'} \lesssim \phi_n(\delta)\sqrt{\mathrm{Var}(\|Z_Q\|_{\mathcal{G}_n})},
	\end{align}
	where $\mathcal{F}_{n, \delta}' = \{ f - g: f, g \in \mathcal{F}_n, \: \rho_1(f, g) < \delta\|F_n\|_{P,2}\}$ and $\mathcal{G}_{n,\delta}' = \{ f - g: f, g \in \mathcal{G}_n, \: \rho_2(f, g) < \delta\|G_n\|_{Q,2}\}$. Let $r_n = \inf\big\{  \psi_n(\delta) \vee \phi_n(\delta): \delta > 0\big\}$. Then, for each $M \geq 0$,
	\begin{align*}
		&\sup_{s \geq 0} \Big|\mathbb{P}\left\{\|\mathbb{G}_n\|_{\mathcal{F}_n} \leq s \right\} -\mathbb{P}\left\{\|Z_Q\|_{\mathcal{G}_n} \leq s \right\} \Big|\\
		&\quad{}\quad{}\quad{} \lesssim \frac{\|F_n\|_{P, 3}}{\sqrt{n^{1/3}\mathrm{Var}(\|G_P\|_{\mathcal{F}_n})}} + \frac{\|F_n\mathbf{1}\{F_n > M\}\|_{P,3}^3}{\|F_n\|_{P,3}^3} + \frac{ \mathrm{E} \|G_P\|_{\mathcal{F}_n} + M}{\sqrt{n\mathrm{Var}(\|G_P\|_{\mathcal{F}_n})}}\\
		&\quad{}\quad{}\quad{}\quad{} \quad{} \quad{} + \left( \frac{\Delta_{P,Q}(\mathcal{F}_n, \mathcal{G}_n)}{\mathrm{Var}(\|G_P\|_{\mathcal{F}_n} ) \vee \mathrm{Var}(\|Z_Q\|_{\mathcal{G}_n} )} \right)^{1/3} + \sqrt{r_n},
	\end{align*}
	where $\Delta_{P,Q}(\mathcal{F}_n, \mathcal{G}_n)$ is given in eq.~\eqref{eq:theorem:GaussianComparison-PQ-Simultaneous-0} and $\lesssim$ hides an absolute constant independent of $n, r_n, M,$ $\mathcal{F}_n,$ $\mathcal{G}_n, F_n, G_n, \psi_n, \phi_n, P$, and $Q$.
\end{theorem}
\begin{remark}[Bootstrap approximation with Gaussian $Q$-bridge process]\label{remark:theorem:Abstract-Bootstrap-Sup-Empirical-Process}
	Theorem~\ref{theorem:Abstract-Bootstrap-Sup-Empirical-Process} holds without changes (and identical proof) with the Gaussian $Q$-bridge process $\{G_Q(f) : f \in \mathcal{F}_n\}$ substituted for the Gaussian $Q$-motion $\{Z_Q(f) : f \in \mathcal{F}_n\}$ (see also Remark~\ref{remark:theorem:GaussianComparison-PQ-Simultaneous}).
\end{remark}

Theorem~\ref{theorem:Abstract-Bootstrap-Sup-Empirical-Process} is an immediate consequence of Theorems~\ref{theorem:CLT-Max-Norm-Simultaneous} and~\ref{theorem:GaussianComparison-PQ-Simultaneous} and the triangle inequality, and thus a mathematical triviality. What matters is its statistical interpretation: It implies that we can approximate the sampling distribution of the supremum of an empirical process $\{\mathbb{G}_n (f) : f \in \mathcal{F}_n\}$ with the distribution of the supremum of a Gaussian $Q$-motion $\{Z_Q(g) : g\in \mathcal{G}_n \}$ provided that their covariance functions $\mathcal{C}_P$ and $\mathcal{E}_Q$ do not differ by too much. Importantly, up to the smoothness condition in~\eqref{eq:theorem:Abstract-Bootstrap-Sup-Empirical-Process-0}, we are completely free to choose whatever $Q$-measure and function class $\mathcal{G}_n$ suit us best. 

This interpretation of Theorem~\ref{theorem:Abstract-Bootstrap-Sup-Empirical-Process} motivates the following bootstrap procedure for estimating the distribution of the supremum of an empirical process $\{\mathbb{G}_n (f) : f \in \mathcal{F}_n\}$ when the $P$-measure is unknown.

\begin{algorithm}[Gaussian process bootstrap]\label{algorithm:GaussianProcessBootstrap} Let $X_1,\ldots, X_n$ be a simple random sample drawn from distribution $P$.
	\begin{itemize}[leftmargin=60pt]
		\item[Step 1:] Construct a positive semi-definite estimate $\widehat{\mathcal{C}}_n$ of the covariance function $\mathcal{C}_P$ of the empirical process $\{\mathbb{G}_n(f) : f \in \mathcal{F}_n\}$. 
		
		\item[Step 2:] Construct a Gaussian process $\widehat{Z}_n = \{\widehat{Z}_n(f): f \in \mathcal{F}_n\}$ such that for all $f,g \in \mathcal{F}_n$,
		\begin{align}\label{eq:theorem:Abstract-Bootstrap-Sup-Empirical-Process-1}
			\mathrm{E}\big[\widehat{Z}_n(f)\mid X_1, \ldots, X_n\big] = 0 \quad \mathrm{and} \quad
			\mathrm{E}\big[\widehat{Z}_n(f)\widehat{Z}_n(g) \mid X_1, \ldots, X_n\big] = \widehat{\mathcal{C}}_n(f,g).
		\end{align}
		
		\item[Step 3:] Approximate the distribution of $\|\mathbb{G}_n\|_{\mathcal{F}_n}$ by drawing Monte Carlo samples from $\|\widehat{Z}_n\|_{\mathcal{F}_n}$. 
	\end{itemize}
\end{algorithm}
By Kolmogorov's consistency theorem the Gaussian process $\widehat{Z}_n$ will always exist. If the estimate $\widehat{\mathcal{C}}_n$ is uniformly consistent and if the envelope function $F_n$ and the supremum of the Gaussian $P$-bridge process satisfy certain moment conditions, Theorem~\ref{theorem:Abstract-Bootstrap-Sup-Empirical-Process} readily implies uniform consistency of the Gaussian process bootstrap as $n \rightarrow \infty$. The challenge is to actually construct a (measurable) version of $\widehat{Z}_n$ from which we can draw Monte Carlo samples. This is the content the next section.


\begin{remark}[Relation to the \emph{Gaussian multiplier bootstrap} by~\cite{chernozhukov2014GaussianApproxEmp}]
	The Gaussian multiplier bootstrap is a special case of the Gaussian process bootstrap. To see this, note that if we estimate the covariance function $\mathcal{C}_P$ nonparametrically via the sample covariance function $\mathcal{C}_{P_n}(f,g) = P_n fg - (P_nf)(P_ng)$, $f,g \in \mathcal{F}_n$, then the corresponding Gaussian $P_n$-Bridge process $\{G_{P_n}(f) : f \in \mathcal{F}_n\}$ can be expanded into a finite sum of non-orthogonal functions
	\begin{align*}
		G_{P_n}(f) = \frac{1}{\sqrt{n}}\sum_{i=1}^n \xi_i \big(f(X_i) - P_n f\big), \quad{} f \in \mathcal{F}_n,
	\end{align*}
	where $\xi_1, \ldots, \xi_n$ are i.i.d. standard normal random variables independent of $X_1, \ldots, X_n$. While this representation is extremely simple (in theory and practice!), unlike the Gaussian process bootstrap it does not allow integrating additional structural information about the covariance function. Empirically, this often results in less accurate bootstrap estimates~\citep{giessing2023bootstrap}.
\end{remark}

\subsection{A practical guide to the Gaussian process bootstrap}\label{subsec:GaussianProcessBootstrap-Implementation-Consistency}
To implement the Gaussian process bootstrap we need two things: a uniformly consistent estimate of the covariance function and a systematic way of constructing Gaussian processes for given covariance functions. 

Finding consistent estimates of the covariance function is relatively straightforward. For example, a natural estimate is the nonparametric sample covariance function defined as $\mathcal{C}_{P_n}(f,g) := P_n(fg) - (P_nf)(P_n g)$, $f, g \in \mathcal{F}_n$. Under suitable conditions on the probability measure $P$ and the function class $\mathcal{F}_n$ this estimate is uniformly consistent in high dimensions~\cite[e.g.][]{koltchinskiiConcentration2017}. Of course, an important feature of the Gaussian process bootstrap is that it can be combined with any positive semi-definite estimate of the covariance function. In particular, we can use (semi-)parametric estimators to exploit structural constraints induced by $P$ and $\mathcal{F}_n$. For concrete examples we refer to Section~\ref{sec:Applications} and also~\cite{giessing2023bootstrap}.

Constructing Gaussian processes for given covariance functions and defined on potentially arbitrary index sets is a more challenging problem. We propose to base the construction on an almost sure version of the classical Karhunen-Lo{\`e}ve decomposition.

To develop this decomposition in the framework of this paper, we need new notation and concepts. In the following, $(\mathcal{F}_n, \mathcal{B}_n, \mu)$ denotes a measurable space for some finite measure $\mu$ with support $\mathcal{F}_n$. Typically, $\mathcal{B}_n$ is the Borel $\sigma$-algebra on $\mathcal{F}_n$ and the measure $\mu$ is chosen for convenience. For example, $\mu$ can be set to the Lebesgue measure when $\mathcal{F}$ an interval, or the counting measure when $\mathcal{F}_n$ is a finite (discrete) set. The space $L_2(\mathcal{F}_n, \mathcal{B}_n, \mu)$ equipped with the inner product $\langle \psi, \phi \rangle := \int_{\mathcal{F}_n} \psi(f) \phi(f) d\mu (f)$ is a Hilbert space. Given a positive semi-definite and continuous kernel  $K : \mathcal{F}_n \times \mathcal{F}_n \rightarrow \mathbb{R}$ we define the linear operator $T_K$ on $L_2(\mathcal{F}_n, \mathcal{B}_n, \mu)$ via
\begin{align}\label{eq:subsec:GaussianProcessBootstrap-Implementation-Consistency-1}
	(T_K\psi)(g) := \langle K( \cdot, g), \psi \rangle = \int_{\mathcal{F}_n} K(f,g ) \psi(f) d\mu(f) , \quad \psi \in  L^2(\mathcal{F}_n, \mathcal{B}_n, \mu), \quad g \in \mathcal{F}_n.
\end{align}
If $\int_{\mathcal{F}_n \times \mathcal{F}_n} K^2(f, g) d\mu(f) d\mu(g) < \infty$, then $T_K$ is a bounded linear operator. If $(\mathcal{F}_n, d_Q)$ is a compact metric space, then $T_K$ is a compact operator. 
In this case, $T_K$ has a spectral decomposition which depends on the kernel $K$ alone; the measure $\mu$ is exogenous. For further details we refer the reader to Chapters 4 and 7 in~\cite{hsing2015theoretical}.

The following result is well-known~\citep{jain1970note}; since it is slightly adapted to our setting we have included its proof in the appendix.

\begin{proposition}[Karhunen-Lo{\`e}ve decomposition of Gaussian processes]\label{lemma:KarhunenLoeve-GP} Let $(\mathcal{F}_n, e_Q)$ be a compact pseudo-metric space and $Z_Q = \{Z_Q(f) : f \in \mathcal{F}_n\}$ be a Gaussian $Q$-motion. Suppose that $Z_Q$ has a continuous covariance function $\mathcal{E}_Q$ and $\mathrm{E}\|Z_Q\|_{\mathcal{F}_n} < \infty$. Let $\{(\lambda_k,\varphi_k )\}_{k=1}^\infty$ be the eigenvalue and eigenfunction pairs of the linear operator $T_{\mathcal{E}_Q}$ corresponding to $\mathcal{E}_Q$, and $\{\xi_k\}_{k=1}^\infty$ be a sequence of i.i.d. standard normal random variables. Then, with probability one,
	\begin{align*}
		\lim_{m \rightarrow \infty}\|Z_Q^m - Z_Q\|_{\mathcal{F}_n} = 0,
	\end{align*}
	where $Z_Q^m(f) :=  \sum_{k=1}^m \xi_k \sqrt{\lambda_k} \varphi_k(f)$, $f \in \mathcal{F}_n$, is an almost surely bounded and continuous Gaussian process on $\mathcal{F}_n$ for all $m \in \mathbb{N}$. 
\end{proposition}

This proposition justifies the following constructive approximation of the Gaussian proxy process $\widehat{Z}_n$ defined in~\eqref{eq:theorem:Abstract-Bootstrap-Sup-Empirical-Process-1}: Let $\widehat{\mathcal{C}}_n$ be a generic estimate of the covariance function $\mathcal{C}_P$ of the empirical process $\{\mathbb{G}_n(f) : f \in \mathcal{F}_n\}$. Suppose that $\widehat{\mathcal{C}}_n$ is such that the associated integral operator $T_{\widehat{\mathcal{C}}_n}$ defined in~\eqref{eq:subsec:GaussianProcessBootstrap-Implementation-Consistency-1} admits a spectral decomposition. Denote by $\big\{(\widehat{\lambda}_k, \widehat{\varphi}_k)\big\}_{k=1}^\infty$ its eigenvalue and eigenfunction pairs, and, for concreteness, assume that the $\widehat{\lambda}_k$'s are sorted in nonincreasing order.  Given a sequence of i.i.d. standard normal random variables $\{\xi_k\}_{k=1}^\infty$ and truncation level $m \geq 1$, define the Gaussian process and associated covariance function
\begin{align}\label{eq:subsec:GaussianProcessBootstrap-Implementation-Consistency-2}
	\widehat{Z}_n^m(f) := \sum_{k=1}^m \xi_k \sqrt{\widehat{\lambda}_k} \widehat{\varphi}_k(f), \quad \quad \widehat{\mathcal{C}}_n^m(f, g) :=  \sum_{k=1}^m \widehat{\lambda}_k \widehat{\varphi}_k(f)\widehat{\varphi}_k(g),  \quad f, g \in \mathcal{F}_n. 
\end{align}
While there are only a few situations in which the eigenvalues and eigenfunctions of $T_{\widehat{\mathcal{C}}_n}$ can be found analytically, from a computational perspective this is a standard problem and efficient numerical solvers exist~\citep[][]{ghanem2003stochastic, berlinet2011reproducing}. Thus, constructing the process $\widehat{Z}_n^m$ in~\eqref{eq:subsec:GaussianProcessBootstrap-Implementation-Consistency-2} does not pose any practical challenges.  Proposition~\ref{lemma:KarhunenLoeve-GP} now guarantees (under appropriate boundedness and smoothness conditions on $\widehat{Z}_n$) that the approximation error between $\widehat{Z}_n$ and $\widehat{Z}_n^m$ can be made arbitrarily small by choosing $m \geq 1$ sufficiently large. In fact, by Mercer's theorem the error can be quantified in terms of the operator norm of the difference between the covariance functions $\widehat{\mathcal{C}}_n$ and $\widehat{\mathcal{C}}_n^m$~\citep[e.g.][Lemma 4.6.6 and Corollary~\ref{corollary:theorem:Abstract-Bootstrap-Sup-Empirical-Process-1}]{hsing2015theoretical}.

In above discussion we have implicitly imposed several assumptions on the estimate $\widehat{\mathcal{C}}_n$. For future reference, we summarize these assumptions in a single definition.
\begin{definition}[Admissibility of $\widehat{\mathcal{C}}_n$]
	We say that an estimate $\widehat{\mathcal{C}}_n : \mathcal{F}_n \times \mathcal{F} _n\rightarrow \mathbb{R}$ of the covariance function $\mathcal{C}_P$ is admissible if it is continuous, symmetric, positive semi-definite, and its associated integral operator $T_{\widehat{\mathcal{C}}_n}$ is a bounded linear operator on $L_2(\mathcal{F}_n, \mathcal{B}_n, \mu)$ for some finite measure $\mu$ with support $\mathcal{F}_n$. 
\end{definition}
\begin{remark}[Admissible estimates exist]
	Under the assumption that $\mathcal{F}_n$ has an envelope function, the nonparametric sample covariance $\mathcal{C}_{P_n}(f,g) = P_n(fg) - (P_nf)(P_n g)$ is admissible. Indeed, by definition, $\mathcal{C}_{P_n}$ is continuous (w.r.t. $e_{P_n}$), symmetric, positive semi-definite, and, by existence of the envelope function, $\int_{\mathcal{F}_n \times \mathcal{F}_n} \mathcal{C}_{P_n}^2(f, g) dP_n(f) dP_n(g) < \infty$. Also, the existence of an envelop function implies that $(\mathcal{F}_n, e_{P_n})$ is a compact metric space. See also Section~\ref{sec:Applications}.
\end{remark}

The next result establishes consistency of the Gaussian process bootstrap based on the truncated Karhunen-Lo{\`e}ve expansion in~\eqref{eq:subsec:GaussianProcessBootstrap-Implementation-Consistency-2}. It is a simple corollary of Theorem~\ref{theorem:Abstract-Bootstrap-Sup-Empirical-Process}. Note that the intrinsic standard deviation metric associated with (the pushforward probability measure induced by) $\widehat{Z}_n^m$ can be expressed in terms of its covariance function as $d_{\widehat{\mathcal{C}}_n^m}^2(f,g) = \widehat{\mathcal{C}}_n^m(f,f) +\widehat{\mathcal{C}}_n^m(g,g) - 2 \widehat{\mathcal{C}}_n^m(f,g)$. 

\begin{corollary}[Consistency of the Gaussian process bootstrap]\label{corollary:theorem:Abstract-Bootstrap-Sup-Empirical-Process-1}
	Let $\widehat{\mathcal{C}}_n$ be an admissible estimate of $\mathcal{C}_P$ and  $\widehat{\mathcal{C}}_n^m$ its best rank-$m$ approximation. Let $\mathcal{F}_n \subset L_2(S, \mathcal{S}, P)$ have envelope $F_n \in  L_3(S, \mathcal{S}, P)$. If the Gaussian processes $\{G_P(f) : f \in \mathcal{F}_n\}$ and $\{\widehat{Z}_n^m(f) : f \in \mathcal{F}_n\}$ have almost surely uniformly continuous sample paths w.r.t. their respective standard deviation metrics $d_P$ and $d_{\widehat{\mathcal{C}}_n^m}$, and $\mathrm{E}\|G_P\|_{\mathcal{F}_n} \vee \mathrm{E}\big[\|\widehat{Z}_n^m\|_{\mathcal{F}_n} \mid X_1, \ldots, X_n\big]< \infty$, then, for each $M \geq 0$,
	\begin{align*}
		&\sup_{s \geq 0} \Big|\mathbb{P}\left\{\|\mathbb{G}_n\|_{\mathcal{F}_n} \leq s \right\} -\mathbb{P}\left\{\|\widehat{Z}_n^m\|_{\mathcal{F}_n} \leq s \mid X_1, \ldots, X_n \right\}\Big| \\
		&\quad{}\quad{}\quad{} \lesssim \frac{\|F_n\|_{P, 3}}{\sqrt{n^{1/3}\mathrm{Var}(\|G_P\|_{\mathcal{F}_n})}} + \frac{\|F_n\mathbf{1}\{F_n > M\}\|_{P,3}^3}{\|F_n\|_{P,3}^3} + \frac{ \mathrm{E} \|G_P\|_{\mathcal{F}_n} + M}{\sqrt{n\mathrm{Var}(\|G_P\|_{\mathcal{F}_n})}}\\
		&\quad{}\quad{}\quad{}\quad{} \quad{} \quad{} + \left( \frac{\sup_{f, g \in \mathcal{F}_n} \big| \mathcal{C}_P(f,g) - \widehat{\mathcal{C}}_n(f, g) \big|}{ \mathrm{Var}(\|G_P\|_{\mathcal{F}_n} )} \right)^{1/3}  + \left( \frac{\sup_{f, g \in \mathcal{F}_n} \big| \widehat{\mathcal{C}}_n(f,g) - \widehat{\mathcal{C}}_n^m(f, g) \big|}{ \mathrm{Var}(\|G_P\|_{\mathcal{F}_n} )} \right)^{1/3},
	\end{align*}
	where $\lesssim$ hides an absolute constant independent of $n, m, M, \mathcal{F}_n, F_n, P_n$, and $P$.
	
	Moreover, if $\mathcal{F}_n$ is compact w.r.t. $d_{\widehat{\mathcal{C}}_n}$, then the last term on the right hand side in above display vanishes as $m \rightarrow \infty$.
\end{corollary}

In general, the strong variance $\mathrm{Var}(\|G_P\|_{\mathcal{F}_n} )$ depends on the dimension of the function class $\mathcal{F}_n$. Hence, the truncation level $m$ has to be chosen (inverse) proportionate to $\mathrm{Var}(\|G_P\|_{\mathcal{F}_n})$ to ensure that in the upper bound of Corollary~\ref{corollary:theorem:Abstract-Bootstrap-Sup-Empirical-Process-1} the deterministic approximation error $\sup_{f, g \in \mathcal{F}_n} \big| \widehat{\mathcal{C}}_n(f,g) - \widehat{\mathcal{C}}_n^m(f, g) \big|$ is negligible compared to the stochastic estimation errors.

We conclude this section with a consistency result on the bootstrap approximation of quantiles of suprema of empirical processes. This result is the empirical process analogue to Proposition~\ref{lemma:Bootstrap-Max-Norm-Asymptotic-Size}. It is relevant in the context of hypothesis testing and construction of confidence regions. For $\alpha \in (0,1)$ we denote the conditional $\alpha$-quantile of the supremum of $\{\widehat{Z}_n^m(f) : f \in \mathcal{F}_n\}$ by
\begin{align*}
	c_n(\alpha; \widehat{\mathcal{C}}_n^m) &:= \inf \left\{s \geq 0: \mathbb{P}\left\{  \|\widehat{Z}_n^m\|_{\mathcal{F}_n} \leq s \mid X_1, \ldots, X_n\right\} \geq \alpha \right\}.
\end{align*}

We have the following result:

\begin{theorem}[Quantiles of the Gaussian process bootstrap]\label{theorem:Abstract-Bootstrap-Sup-Empirical-Process-Quantiles}
	Consider the setup of Corollary~\ref{corollary:theorem:Abstract-Bootstrap-Sup-Empirical-Process-1}. Let $(\Theta_n)_{n \geq 1} \in \mathbb{R}$ be a sequence of arbitrary random variables, not necessarily independent of $X_1, \ldots, X_n$. Then, for each $M \geq 0$,
	\begin{align*}
		&\sup_{\alpha \in (0,1)} \Big|\mathbb{P}\left\{ \|\mathbb{G}_n\|_{\mathcal{F}_n} + \Theta_n \leq 	c_n(\alpha; \widehat{\mathcal{C}}_n^m) \right\}  - \alpha \Big| \\
		&\quad{}\quad{}\quad{} \lesssim \frac{\|F_n\|_{P, 3}}{\sqrt{n^{1/3}\mathrm{Var}(\|G_P\|_{\mathcal{F}_n})}} + \frac{\|F_n\mathbf{1}\{F_n > M\}\|_{P,3}^3}{\|F_n\|_{P,3}^3} + \frac{ \mathrm{E} \|G_P\|_{\mathcal{F}_n} + M}{\sqrt{n\mathrm{Var}(\|G_P\|_{\mathcal{F}_n})}}\\
		& \quad{}\quad{}\quad{}\quad{} \quad{} \quad{}  + \inf_{\delta > 0}\left\{ \left(\frac{\delta }{\mathrm{Var}(\|G_P\|_{\mathcal{F}_n})} \right)^{1/3}  + \mathbb{P}\left\{\sup_{f, g \in \mathcal{F}_n} \big| \mathcal{C}_P(f,g) - \widehat{\mathcal{C}}_n^m(f, g) \big| > \delta\right\}\right\}\\
		&\quad{}\quad{}\quad{}\quad{} \quad{} \quad{} \quad{}\quad{}\quad{}\quad{} + \inf_{\eta > 0} \left\{ \frac{\eta }{ \sqrt{\mathrm{Var}(\|G_P\|_{\mathcal{F}_n})} } + \mathbb{P}\left\{|\Theta_n| > \eta \right\} \right\},
	\end{align*}
	where $\lesssim$ hides an absolute constant independent of $n, m, M, \mathcal{F}_n, F_n, P_n$, and $P$. 
\end{theorem}
In statistical applications, the statistic of interest is rarely a simple empirical process. Instead, the empirical process usually arises as the leading term of a (functional) Taylor expansion. The random sequence $(\Theta_n)_{n\geq 1}$ in above theorem can be used to capture the higher-order approximation errors of such a expansion. 

\begin{remark}[Additional practical considerations]
	It is often infeasible to draw Monte Carlo samples directly from $\|Z_n^m\|_{\mathcal{F}_n}$. In practice, we suggest approximating $\mathcal{F}_n$ via a finite $\delta \|F\|_{P,2}$-net $\mathcal{G}_n \subseteq \mathcal{F}_n$ with respect to $d_P$. With this additional approximation step and by Corollary~\ref{corollary:theorem:GaussianComparison-PQ-Simultaneous-2}, the conclusion of the Corollary~\ref{corollary:theorem:Abstract-Bootstrap-Sup-Empirical-Process-1} holds with $\|\widehat{Z}_n^m\|_{\mathcal{G}_n}$ substituted for $\|\widehat{Z}_n^m\|_{\mathcal{F}_n}$ and the additional quantity $\big(\delta\|F\|_{P,2} \sup_{f \in \mathcal{F}_n} \sqrt{Pf^2}\big)^{1/3} \mathrm{Var}(\|G_P\|_{\mathcal{F}_n})^{-1/3}$ in the upper bound on the Kolmogorov distance. 
	This shows that the level of discretization $\delta > 0$ should be chosen proportional to $\mathrm{Var}(\|G_P\|_{\mathcal{F}_n})$.
\end{remark}

\subsection{Relation to previous work}\label{subsec:RelationPreviousWork-Stat}
The papers by~\cite{chernozhukov2014AntiConfidenceBands,chernozhukov2014GaussianApproxEmp, chernozhukov2016empirical} are currently the only other existing works on Gaussian and bootstrap approximations of suprema of empirical processes indexed by (potentially) non-Donsker function classes. In this section, we compare their results to our Theorems~\ref{theorem:CLT-Max-Norm-Simultaneous}--\ref{theorem:Abstract-Bootstrap-Sup-Empirical-Process-Quantiles}. To keep the comparison we focus on the key aspects that motivated us to write this paper.

It is important to note that the results presented in~\cite{chernozhukov2014AntiConfidenceBands,chernozhukov2014GaussianApproxEmp, chernozhukov2016empirical} differ slightly in nature from ours. Instead of establishing bounds on Kolmogorov distances, Chernozhukov and his co-authors derive coupling inequalities. However, through standard arguments (Strassen's theorem and anti-concentration) these coupling inequalities indirectly yield bounds on Kolmogorov distances. These implied bounds on the Kolmogorov distances are at most as sharp as the ones in the coupling inequalities, up to multiplicative constants. Since we do not care about absolute constants, but only the dependence of the upper bounds on characteristics of the function class $\mathcal{F}_n$ and the law $P$, we can meaningfully compare their findings and ours.

\begin{itemize}
	\item \emph{Unbounded function classes.}
	In classical empirical process theory the function class $\mathcal{F}_n$ is typically assumed to be uniformly bounded, i.e. $\sup_{x \in S} |f(x)| < \infty$ for all $f \in \mathcal{F}_n$~\citep[e.g.][]{vandervaart1996weak}.
	A key feature of Theorems~\ref{theorem:CLT-Max-Norm-Simultaneous}--\ref{theorem:Abstract-Bootstrap-Sup-Empirical-Process-Quantiles} as well as Theorems A.2, 2.1, and 2.1-2.3 in~\cite{chernozhukov2014AntiConfidenceBands, chernozhukov2014GaussianApproxEmp, chernozhukov2016empirical} is that they hold for unbounded function classes and only require the envelope function $F_n$ to have finite moments. Relaxing the uniform boundedness of the function class is useful, among other things, for inference on high-dimensional statistical models, functional data analysis, nonparametric regression, and series/ sieve estimation~\citep{chernozhukov2014AntiConfidenceBands, chernozhukov2014GaussianApproxEmp, giessing2023bootstrap}.
		
	\item \emph{Entropy conditions.} The upper bounds provided in Theorems A.2, 2.1, and 2.1-2.3 in~\cite{chernozhukov2014AntiConfidenceBands, chernozhukov2014GaussianApproxEmp, chernozhukov2016empirical} depend on a combination of (truncated) second and $q$th ($q \geq 3$) moments of $\max_{1 \leq i \leq n} F_n$ $(X_i)$, ``local quantities'' of order $\sup_{f \in \mathcal{F}_n} P|f|^3$ and $\sup_{f \in \mathcal{F}_n} \sqrt{Pf^4}$ (disregarding $(\log n)$-factors), $\mathrm{E} \|G_P\|_{\mathcal{F}_{n,\delta}'}  \vee  \big\|\|\mathbb{G}_n\|_{\mathcal{F}_{n,\delta}'} \big\|_{P,1}$, and the entropy number $\log N(\mathcal{F}, e_P,\\ \delta \|F_n\|_{P, 2})$, $\delta > 0$ arbitrary. These upper bounds are not only more complex than ours but also weaker in terms of their implications: 
	Since metric entropy with respect to intrinsic standard deviation metric $e_P$ typically scales linearly in the (VC-)dimension of the statistical model (see Appendix~\ref{sec:DerivationIntro} for two (counter-)examples), the upper bounds in~\cite{chernozhukov2014AntiConfidenceBands, chernozhukov2014GaussianApproxEmp, chernozhukov2016empirical} are vacuous in high-dimensional situations without sparsity or when the (VC-)dimension exceeds the sample size. In contrast, the upper bounds in our Theorems~\ref{theorem:CLT-Max-Norm-Simultaneous}--\ref{theorem:Abstract-Bootstrap-Sup-Empirical-Process-Quantiles} depend only on the expected value and standard deviation of the supremum of the Gaussian $P$-bridge process and the (truncated) third moments of the envelope. Under mild assumptions, these quantities can be upper bounded independently of the (VC-)dimension, thus offering useful bounds even in high-dimensional problems (see Section~\ref{sec:Applications}).
	
	The entropy term in~\cite{chernozhukov2014AntiConfidenceBands, chernozhukov2014GaussianApproxEmp, chernozhukov2016empirical} is due to their proofs relying on a dimension-dependent version of Proposition~\ref{lemma:CLT-Max-Norm}. 
	
	\item \emph{Lower bounds on weak variances.} Lemmas 2.3 and 2.4 in~\cite{chernozhukov2014GaussianApproxEmp} and all of the results in~\cite{chernozhukov2014AntiConfidenceBands, chernozhukov2016empirical} require a strictly positive lower bound on the weak variance of the Gaussian $P$-bridge process, i.e. $\inf_{f \in \mathcal{F}_n}\mathrm{Var}\big(G_P(f)\big) \geq \underline{\sigma}^2 > 0$. This assumption automatically limits the applicability of these lemmas to studentized statistics/ standardized function classes~\citep{chernozhukov2014AntiConfidenceBands, chernozhukov2014GaussianApproxEmp} and excludes relevant scenarios with variance decay~\citep{lopes2020bootstrapping, lopes2022improved, lopes2023bootstrapping}. In contrast, our Theorems~\ref{theorem:CLT-Max-Norm-Simultaneous}--\ref{theorem:Abstract-Bootstrap-Sup-Empirical-Process-Quantiles} apply to all function classes for which the strong variance of the Gaussian $P$-bridge process, $\mathrm{Var}(\|G_P\|_{\mathcal{F}_n})$, does not vanish ``too fast''. A slowly vanishing strong variance is a weaker requirement than a strictly positive lower bound on the weak variance (see~\cite{giessing2022anticoncentration} and also Lemmas~\ref{lemma:AntiConcentration-SeparableProcess} and~\ref{lemma:BoundsVariance-SeparableProcess}).
	
	The lower bound on the weak variance is an artifact of~\possessivecites{chernozhukov2014AntiConfidenceBands}{chernozhukov2014AntiConfidenceBands, chernozhukov2014GaussianApproxEmp, chernozhukov2016empirical} proof technique which requires the joint distribution of the finite-dimensional marginals of the Gaussian $P$-bridge process to be non-degenerate. 
\end{itemize}

\section{Applications}\label{sec:Applications}

\subsection{Confidence ellipsoids for high-dimensional parameter vectors}\label{subsec:Application-Vector}
Confidence regions are fundamental to uncertainty quantification in multivariate statistics. Typically, their asymptotic validity relies on multivariate CLTs. However, in high dimensions, when the number of parameters exceeds the sample size, validity of confidence regions needs to be justified differently. In this section, we show how the Gaussian process bootstrap offers a practical solution to this problem. Existing work on bootstrap confidence regions in high dimensions has focused exclusively on conservative, rectangular confidence regions~\citep[e.g.][]{chernozhukov2013GaussianApproxVec,chernozhukov2014GaussianApproxEmp, chernozhukov2023high-dimensional}. For the first time, our results allow construction of tighter, elliptical confidence regions and without sparsity assumptions.

Let $\theta_0 \in \mathbb{R}^d$ be an unknown parameter and $\hat{\theta}_n \equiv \hat{\theta}_n(X_1, \ldots, X_n)$ an estimator for $\theta_0$ based on the simple random sample $X_1, \ldots, X_n$. Then, an asymptotic $1-\alpha$ confidence ellipsoid for $\theta_0$ can be constructed as
\begin{align}\label{eq:Application-Vector-0}
	\mathcal{E}_n(c) :=  \left\{ \theta \in \mathbb{R}^d : \sqrt{n}\| \hat{\theta}_n - \theta\|_2 \leq c \right\},
\end{align}
where $\|\cdot\|_2$ denotes the Euclidean norm and $c > 0$ solves
\begin{align*}
	\lim_{n \rightarrow \infty} \mathbb{P} \left\{ \sqrt{n} \| \hat{\theta}_n - \theta_0 \|_2 \leq c \right\} = 1 - \alpha.
\end{align*}
In general, it is difficult to determine $c$ from above equation because the sampling distribution of $\sqrt{n} \| \hat{\theta}_n - \theta_0 \|_2 $ often does not have an explicit form. However, in many cases, $\hat{\theta}_n$ admits an asymptotically linear expansion of the form
\begin{align*} 
	\hat{\theta}_n - \theta_0 = \frac{1}{n} \sum_{i=1}^n \psi_i+ R_n,
\end{align*}
where $\psi_1, \ldots, \psi_n$ are centered i.i.d. random vectors (influence functions) and $R_n$ is a remainder term.
Thus, in these cases, we have
\begin{align*}
	\sqrt{n} \|\hat{\theta}_n - \theta_0\|_2 = \sup_{u \in S^{d-1}}\left| \frac{1}{\sqrt{n}} \sum_{i=1}^n \psi_i'u \right|  + \Theta_n, \quad \quad \text{where} \quad \quad  |\Theta_n| \leq \sqrt{n} \|R_n\|_2.
  \end{align*}
Using this formulation, we can now apply Theorem~\ref{theorem:Abstract-Bootstrap-Sup-Empirical-Process-Quantiles} to approximate the quantile $c$ and complete the construction of the asymptotic confidence ellipsoid. By Algorithm~\ref{algorithm:GaussianProcessBootstrap} we need to construct a Gaussian process with index set $S^{d-1}$ and whose covariance function approximates the bi-linear map $(u, v) \mapsto u'\Omega_\psi v$, where $\Omega_\psi= \mathrm{E}[\psi_1\psi_1']$. Clearly, the following process does the job:
\begin{align*}
	\{ \widehat{Z}_\psi'u : u \in S^{d-1}\}, \quad \quad \text{where} \quad \quad \widehat{Z}_\psi \mid X_1, \ldots, X_n \sim N(0, \widehat{\Omega}_\psi ),
\end{align*}
and $\widehat{\Omega}_\psi = n^{-1} \sum_{i=1}^n (\psi_i - \overline{\psi}_n)(\psi_i - \overline{\psi}_n)'$ and $\overline{\psi}_n = n^{-1} \sum_{i=1}^n \psi_i$. We denote the $\alpha$-quantile of the supremum of this process by
\begin{align}\label{eq:Application-Vector-2}
	c_n(\alpha) :=  \inf\left\{ c \geq 0 :\mathbb{P} \left\{\|\widehat{Z}_\psi\|_2 \leq c\mid X_1, \ldots, X_n \right\} \geq \alpha \right\}.
\end{align}
To show that these quantiles uniformly approximate the quantiles of the (asymptotic) distribution of $\sqrt{n} \|\hat{\theta}_n - \theta_0\|_2$ introduce the following assumptions:

\begin{definition}[Sub-Gaussian random vector]\label{defintion:SubGaussian-RV}
	We call a centered random vector $X \in \mathbb{R}^d$ sub-Gaussian if $\|X'u\|_{\psi_2}^2 \lesssim \mathrm{E}[(X'u)^2]$ for all $u \in \mathbb{R}^d$.
\end{definition}

\begin{assumption}[Sub-Gaussian influence functions]\label{assumption:SubGaussianInfluenceFunction}
	The influence functions $\psi, \psi_1, \ldots, \psi_n \in \mathbb{R}^d$ are i.i.d. sub-Gaussian random vectors with covariance matrix $\Omega_\psi$ and (i) $ r(\Omega_\psi)/r(\Omega_\psi^2) = O(1)$, (ii) $r(\Omega_\psi) = o(n^{1/6})$, and (iii) $\sqrt{n}\|R_n\|_2  = o_p(1)$, where $r(A) = \mathrm{tr}(A)/\|A\|_{op}$ is the effective rank matrix $A$.
\end{assumption}

\begin{assumption}[Heavy-tailed influence functions]\label{assumption:HeavyTailedInfluenceFunction}
	The influence functions $\psi, \psi_1, \ldots, \psi_n\in \mathbb{R}^d$ are centered i.i.d. random vectors with covariance matrix $\Omega_\psi$ such that $M_n^2:= \mathrm{E}[\max_{1 \leq i \leq n}\|\psi_i\|_2^2] / \|\Omega_\psi\|_{op}$ $< \infty$ and $m_{2,s}^s := \mathrm{E}[\|\psi_1\|_2^s] / \|\Omega_\psi\|_{op}^{s/2} < \infty$ for some $s > 3$. Furthermore, (i) $r(\Omega_\psi)/r(\Omega_\psi^2) = O(1)$, (ii) $m_{2, s}/m_{2,3} = o(n^{1/3-1/s})$, (iii) $m_{2, 3} = o(n^{1/6})$, (iv) $ M_n^2 (\log n)  = o(n)$, and (v) $\sqrt{n}\|R_n\|_2 = o_p(1)$,
	where $r(A) = \mathrm{tr}(A)/\|A\|_{op}$ is the effective rank matrix $A$.
\end{assumption}

Under either assumption, we have the following result:

\begin{proposition}[Bootstrap confidence ellipsoid]\label{theorem:ConfidenceEllipsoid}
	If either Assumption~\ref{assumption:SubGaussianInfluenceFunction} or Assumption~\ref{assumption:HeavyTailedInfluenceFunction} holds, then
	\begin{align*}
		\lim_{n \rightarrow \infty} \sup_{\alpha \in (0,1)} \left| \mathbb{P}\Big\{ \theta_0 \in \mathcal{E}_n\big(c_n(1-\alpha)\big) \Big\} - (1- \alpha)\right| = 0,
	\end{align*}
	with the confidence ellipsoid $\mathcal{E}_n$ as defined in~\eqref{eq:Application-Vector-0} and the conditional quantile $c_n(1-\alpha)$ as in~\eqref{eq:Application-Vector-2}.
\end{proposition}

It is worth noticing that neither Assumption~\ref{assumption:SubGaussianInfluenceFunction} nor~\ref{assumption:HeavyTailedInfluenceFunction} require sparsity of the parameter $\theta_0$, the estimator $\hat{\theta}_n$, or the influence functions $\psi_1, \ldots, \psi_n$. Instead, Assumption~\ref{assumption:SubGaussianInfluenceFunction} and~\ref{assumption:HeavyTailedInfluenceFunction} already hold if the covariance matrix $\Omega_\psi$ has bounded effective rank~\citep[see][Section 2.1 and Appendices A.1 and A.2]{giessing2023bootstrap}. In this context it is important to notice that, unlike~\cite{chernozhukov2014GaussianApproxEmp}, we do not require a strictly positive lower bound on $\inf_{\|u\|_2=1} \mathrm{E}[(Z_\psi'u)^2] \equiv \lambda_{\min}(\Omega_\psi)$, the smallest eigenvalue of the covariance matrix $\Omega_\psi$. If such a lower bound was needed, the effective rank $r(\Omega_\psi)$ would grow linearly in the dimension $d$ and Assumptions~\ref{assumption:SubGaussianInfluenceFunction} and~\ref{assumption:HeavyTailedInfluenceFunction} would be violated if $d \gg n$ (see Appendix~\ref{sec:DerivationIntro}). Moreover, Proposition~\ref{theorem:ConfidenceEllipsoid} cannot be deduced from Theorem 2.1 in~\cite{chernozhukov2014GaussianApproxEmp} because the upper bound in their coupling inequality would feature the term $\sqrt{d \log (1 + 1/\epsilon)}/n^{1/4}$ (and other terms as well), which is a remnant of the entropy number of the $\epsilon$-net discretization of the $d$-dimensional Euclidean unit-ball.

\subsection{Inference on the spectral norm of high-dimensional covariance matrices}\label{subsec:Application-Matrix}
Spectral statistics of random matrices play an important role in multivariate statistical analysis. The asymptotic distributions of spectral statistics are well established in low dimensions~\citep[e.g.][]{anderson1963asymptotic, waternaux1976asymptotic, fujikoshi1980asymptotic} and when the dimension is comparable to the sample size~\citep[e.g.][]{johnstone2001on, elKaroui2007tracy-widom, peche2009universality, bao2015universality}. 
In the high-dimensional case, when asymptotic arguments do not apply, bootstrap procedures have proved to be effective in approximating distributions of certain maximum-type spectral statistics~\citep[e.g.][]{han2018on, naumov2019bootstrap, lopes2019BootstrappingSpectral, lopes2020bootstrapping, lopes2022improved}. Here, we demonstrate that the Gaussian process bootstrap is a viable alternative to these bootstrap procedures in approximating the distribution of the spectral norm of a high-dimensional sample covariance matrix. 

Let $X_1, \ldots, X_n \in \mathbb{R}^d$ be i.i.d. random vectors with law $P$, mean zero, and covariance matrix $\Sigma \in \mathbb{R}^{d \times d}$. Consider the spectral statistic
\begin{align*}
	T_n := \sqrt{n} \| \widehat{\Sigma}_n - \Sigma \|_{op}, \quad\quad{}  \widehat{\Sigma}_n = n^{-1} \sum_{i=1}^n X_iX_i'.
\end{align*}
Since there does not exist a closed form expression of the sampling distribution of $T_n$ when $d \gg n$, we apply Algorithm~\ref{algorithm:GaussianProcessBootstrap} to obtain an approximation based on a Gaussian proxy process. We introduce the following notation: Let $x, u, v \in \mathbb{R}^d$ and note that $(x'u)(x'v) = \mathrm{vech}'(xx') H_d'(v \otimes u)$, where $\otimes$ denotes the Kronecker product, $\mathrm{vech}(\cdot)$ the half-vectorization operator which turns a symmetric matrix $A \in \mathbb{R}^{d \times d}$ into a $d(d+1)/2$ column vector (of unique entries), and $H_d $ the duplication matrix such that $H_d\: \mathrm{vech}(A)=\mathrm{vec}(A)$ where $\mathrm{vec}(\cdot)$ is the ordinary vectorization operator. Whence, $T_n \equiv \sqrt{n}\|Q_n - Q\|_{{\mathcal{F}_n}}$ where $Q_n$  is the empirical measure of the collection of $Y_i = \mathrm{vech}'(X_iX_i') H_d'$, $1 \leq i \leq n$, $Q$ the pushforward of $P$ under the map $X \mapsto Y = \mathrm{vech}'(XX') H_d'$, and $\mathcal{F}_n =\{ y \mapsto  y(v \otimes u) : u, v \in S^{d-1}\}$. Since each $f \in \mathcal{F}_n$ has a (not necessarily unique) representation in terms of $u, v \in S^{d-1}$, in the following we identify $f \in \mathcal{F}_n$ with pairs $(u,v) \in S^{d-1} \times S^{d-1}$. The covariance function associated with the empirical process $\{\sqrt{n}(Q_n - Q)(f): f \in \mathcal{F}_n\}$ is thus given by
\begin{align*}
	\big((u_1, v_1), (u_2, v_2)\big) &\mapsto 
						     (v_1 \otimes u_1)' H_d \Omega H_d'(v_2 \otimes u_2),
\end{align*}
where $\Omega = \mathrm{E}\left[ \mathrm{vech}(XX'-\Sigma) \mathrm{vech}'(XX'-\Sigma) \right] \in  \mathbb{R}^{d(d+1)/2 \times d(2+1)/2}$. Thus, in the light of Algorithm~\ref{algorithm:GaussianProcessBootstrap} the natural choice for the Gaussian bootstrap process is
\begin{align}\label{eq:subsec:Application-Matrix-2}
	\big\{\widehat{Z}_n'H_d'(v \otimes u) : u, v \in S^{d-1}\big\}, \quad \quad \text{where} \quad \quad \widehat{Z}_n \mid X_1, \ldots, X_n \sim N(0, \widehat{\Omega}_n),
\end{align}
and $\widehat{\Omega}_n = n^{-1} \sum_{i=1}^n \mathrm{vech}(X_iX_i'-\widehat{\Sigma}_n) \mathrm{vech}'(X_iX_i'-\widehat{\Sigma}_n)$ is the sample analogue of $\Omega$.

We make the following assumption:

\begin{assumption}[Sub-Gaussian data]\label{assumption:SubGaussianData}
	The data $X, X_1, \ldots, X_n \in \mathbb{R}^d$ are i.i.d. sub-Gaussian random vectors with covariance matrix $\Sigma$ and $\inf_{ \|u\|_2 = 1} \mathrm{Var}((X'u)^2) \geq \kappa > 0$.
\end{assumption}

\begin{remark}[On the lower bound on the variances]
	The strictly positive lower bound $\kappa > 0$ on the variance of the quadratic from $u'X'Xu$ is mild. The existence of the lower bound $\kappa > 0$ is equivalent to $\mathrm{E}[(X'u)^4]^{1/4} > \mathrm{E}[(X'u)^2]^{1/2}$ for all $u \in S^{d-1}$. The latter inequality holds if the law of $X$ does not concentrate on a lower dimensional subspace of $\mathbb{R}^d$. Since the bounds in Proposition~\ref{theorem:ConsistencyOperatorNorm} are explicit in $\kappa > 0$, Proposition~\ref{theorem:ConsistencyOperatorNorm} also applies to scenarios in which $\kappa \equiv \kappa(n, d) \rightarrow 0$ as $n,d \rightarrow \infty$. Similar lower bounds on the variance of $u'X'Xu$ appear in~\cite{lopes2022improved, lopes2023bootstrapping}.
\end{remark}

We have the following result:

\begin{proposition}[Bootstrap approximation of the distribution of spectral norms of covariance matrices]\label{theorem:ConsistencyOperatorNorm}
	 Suppose that Assumption~\ref{assumption:SubGaussianData} holds. Let $\widehat{S}_n = \mathrm{vec}^{-1}(H_d\widehat{Z}_n) \in \mathbb{R}^{d \times d}$ with $\widehat{Z}_n$ as given in~\eqref{eq:subsec:Application-Matrix-2}. Then,
	\begin{align*}
		&\sup_{s \geq 0 } \Big|\mathbb{P} \left\{ \sqrt{n} \| \widehat{\Sigma}_n - \Sigma \|_{op}  \leq s\right\} - \mathbb{P} \left\{\|\widehat{S}_n\|_{op} \leq s \mid X_1, \ldots, X_n \right\}  \Big| \\
		&\quad\quad \quad\lesssim  O_p\left(\left( \frac{\|\Sigma\|_{op}^{2/3}}{\kappa^{2/3}} \vee \frac{\|\Sigma\|_{op}^{4/3}}{\kappa^{4/3}}  \right) \left( \sqrt{\frac{(\log en)^2  \mathrm{r}^7(\Sigma)}{n}} \vee \frac{(\log en)^2  \mathrm{r}^4(\Sigma)}{n} \right)^{1/3} \right)\\
		&\quad{}\quad{}\quad{}\quad{} \quad{} \quad{} + \left( \frac{\|\Sigma\|_{op}}{\kappa }  \vee \frac{\|\Sigma\|_{op}^2}{\kappa^2 } \right) \frac{r^2(\Sigma)}{ n^{1/6} },
	\end{align*}
	where $r(\Sigma) = \mathrm{tr}(\Sigma)/\|\Sigma\|_{op}$ is the effective rank of $\Sigma$ and $\lesssim$ hides an absolute constant independent of $n, d, \Sigma$, and $\kappa$.
\end{proposition}
\begin{remark}
	Note that the matrix $\widehat{S}_n$ is symmetric just as the target matrix $ \widehat{\Sigma}_n - \Sigma$. 
\end{remark}
We conclude this section with a comparison of Proposition~\ref{theorem:ConsistencyOperatorNorm} to existing results in the literature: First, if $(\|\Sigma\|_{op}/\kappa  \vee \|\Sigma\|_{op}^2 /\kappa^2 ) r^2(\Sigma) = o(n^{1/6})$, then the upper bound in above theorem is asymptotically negligible. In this case, the bootstrapped distribution of the Gaussian proxy statistic $\|\widehat{S}_n\|_{op}$ consistently approximates the distribution of $\| \widehat{\Sigma}_n - \Sigma \|_{op}$. Since this rate depends only on the effective rank, it is dimension-free and cannot be derived through the results in~\cite{chernozhukov2014GaussianApproxEmp}.
Second, unlike the results in~\cite{lopes2022improved} and~\cite{lopes2023bootstrapping}, Proposition~\ref{theorem:ConsistencyOperatorNorm} does not rely on specific assumptions about the decay of the eigenvalues of $\Sigma$. And yet, under certain circumstances, the consistency rate provided by Proposition~\ref{theorem:ConsistencyOperatorNorm} can be faster than theirs. Specifically, the bootstrap procedure described in~\cite{lopes2022improved} achieves consistency at a rate of $n^{-\frac{\beta - 1/2}{2\beta + 4 + \epsilon}}$ for $\epsilon > 0$ and $\beta > 1/2$. In our context, the parameter $\beta$ determines the rate at which the eigenvalues of $\Sigma$ decrease, i.e.  $\lambda_k(\Sigma) \asymp k^{-2\beta}$, $k = 1, \ldots, d$. To achieve a rate faster than $n^{-1/6}$, $\beta$ must be greater than $(7 + \epsilon)/4$ which requires an extremely fast decay of the eigenvalues.
Third,~\cite{lopes2023bootstrapping} conduct extensive numerical experiments and observe that the bootstrap approximation exhibits a sharp phase transition from accurate to inaccurate when $\beta$ switches from greater than $1/2$ to less than $1/2$. This observation aligns not only with their own theoretical findings but also with the upper bound presented in Proposition~\ref{theorem:ConsistencyOperatorNorm}, since, under their modeling assumptions, the effective rank $r(\Sigma)$ remains bounded if $\beta > 1/2$ but diverges if $\beta \leq 1/2$.

\subsection{Simultaneous confidence bands for functions in reproducing kernel Hilbert spaces}\label{subsec:Application-FDA}
Reproducing kernel Hilbert spaces (RKHS) are an integral part of statistics, with applications in classical non-parametric statistics~\citep{wahba1990spline}, machine learning~\citep{schoelkopf2002learning,steinwart2008support} and, most recently, (deep) neural nets~\citep{belkin2018understand, jacot2018neural, bohn2019representer,unser2019representer, chen2020deep}. In this section, we consider constructing simultaneous confidence bands for functions in RKHS by bootstrapping the distribution of a bias-corrected kernel ridge regression estimator. Recently, this problem has been addressed by~\cite{singh2023kernel} with a symmetrizied multiplier bootstrap. Here, we propose an alternative based on the Gaussian process bootstrap using a truncated Karhunen-Lo{\`e}ve decomposition. We point out several commonalities and differences between the two procedures.

In the following, $\mathcal{H}$ denotes the RKHS of continuous functions from $S$ to $\mathbb{R}$ associated to the symmetric and positive definite kernel $k: S \times S \rightarrow \mathbb{R}$. We allow $\mathcal{H}$, $S$, and $k$ to change with the sample size $n$, but we do not make this dependence explicit. We write $\| \cdot\|_{\mathcal{H}}$ for the norm induced by the inner product $\langle \cdot,\cdot \rangle_\mathcal{H}$ and $\| \cdot\|_\infty$ for the supremum norm. For $x \in S$, we let $k_x: S \rightarrow \mathbb{R}$ be the function $y \mapsto k(x, y)$. Then, $k_x \in \mathcal{H}$ for all $x \in S$ and, by the reproducing property,  $f(x) = \langle f, k_x\rangle_{\mathcal{H}}$ for all $f \in \mathcal{H}$ and $x \in S$. The kernel induces the so-called kernel metric $d_k(x,y) := \|k_x - k_y\|_{\mathcal{H}}$, $x, y \in S$.
Given $f \in \mathcal{H}$ ($z \in S$) we denote its dual by $f^* \in \mathcal{H}^*$ ($z^* \in S^* \subset \mathcal{H}$). For $f,g \in \mathcal{H}$ we define the tensor product $f \otimes g^* : \mathcal{H} \rightarrow \mathcal{H}$ by $h \mapsto (f \otimes g^*)(h)  := \langle g, h\rangle_{\mathcal{H}} f$. For operators on $\mathcal{H}$ we use $\|\cdot\|_{op}$, $\|\cdot\|_{HS}$, and $\mathrm{tr}(\cdot)$ to denote operator norm, Hilbert-Schmidt norm, and trace, respectively. For further details on RKHSs we refer to~\cite{berlinet2011reproducing}. 

Let $Y \in \mathbb{R}$ and $X \in S$ have joint law $P$. Given a simple random sample $(Y_1, X_1), \ldots, (Y_n,$ $X_n)$ our goal is to construct uniform confidence bands for the conditional mean function $x \mapsto \mathrm{E}[Y \mid X = x]$ or, rather, its best approximation in hypothesis space $\mathcal{H}$, i.e.
\begin{align*}
	f_0 \in \argmin_{f \in \mathcal{H}} \mathrm{E}\left[\big(Y - f(X)\big)^2 \right].
\end{align*}
To this end, consider the classical kernel ridge regression estimator
\begin{align*}
	\widehat{f}_n \in \argmin_{f \in \mathcal{H}}  \left\{ \frac{1}{n}\sum_{i=1}^n \big(Y_i - f(X_i)\big)^2 + \lambda\|f\|_{\mathcal{H}}^2 \right\}, \quad  \lambda > 0,
\end{align*}
and define its bias-corrected version as
\begin{align}\label{eq:subsec:Application-FDA-1}
	\widehat{f}^{\mathrm{bc}}_n : = \widehat{f}_n  + \lambda(\widehat{T}_n + \lambda)^{-1} \widehat{f}_n ,  \quad \quad \widehat{T}_n =  \frac{1}{n} \sum_{i=1}^n \left(k_{X_i} \otimes k_{X_i}^* \right).
\end{align}
We propose to construct simultaneous $1-\alpha$ confidence bands for $f_0$ based on $\widehat{f}^{\mathrm{bc}}_n$ via the rectangle
\begin{align}\label{eq:subsec:Application-FDA-0}
	\mathcal{R}_n(c) :=  \left\{ f \in \mathcal{F} : \sqrt{n}\| \widehat{f}^{\mathrm{bc}}_n - f\|_\infty \leq c \right\}, 
\end{align}
where $c > 0$ approximates the (asymptotic) $1-\alpha$ quantile of the law of $\sqrt{n}\| \widehat{f}^{\mathrm{bc}}_n - f_0\|_\infty$. To compute $c> 0$ we proceed in two steps: First, we show that $\sqrt{n}\| \widehat{f}^{\mathrm{bc}}_n - f_0\|_\infty$ can be written as the sum of the supremum of an empirical process and a negligible remainder term. Then, we apply the strategy developed in Section~\ref{subsec:GaussianProcessBootstrap-Implementation-Consistency} to bootstrap the supremum of the empirical process.

By Lemma~\ref{lemma:BC-KRR-Expansion} in Appendix~\ref{subsec:AuxiliaryResults-Applications},
\begin{align*} 
\widehat{f}^{\mathrm{bc}}_n - f_0 = (T + \lambda)^{-2} T \left( \frac{1}{n} \sum_{i=1}^n \big(Y_i - f_0(X_i)\big) k_{X_i} \right) + R_n,
\end{align*}
where $T =  \mathrm{E}[k_X \otimes k_X^*]$ and $R_n$ is a higher-order remainder term. 
Since $k$ is a reproducing kernel, we have $k_x(z) = \langle k_x, k_z\rangle_\mathcal{H} = \langle k_x, z^* \rangle_\mathcal{H}$ for all $x, z \in S$. Hence, above expansion implies 
\begin{align*}
	\sqrt{n} \|\widehat{f}^{\mathrm{bc}}_n - f_0\|_{\infty} = \sup_{u \in S^*} \left| \left \langle \frac{1}{\sqrt{n}} \sum_{i=1}^n V_i, u \right\rangle_\mathcal{H}  \right| + \Theta_n,\quad \quad \text{where} \quad \quad  |\Theta_n| \leq \sqrt{n} \|R_n\|_\infty,
\end{align*}
and $V_i = (T + \lambda)^{-2} T \big(Y_i - f_0(X_i)\big) k_{X_i}$, $ 1 \leq i \leq n$. 
By Lemma~\ref{lemma:BC-KRR-Remainder} in Appendix~\ref{subsec:AuxiliaryResults-Applications}, $\sqrt{n}\|R_n\|_\infty$ is negligible with high probability. Moreover, since $f_0$ is the best approximation in square loss, the random elements $V_i$'s have mean zero. Thus, $\sup_{u \in S^*} | \langle n^{-1/2} \sum_{i=1}^n V_i, u\rangle_\mathcal{H}| \equiv \sqrt{n}\|Q_n - Q\|_{\mathcal{F}_n}$, where $Q_n$ is the empirical measure of the $V_i$'s, $Q$ the pushforward of $P$ under the map $(Y, X) \mapsto V =  (T + \lambda)^{-2} T \big(Y - f_0(X)\big) k_{X}$, and $\mathcal{F}_n = \{ v \mapsto \langle  v, u\rangle_\mathcal{H} : u \in S^* \}$. Since the functions $ f \in \mathcal{F}_n$ are just the evaluation functionals of $u \in S^*$, in the following we identify $f \in \mathcal{F}_n$ with its corresponding $u \in S^*$. The covariance function $\mathcal{C}: S^* \times S^* \rightarrow \mathbb{R}$ associated with the empirical process $\{ \sqrt{n}(Q_n - Q)(f) : f \in \mathcal{F}_n \}$ is thus given by
\begin{align*}
	(u_1, u_2) \mapsto \mathrm{E}[\langle V, u_1\rangle_{\mathcal{H}} \langle V, u_2\rangle_{\mathcal{H}}] = \mathrm{E}\big[\big\langle (V \otimes V^*)u_1, u_2\big\rangle_{\mathcal{H}}\big] 
	= \big\langle \Omega u_1, u_2\big\rangle_{\mathcal{H}},
\end{align*}
with covariance operator
\begin{align*}
	\Omega :=   \sigma_0^2 T (T + \lambda)^{-2} T (T + \lambda)^{-2} T, \quad \quad \text{where} \quad \quad \sigma_0^2 := \mathrm{E}\big[\big(Y - f_0(X)\big)^2\big], 
\end{align*}
where we have used that the operators $T$ and $(T + \lambda)^{-1}$ commute.

We proceed to construct a Gaussian proxy process as outlined in Section~\ref{subsec:GaussianProcessBootstrap-Implementation-Consistency}: Let $\widehat{\Omega}_n = \widehat{\sigma}_n^2\widehat{T}_n(\widehat{T}_n + \lambda)^{-2}\widehat{T}_n (\widehat{T}_n + \lambda)^{-2}\widehat{T}_n$ and $\widehat{\sigma}_n^2 = n^{-1} \sum_{i=1}^n\big(Y_i - \widehat{f}_n(X_i)\big)^2$ be the plug-in estimates of the covariance operator $\Omega$ and the variance $\sigma_0^2$, respectively. Define $\widehat{\mathcal{C}}_n : S^* \times S^* \rightarrow \mathbb{R}$ by $(u_1, u_2) \mapsto \big\langle \widehat{\Omega}_n u_1, u_2\big\rangle_{\mathcal{H}}$. Recall definition~\eqref{eq:subsec:GaussianProcessBootstrap-Implementation-Consistency-1} of the integral operator $T_K$. In the present setup, $T_{\widehat{\mathcal{C}}_n} = \widehat{\Omega}_n$ by Fubini's theorem. 
Denote by $\big\{(\widehat{\lambda}_k, \widehat{\varphi}_k)\big\}_{k=1}^\infty$ the eigenvalue and eigenfunction pairs of $\widehat{\Omega}_n$. Further, let $\{\xi_k\}_{k=1}^\infty$ be a sequence of i.i.d. standard normal random variables. Then, for $m \geq 1$ and $u, u_1, u_2 \in S^*$ define
\begin{align}\label{eq:subsec:Application-FDA-2}
	&\widehat{Z}_n^m(u) := \sum_{k=1}^m \xi_k \sqrt{\widehat{\lambda}_k} \widehat{\varphi}_k(u), \quad \quad  \widehat{\mathcal{C}}_n^m(u_1, u_2) := \sum_{k=1}^m \widehat{\lambda}_k \widehat{\varphi}_k(u_1)\widehat{\varphi}_k(u_2) = \big\langle \widehat{\Omega}_n^m u_1, u_2  \big\rangle_\mathcal{H}.
\end{align}
where $\widehat{\Omega}_n^m$ is the best rank-$m$ approximation of $\widehat{\Omega}_n$.

Given Proposition~\ref{lemma:KarhunenLoeve-GP} we postulate that the process $\widehat{Z}_n^\infty$ is an almost sure version of a Gaussian process on $S^*$ with covariance function $\widehat{\mathcal{C}}_n$ (or, equivalently, with covariance operator $\widehat{\Omega}_n$. Consequently, Theorem~\ref{theorem:Abstract-Bootstrap-Sup-Empirical-Process-Quantiles} guarantees validity of the Gaussian process bootstrap based on $\widehat{Z}_n^m$. To make all these claims rigorous, consider the following assumptions:

\begin{assumption}[On the kernel]\label{assumption:kernel}
	The kernel $k : S \times S \rightarrow \mathbb{R}$ is symmetric, positive semi-definite, continuous, and bounded, i.e. $\sup_{x \in S} \sqrt{|k(x,x)|} =: \kappa < \infty$.
\end{assumption}
\begin{remark}
	The assumptions on the kernel are standard and important~\citep[][]{berlinet2011reproducing}. 
	The continuity of $k$ guarantees that $k_X$ is a random element on $\mathcal{H}$ whenever $X \in S$ is a random variable. It also implies that the RKHS $\mathcal{H}$ is separable whenever $S$ is separable. The boundedness and the reproducing property of $k$ imply that $\|f\|_\infty \leq  \kappa \|f\|_{\mathcal{H}}$ for all $f \in \mathcal{H}$. 
\end{remark}

\begin{assumption}[On the data]\label{assumption:data}
	The data $(Y_1, X_1), \ldots, (Y_n, X_n) \in \mathbb{R} \times S$ are i.i.d. random elements defined on an abstract product probability space $(\Omega, \mathcal{A}, \mathbb{P})$. The $Y_i$'s are almost surely bounded, i.e. there exists an absolute constant $B > 0$ such that $\max_{1 \leq i \leq n}|Y_i| \leq B$ almost surely.
\end{assumption}
\begin{remark} 
	The almost sure boundedness of the $Y_i$'s is a strong assumption. We introduce this assumption to keep technical arguments at a minimum.~\cite{singh2023kernel} impose an equivalent boundedness condition on the pseudo-residuals $\varepsilon_i = Y_i - f_0(X_i)$, $1 \leq i \leq n$. 
\end{remark}

\begin{assumption}[On the population and sample covariance operators]\label{assumption:cov-operator}
	For all $n \geq 1$ (i) there exists $\omega_S > 0$ such that $\inf_{u \in S^*} \langle \Omega u, u \rangle_{\mathcal{H}} \geq \omega_S$ and (ii) $\mathrm{tr}(\Omega) \vee \mathrm{tr}(\widehat{\Omega}_n) < \infty$.
\end{assumption}
\begin{remark}
	Condition (i) is the Hilbert space equivalent to the lower bound on the variance in Assumption~\ref{assumption:SubGaussianData}. While~\cite{singh2023kernel} do not explicitly impose a lower bound on the covariance operator, such a lower bound is implied by their Assumption 5.2~\cite[see][for details]{giessing2022anticoncentration}. In eq.~\eqref{eq:theorem:Uniform-CI-Bands-RKHS-8} in Appendix~\ref{sec:Applications} we provide a general non-asymptotic complement to below Proposition~\ref{theorem:Uniform-CI-Bands-RKHS} which applies even if $\omega_S \equiv \omega_S(n) \rightarrow 0$ as $n \rightarrow \infty$. Condition (ii) is a classical assumption in learning theory on RKHS~\citep{mendelson2002geometric}. Together with Conditions (i) it implies that $\Omega$ and $\widehat{\Omega}_n$ are finite rank and trace class operators, respectively.
\end{remark}

Denote the $\alpha$-quantile of the supremum of the Gaussian proxy process $\widehat{Z}_n^m$ in~\eqref{eq:subsec:Application-FDA-2} by
\begin{align}\label{eq:subsec:Application-FDA-3}
	c_n^m(\alpha) &:= \inf \left\{s \geq 0: \mathbb{P}\left\{ \sup_{u \in S^*} \left| \big \langle\widehat{Z}_n^m, u \big\rangle_\mathcal{H}  \right| \leq s \mid (Y_1, X_1), \ldots, (Y_n, X_n) \right\} \geq \alpha \right\}.
\end{align}

An application of Theorem~\ref{theorem:Abstract-Bootstrap-Sup-Empirical-Process-Quantiles} yields:

\begin{proposition}[Bootstrap quantiles]\label{theorem:Uniform-CI-Bands-RKHS}
	Let $(S, d_k)$ be a compact metric space such that $\int_0^\infty \sqrt{N(S, d_k, \varepsilon)} \\ d\varepsilon < \infty$. If Assumptions~\ref{assumption:kernel},~\ref{assumption:data}, and~\ref{assumption:cov-operator}, and the rates in eq.~\eqref{eq:theorem:Uniform-CI-Bands-RKHS-7} in Appendix~\ref{sec:Proofs-Applications} hold, then
	\begin{align*}
		& \sup_{\alpha \in (0,1)} \left| \mathbb{P}\left\{ 	\sqrt{n} \|\widehat{f}^{\mathrm{bc}}_n - f_0\|_{\infty} \leq c_n^m(\alpha) \right\} - \alpha \right|\\
		& \quad  = o(1) + \mathbb{P}\left\{ \big\| \widehat{\Omega}_n - \widehat{\Omega}_n^m\big\|_{op} > \big( \kappa  \mathfrak{n}_1(\lambda) + \bar{\sigma}^2 \big) \|T^3(T + \lambda)^{-4}\|_{op}   \sqrt{\frac{\log n}{n\lambda^2}}\right\},
	\end{align*}
	where $\bar{\sigma}^2 \geq \sigma_0^2 \vee \kappa^2(B + \kappa\|f_0\|_{\mathcal{H}})^2 \vee 1$ and $\mathfrak{n}_1^2(\lambda) = \mathrm{tr}\left((T+ \lambda)^{-2}T\right)$.
\end{proposition}

The finite metric entropy condition on the set $S$ ensures that the Gaussian bootstrap process~\eqref{eq:subsec:Application-FDA-2} is almost surely bounded and uniformly continuous on $S^* \subset \mathcal{H}$ (as required by Theorem~\ref{theorem:Abstract-Bootstrap-Sup-Empirical-Process-Quantiles}). This condition is not merely technical but also intuitive: Since the RKHS $\mathcal{H}$ is the completion of $\mathrm{span}(\{k_x : x \in S\})$ and $k$ is bounded and continuous, conditions that guarantee the continuity of Gaussian bootstrap processes on (a subset of) $\mathcal{H}$ should indeed be attributable to properties of $S$. Importantly, the metric entropy condition on $S$ does not impose restrictions on the dimension of $\mathcal{H}$. Only Assumption~\ref{assumption:cov-operator} (ii) implicitly imposes restrictions on the dimension of $\mathcal{H}$. 

Since under the conditions of Proposition~\ref{theorem:Uniform-CI-Bands-RKHS}, $\lim_{m \rightarrow \infty} \| \widehat{\Omega}_n - \widehat{\Omega}_n^m\|_{op} = 0$ almost surely for all $n \geq 1$, it follows that the bootstrap confidence band proposed in~\eqref{eq:subsec:Application-FDA-0} is asymptotically valid:

\begin{corollary}[Simultaneous bootstrap confidence bands]\label{corollary:theorem:Uniform-CI-Bands-RKHS}
	Under the setup of Proposition~\ref{theorem:Uniform-CI-Bands-RKHS},
	\begin{align*}
			\lim_{n \rightarrow \infty} \lim_{m\rightarrow \infty}  \sup_{\alpha \in (0,1)} \left| \mathbb{P}\Big\{ f_0 \in \mathcal{R}_n\big(c_n^m(1-\alpha)\big) \Big\} - (1- \alpha)\right| = 0,
	\end{align*}
	with the uniform confidence band $\mathcal{R}_n$ as defined in~\eqref{eq:subsec:Application-FDA-0} and quantile $c_n(1-\alpha)$ as in~\eqref{eq:subsec:Application-FDA-3}.
\end{corollary}

A thorough comparison of the Gaussian process bootstrap and~\possessivecite{singh2023kernel} symmetrized multiplier bootstrap is beyond the scope of this paper. In practice, both methods yield biased confidence bands for $f_0$, albeit for different reasons:~\possessivecite{singh2023kernel} bias stems from constructing a confidence band for the pseudo-true regression function $f_\lambda = (T + \lambda)^{-1} T f_0$ without correcting the regularization bias induced by $\lambda >0 $, ours is due to using an $m$-truncated Karhunen-Lo{\`e}ve decomposition based on a finite number of eigenfunctions. In future work we will explore ways to mitigate these biases by judiciously choosing $\lambda \equiv \lambda(n) \rightarrow 0$ and $m \equiv m(d) \rightarrow \infty$.

\section{Conclusion}\label{sec:Conclusion}

In this paper we have developed a new approach to theory and practice of Gaussian and bootstrap approximations of the sampling distribution of suprema of empirical processes. We have put special emphasize on non-asymptotic approximations that are entropy- and weak variance-free, and have allowed the function class $\mathcal{F}_n$ to vary with the sample size $n$ and to be non-Donsker. We have shown that such general approximation results are useful, among other things, for inference on high-dimensional statistical models and reproducing kernel Hilbert spaces. However, theory and methodology in this paper have three limitations that need to be addressed in future work:
\begin{itemize}
	\item \emph{Reliance on independent and identically distributed data.} All statistically relevant results in this paper depend on Proposition~\ref{lemma:CLT-Max-Norm}, which heavily relies on the assumption of independent and identically distributed data. Expanding Proposition~\ref{lemma:CLT-Max-Norm} to accommodate non-identical distributed data would be a first step towards solving simultaneous and large-scale two-sample testing problems and conducting inference in high-dimensional fixed design settings. Currently, the results in this paper are exclusively applicable to one-sample testing and unconditional inference.
	
	\item \emph{Lack of tight lower bounds on the strong variances of most Gaussian processes.} One the most notable features of the results in this paper is the fact that all upper bounds depend on the inverse of the strong variance of some Gaussian proxy process. Unfortunately, in statistical applications this poses a formidable challenge since up until now there exist only few techniques to derive tight lower bounds on these strong variances~\citep[][]{giessing2023bootstrap, giessing2022anticoncentration}. We either need new tools or we need to develop Gaussian and bootstrap approximations for statistics other than maxima/ suprema. The latter will require new anti-concentration inequalities.
	
	\item \emph{Biased quantile estimates due to bootstrapping non-pivotal statistic.} The Gaussian process bootstrap is based on a non-pivotal statistic, i.e. the sampling distribution of the supremum depends on the unknown population covariance function. In practice, when bootstrapping non-pivotal statistics the estimated quantiles often differ substantially from the true quantiles. Several bias correction schemes have been proposed in the classical setting~\citep{davison1986bootstrap, beran1987prepivoting, hall1988resampling, shi1992accurate}. Since the Gaussian process bootstrap is not a re-sampling procedure in the classical sense, these techniques do not apply. In~\cite{giessing2023bootstrap} we therefore develop the spherical bootstrap to improve accuracy and efficiency when bootstrapping $\ell_p$-statistics. However, this approach does not generalize to arbitrary empirical processes as the one in Section~\ref{subsec:Application-FDA}. Thus, there is an urgent need for new bias correction schemes.
\end{itemize}

\newpage
\section*{Acknowledgement}
Alexander Giessing is supported by NSF grant DMS-2310578.

\newpage
\normalsize
\setcounter{page}{1}

\bibliography{GBA_51_ref}

\newpage

\title{\textsc{Supplementary Materials for ``Gaussian and Bootstrap Approximations for Empirical Processes''}}
\author{Alexander Giessing\protect\footnotemark[1]}

\date{\today}

\maketitle

\appendix

\setcounter{page}{1}

\noindent{\bf\LARGE Contents}

\startcontents[sections]
\printcontents[sections]{l}{1}{\setcounter{tocdepth}{2}}

\section{Two toy examples}\label{sec:DerivationIntro}

In this section we present two examples to illustrate the limitations of the existing and the advantage of the new Gaussian approximation results. The presentation is deliberately expository; we omit proofs as much as possible. In both examples, we take $\mathcal{F}_n = \{x \mapsto f(x) = x'u : u \in S^{d-1}\}$ with $S^{d-1} = \{ u \in \mathbb{R}^d : \|u \|_2 = 1\}$, which plays a role in the construction of high-dimensional confidence regions and multiple testing problems (see Section~\ref{subsec:Application-Vector}). 


Consider a simple random sample $X_1, \ldots, X_n \in \mathbb{R}^d$ from the law of $X = a \cdot \xi$, where $a \in \mathbb{R}^d$ is a fixed vector and $\xi \in \mathbb{R}$ a centered random variable with finite third moment.  Then, $\sup_{f \in \mathcal{F}_n} |\mathbb{G}_n(f)| = \|n^{-1/2} \sum_{i=1}^n X_i\|_2 = \|a\|_2 |n^{-1/2} \sum_{i=1}^n \xi_i| $. Whence, the supremum of the empirical process reduces to the average of i.i.d. scalar-valued random variables and the classical univariate Berry-Essen theorem yields the following bound on the Kolmogorov distance:
\begin{align}\label{eq:sec:Intro-3}
	\varrho_n \lesssim \frac{\mathrm{E}[|\xi|^3]}{\sqrt{n} \: \mathrm{E}[|\xi|^2]^{3/2}} \rightarrow 0 \quad \mathrm{as} \quad n \rightarrow \infty.
\end{align} 
This non-asymptotic bound is independent of the dimension $d$ and the weak variance of the Gaussian proxy process. If we ignore the low-rank structure of the data and instead use the Gaussian approximation results in~\cite{chernozhukov2014GaussianApproxEmp} (Theorem 2.1 combined with Lemma 2.3) we obtain (qualitatively)
\begin{align}\label{eq:sec:Intro-4}
	\varrho_n \lesssim_\star \frac{\mathrm{polylog}\left(N \left(\mathcal{F}_n, e_P, \varepsilon \right) \vee n \right)}{ n^{1/6} \sqrt{\inf_{ u \in S^{d-1} } \mathrm{Var}(Z'u)}},
\end{align} 
where $Z \sim N\left(0,  \mathrm{E}[\xi^2] \cdot aa'  \right)$ and $\lesssim_\star$ hides a complicated multiplicative factor which, among other things, depends on the sample size $n$, third and higher moments of $\|Z\|_2$, and the inverse of the discretization level $\varepsilon > 0$ (see Section~\ref{subsec:RelationPreviousWork-Stat}). Since the covariance matrix $\mathrm{E}[\xi^2] \cdot aa' $ has rank one, the unit ball in $\mathcal{F}_n$ with respect to the intrinsic standard deviation metric $e_P$ is isometrically isomorphic to the interval $[-1,1] \subset \mathbb{R}$. Hence, the metric entropy can be upper bounded independently of the dimension $d$; in particular, we have $ \log N \left(\mathcal{F}_n, e_P, \varepsilon \right) \leq \log (1 + 2/\varepsilon)$. However, the low-rank structure of the data also implies that the weak variance of the associated Gaussian proxy process vanishes; indeed, $\inf_{ u \in S^{d-1} } \mathrm{Var}(Z'u) =  \mathrm{E}[\xi^2] \cdot \inf_{u \in S^{d-1}} (a'u)^2 =0$. Thus, the upper bound in~\eqref{eq:sec:Intro-4} is in fact invalid (or trivial if we interpret $1/0 = \infty$) and fails to replicate the univariate Berry-Esseen bound in~\eqref{eq:sec:Intro-3}.

In contrast, the new Gaussian approximation inequality (Theorem~\ref{theorem:CLT-Max-Norm-Simultaneous} in Section~\ref{subsec:GaussianApproximation}; Theorem A.1 in~\cite{giessing2023bootstrap}) is agnostic to the covariance structure of the data and yields
\begin{align}\label{eq:sec:Intro-5}
	\varrho_n \lesssim \frac{(\mathrm{E}[\|X \|_2^3])^{1/3}}{n^{1/6}\sqrt{\mathrm{Var}(\|Z\|_2)}} + \frac{\mathrm{E} \left[ \|X\|_2^3 \mathbf{1}\{\|X\|_2^3 > n\: \mathrm{E}[\|X \|_2^3]\}\right]}{\mathrm{E} \left[ \|X\|_2^3\right]}  + \frac{ \mathrm{E}[\|Z\|_2]}{\sqrt{n\mathrm{Var}(\|Z\|_2)}},
\end{align}
where $\lesssim$ hides an absolute constant independent of $n, d$, and the distribution of the $X_i$'s and $Z \sim N\left(0,  \mathrm{E}[\xi^2] \cdot aa'  \right)$.
The upper bound in this inequality depends only on the third moment of the envelop function $x \mapsto \|x\|_2$ and the strong variance of the Gaussian proxy process $\mathrm{Var}(\|Z\|_2)$. If we use that $\|X\|_2 = \|a\|_2 \cdot|\xi|$ and $\|Z\|_2 = \|a\|_2 \cdot |g|$ for $g \sim N(0,1)$, we obtain
\begin{align}\label{eq:sec:Intro-6}
	\varrho_n \lesssim \frac{\mathrm{E}[|\xi|^3]^{1/3}}{n^{1/6} \: \mathrm{E}[|\xi|^2]^{1/2} } + \frac{\mathrm{E} \left[ |\xi|^3 \mathbf{1}\{|\xi|^3 > n\: \mathrm{E}[|\xi|^3]\}\right]}{\mathrm{E} \left[ |\xi|^3\right]} \rightarrow 0 \quad \mathrm{as} \quad n \rightarrow \infty.
\end{align}
This inequality is obviously dimension- and weak variance-free. In this sense, it recovers the essential feature of the univariate Berry-Esseen bound in~\eqref{eq:sec:Intro-3} and improves over the bound in~\eqref{eq:sec:Intro-4}. The dependence on the sample size $n$ and the moments of $\xi$ is still sub-optimal; but refinements in this direction are beyond the scope of this paper. Obviously, in this example, we have chosen a rank one covariance matrix only to be able to compare the Gaussian approximation result with the Berry-Esseen theorem. Any low-rank structure implies a vanishing weak variance and, hence, a breakdown of the results in~\cite{chernozhukov2014GaussianApproxEmp}.

Next, suppose that the data $X_1, \ldots, X_n \in \mathbb{R}^d$ are a simple random sample drawn from the law of a random vector $X = (X^{(1)}, \ldots, X^{(d)})'$ with mean zero, element-wise bounded entries $\max_{1 \leq k \leq d}|X^{(k)}| < B$ almost surely for some $B > 0$, and equi-correlated covariance matrix $\Sigma = (1-\rho) I_d + \rho \mathbf{1}_d \mathbf{1}_d'$ for some $\rho \in (-1/(d-1), 1)$. The constraints on $\rho$ guarantee that $\Sigma$ has full rank. Therefore, the multivariate Berry-Esseen theorem by~\cite{bentkus2003DependenceBerryEsseen} (Theorem 1.1) implies
\begin{align}\label{eq:sec:Intro-7}
	\varrho_n \lesssim \frac{d^{1/4} \: \mathrm{E}[\|\Sigma^{-1/2}X\|_2^3]}{\sqrt{n} }.
\end{align}
This upper bound is not useful in high-dimensional settings because the expected value is polynomial in the dimension $d$, i.e. $\mathrm{E}[\|\Sigma^{-1/2}X\|_2^3] \geq \mathrm{E}[\|\Sigma^{-1/2}X\|_2^2]^{3/2} = d^{3/2}$. 

The Gaussian approximation results by~\cite{chernozhukov2014GaussianApproxEmp} yield (again) inequality~\eqref{eq:sec:Intro-4}. Since the covariance matrix $\Sigma$ has full rank, the weak variance is now strictly positive and equal to the smallest eigenvalue of $\Sigma$, i.e. $\inf_{ u \in S^{d-1} } \mathrm{Var}(Z'u) =  \inf_{u \in S^{d-1}} u'\Sigma u = 1 - \rho > 0$. However, in this example, the metric entropy with respect to the intrinsic standard deviation metric $e_P$ poses a problem. Indeed, let $\lambda_1 = 1 +(d-1) \rho$ and $\lambda_2 = \ldots = \lambda_d = 1-\rho$ be the eigenvalues of $\Sigma$. Then, the unit ball in $\mathcal{F}_n$ with respect to $e_P$ can be identified with the weighted Euclidean ball $B_\lambda^d(0,\varepsilon) := \{u \in \mathbb{R}^d: \sum_{k=1}^d \lambda_k u_k \leq \varepsilon^2 \}$. Let $B^d(0,1)$ be the standard unit ball in $\mathbb{R}^d$ with respect to the Euclidean distance. Then, $B^d(0,\varepsilon) \subseteq B_\lambda^d(0, \varepsilon)$ whenever $\varepsilon < 1-\rho$.  Therefore, standard arguments $N \left(\mathcal{F}_n, e_P, \varepsilon \right) \geq \mathrm{vol}\left(B^d(0,1)\right)/ \mathrm{vol}\left(B_\lambda^d(0,\varepsilon)\right) \geq \mathrm{vol}\left(B^d(0,1)\right)/ \mathrm{vol}\left(B^d(0,\varepsilon)\right) \geq \left((1-\rho)/\varepsilon\right)^d$. Thus, the metric entropy grows linear in the dimension $d$, i.e. $\log N \left(\mathcal{F}_n, e_P, \varepsilon \right) \geq d \log\left((1-\rho)/\varepsilon\right)$ for $\varepsilon \downarrow 0$. We conclude that the results in~\cite{chernozhukov2014GaussianApproxEmp} are again not useful in high dimensions with $d \gg n$.

The new Gaussian approximation inequality implies (again)~\eqref{eq:sec:Intro-5}. Since $\max_{1 \leq k \leq d}|X^{(k)}| < B$ almost surely and $\mathrm{tr}(\Sigma) = d$ it follows that $\sqrt{d} = (\mathrm{E}[\|Z\|_2^2])^{1/2} = (\mathrm{E}[\|X \|_2^2])^{1/2} \leq (\mathrm{E}[\|X \|_2^3])^{1/3} \leq B \sqrt{d}$. Moreover, by Theorem A.6 in~\cite{giessing2023bootstrap} and since $\mathrm{tr}(\Sigma^2) = d(1-\rho)^2 + \rho^2 d^2$, we have $\mathrm{Var}(\|Z\|_2) \geq \mathrm{tr}(\Sigma^2)/ \mathrm{tr}(\Sigma) =  (1-\rho)^2 + \rho^2 d$. Therefore, inequality~\eqref{eq:sec:Intro-5} simplifies to
\begin{align}\label{eq:sec:Intro-8}
	\varrho_n 
	\lesssim \frac{B }{n^{1/6} \rho }  + B^3  \mathbf{1}\{ B  > n \} \rightarrow 0 \quad \mathrm{as} \quad n \rightarrow \infty.
\end{align}
This inequality is not only  dimension- and weak variance-free but also improves qualitatively over, both, the results by~\cite{chernozhukov2014GaussianApproxEmp} and~\cite{bentkus2003DependenceBerryEsseen}. Intuitively, the reason why we are able to shed the $d^{1/4}$-factor from the upper bound compared to the results in~\cite{bentkus2003DependenceBerryEsseen} is that we only take the supremum over all Euclidean balls with center at the origin whereas he takes the supremum over all convex sets in $\mathbb{R}^d$. (Note that to apply his bound in the context of this example, we need to take the supremum over at least all weighted Euclidean balls $B_\lambda^d(0,\varepsilon) := \{u \in \mathbb{R}^d: \sum_{k=1}^d \lambda_k u_k \leq \varepsilon^2 \}$.) This second example is related to the first one in so far as the covariance matrix has ``approximately'' rank one. Indeed, as the dimension increases the law of the random vector $X$ concentrates in the neighborhood of the one-dimensional subspace spanned by eigenvector associated to the largest eigenvalue of $\Sigma$. This becomes even more obvious if we consider the standardized covariance $\Sigma/ \|\Sigma\|_{op}$ with eigenvalues $\lambda_1 = \rho + (1-\rho)/d \rightarrow \rho$ and $\lambda_2 = \ldots \lambda_d = (1 - \rho)/d \rightarrow 0$ as $d \rightarrow \infty$.  

In both examples, the exact and approximate low-rank structures of the data are crucial in order to go from the abstract bound~\eqref{eq:sec:Intro-5} to the dimension-/ entropy-free bounds~\eqref{eq:sec:Intro-6} and~\eqref{eq:sec:Intro-8}, respectively. This is not coincidental and is in fact representative for the entire theory that we develop in this paper: While our Gaussian and bootstrap approximation inequalities hold without assumptions on the metric entropy, in concrete examples they often only yield entropy-free upper bounds if the trace of the covariance operator is bounded or grows at a much slower rate than the sample size $n$ (e.g. low-rank, bounded effective rank, or variance decay, see Section~\ref{sec:Applications}). This is certainly a limitation of our theory, but a low-rank covariance (function) is an empirically well-documented property of many data sets and therefore a common assumption in multivariate and high-dimensional statistics~\citep[e.g.][]{anderson2003introduction, vershynin2018HighDimProb}. 

\newpage

\section{Auxiliary results}\label{sec:AuxiliaryResults}

\subsection{Smoothing inequalities and partial derivatives}\label{subsec:AuxiliaryResults-Smoothing}

\begin{lemma}\label{lemma:Smooth-Lipschitz-Approx} Let $X, Z \in \mathbb{R}$ be arbitrary random variables. There exists a map $h_{s, \lambda} \in C^\infty_c(\mathbb{R})$ such that
	\begin{itemize}
		\item[(i)] for all $s, t \in \mathbb{R}$ and $\lambda >0$,
		\begin{align*}
			|D^k h_{s, \lambda}(t)| \leq C_k \lambda^{-k} \mathbf{1}\{ s \leq t \leq s + 3 \lambda \},
		\end{align*}
		where $C_k > 0$ is a constant depending only on $k \in \mathbb{N}_0$; and
		\item[(ii)] for all $\lambda > 0$,
		\begin{align*}
			\left| \sup_{s \in \mathbb{R}} \left|\mathbb{P}\left\{ X \leq s \right\} - \mathbb{P}\left\{Z \leq s\right\}\right| - \sup_{s \in \mathbb{R}} \big|\mathrm{E}[h_{s, \lambda}(X) - h_{s, \lambda}(Z)] \big| \right| \leq \zeta_{3\lambda}(X) \wedge \zeta_{3\lambda}(Z),
		\end{align*}
		where $\zeta_\lambda(V) := \sup_{s \in \mathbb{R}} \mathbb{P}\{s \leq V \leq s + \lambda \}$ for real-valued $V \in \mathbb{R}$.
	\end{itemize}
\end{lemma}
\begin{remark}
	We can take $C_0 = C_1 =1$ (see proof).
\end{remark}

\begin{lemma}\label{lemma:DifferentiatingUnderIntegral}
	Let $h_{s, \lambda} \in C^\infty_0(\mathbb{R})$ be the map from Lemma~\ref{lemma:Smooth-Lipschitz-Approx} and define $x \mapsto h(x) := h_{s, \lambda}(\|x\|_\infty)$ for $x \in \mathbb{R}^d$. Let $(P_t)_{t \geq 0}$ be the Ornstein-Uhlenbeck semi-group with stationary measure $N(0, \Sigma)$ and positive definite covariance matrix $\Sigma$.
	\begin{itemize}
		\item[(i)] For arbitrary indices $1 \leq i_1, \ldots, i_k \leq d$, $k \geq 1$, and all $x_0 \in \mathbb{R}^d$,
		\begin{align*}
			\frac{\partial^k}{\partial x_{i_1} \cdots \partial x_{i_k}} \left( \int_0^\infty P_t h(x) dt\right) \Big|_{x=x_0} = \int_0^\infty e^{-kt}P_t \left(\frac{\partial^k h}{\partial x_{i_1} \cdots \partial x_{i_k}} \right) (x_0)dt;
		\end{align*}
		\item[(ii)] for almost every $x_0 \in \mathbb{R}^d$ the absolute value of the derivative in (i) can be upper bounded by
		\begin{align*}
			C_k\lambda^{-k} \int_0^\infty e^{-kt}\mathrm{E}\left[\mathbf{1}\left\{ s \leq \|V_0^t\|_\infty \leq s + 3 \lambda \right\} \mathbf{1}\left\{|V_{0{i_1}}^t| \geq |V_{0\ell}^t|, \: \ell \neq {i_1} \right\} \right] dt \:  \mathbf{1}\{i_1 = \ldots = i_k\},
		\end{align*}
		where $V_0^t = e^{-t}x_0 + \sqrt{1- e^{-2t}}Z$, $ Z\sim N(0, \Sigma)$, and $C_k  >0$ is the absolute constant from Lemma~\ref{lemma:Smooth-Lipschitz-Approx}.
	\end{itemize}
\end{lemma}
\begin{remark}
	While claim (i) looks a lot like a ``differentiating under the integral sign" type result, it is more accurate to think of it as a specific smoothing property of the Ornstein-Uhlenbeck semigroup with stationary measure $N(0, \Sigma)$ when applied to (compositions of Lipschitz continuous functions with) the map $x \mapsto \|x\|_\infty$.
\end{remark}

\begin{lemma}\label{lemma:L1ConvergencePartialRegularization}
	Let $\varrho \in L^1(\mathbb{R})$ with $\int \varrho(r) dr = 1$. For $\eta > 0$ and a map $h$ on $\mathbb{R}^d$ set
	\begin{align*}
		(\varrho_\eta \ast_{(i)} h) (x) := \int \varrho(r) h(x - r\eta e_i) dr,
	\end{align*}
	where $e_i$ denotes the $i$th standard unit vector in $\mathbb{R}^d$.
	For $k \in \mathbb{N}$, $f_1, \ldots, f_k \in B(\mathbb{R}^d)$, $g \in L^1(\mathbb{R}^d)$, and arbitrary indices $1 \leq i_1, \ldots, i_k \leq d$,
	\begin{align*}
		\left\| \prod_{j=1}^k(\varrho_\eta \ast_{(i_j)} f_j) g - \prod_{j=1}^kf_jg \right\|_1 \rightarrow 0 \quad{} \mathrm{as} \quad{} \eta \rightarrow 0.
	\end{align*}
\end{lemma}
\begin{remark}
	The proof of this result on ``partially regularized'' functions is similar to the one on ``fully regularized'' functions~\citep[][Theorem 8.14 (a)]{folland1999real}. The conditions on the functions $f_1, \ldots, f_k, g, \varrho$ can probably be relaxed, but they are sufficiently general to apply to the situations that we encounter.
\end{remark}

\begin{lemma}[Partial derivatives of compositions of almost everywhere diff'able functions]\label{lemma:ChainRule-SecondOrder} 
	Let $g \in C^k(\mathbb{R})$ and $f \in C^k(\mathbb{R}^d\setminus \mathcal{N})$ and $ \mathcal{N}$ is a null set with respect to the Lebesgue measure on $\mathbb{R}^d$. Then, $g \circ f \in C^k(\mathbb{R}^d\setminus  \mathcal{N})$ and, for arbitrary indices $1 \leq i_1, \ldots, i_k \leq d$, $k \geq 1$, and all $x \in \mathbb{R}^d \setminus  \mathcal{N}$,
	\begin{align*}
		\frac{\partial^k (g \circ f)}{\partial x_{i_1} \cdots \partial x_{i_k}}(x)= \sum_{\pi \in \Pi} (D^{|\pi|} g \circ f) \prod_{B \in \pi} \frac{\partial^{|B|}f}{\prod_{j \in B} \partial x_j}(x),
	\end{align*}
	where $\Pi$ is the set of all partitions of $\{i_1, \ldots, i_k\}$, $|\pi|$ denotes the number of ``blocks of indices'' in partition $\pi \in \Pi$, and $|B|$ is the number of indices in block $B \in \pi$.
\end{lemma}
\begin{remark}
	This result is a trivial modification of the multivariate version of Fa{\`a} di Bruno's formula due to~\cite{hardy2006combinatorics}. The modification is that we only require $f$ to be $k$-times differentiable almost everywhere on $\mathbb{R}^d$. The formula is generally false if $g$ is only $k$-times differentiable almost everywhere on $\mathbb{R}$.
\end{remark}
\begin{remark}
	The first three partial derivatives are of particular interest to us. They are given by (whenever they exist)
	\begin{align*}
		\frac{\partial (g \circ f)}{\partial x_i}(x) &= \left( Dg \circ f\right) (x) \frac{\partial f}{\partial x_i}(x),\\
		\frac{\partial^2 (g \circ f)}{\partial x_i \partial x_j}(x) &= \left( D^2 g \circ f\right) (x) \frac{\partial f}{\partial x_i}(x) \frac{\partial f}{\partial x_j}(x) + \left( D g \circ f\right) (x) \frac{\partial^2 f}{\partial x_i \partial x_j}(x),\\
		\frac{\partial^3 (g \circ f)}{\partial x_i \partial x_j \partial x_k}(x) &= \left( D^3 g \circ f\right) (x) \frac{\partial f}{\partial x_i}(x) \frac{\partial f}{\partial x_j}(x)  \frac{\partial f}{\partial x_k}(x)\\
		&\quad{} + \left( D^2 g \circ f\right) (x) \left[ \frac{\partial^2 f}{\partial x_i \partial x_j}(x)\frac{\partial f}{\partial x_k}(x) +   \frac{\partial^2 f}{\partial x_i \partial x_k}(x)\frac{\partial f}{\partial x_j}(x) + \frac{\partial^2 f}{\partial x_k \partial x_j}(x)\frac{\partial f}{\partial x_i}(x)\right]\\
		&\quad{} + \left( D g \circ f\right) (x) \frac{\partial^3 f}{\partial x_i \partial x_j \partial x_k}(x).
	\end{align*}
\end{remark}

\begin{lemma}[Partial and total derivatives of $\ell_\infty$-norms]\label{lemma:DerivativesMaxNorm}
	The map $f(x) = \max_{1 \leq k \leq d} |x_k|$ is partially differentiable of any order almost everywhere on $\mathbb{R}^d$ with partial derivatives (whenever they exist)
	\begin{align*}
		\frac{\partial f}{\partial x_i}(x) = \mathrm{sign}(x_i) \mathbf{1}\{ |x_i| \geq |x_k|, \: \forall k \} \quad{} \quad{} \text{and} \quad{}\quad{}  \frac{\partial^k f}{\partial x_{i_1} \cdots \partial x_{i_k}}(x) = 0,
	\end{align*}
	for $1 \leq i, i_1, \ldots, i_k \leq d$ and $k \geq 2$. Moreover, $f$ is twice totally differentiable almost everywhere on $\mathbb{R}^d$ with Jacobian and Hessian matrices (whenever they exist) 
	\begin{align*}
		Df(x) = \left[ \frac{\partial f}{\partial x_1}(x), \ldots, \frac{\partial f}{\partial x_d}(x)\right] \quad{}\quad{} \text{and} \quad{}\quad{} D^2f(x) = \mathbf{0} \in \mathbb{R}^{d \times d}.
	\end{align*}
\end{lemma}
\begin{remark}
	Note that the first partial derivative can be re-written (less compact) as a piece-wise linear function.
\end{remark}

\subsection{Anti-concentration inequalities and lower bounds on variances}\label{subsec:AuxiliaryResults-AntiConcentration}

\begin{lemma}[\citeauthor{giessing2022anticoncentration},~\citeyear{giessing2022anticoncentration}]\label{lemma:AntiConcentration-SeparableProcess}
	Let $X = (X_u)_{u \in U}$ be a centered separable Gaussian process indexed by a semi-metric space $U$. Set $Z = \sup_{u \in U}X_u$ and assume that $0 \leq Z < \infty$ a.s. For all $\varepsilon \geq 0$,
	\begin{align*}
		\frac{\varepsilon/\sqrt{12}}{ \sqrt{\mathrm{Var}(Z) + \varepsilon^2/ 12}} \leq \sup_{t \geq 0} \mathbb{P}\left\{ t \leq  Z \leq t + \varepsilon  \right\} \leq \frac{\varepsilon\sqrt{12}}{ \sqrt{\mathrm{Var}(Z) + \varepsilon^2/ 12}}.
	\end{align*}
	The result remains true if $Z$ is replaced by $\widetilde{Z} = \sup_{u \in U}|X_u|$.
\end{lemma}
\begin{remark}
	If the covariance function of $X$ is positive definite, then above inequalities hold even for uncentered $X = (X_u)_{u \in U}$ and $Z \in [-\infty, \infty)$ a.s.
\end{remark}

\begin{lemma}[\citeauthor{giessing2022anticoncentration},~\citeyear{giessing2022anticoncentration}]\label{lemma:BoundsVariance-SeparableProcess}
	Let $X = (X_u)_{u \in U}$ be a separable Gaussian process indexed by a semi-metric space $U$ such that $\mathrm{E}[X_u] = 0$, $0 < \underline{\sigma}^2 \leq \mathrm{E}[X_u^2] \leq \bar{\sigma}^2 < \infty$, and $|\mathrm{Corr}(X_u, X_v)| \leq \rho$ for all $u, v \in U$. Set $Z = \sup_{u \in U}X_u$ and assume that $Z < \infty$ a.s. Then, $ 0 \leq \mathrm{E}[Z] < \infty$ and there exist absolute constants $c, C > 0$ such that
	\begin{align*}
		\frac{1}{C} \left(\frac{\underline{\sigma}}{1 + \mathrm{E}[Z/\underline{\sigma}]}\right)^2 \leq \mathrm{Var}(Z) \leq C \left[ \bar{\sigma}^2 \wedge \left( \bar{\sigma}^2 \rho +  \left(\frac{\bar{\sigma}}{ (\mathrm{E}[Z/ \bar{\sigma}] - c)_+}\right)^2 \right) \right],
	\end{align*}
	with the convention that ``$1/0 = \infty$''. The result remains true if $Z$ is replaced by $\widetilde{Z} = \sup_{u \in U}|X_u|$.
\end{lemma}

\begin{lemma}[\citeauthor{lecam1986asymptotic},~\citeyear{lecam1986asymptotic}, p. 402]\label{lemma:Kolmogorov-Coupling-AntiConcentration} For $X, Z \in \mathbb{R}$ arbitrary random variables and $\lambda > 0$,
	\begin{align*}
		\sup_{s \geq 0} \Big| \mathbb{P} \left\{ X \leq s \right\}  - \mathbb{P} \left\{ Z \leq s \right\}  \Big| &\leq \mathbb{P}\left\{  |X - Z| > \lambda \right\} + \zeta_{\lambda}(X) \wedge \zeta_{\lambda}(Z),
	\end{align*}	
	where $\zeta_\lambda(V) := \sup_{s \in \mathbb{R}} \mathbb{P}\{s \leq V \leq s + \lambda\}$ for real-valued $V \in \mathbb{R}$.
\end{lemma}


\begin{lemma}\label{lemma:VarianceMaximum}
	Let $X, Z \in \mathbb{R}$ be arbitrary random variables. Then, 
	\begin{align*}
		\left| \sqrt{\mathrm{Var}(X)} - \sqrt{\mathrm{Var}(X \vee Z)}\right| \leq \sqrt{\mathrm{Var}\big( (Z-X)_+\big)},
	\end{align*}
	where $(a)_+ = \max(a, 0)$ for $a \in \mathbb{R}$. Moreover, if $\mathrm{E}[Z] \geq \mathrm{E}[X]$, then
	\begin{align*}
		\mathrm{Var}\big((Z-X)_+\big) \leq 	\mathrm{Var}(Z-X).
	\end{align*}
\end{lemma}

\subsection{Quantile comparison lemmas}\label{subsec:AuxiliaryResults-QuantileComparison}
Throughout this section, $Z_n \mid X_1, \ldots, X_n \sim N(0, \widehat{\Sigma}_n)$ and $Z \sim N(0, \Sigma)$. For $\alpha \in (0,1)$ we define the $\alpha$th quantile of $\|Z_n\|_\infty$ and $\|Z\|_\infty$, respectively, by
\begin{align}\label{eq:subsec:QuantileComparison-0}
	\begin{split}
	c_n(\alpha; \widehat{\Sigma}_n) &:= \inf \left\{s \geq 0: \mathbb{P}\left\{\|Z_n\|_\infty \leq s \mid X_1, \ldots, X_n \right\}  \geq \alpha \right\},\\
	c_n(\alpha; \Sigma) &:= \inf \left\{s \geq 0: \mathbb{P}\left\{\|Z\|_\infty \leq s \right\}  \geq \alpha \right\}.
	\end{split}
\end{align}

\begin{lemma}\label{lemma:Bootstrap-Max-Norm-Quantil-Comparison}
	For all $\delta > 0$,
	\begin{align*}
		&\inf_{\alpha \in (0,1)} \mathbb{P}\Big\{c_n(\alpha; \widehat{\Sigma}_n) \leq c_n(\pi_n(\delta) + \alpha; \Sigma) \Big\} \geq 1 - \mathbb{P}\left\{ \max_{j,k} |\widehat{\Sigma}_{n,jk} - \Sigma_{jk}|  > \delta \right\}, \quad{} \text{and}\\
		&\inf_{\alpha \in (0,1)} \mathbb{P}\Big\{c_n(\alpha; \Sigma) \leq c_n(\pi_n(\delta) + \alpha; \widehat{\Sigma}_n) \Big\} \geq 1 - \mathbb{P}\left\{ \max_{j,k} |\widehat{\Sigma}_{n,jk} - \Sigma_{jk}|  > \delta \right\},
	\end{align*}
	where $\pi_n(\delta) = K\delta^{1/3} \left(\mathrm{Var}(\|\Sigma^{1/2} Z\|_\infty)\right)^{-1/3}$ and $K > 0$ is an absolute constant.
\end{lemma}

If we use non-Gaussian proxy statistics to test hypothesis (such as Efron's empirical bootstrap), we replace Lemma~\ref{lemma:Bootstrap-Max-Norm-Quantil-Comparison} by the following result:

\begin{lemma}\label{lemma:Bootstrap-Max-Norm-Quantil-Comparison-Non-Gaussian}
	Let $T_n = S_n + R_n \geq 0$ be a statistic; $S_n$ and $R_n$ need not be independent. For $\alpha \in (0,1)$ arbitrary define $c_n^T(\alpha) := \inf \{s \geq 0: \mathrm{P}(T_n \leq s \mid X_1, \ldots, X_n)  \geq \alpha \}$. Then, for all $\delta, \eta > 0$,
	\begin{align*}
		& \inf_{\alpha \in (0,1)} \mathbb{P} \left\{ c_n^T(\alpha) \leq c_n( \kappa_n(\delta) + \eta + \alpha; \Sigma) \right\} \geq 1 - \mathbb{P}\left\{\gamma_n + \rho_n(\delta) > \eta \right\} \quad{} \text{and} \\
		&\inf_{\alpha \in (0,1)} \mathbb{P} \left\{ c_n(\alpha; \Sigma) \leq c_n^T( \kappa_n(\delta) + \eta + \alpha) \right\} \geq 1 - \mathbb{P}\left\{\gamma_n + \rho_n(\delta) > \eta  \right\},
	\end{align*}
	where $\gamma_n = \sup_{s \geq 0} | \mathbb{P}\{S_n \leq s \mid X_1, \ldots, X_n\} - \mathbb{P}\{\|\Sigma_n^{1/2}Z\|_\infty \leq s \}|$, $\kappa_n(\delta) =  K \delta (\mathrm{Var}(\|\Sigma_n^{1/2} Z\|_\infty))^{-1/2}$, $\rho_n(\delta)= \mathbb{P}\{|R_n| > \delta \mid X_1, \ldots, X_n\}$, and $K > 0$ is an absolute constant.
\end{lemma}

For $\alpha \in (0,1)$ we define the (conditional) $\alpha$th quantile of the supremum of $\{\widehat{Z}_n^m(f) : f \in \mathcal{F}_n\}$ by
\begin{align*}
	c_n(\alpha; \widehat{\mathcal{C}}_n^m) &:= \inf \left\{s \geq 0: \mathbb{P}\left\{  \|\widehat{Z}_n^m\|_{\mathcal{F}_n} \leq s \mid X_1, \ldots, X_n\right\} \geq \alpha \right\},
\end{align*}
and the $\alpha$th quantile of the supremum of the Gaussian $P$-bridge process $\{G_P(f) : f \in \mathcal{F}_n\}$ by
\begin{align*}
	c_n(\alpha; \mathcal{C}_P) &:= \inf \left\{s \geq 0: \mathbb{P}\left\{  \|G_P\|_{\mathcal{F}_n} \leq s \right\} \geq \alpha \right\},
\end{align*}

The following two lemmas are straightforward generalizations of the preceding lemmas to empirical processes. The proofs of Lemmas~\ref{lemma:Bootstrap-Sup-Empirical-Process-Quantil-Comparison} and~\ref{lemma:Bootstrap-Sup-Empirical-Process-Quantil-Comparison-Non-Gaussian} are identical to the ones of Lemmas~\ref{lemma:Bootstrap-Max-Norm-Quantil-Comparison} and~\ref{lemma:Bootstrap-Max-Norm-Quantil-Comparison-Non-Gaussian}, respectively. We therefore omit them.

\begin{lemma}\label{lemma:Bootstrap-Sup-Empirical-Process-Quantil-Comparison}
	For all $\delta > 0$,
	\begin{align*}
		&\inf_{\alpha \in (0,1)} \mathbb{P}\left\{ 	c_n(\alpha; \widehat{\mathcal{C}}_n^m) \leq c_n(\pi_n(\delta; \mathcal{C}_P) + \alpha) \right\} \geq 1 - \mathbb{P}\left\{  \sup_{f, g \in \mathcal{F}_n} \big| \mathcal{C}_P(f,g) - \widehat{\mathcal{C}}_n^m(f, g) \big| > \delta \right\}, \quad{} \text{and}\\
		&\inf_{\alpha \in (0,1)} \mathbb{P}\left\{ c_n(\alpha; \mathcal{C}_P) \leq c_n(\pi_n(\delta) + \alpha; \widehat{\mathcal{C}}_n^m) \right\} \geq 1 - \mathbb{P}\left\{ \sup_{f, g \in \mathcal{F}_n} \big| \mathcal{C}_P(f,g) - \widehat{\mathcal{C}}_n^m(f, g) \big| > \delta \right\},
	\end{align*}
	where $\pi_n(\delta) = K\delta^{1/3} \left(\mathrm{Var}(\|G_P\|_{\mathcal{F}_n})\right)^{-1/3}$ and $K > 0$ is an absolute constant.
\end{lemma}

\begin{lemma}\label{lemma:Bootstrap-Sup-Empirical-Process-Quantil-Comparison-Non-Gaussian}
	Let $\{Z_n(f) = Z_n^1(f) + Z_n^2(f) : f \in \mathcal{F}_n\}$ be an arbitrary stochastic process. For $\alpha \in (0,1)$ arbitrary define $c_{n, Z_n}(\alpha) := \inf \{s \geq 0: \mathbb{P}\{\|Z_n\|_{\mathcal{F}_n} \leq s \mid X_1, \ldots, X_n\} \geq \alpha \}$. Then, for all $\delta, \eta > 0$,
	\begin{align*}
		&\inf_{\alpha \in (0,1)} \mathbb{P}\left\{ c_{n, Z_n}(\alpha) \leq c_n(\kappa_n(\delta) +  \eta + \alpha; \mathcal{C}_P) \right\} \geq 1 - \mathbb{P} \left\{ \gamma_{n, Z_n^1} + \rho_{n,Z_n^2}(\delta) > \eta \right\},  \quad{} \text{and}\\
		&\inf_{\alpha \in (0,1)} \mathbb{P}\left\{ c_n(\alpha; \mathcal{C}_P) \leq c_{n, Z_n}(\kappa_n(\delta) +  \eta + \alpha) \right\} \geq 1 - \mathbb{P} \left\{ \gamma_{n, Z_n^1} + \rho_{n,Z_n^2}(\delta) > \eta \right\},
	\end{align*}
	where $\gamma_{n,Z_n^1} = \sup_{s \geq 0} | \mathbb{P}\{\|Z_n^1\|_{\mathcal{F}_n} \leq s \mid X_1, \ldots, X_n\} - \mathbb{P}\{\|G_P\|_{\mathcal{F}_n} \leq s\}  |$, $\kappa_n(\delta) =  K \delta /\sqrt{\mathrm{Var}(\|G_P\|_{\mathcal{F}_n})}$, $\rho_{n,Z_n^2}(\delta)= \mathbb{P}\left\{ \|Z_n^2\|_{\mathcal{F}_n} > \delta \mid X_1, \ldots, X_n \right\}$, and $K > 0$ is an absolute constant.
\end{lemma}

\subsection{Boundedness and continuity of centered Gaussian processes}\label{subsec:GaussianProcesses}
This section contains classical results on boundeness and continuity of the sample paths of centered Gaussian processes. We provide a proof for Lemma~\ref{lemma:ModulusContinuity} in Appendix~\ref{subsec:ProofModulusContinuity}; all other results (with proofs or references to proofs) can be found in Appendix A of~\cite{vandervaart1996weak}.

Throughout this section $X = (X_u)_{u \in U}$ denotes a centered separable Gaussian process indexed by a semi-metric space $U$, $Z = \sup_{u \in U} |X_u|$, $\sigma^2 = \sup_{u \in U} \mathrm{E}[X_u^2]$, and $d_X$ the intrinsic standard deviation metric associated with $X$.


\begin{lemma}[Equivalence of bounded sample path and finite expectation]\label{lemma:AS-Bounded}
	$X$ is almost surely bounded on $U$ if and only if $\mathrm{E}[Z] < \infty$.
\end{lemma}

\begin{lemma}[Reverse Liapunov inequality]\label{lemma:ReverseLiapunov}
	If $X$ is almost surely bounded on $U$, then there exist constants $K_{p, q} > 0$ depending on $0 < p \leq q < \infty$ only such that 
	\begin{align*}
		\left(\mathrm{E}[Z^q]\right)^{1/q} \leq K_{p,q} 	\left(\mathrm{E}[Z^p]\right)^{1/p} .
	\end{align*}
\end{lemma}


\begin{lemma}[Sudakov's lower bound]\label{lemma:Sudakov}
	Let $N(U, d_X, \varepsilon)$ be the $\varepsilon$-covering number of $U$ w.r.t. $d_X$. Then, there exists an absolute constant $K > 0$ such that
	\begin{align*}
		\sup_{\varepsilon > 0} \varepsilon \sqrt{N(U, d_X, \varepsilon)} \leq K\mathrm{E}[Z].
	\end{align*}
	Consequently, if $\mathrm{E}[Z] < \infty$, then $U$ is totally bounded w.r.t. $d_X$.
\end{lemma}

\begin{lemma}[Metric entropy condition for bounded and continuous sample paths]\label{lemma:MetricEntropCondition}
	Let $N(U, d_X, \varepsilon)$ be the $\varepsilon$-covering number of $U$ w.r.t. $d_X$. If $\int_0^\infty \sqrt{N(U, d_X, \varepsilon)}d\varepsilon < \infty$, then there exists a version of $X$ that is almost surely bounded and has almost surely uniformly $d_X$-continuous sample paths.
\end{lemma}

\begin{lemma}[Continuous sample paths and modulus of continuity]\label{lemma:ModulusContinuity}
	If $X$ is almost surely bounded on $U$, then the sample paths of $X$ on $U$ are almost surely uniformly $d_X$-continuous if and only if 
	\begin{align}\label{eq:lemma:ModulusContinuity-0}
		\lim_{\delta \rightarrow 0} \mathrm{E} \left[ \sup_{d_X(u, v) < \delta} |X_u - X_v| \right] = 0.
	\end{align}
\end{lemma}
\begin{remark}
	Sufficiency of~\eqref{eq:lemma:ModulusContinuity-0} holds for arbitrary stochastic processes on general metric spaces. Necessity of~\eqref{eq:lemma:ModulusContinuity-0} holds only for Gaussian processes. See comment in proof.
\end{remark}

\subsection{Auxiliary results for applications}\label{subsec:AuxiliaryResults-Applications}
In this section we collect several technical results needed for the applications in Section~\ref{sec:Applications}.

\begin{lemma}\label{lemma:Tesseract}
	Let $X, X_1, \ldots, X_n \in \mathbb{R}^d$ be i.i.d. sub-Gaussian random vectors with mean zero and covariance matrix $\Sigma$. Define $T_4: \mathbb{R}^d \times  \mathbb{R}^d \times  \mathbb{R}^d \times  \mathbb{R}^d \rightarrow \mathbb{R}$ by $T_4(t,u,v,w) := \mathrm{E}[ (X't) (X'u) (X'v) (X'w)]$. Then,
	\begin{align*}
		\mathrm{E}  \left\| \frac{1}{n}\sum_{i=1}^n X_i \otimes X_i \otimes X_i \otimes X_i - T_4 \right\|_{op}  \lesssim  \mathrm{r}(\Sigma) \|\Sigma\|_{op}^2  \left( \sqrt{\frac{(\log en)^2  \mathrm{r}(\Sigma)}{n}} \vee \frac{(\log en)^2  \mathrm{r}(\Sigma)}{n} \right),
	\end{align*}
	where $\mathrm{r}(\Sigma) = \mathrm{tr}(\Sigma)/ \|\Sigma\|_{op}$ and $\lesssim$ hides an absolute constant independent of $n$, $d$, and $\Sigma$. (Here, we tacitly identify the $X_i$'s with linear maps $\mathbb{R}^d \rightarrow \mathbb{R}$ and $\otimes$ denotes the tensor product between these linear maps.) 
\end{lemma}
\begin{remark}
	This result is useful because the upper bound is dimension-free in the sense that it only depends on the effective rank $\mathrm{r}(\Sigma)$ and the operator norm $\|\Sigma\|_{op}$. However, the dependence on the sample size $n$ is sub-optimal. 
\end{remark}

\begin{lemma}\label{lemma:BC-KRR-Expansion}
	The bias-corrected kernel ridge regression estimator defined in~\eqref{eq:subsec:Application-FDA-1} satisfies
	\begin{align*}
		\widehat{f}^{\mathrm{bc}}_n - f_0 = (T + \lambda)^{-2} T \left( \frac{1}{n} \sum_{i=1}^n \big(Y_i - f_0(X_i)\big) k_{X_i} \right) + R_n,
	\end{align*}
	where $R_n$ is a higher-order remainder term and
	\begin{align*}
		\sqrt{n}\|R_n\|_\infty &\lesssim  \kappa \left( \lambda^{-1} \|\widehat{T}_n - T\|_{op} + \lambda^{-2}\|\widehat{T}_n - T\|_{op}^2 \right)  \left\| \frac{1}{n}\sum_{i=1}^n (T + \lambda)^{-1} \varepsilon_i k_{X_i}\right\|_{\mathcal{H}}\\
		&\quad\quad + \kappa \left( \sqrt{n}\|\widehat{T}_n - T\|_{op}^2  + \sqrt{n}\lambda^2 \right) \| (T + \lambda)^{-2}f_0\|_{\mathcal{H}},
	\end{align*}
	where $\lesssim$ hides an absolute constant.
\end{lemma}

\begin{lemma}\label{lemma:BC-KRR-Remainder}
	Let $\delta \in (0,1)$ and $R_n$ be the remainder term in~\eqref{eq:subsec:Application-FDA-1}. Let $S$ be a separable metric space (w.r.t. some metric). If Assumptions~\ref{assumption:kernel} and~\ref{assumption:data} hold, then, with probability at least $1- \delta$,
	\begin{align*}
	\sqrt{n}\|R_n\|_\infty& \lesssim \left(\sqrt{\frac{ \bar{\sigma}^2  \mathfrak{n}_1^2(\lambda)  }{n\lambda^2 } } \vee \frac{ \bar{\sigma}^2 }{n \lambda^2}\right)\left(\frac{\kappa^4\log^3(2/\delta)}{\sqrt{n} \lambda} \vee  \kappa^2 \log^2(2/\delta) \right) \\
			& \quad\quad +  \left( \frac{\kappa^4\log^2(2/\delta)}{\sqrt{n}} \vee \sqrt{n} \lambda^2 \right) \kappa \| (T + \lambda)^{-2}f_0\|_{\mathcal{H}},
	\end{align*}
	where $\bar{\sigma}^2 \geq \sigma_0^2 \vee \kappa^2(B + \kappa\|f_0\|_{\mathcal{H}})^2 \vee 1$, $\mathfrak{n}_1^2(\lambda) = \mathrm{tr}\left((T+ \lambda)^{-2}T\right)$, and $\lesssim$ hides an absolute constant.
\end{lemma}
\begin{remark}
	The quantity $\mathfrak{n}_1(\lambda)$ also appears in~\cite{singh2023kernel}. For an interpretation and its relation to the effective rank of the operator $T$ we refer to Section H (ibid., pp. 80ff).
\end{remark}
\begin{remark}
	 Above upper bound is $o_p(1)$ if (i) $\bar{\sigma} \kappa^3 (\log n)^3 \vee \bar{\sigma} \kappa^2 (\log n)^2 \mathfrak{n}_1(\lambda)= o( \sqrt{n} \lambda)$,  (ii) $\kappa^5 (\log n)^2 \|(T + \lambda)^{-2}f_0\|_{\mathcal{H}} = o(\sqrt{n})$, and (iii) $\lambda^2 \kappa \|(T + \lambda)^{-2}f_0\|_{\mathcal{H}} = o(1)$.
\end{remark}

\begin{lemma}\label{lemma:BC-KRR-Consistency-Operator}
	Let $\delta \in (0,1)$ and $\Omega$ and $\widehat{\Omega}_n$ be the covariance operator and its sample analogue as define in Section~\ref{subsec:Application-FDA}. If Assumptions~\ref{assumption:kernel} and~\ref{assumption:data} hold, then, with probability at least $1- \delta$,
	\begin{align*}
		\|\widehat{\Omega}_n - \Omega\|_{op} \lesssim \|T^3(T + \lambda)^{-4}\|_{op} \sqrt{\frac{ (\kappa^4 +\kappa^2  \mathfrak{n}_1^2(\lambda) + \bar{\sigma}^4)\log(2/\delta)}{n\lambda^2}},
	\end{align*}
	where $\bar{\sigma}^2 \geq \sigma_0^2 \vee \kappa^2(B + \kappa\|f_0\|_{\mathcal{H}})^2 \vee 1$, $\mathfrak{n}_1^2(\lambda) = \mathrm{tr}\left((T+ \lambda)^{-2}T\right)$, and $\lesssim$ hides an absolute constant.
\end{lemma}

\begin{remark}
	Above upper bound is $o_p(1)$ if $\|T^3(T + \lambda)^{-4}\|_{op} (\kappa^2 \vee \kappa \mathfrak{n}_1(\lambda) \vee \bar{\sigma}^2 ) \sqrt{\log n} = o(\sqrt{n} \lambda)$.
\end{remark}

\begin{lemma}\label{lemma:Consistency-Gram-Matrix}
	Let $\delta \in (0,1)$ and $\alpha \in \mathbb{N}_0$. Let $S$ be a separable metric space (w.r.t. some metric). If Assumptions~\ref{assumption:kernel} and~\ref{assumption:data} hold, then, with probability at least $1- \delta$,
	\begin{align*}
		(i)&\quad \quad  \left\| \frac{1}{n} \sum_{i=1}^n (T + \lambda)^{-\alpha}\big(Y_i - f_0(X_i)\big) k_{X_i}\right\|_{\mathcal{H}}  \lesssim \sqrt{ \frac{\sigma_0^2 \mathfrak{n}_\alpha^2(\lambda) \log(2/\delta)}{n}} \vee \frac{ \lambda^{-\alpha}\kappa (B + \kappa\|f_0\|_{\mathcal{H}}) \log(2/\delta) }{n},\\
		(ii)&\quad \quad \left\|\frac{1}{n}\sum_{i=1}^n (T + \lambda)^{-\alpha} \left(( k_{X_i} \otimes k_{X_i}^*) - T \right)\right\|_{HS} \lesssim   \sqrt{ \frac{\kappa^2  \mathfrak{n}_\alpha^2(\lambda) \log(2/\delta)}{n}}  \vee \frac{ \lambda^{-\alpha}\kappa^2 \log(2/\delta) }{n},
	\end{align*}
	where $\bar{\sigma}^2 \geq \sigma_0^2 \vee \kappa^2(B + \kappa\|f_0\|_{\mathcal{H}})^2 \vee 1$, $\mathfrak{n}_\alpha^2(\lambda) = \mathrm{tr}\left((T+ \lambda)^{-2\alpha}T\right)$ and $\lesssim$ hides an absolute constant.
\end{lemma}
\begin{remark}
	If $\mathcal{H}$ is pre-Gaussian, then weaker conditions and generic chaining arguments yield tighter bounds~\cite[e.g.][Theorem 9]{koltchinskiiConcentration2017}. 
\end{remark}
The next lemma is a version of Bernstein's exponential tail bound for random elements on separable Hilbert spaces.

\begin{lemma}[Theorem 3.3.4,~\citeauthor{yurinsky2006sums},~\citeyear{yurinsky2006sums}; Lemma G.2,~\citeauthor{singh2023kernel},~\citeyear{singh2023kernel}]\label{lemma:Bernstein-RE}
	Let $X, X_1, \ldots, X_n$ be i.i.d. centered random elements on a separable Hilbert space $\mathcal{H}$ with induced norm $\| \cdot\|_{\mathcal{H}}$. Suppose that there exist absolute constants $\nu, \sigma > 0$ such that $\sum_{i=1}^n \mathrm{E}\|X_i\|_{\mathcal{H}}^k \leq (k!/2) \sigma^2 \nu^{k-2}$ for all $k \geq 2$. Then, for $t > 0$ arbitrary,
	\begin{align*}
		\mathbb{P} \left\{  \max_{1 \leq m \leq n} \left\| \sum_{i=1}^m X_i \right\|_{\mathcal{H}} > t \sigma \right\} \leq 2\exp \left( \frac{-t^2/2}{1 + t \nu/\sigma} \right).
	\end{align*}
	In particular, for $\delta \in (0,1)$ arbitrary, with probability at least $1-\delta$,
	\begin{align*}
		\left\| \frac{1}{n}\sum_{i=1}^n X_i \right\|_{\mathcal{H}} \lesssim \sqrt{\frac{\sigma^2 \log(2/\delta)}{n}}  \vee \frac{\nu \log(2/\delta)}{n},
	\end{align*}
	where $\lesssim$ hides an absolute constant independent of $\delta, n, \nu, \sigma$, and $\mathcal{H}$.
\end{lemma}

\newpage
\section{Proofs of the results in Section~\ref{sec:ResultsMaximaVec}}\label{sec:Proofs-ResultsMaximaVec}

\subsection{Proof of Proposition~\ref{lemma:CLT-Max-Norm}}
\begin{proof}[\textbf{Proof of Proposition~\ref{lemma:CLT-Max-Norm}}]
	Our proof is inspired by~\cite{nourdin2012normal} (Theorem 3.7.1) who establish a Berry-Esseen type bound for the uni-variate case using an inductive argument~\citep[attributed to][]{bolthausen1984estimate}. The multi-variate case requires several modifications some of which we take from~\cite{goetze1991RateConvMultCLT},~\cite{bhattacharyaExpositionGoetze2010}, and~\cite{fang2021HighDimCLT}. We also borrow a truncation argument from~\cite{chernozhukov2017CLTHighDim}, which explains the qualitative similarity between their bound and ours. Original ideas in our proof are mostly those related to the way in which we use our Gaussian anti-concentration inequality (Lemma~\ref{lemma:AntiConcentration-SeparableProcess}) and exploit the mollifying properties of the Ornstein-Uhlenbeck semi-group operator to by-pass dimension dependent smoothing inequalities (Lemmas~\ref{lemma:Smooth-Lipschitz-Approx} and~\ref{lemma:DifferentiatingUnderIntegral}).
	
	Our proof strategy has two drawbacks: First, the inductive argument relies substantially on the i.i.d. assumption of the data. Generalizing this argument to independent but non-identical data requires additional assumptions on the variances similar to the uniform asymptotic negligibility condition in the classical Lindeberg-Feller CLT. We leave this generalization to future research. Second, the recent results by~\cite{chernozhukov2021NearlyOptimalCLT} suggest that our Berry-Esseen type bound is not sharp. Unfortunately, their proof technique (based on delicate estimates of Hermite polynomials) is inherently dimension dependent. Extending their approach to the coordinate-free Wiener chaos decomposition is another interesting research task.
		
	\vspace{10pt}
	\noindent
	\textbf{The case of positive definite $\Sigma$.} 
	\vspace{10pt}
	
	\noindent
	Suppose that $\Sigma \in \mathbb{R}^{d \times d}$ is positive definite. Let $Z \sim N(0, \Sigma)$ be independent of $X, X_1, \ldots, X_n \in \mathbb{R}^d$ and define, for each $n \geq 1$,
	\begin{align*}
		\Delta_n : = \sup_{s, t \geq 0} \Big|\mathbb{P}\left\{\|e^{-t} S_n + \sqrt{1 -e^{-2t}} Z\|_\infty \leq s \right\} -\mathbb{P}\left\{\|Z\|_\infty \leq s \right\} \Big|.
	\end{align*}
	Further, for each $n \geq 1$, let $C_{n,d} \geq 1$ be the smallest constant greater than or equal to one such that, for all i.i.d. random variables $X, X_1, \ldots, X_n \in \mathbb{R}^d$ with $\mathrm{E}[\|X\|_\infty^3] < \infty$, $\mathrm{E}[X] = 0$, and $\mathrm{E}[XX'] = \Sigma$,
	\begin{align*}
		\Delta_n \leq C_{n,d} B_n,
	\end{align*}
	where
	\begin{align*}
		B_n := \frac{(\mathrm{E}[\|X \|_\infty^3])^{1/3}}{n^{1/6}\sqrt{\mathrm{Var}(\|Z\|_\infty)}} + \frac{\mathrm{E} \left[ \|X\|_\infty^3 \mathbf{1}\{\|X\|_\infty > M\}\right]}{\mathrm{E} \left[ \|X\|_\infty^3\right]} + \frac{ 12\mathrm{E}[\|Z\|_\infty] + M}{\sqrt{n\mathrm{Var}(\|Z\|_\infty)}} .
	\end{align*}
	The factor 12 in front of $\mathrm{E}[\|Z\|_\infty]$ ensures that $B_n \sqrt{n} \geq 1$ so that $C_{n,d} \leq \sqrt{n}$. Indeed, one easily computes $(\mathrm{E}[\|Z\|_\infty^2])^{1/2} \leq 2(2 \sqrt{\pi} + 1) \mathrm{E}[\|Z\|_\infty] \leq 12 \mathrm{E}[\|Z\|_\infty]$~\citep[by the equivalence of moments of suprema of Gaussian processes, e.g.][Proposition A.2.4]{vandervaart1996weak} and, hence, $B_n \sqrt{n} \geq \sqrt{\mathrm{E}[\|Z\|_\infty^2]/\mathrm{Var}(\|Z\|_\infty)} \geq 1$.
	
	While the upper bound $C_{n,d} \leq \sqrt{n}$ is too loose to conclude the proof, it is nonetheless an important first step towards a tighter bound. For the moment, assume that there exists an absolute constant $K \geq 1$ such that
	\begin{align}\label{eq:theorem:CLT-Max-Norm-1}
		\Delta_n \leq B_n \left[ \left(1 + 2K^2\right)\sqrt{C_{n-1,d}} + 1 \right] \quad{} \forall n \geq 2.
	\end{align}
	Then, by construction of $C_{n,d}$,
	\begin{align}\label{eq:theorem:CLT-Max-Norm-2}
		C_{1, d} = 1, \quad{} C_{n,d} \leq \left[\left(1 + 2K^2\right)\sqrt{C_{n-1,d}} + 1\right] \wedge \sqrt{n} \quad{}\forall n \geq 2.
	\end{align}
	We shall now show that this difference inequality implies that $\sup_{d \geq 1}\sup_{n \geq 1} C_{n,d} < \infty$ independent of the distribution of the $X_i$'s: Define the map $x \mapsto F(x) := (1 + 2 K^2) \sqrt{x} + 1$ and consider the nonlinear first-order difference equation
	\begin{align*}
		x_1 = 1, \quad{} x_n = F(x_{n-1}) \quad{} \forall n \geq 2.
	\end{align*}
	We easily verify that the fixed point $x^* > 0$ solving $x = F(x)$ satisfies
	\begin{align*}
		x^* = \frac{1}{2} \left( 4K^4 + 4 K^2 + \sqrt{(2K^2 +1)^2 (4K^4 + 4 K^2 + 5)} + 3\right),
	\end{align*}
	and that $F(x) > x$ for all $x \in (0, x^*)$ and $F(x) < x$ for all $(x^*, \infty)$. We also notice that $F$ is monotone increasing on $\mathbb{R}_+$. Thus,
	\begin{align*}
		[(F \circ F)(x) - x](x - x^*) < 0 \quad{}\quad{} \forall x \in \mathbb{R}_+ \setminus \{x^*\}.
	\end{align*}
	Hence, by Theorem 2.1.2 in~\cite{sedaghat2003nonlinear} every trajectory $\{F^n(x_1)\}_{n \geq 1}$ with $x_1 \in \mathbb{R}_+$ converges to $x^*$. In particular, $\lim_{n \rightarrow \infty} x_n = \lim_{n \rightarrow \infty} F^n(1) = x^*$. Returning to the inequality~\eqref{eq:theorem:CLT-Max-Norm-2} we conclude that there exists $N_0 \geq 2$ such that $C_{n,d} \leq x^* + 1$ for all $n > N_0$ and all $d \geq 1$. Since $C_{n,d} \leq \sqrt{n}$ for all $n \geq 1$, it follows that $C_{n,d} \leq (x^* + 1) \vee \sqrt{N_0} < \infty$ for all $n \geq 1$ and $d \geq 1$.
	
	To complete the proof of the theorem, it remains to show that eq.~\eqref{eq:theorem:CLT-Max-Norm-1} holds. Let $s \in \mathbb{R}$, $t, \lambda \geq 0$ be arbitrary and $h_{s, \lambda}$ be the map from Lemma~\ref{lemma:Smooth-Lipschitz-Approx}. Define $x \mapsto h(x) := h_{s, \lambda}(\|x\|_\infty)$.
	By Lemma~\ref{lemma:AntiConcentration-SeparableProcess} and Lemma~\ref{lemma:Smooth-Lipschitz-Approx} (ii),
	\begin{align}\label{eq:theorem:CLT-Max-Norm-3}
		\Delta_n \leq \sup_{s \in \mathbb{R}, t \geq 0} \big|\mathrm{E}[P_th(S_n) - h(Z)] \big| + \frac{2 \sqrt{3}\lambda }{ \sqrt{\mathrm{Var}(\|Z\|_\infty) + \lambda^2/12}},
	\end{align} 
	where $P_t h$ denotes the Ornstein-Uhlenbeck semi-group with stationary measure $N(0, \Sigma)$, i.e.
	\begin{align*}
		P_th(x) := \mathrm{E}\left[h\left( e^{-t}x + \sqrt{1- e^{-2t}}Z \right) \right] \quad{}\quad{} \forall x \in \mathbb{R}^d.
	\end{align*}
	Since $x \mapsto P_th(x) - \mathrm{E}[h(Z)]$ is Lipschitz continuous (with constant $\lambda^{-1} e^{-t}$) and $\Sigma$ positive definite, Proposition 4.3.2 in~\cite{nourdin2012normal} implies that
	\begin{align}\label{eq:theorem:CLT-Max-Norm-4}
		\mathrm{E}\left[P_th(S_n) - h(Z)\right] =  \mathrm{E}\left[\mathrm{tr}\left(\Sigma D^2 G_h(S_n) \right) - S_n' D G_h(S_n) \right],
	\end{align}
	where $G_h \in C^\infty(\mathbb{R}^d)$ and
	\begin{align*}
		G_h(x) := \int_0^\infty \Big(\mathrm{E}[h(Z)]  - P_uP_th(x)  \Big) du \quad{}\quad{} \forall x \in \mathbb{R}^d. 
	\end{align*}
	Since $P_u P_t f = P_{u +t} f$ almost surly for all integrable maps $f$ (semi-group property!), we have
	\begin{align*}
		G_h(x) = \int_t^\infty \mathrm{E}\left[h(Z) - h\left( e^{-u}x + \sqrt{1- e^{-2u}}Z \right) \right]du.
	\end{align*} 
	We proceed by re-writing eq.~\eqref{eq:theorem:CLT-Max-Norm-4} in multi-index notation as
	\begin{align*}
		\big|\mathrm{E}[P_th(S_n) - h(Z)] \big|&= \left| \mathrm{E}\left[\mathrm{tr}\left(\Sigma D^2 G_h(S_n) \right) - S_n' D G_h(S_n) \right]\right| \nonumber\\
		&= \left|\sum_{i=1}^n  \mathrm{E}\left[  \frac{1}{n} \mathrm{tr}\left(\widetilde{X}_i \widetilde{X}_iD^2 G_h(S_n) \right)  - \frac{X_i'}{\sqrt{n}} D G_h(S_n) \right]\right| \nonumber\\
		&= \left|\sum_{i=1}^n \mathrm{E}\left[ \sum_{|\alpha| = 2}  D^\alpha G_h(S_n) \left(\frac{\widetilde{X}_i}{\sqrt{n}}\right)^\alpha  - \sum_{|\alpha|=1} D^\alpha G_h\left(S_n\right)\left(\frac{X_i}{\sqrt{n}}\right)^\alpha  \right] \right|,
	\end{align*}
	where the $\widetilde{X}_i$'s are independent copies of the $X_i$'s. A Taylor expansion around $S_n^i : = S_n - n^{-1/2} X_i$ with exact integral remainder term yields
	\begin{align}\label{eq:theorem:CLT-Max-Norm-5}
		&\big|\mathrm{E}[P_th(S_n) - h(Z)] \big| \nonumber\\
		&\quad{}= \left| \mathrm{E}\left[ \sum_{i=1}^n\sum_{|\alpha| = 2}  D^\alpha G_h(S_n^i) \left(\frac{\widetilde{X}_i}{\sqrt{n}}\right)^\alpha  \right]  + \mathrm{E}\left[\sum_{i=1}^n \sum_{|\alpha| = 2}\sum_{|\beta| = 1} D^{\alpha + \beta} G_h\left(S_n^i + \theta \frac{X_i}{\sqrt{n}}\right) \left(\frac{\widetilde{X}_i}{\sqrt{n}}\right)^\alpha \left(\frac{ X_i}{\sqrt{n}}\right)^\beta \right] \right.\nonumber\\
		&\quad{} \left. \quad{}- \mathrm{E}\left[\sum_{i=1}^n\sum_{|\alpha|=1} D^\alpha G_h\left(S_n^i\right)\left(\frac{X_i}{\sqrt{n}}\right)^\alpha  \right] - \mathrm{E}\left[ \sum_{i=1}^n\sum_{|\alpha| = 1} \sum_{|\beta| = 1}   D^{\alpha + \beta} G_h(S_n^i) \left(\frac{X_i}{\sqrt{n}}\right)^{\alpha + \beta}  \right]  \right.\nonumber\\
		& \quad{}\left. \quad{}  -  \mathrm{E}\left[\sum_{i=1}^n \sum_{|\alpha| = 1} \sum_{|\beta| = 2}  D^{\alpha + \beta} G_h\left(S_n^i + \theta \frac{X_i}{\sqrt{n}}\right) \left(\frac{X_i}{\sqrt{n}}\right)^{\alpha + \beta} \frac{2(1 - \theta)}{\beta!} \right]  \right|\nonumber\\
		&\quad{}= \left|  \mathrm{E}\left[\sum_{i=1}^n \sum_{|\alpha| = 2} \sum_{|\beta|=1} D^{\alpha + \beta} G_h\left(S_n^i + \theta \frac{X_i}{\sqrt{n}}\right) \left\{ \left(\frac{\widetilde{X}_i}{\sqrt{n}}\right)^\alpha \left(\frac{ X_i}{\sqrt{n}}\right)^\beta - \left(\frac{X_i}{\sqrt{n}}\right)^{\alpha + \beta} \frac{2(1 -\theta)}{\alpha!}  \right\} \right]   \right|,
	\end{align}
	where $\theta \sim Unif(0,1)$ is independent of the $X_i$'s and $\widetilde{X}_i$'s. The first and fourth term cancel out and the third term vanishes because $S_n^i$ and $X_i$ are independent and mean zero. (Eq.~\eqref{eq:theorem:CLT-Max-Norm-5} is essentially a re-statement of Lemmas 2.9 and 2.4 in~\cite{goetze1991RateConvMultCLT} and~\cite{raic2019MultBerryEsseen}, respectively.)
	
	Notice that $G_h \in C^\infty(\mathbb{R}^d)$ because it is the convolution of a bounded Lipschitz map with the density of $N(0, \Sigma)$~\citep[e.g.][Proposition 4.2.2]{nourdin2012normal}. Derivatives on $G_h$ are usually obtained by differentiating the density of $N(0,\Sigma)$ (e.g.~\citeauthor{goetze1991RateConvMultCLT},~\citeyear{goetze1991RateConvMultCLT}, Lemma 2.1;~\citeauthor{raic2019MultBerryEsseen}, \citeyear{raic2019MultBerryEsseen}, Lemma 2.5 and 2.6;~\citeauthor{fang2021HighDimCLT},~\citeyear{fang2021HighDimCLT}, Lemma 2.2; and~\citeauthor{chernozhukov2021NearlyOptimalCLT},~\citeyear{chernozhukov2021NearlyOptimalCLT}, Lemmas 6.1, 6.2, 6.3). Here, we proceed differently. Let $1 \leq j, k, \ell \leq d$ be the indices corresponding to the multi-indices $|\alpha| = 2$ and $|\beta| = 1$. By Lemma~\ref{lemma:DifferentiatingUnderIntegral} (i) we have
	\begin{align*}
		D^{\alpha + \beta} G_h\left(S_n^i + \theta \frac{X_i}{\sqrt{n}}\right) &= \frac{\partial^3 G_h}{\partial x_j \partial x_k \partial x_\ell} \left(S_n^i + \theta \frac{X_i}{\sqrt{n}}\right) \\
		& = - \int_t^\infty e^{-3u} \mathrm{E}_Z\left[\frac{\partial^3 h}{\partial x_j \partial x_k \partial x_\ell} \left( e^{-u} x + \sqrt{1- e^{-2u}}Z \right) \Big|_{ x = S_n^i + \theta \frac{X_i}{\sqrt{n}} } \right]du,
	\end{align*}
	where $\mathrm{E}_Z$ denotes the expectation taken with respect to the law of $Z$ only. And by Lemma~\ref{lemma:DifferentiatingUnderIntegral} (ii),
	\begin{align*}
		&\sum_{j,k, \ell}  e^{-3u}\mathrm{E}_Z\left[\left|\frac{\partial^3 h}{\partial x_j \partial x_k \partial x_\ell} \left( e^{-u} x + \sqrt{1- e^{-2u}}Z \right) \Big|_{ x = S_n^i + \theta \frac{X_i}{\sqrt{n}} } \right|\right]\\
		&\quad{} \leq C_3\lambda^{-3}e^{-3u}\sum_{j=1}^d	\mathrm{E}_Z\left[\mathbf{1}_{[s, s+ 3\lambda]} \left(\left\|V_i\right\|_\infty \right) \mathbf{1}\left\{\left|V_{ij}\right| \geq \left| V_{im}\right|, \: \forall m \right\} \right],
	\end{align*}
	where
	\begin{align*}
		V_i :=	e^{-u} \left(S_n^i + \theta \frac{X_i}{\sqrt{n}}\right) + \sqrt{1- e^{-2u}}Z, \quad\quad 1 \leq i \leq n.
	\end{align*}
	Let $\alpha_{k, \ell} \in \{1, 2\}$ be the value of $\alpha!$ for $|\alpha|=2$ corresponding to the indices $k, \ell$. An application of H{\"o}lder's inequality yields
	\begin{align}\label{eq:theorem:CLT-Max-Norm-6}
		&\frac{ e^{-3u}}{n^{3/2}}\sum_{i=1}^n\sum_{j,k, \ell} \mathrm{E}_Z\left[\left|\frac{\partial^3 h}{\partial x_j \partial x_k \partial x_\ell} \left( e^{-u} x + \sqrt{1- e^{-2u}}Z \right) \Big|_{ x = S_n^i + \theta \frac{X_i}{\sqrt{n}} }\right|\right]\left| \widetilde{X}_{ij} \widetilde{X}_{ik} X_{i\ell} -  X_{ij}X_{ik}X_{i\ell}\frac{2(1-\theta)}{\alpha_{k, \ell}}\right| \nonumber\\
		&\leq \frac{C_3e^{-3u}}{n^{3/2}\lambda^3} \mathrm{E}_Z\left[\sum_{i=1}^n \sum_{j=1}^d \mathbf{1}\left\{\left|V_{ij}\right| \geq \left| V_{im}\right|,\: m \neq j\right\}\mathbf{1}_{[s, s+ 3\lambda]} \left(\left\|V_i\right\|_\infty \right) \left|\widetilde{X}_{ij}^2 X_{ij} -  X_{ij}^3\frac{2(1-\theta)}{\alpha_{j, j}}\right| \right]\nonumber\\
		&\leq \frac{C_3e^{-3u}}{n^{3/2}\lambda^3} \mathrm{E}_Z\left[\sum_{i=1}^n  \mathbf{1}_{[s, s+ 3\lambda]} \left(\left\|V_i\right\|_\infty \right)\left(\sum_{j=1}^d\mathbf{1} \left\{\left|V_{ij}\right| \geq \left| V_{im}\right|,\: m \neq j\right\}  \right) \max_{1 \leq j \leq d}\left( \widetilde{X}_{ij}^2| X_{ij}| + | X_{ij}|^3\right) \right]\nonumber\\
		&= \frac{C_3e^{-3u}}{n^{3/2}\lambda^3} \mathrm{E}_Z\left[\sum_{i=1}^n  \mathbf{1}_{[s, s+ 3\lambda ]} \left(\left\|V_i\right\|_\infty \right)  \left(\|\widetilde{X}_i \|_\infty^2 \| X_i \|_\infty + \|X_i \|_\infty^3 \right) \right],
	\end{align}
	where in the last line we have used that $ \sum_{j=1}^d\mathbf{1} \left\{\left|V_{ij}\right| \geq \left| V_{im}\right|,\: \forall m \right\} = 1$ almost surely because $|Corr(Z_j, Z_k)| = |Corr(X_{1j}, X_{1k})| < 1$ for all $j \neq k$ and $1 \leq i \leq n$ (since $\Sigma$ is positive definite no pair of entries in $X_i$ and $Z$ can be perfectly (positively or negatively) correlated!). Taking expectation with respect to the $X_i$'s and $\theta$ over above inequality, we obtain
	\begin{align}\label{eq:theorem:CLT-Max-Norm-7}
		& \frac{C_3e^{-3u}}{n^{3/2}\lambda^3}\mathrm{E}\left[\sum_{i=1}^n  \mathbf{1}_{[s, s+ 3\lambda]} \left(\left\|V_i\right\|_\infty \right)  \left(\|\widetilde{X}_i \|_\infty^2 \| X_i \|_\infty + \|X_i \|_\infty^3 \right) \right] \nonumber\\
		&\quad{} = \frac{C_3e^{-3u}}{n^{3/2}\lambda^3} \mathrm{E}\left[\sum_{i=1}^n \mathrm{E}\left[ \mathbf{1}_{[s, s+ 3\lambda]} \left(\left\|V_i\right\|_\infty \right) \mid X_i, \theta \right]  \left(\|\widetilde{X}_i \|_\infty^2 \| X_i \|_\infty + \|X_i \|_\infty^3 \right) \right] \nonumber\\
		&\quad{} =  \frac{C_3e^{-3u}}{n^{3/2}\lambda^3}\mathrm{E}\left[\sum_{i=1}^n \mathbb{P}\left\{  s \leq \|V_i\|_\infty \leq s + 3\lambda \mid X_i, \theta \right\} \left(\|\widetilde{X}_i \|_\infty^2 \| X_i \|_\infty + \|X_i \|_\infty^3 \right) \right].
	\end{align}
	Notice that $\mathrm{E}[S_n^i] = 0$ and $\mathrm{E}[S_n^iS_n^{i'}] =\frac{n-1}{n} \mathrm{E}[S_n S_n']$, and $Z \overset{d}{=}\sqrt{\frac{n-1}{n}}Z + \frac{1}{\sqrt{n}}\widetilde{Z}$, where $\widetilde{Z}$ is an independent copy of $Z \sim N(0, \Sigma)$. Set $Z_n :=\sqrt{\frac{n-1}{n}}Z$ and bound the probability in line~\eqref{eq:theorem:CLT-Max-Norm-7} by
	\begin{align}\label{eq:theorem:CLT-Max-Norm-8}
		&\mathbb{P}\left\{  s \leq \|V_i\|_\infty \leq s + 3\lambda \mid X_i, \theta \right\} \nonumber\\
		&\quad{}=  \mathbb{P}\left\{ \|V_i\|_\infty  > s \mid X_i, \theta \right\} - \mathbb{P}\left\{\|V_i\|_\infty  > s + 3\lambda \mid X_i, \theta \right\}  \nonumber\\
		&\quad{}\leq \mathbb{P}\left\{ \| e^{-u} S_n^i+ \sqrt{1- e^{-2u}} Z_n \|_\infty  + \sqrt{1 - e^{-2u}} n^{-1/2} \|\widetilde{Z}\|_\infty +  e^{-u} n^{-1/2} \|\theta  X_i\|_\infty  > s \mid X_i, \theta \right\} \nonumber\\
		&\quad{}\quad{} - \mathbb{P}\left\{ \| e^{-u} S_n^i+ \sqrt{1- e^{-2u}} Z_n \|_\infty - \sqrt{1 - e^{-2u}} n^{-1/2} \|\widetilde{Z}\|_\infty - e^{-u} n^{-1/2} \| \theta X_i\|_\infty  > s + 3\lambda \mid X_i, \theta \right\} \nonumber\\
		&\quad{}= \mathbb{P}\left\{ \| e^{-u} S_n^i+ \sqrt{1- e^{-2u}} Z_n\|_\infty + \sqrt{1 - e^{-2u}} n^{-1/2} \|\widetilde{Z}\|_\infty + e^{-u} n^{-1/2} \|\theta  X_i\|_\infty  > s \mid X_i, \theta \right\} \nonumber\\
		&\quad{}\quad{} \pm \mathbb{P}\left\{ \|  Z_n\|_\infty + \sqrt{1 - e^{-2u}} n^{-1/2} \|\widetilde{Z}\|_\infty +  e^{-u} n^{-1/2} \|\theta  X_i\|_\infty  > s \mid X_i, \theta \right\} \nonumber\\
		&\quad{}\quad{} \pm \mathbb{P}\left\{ \| Z_n\|_\infty - \sqrt{1 - e^{-2u}} n^{-1/2} \|\widetilde{Z}\|_\infty -  e^{-u} n^{-1/2} \|\theta  X_i\|_\infty  > s + 3 \lambda \mid X_i, \theta \right\} \nonumber\\
		&\quad{}\quad{} - \mathbb{P}\left\{ \| e^{-u} S_n^i+ \sqrt{1- e^{-2u}} Z_n\|_\infty - \sqrt{1 - e^{-2u}} n^{-1/2} \|\widetilde{Z}\|_\infty -  e^{-u} n^{-1/2} \|\theta  X_i\|_\infty  > s + 3\lambda  \mid X_i, \theta \right\}  \nonumber\\
		\begin{split}
			&\quad{} \leq 2 \Delta_{n-1} + \mathbb{P}\left\{ s - \sqrt{1 - e^{-2u}} n^{-1/2} \|\widetilde{Z}\|_\infty - e^{-u} n^{-1/2} \| \theta  X_i\|_\infty \leq \|Z_n\|_\infty \leq   s + 3\lambda \right.\\
			& \left. \quad{}\quad{}\quad{}\quad{}\quad{}\quad{}\quad{}\quad{}\quad{}\quad{}\quad{}\quad{}\quad{}\quad{}\quad{} + \sqrt{1 - e^{-2u}} n^{-1/2} \|\widetilde{Z}\|_\infty + e^{-u} n^{-1/2} \| \theta  X_i\|_\infty \mid X_i, \theta \right\},
		\end{split}
	\end{align}
	where we have used that under the i.i.d. assumption
	\begin{align*}
		\Delta_{n-1} 
		&\equiv \sup_{s \in \mathbb{R}, t \geq 0} \Big|\mathbb{P}\left\{\| e^{-t} S_n^i+ \sqrt{1- e^{-2t}} Z_n\|_\infty \leq s \right\} -\mathbb{P}\left\{\|Z_n\|_\infty \leq s \right\} \Big|.
	\end{align*}
	By Lemma~\ref{lemma:AntiConcentration-SeparableProcess}, Lemma~\ref{lemma:Smooth-Lipschitz-Approx} (ii), and  monotonicity and concavity of the map $x \mapsto x/\sqrt{a + x^2}$, $a > 0$, we have, for arbitrary $M > 0$,
	\begin{align}\label{eq:theorem:CLT-Max-Norm-9}
		&\mathbb{P}\left\{ s - \sqrt{1 - e^{-2u}} n^{-1/2} \|\widetilde{Z}\|_\infty - e^{-u} n^{-1/2} \| \theta  X_i\|_\infty \leq \|Z_n\|_\infty \leq   s + 3\lambda \right. \nonumber\\
		& \left. \quad{}\quad{}\quad{}\quad{}\quad{}\quad{}\quad{}\quad{}\quad{}\quad{}\quad{} + \sqrt{1 - e^{-2u}} n^{-1/2} \|\widetilde{Z}\|_\infty + e^{-u} n^{-1/2} \| \theta  X_i\|_\infty \mid X_i, \theta \right\} \nonumber\\
		&\quad{}\leq \frac{6 \sqrt{3} \lambda + 2\sqrt{3}\sqrt{1 - e^{-2u}} n^{-1/2}\mathrm{E}[\|\widetilde{Z}\|_\infty] + 4\sqrt{3} e^{-u} n^{-1/2} \|\theta X_i\|_\infty}{\sqrt{\mathrm{Var}(\|Z\|_\infty) + (1 - e^{-2u}) n^{-1}\mathrm{E}[\|\widetilde{Z}\|_\infty]^2/12 + 3\lambda^2/4 + e^{-2u}n^{-1}\|\theta X_i\|_\infty^2/3 }}\nonumber\\
		\begin{split}
			&\quad{} \leq  \frac{ 4\sqrt{3} n^{-1/2}M}{\sqrt{\mathrm{Var}(\|Z\|_\infty) + n^{-1} M^2/3 }} \mathbf{1}\{\| X_i\|_\infty \leq M\} + 12 \cdot \mathbf{1}\{\|X_i\|_\infty > M\} \\
			&\quad{} \quad{} +\frac{6 \sqrt{3} \lambda}{\sqrt{\mathrm{Var}(\|Z\|_\infty) + 3\lambda^2/4 }}  +  \frac{2\sqrt{3}n^{-1/2}\mathrm{E}[\|Z\|_\infty] }{\sqrt{\mathrm{Var}(\|Z\|_\infty) + n^{-1}\mathrm{E}[\|Z\|_\infty]^2/12}}.
		\end{split}
	\end{align}	
	Combine eq.~\eqref{eq:theorem:CLT-Max-Norm-8}--\eqref{eq:theorem:CLT-Max-Norm-9} with~\eqref{eq:theorem:CLT-Max-Norm-7} and integrate over $u \in (t, \infty)$ to conclude via eq.~\eqref{eq:theorem:CLT-Max-Norm-3}--\eqref{eq:theorem:CLT-Max-Norm-6} and the i.i.d. assumption that there exists an absolute constant $K \geq 1$ such that
	\begin{align}
		\begin{split}\label{eq:theorem:CLT-Max-Norm-10}
			\Delta_n &\leq  \frac{K \lambda}{\sqrt{\mathrm{Var}(\|Z\|_\infty)}} + \frac{K}{n^{1/2}\lambda^3}\mathrm{E}\left[\|X \|_\infty^3\right]\Delta_{n-1}\\
			&\quad{}+ \frac{K}{n^{1/2}\lambda^2}\frac{\mathrm{E}[\|X \|_\infty^3]}{\sqrt{\mathrm{Var}(\|Z\|_\infty)}} + \frac{K}{n\lambda^3} \frac{\mathrm{E}[\|X\|_\infty^3] \mathrm{E}[\|Z\|_\infty]}{\sqrt{\mathrm{Var}(\|Z\|_\infty)}} \\
			&\quad{}+\frac{K}{n\lambda^3} \frac{\mathrm{E}[\|X\|_\infty^3] M}{\sqrt{\mathrm{Var}(\|Z\|_\infty)}}  +\frac{K}{n^{1/2}\lambda^3} \mathrm{E} \left[ \|X\|_\infty^3 \mathbf{1}\{\|X\|_\infty > M\}\right],
		\end{split}
	\end{align}
	where we have used Harris' association inequality to simplify several summands, i.e.
	\begin{align*}
		\sum_{i=1}^n \mathrm{E}[\|\widetilde{X}_i\|_\infty^2] \mathrm{E}[\|X_i\|_\infty \mathbf{1}\{\|X_i\|_\infty \leq M\}] \leq \sum_{i=1}^n \mathrm{E}[\|X_i\|_\infty^2] \mathrm{E}[\|X_i\|_\infty] \leq \sum_{i=1}^n \mathrm{E}\left[ \|X_i\|_\infty^3\right],
	\end{align*}
	and
	\begin{align*}
		\sum_{i=1}^n \mathrm{E} [ \|\widetilde{X}_i\|_\infty^2] \mathrm{E} \left[ \|X_i\|_\infty \mathbf{1}\{\|X_i\|_\infty > M\}\right] \leq \sum_{i=1}^n \mathrm{E} \left[ \|X_i\|_\infty^3 	\mathbf{1}\{\|X_i\|_\infty > M\}\right].
	\end{align*}
	Observe that eq.~\eqref{eq:theorem:CLT-Max-Norm-10} holds for arbitrary $\lambda > 0$. Setting 
	\begin{align*}
		\lambda = K \left(\frac{ C_{n-1,d}}{n- 1} \right)^{1/6} \left(\mathrm{E}[\|X\|_\infty^3] \right)^{1/3}
	\end{align*}
	we deduce from eq.~\eqref{eq:theorem:CLT-Max-Norm-10}, the definition of $C_{n-1,d}$, and $K C_{n-1,d}^{1/6} \geq 1$ (because $K, C_{n-1,d} \geq 1$!) that, for all $n \geq 2$,
	\begin{align*}
		\Delta_n &\leq  \frac{K \lambda}{\sqrt{\mathrm{Var}(\|Z\|_\infty)}}  + B_{n-1} \sqrt{C_{n-1,d}} + B_n \nonumber \\
		&\leq  B_n \left[\left(\frac{n}{n-1}\right)^{1/6} K^2 C_{n-1,d}^{1/6} +\left(\frac{n}{n-1}\right)^{1/2} \sqrt{C_{n-1,d}} + 1 \right] \nonumber\\
		& \leq B_n \left[ \left(1 + 2K^2\right)\sqrt{C_{n-1,d}} + 1 \right].
	\end{align*}
	We have thus established eq.~\eqref{eq:theorem:CLT-Max-Norm-1}. This concludes the proof of the theorem in the case of a positive definite covariance matrix.
			
	\vspace{10pt}
	\noindent
	\textbf{The case of positive semi-definite $\Sigma \neq \mathbf{0}$.} 
	\vspace{10pt}
	
	\noindent
	Suppose that $\Sigma \in \mathbb{R}^{d \times d}$ is positive semi-definite but not identical to zero.
	
	Take $Y, Y_1, \ldots, Y_n \sim N(0, I_d)$ and $Z \sim N(0, \Sigma)$ such that $X, X_1, \ldots, X_n, Y, Y_1, \ldots, Y_n, Z \in \mathbb{R}^d$ are mutually independent. Let $\eta > 0$ be arbitrary and define $Z^\eta := Z + \eta Y$, $X^\eta := X + \eta Y$, and $S_n^\eta := n^{-1/2} \sum_{i=1}^n X_i^\eta$ with $X_i^\eta := X_i + \eta Y_i$. Clearly, $Z^\eta \sim N(0, \Sigma + \eta^2 I_d)$ and the $X^\eta_i$'s are i.i.d. with mean zero and positive definite covariance $\Sigma + \eta^2 I_d$. Hence, by the first part of the proof there exists an absolute constant $C_* \geq 1$ independent of $n$, $d$, and the distribution of the $X_i^\eta$'s (and hence, independent of $\eta > 0$!) such that for $M \geq 0$ and $n \geq 1$,
	\begin{align}\label{eq:theorem:CLT-Max-Norm-11}
	\Delta_n^\eta \leq C_* B_n^\eta,
	\end{align}
	where
	\begin{align*}
		\Delta_n^\eta := \sup_{s, t \geq 0} \Big|\mathbb{P}\left\{\|e^{-t} S_n^\eta + \sqrt{1 -e^{-2t}} Z^\eta\|_\infty \leq s \right\} -\mathbb{P}\left\{\|Z^\eta\|_\infty \leq s \right\} \Big|,
	\end{align*}
	and
	\begin{align*}
		B_n^\eta := \frac{(\mathrm{E}[\|X^\eta \|_\infty^3])^{1/3}}{n^{1/6}\sqrt{\mathrm{Var}(\|Z^\eta\|_\infty)}} + \frac{\mathrm{E} \left[ \|X^\eta\|_\infty^3 \mathbf{1}\{\|X^\eta\|_\infty > M\}\right]}{\mathrm{E} \left[ \|X^\eta\|_\infty^3\right]} + \frac{ 12\mathrm{E}[\|Z^\eta\|_\infty] + M}{\sqrt{n\mathrm{Var}(\|Z^\eta\|_\infty)}} .
	\end{align*}
	At this point, it is tempting to take $\eta \rightarrow 0$ in eq.~\eqref{eq:theorem:CLT-Max-Norm-11}. However, this would yield the desired result only for the case in which the law of $\|S_n\|_\infty$ is continuous. (Alternatively, we could replace the supremum over $s \in [0, \infty]$ by the supremum over $s \in \mathcal{C}_n$, where $\mathcal{C}_n$ is the set of continuity points of the law of $\|S_n\|_\infty$.) Therefore, we proceed differently. Recall eq.~\eqref{eq:theorem:CLT-Max-Norm-3}, i.e.
	\begin{align*}
		\Delta_n \leq \sup_{s \in \mathbb{R}, t \geq 0} \big|\mathrm{E}[P_th(S_n) - h(Z)] \big| + \frac{2 \sqrt{3}\lambda }{ \sqrt{\mathrm{Var}(\|Z\|_\infty) + \lambda^2/12}},
	\end{align*} 
	where $P_t h$ denotes the Ornstein-Uhlenbeck semi-group with stationary measure $N(0, \Sigma)$,
	\begin{align*}
		P_th(x) := \mathrm{E}\left[h\left( e^{-t}x + \sqrt{1- e^{-2t}}Z \right) \right] \quad{}\quad{} \forall x \in \mathbb{R}^d.
	\end{align*}
	Let $P_t^\eta h(x)$ be the  Ornstein-Uhlenbeck semi-group with stationary measure $N(0, \Sigma + \eta^2 I_d)$ and expand above inequality to obtain the following modified version of eq.~\eqref{eq:theorem:CLT-Max-Norm-3}: 
	\begin{align}
		\Delta_n &\leq \sup_{s \in \mathbb{R}, t \geq 0} \big|\mathrm{E}[P_t^\eta h(S_n^\eta) - h(Z^\eta)] \big| + \frac{2 \sqrt{3}\lambda }{ \sqrt{\mathrm{Var}(\|Z\|_\infty) + \lambda^2/12}}  + \sup_{s \in \mathbb{R}, t \geq 0} \big|\mathrm{E}[P_th(S_n^\eta) - P_th(S_n)] \big| \nonumber\\
		&\quad{} + \sup_{s \in \mathbb{R}, t \geq 0} \big|\mathrm{E}[h(Z^\eta) - h(Z)] \big| + \sup_{s \in \mathbb{R}, t \geq 0} \big|\mathrm{E}[P_t h(S_n^\eta) - P_t^\eta h(S_n^\eta)] \big| \nonumber\\
		&\overset{(a)}{\leq} \sup_{s \in \mathbb{R}, t \geq 0} \big|\mathrm{E}[P_t^\eta h(S_n^\eta) - h(Z^\eta)] \big| + \frac{2 \sqrt{3}\lambda }{ \sqrt{\mathrm{Var}(\|Z\|_\infty) + \lambda^2/12}}  \nonumber\\
		&\quad{} + \frac{\eta}{\lambda} \mathrm{E}\left[\big\|n^{-1/2} \sum_{i=1}^nY_i\big\|_\infty \right] + \frac{\eta}{\lambda} \mathrm{E}[\|Y\|_\infty] + \frac{\eta}{\lambda} \mathrm{E}[\|Y\|_\infty] \nonumber\\
		&\leq \sup_{s \in \mathbb{R}, t \geq 0} \big|\mathrm{E}[P_t^\eta h(S_n^\eta) - h(Z^\eta)] \big| + \frac{2 \sqrt{3}\lambda }{ \sqrt{\mathrm{Var}(\|Z\|_\infty) + \lambda^2/12}} + \frac{3\eta}{\lambda}\sqrt{2\log 2d}, \nonumber
	\end{align}
	where the (a) holds because $h$ is $\lambda^{-1}$-Lipschitz and $P_th$ is $\lambda^{-1} e^{-t}$-Lipschitz w.r.t. the metric induced by the $\ell_\infty$-norm.
	
	Since $\Sigma + \eta^2 I_d$ is positive definite we can proceed as in the first part of the proof and arrive at the following modified version of eq.~\eqref{eq:theorem:CLT-Max-Norm-10}: There exists an absolute constant $K \geq 1$ such that
	\begin{align}
		\begin{split}\label{eq:theorem:CLT-Max-Norm-12}
			\Delta_n &\leq  \frac{K \lambda}{\sqrt{\mathrm{Var}(\|Z^\eta\|_\infty)}} + \frac{K}{n^{1/2}\lambda^3}\mathrm{E}\left[\|X^\eta \|_\infty^3\right]\Delta_{n-1}^\eta\\
			&\quad{}+ \frac{K}{n^{1/2}\lambda^2}\frac{\mathrm{E}[\|X^\eta \|_\infty^3]}{\sqrt{\mathrm{Var}(\|Z^\eta\|_\infty)}} + \frac{K}{n\lambda^3} \frac{\mathrm{E}[\|X^\eta\|_\infty^3] \mathrm{E}[\|Z^\eta\|_\infty]}{\sqrt{\mathrm{Var}(\|Z^\eta\|_\infty)}} \\
			&\quad{}+\frac{K}{n\lambda^3} \frac{\mathrm{E}[\|X^\eta\|_\infty^3] M}{\sqrt{\mathrm{Var}(\|Z^\eta\|_\infty)}}  +\frac{K}{n^{1/2}\lambda^3} \mathrm{E} \left[ \|X^\eta\|_\infty^3 \mathbf{1}\{\|X^\eta\|_\infty > M\}\right]\\
			&\quad{} + \frac{3\eta}{\lambda}\sqrt{2\log 2d}.
		\end{split}
	\end{align}
	Set
	\begin{align*}
		\lambda = K \left(\frac{ C_*}{n- 1} \right)^{1/6} \left(\mathrm{E}[\|X^\eta\|_\infty^3] \right)^{1/3}
	\end{align*}
	and combine eq.~\eqref{eq:theorem:CLT-Max-Norm-11} and eq.~\eqref{eq:theorem:CLT-Max-Norm-12} to conclude (as in the first part of the proof) that for all $n \geq 2$,
	\begin{align}
		\Delta_n &\leq  \frac{K \lambda}{\sqrt{\mathrm{Var}(\|Z^\eta\|_\infty)}}  + B^\eta_{n-1} \sqrt{C_*} + B_n^\eta +  \frac{3\eta}{\lambda}\sqrt{2\log 2d}\nonumber \\
		&\leq  B_n^\eta \left[\left(\frac{n}{n-1}\right)^{1/6} K^2 C_*^{1/6} +\left(\frac{n}{n-1}\right)^{1/2} \sqrt{C_*} + 1 \right] +  \frac{3\eta}{\lambda}\sqrt{2\log 2d} \nonumber\\
		& \leq B_n^\eta \left[ \left(1 + 2K^2\right)\sqrt{C_*} + 1 \right] +\frac{ 3 \eta}{K} \frac{ \sqrt{2\log 2d}}{\left(\mathrm{E}[\|X^\eta\|_\infty^3] \right)^{1/3}} \left(\frac{n-1}{C_*}\right)^{1/6}.\label{eq:theorem:CLT-Max-Norm-13}
	\end{align}
	Since $\eta > 0$ arbitrary, we can take $\eta \rightarrow 0$. To complete the proof, we only need to find the limit of the expression on the right hand side in above display.
	
	Let $(\eta_k)_{k \geq 1}$ be a monotone falling null sequence. For all $0 < \eta_k \leq 1$ and $1 \leq p \leq 3$, $\|X^\eta\|_\infty^p \equiv \|X + \eta Y\|_\infty^p \leq 2^{p-1}\|X\|_\infty^p + 2^{p-1}\|Y\|_\infty^p$ a.s. and $\mathrm{E}[\|X\|_\infty^p] < \infty$ (otherwise the upper bound in Proposition~\ref{lemma:CLT-Max-Norm} is trivial!) and $\mathrm{E}[\|Y\|_\infty^p] < \infty$~\citep[e.g.][Proposition A.2.4]{vandervaart1996weak}. Hence, $(\|X^\eta\|_\infty^p)_{k\geq 1}$ is uniformly integrable for $p \in\{1,2\}$. Since in addition $\|X^\eta\|_\infty^p \rightarrow \|X\|^p$ a.s., the sequence $(\|X^\eta\|_\infty)_{k \geq 1}$ converges in $L^2$. Arguing similarly, we conclude that the sequences $(\|Z^\eta\|_\infty^2)_{k \geq 1}$ converges in $L^2$ as well. Thus,
\begin{align*}
	B_n^\eta \rightarrow \frac{(\mathrm{E}[\|X\|_\infty^3])^{1/3}}{n^{1/6}\sqrt{\mathrm{Var}(\|Z\|_\infty)}} + \frac{\mathrm{E} \left[ \|X\|_\infty^3 \mathbf{1}\{\|X\|_\infty > M\}\right]}{\mathrm{E} \left[ \|X\|_\infty^3\right]}  + \frac{ \mathrm{E}[\|Z\|_\infty] + M}{\sqrt{n\mathrm{Var}(\|Z\|_\infty)}} \equiv B_n \quad{} \mathrm{as} \quad{} k \rightarrow \infty,
\end{align*}
and
\begin{align*}
	\frac{ 3 \eta_k}{K} \frac{ \sqrt{2\log 2d}}{\left(\mathrm{E}[\|X^\eta\|_\infty^3] \right)^{1/3}} \left(\frac{n-1}{C_*}\right)^{1/6} \rightarrow 0\quad{} \mathrm{as} \quad{} k \rightarrow \infty.
\end{align*}
Hence, by eq.~\eqref{eq:theorem:CLT-Max-Norm-13} we have shown that $\Delta_n \lesssim B_n$ for all $M\geq 0$ and $n \geq 1$. This completes the proof of the proposition.
\end{proof}

\subsection{Proof of Corollary~\ref{corollary:lemma:CLT-Max-Norm}}
\begin{proof}[\textbf{Proof of Corollary~\ref{corollary:lemma:CLT-Max-Norm}}]
	Define $\widetilde{X} = (X^{(j)}/\sigma_{(1)})_{j=1}^d$ and $\widetilde{Z} = (Z^{(j)}/\sigma_{(1)})_{j=1}^d$ where $\sigma_{(1)}^2 = \\ \min_{1 \leq j \leq d} \Sigma_{jj}$. Thus, by Lemma~\ref{lemma:BoundsVariance-SeparableProcess},
	\begin{align}\label{eq:corollary:lemma:CLT-Max-Norm-1}
		\mathrm{Var}(\|Z\|_\infty) \gtrsim \left(\frac{\sigma_{(1)}}{1 + \mathrm{E}[\|\widetilde{Z}\|_\infty]}\right)^2.
	\end{align}
	Moreover, for $M, s > 0$ arbitrary,
	\begin{align*}
		\mathrm{E} [ \|\widetilde{X}\|_\infty^3 \mathbf{1}\{\|\widetilde{X}\|_\infty > M\}] \leq \mathrm{E} [ \|\widetilde{X}\|_\infty^{3 + s} \mathbf{1}\{\|\widetilde{X}\|_\infty > M\}] M^{-s} \leq \mathrm{E}[ \|\widetilde{X}\|_\infty^{3 + s}]M^{-s},
	\end{align*}
	and, hence, for $M_{3 + \delta} := \mathrm{E}[\|\widetilde{X} \|_\infty^{3 + \delta}]^{1/(3 + \delta)}$  and $s = \delta$,
	\begin{align}\label{eq:corollary:lemma:CLT-Max-Norm-2}
		\frac{\mathrm{E} [ \|\widetilde{X}\|_\infty^3 \mathbf{1}\{\|\widetilde{X}\|_\infty > n^{1/3}M_{3 + \delta}\}]}{\mathrm{E}[ \|\widetilde{X}\|_\infty^3]} \leq \frac{1}{n^{\delta/3}} \frac{\mathrm{E}[ \|\widetilde{X}\|_\infty^{3 + \delta}]}{M_{3 + \delta}^\delta\mathrm{E}[ \|\widetilde{X}\|_\infty^3] } \leq \frac{1}{n^{\delta/3}} \frac{M_{3 + \delta}^3}{\mathrm{E}[ \|\widetilde{X}\|_\infty^3] }.
	\end{align}
	Combine the upper bound in Proposition~\ref{lemma:CLT-Max-Norm} with eq.~\eqref{eq:corollary:lemma:CLT-Max-Norm-1} and~\eqref{eq:corollary:lemma:CLT-Max-Norm-2} and simplify the expression to conclude the proof of the first claim. (Obviously, this bound is not tight, but it is aesthetically pleasing.) The second claim about equicorrelated coordinates in $X$ follows from the lower bound on the variance of $\|\widetilde{Z}\|_{\infty}$ in Proposition 4.1.1 in~\cite{tanguy2017quelques} combined with the upper bound in Proposition~\ref{lemma:CLT-Max-Norm} and eq.~\eqref{eq:corollary:lemma:CLT-Max-Norm-2}.
\end{proof}

\subsection{Proof of Proposition~\ref{lemma:GaussianComparison}}
\begin{proof}[\textbf{Proof of Proposition~\ref{lemma:GaussianComparison}}]	
	The main proof idea is standard, e.g. similar arguments have been used in proofs by~\cite{fang2021HighDimCLT} (Theorem 1.1) and~\cite{chernozhukov2021NearlyOptimalCLT} (Theorem 3.2). While our bound is dimension-free, it is not sharp~\citep[e.g.][Proposition 2.1]{chernozhukov2021NearlyOptimalCLT}.
	
	\vspace{10pt}
	\noindent
	\textbf{The case of positive definite $\Omega$.} 
	\vspace{10pt}
	
	\noindent
	Suppose that $\Sigma$ is positive semi-definite and $\Omega$ is positive definite. To simplify notation, we set
	\begin{align*}
		\Delta : = \sup_{s \geq 0} \Big|\mathbb{P}\left\{\|Y\|_\infty \leq s \right\}-\mathbb{P}\left\{\|Z\|_\infty \leq s \right\} \Big|.
	\end{align*}
	Moreover, for $s \in \mathbb{R}$, $\lambda \geq 0$ arbitrary denote by $h_{s, \lambda}$ the map from Lemma~\ref{lemma:Smooth-Lipschitz-Approx} and define $y \mapsto h(y) := h_{s, \lambda}(\|y\|_\infty)$. Then, by Lemma~\ref{lemma:AntiConcentration-SeparableProcess} and Lemma~\ref{lemma:Smooth-Lipschitz-Approx} (ii),
	\begin{align}\label{eq:theorem:GaussianComparison-1}
		\Delta \leq \sup_{s \in \mathbb{R}} \big|\mathrm{E}[h(Y) - h(Z)] \big| + \frac{2 \sqrt{3}\lambda }{ \sqrt{\mathrm{Var}(\|Y\|_\infty) \vee \mathrm{Var}(\|Z\|_\infty)  + \lambda^2/12}}.
	\end{align} 
	Since $y \mapsto h(y) - \mathrm{E}[h(Z)]$ is Lipschitz continuous (with constant $\lambda^{-1}$) and $\Omega$ is positive definite, Proposition 4.3.2 in~\cite{nourdin2012normal} implies that
	\begin{align}\label{eq:theorem:GaussianComparison-2}
		\mathrm{E}\left[h(Y) - h(Z)\right] =  \mathrm{E}\left[\mathrm{tr}\left(\Omega D^2 G_h(Y) \right) - Y' D G_h(Y) \right],
	\end{align}
	where $G_h \in C^\infty(\mathbb{R}^d)$ and
	\begin{align*}
		G_h(y) := \int_0^\infty \Big(\mathrm{E}[h(Z)]  - P_th(y)  \Big) dt \quad{}\quad{} \forall y \in \mathbb{R}^d,
	\end{align*}
	and $P_t h$ denotes the Ornstein-Uhlenbeck semi-group with stationary measure $N(0, \Omega)$, i.e.
	\begin{align*}
		P_th(y) := \mathrm{E}\left[h\left( e^{-t}y + \sqrt{1- e^{-2t}}Z \right) \right] \quad{}\quad{} \forall y \in \mathbb{R}^d.
	\end{align*}
	Using Stein's lemma 
	we re-write eq.~\eqref{eq:theorem:GaussianComparison-2} as
	\begin{align*}
		\mathrm{E}\left[h(Y) - h(Z)\right] =  \mathrm{E}\left[\mathrm{tr}\left(\Omega D^2 G_h(Y) \right) - \mathrm{tr}\left(\Sigma D^2 G_h(Y) \right)\right].
	\end{align*}
	Notice that above identity holds even if $\Sigma$ is only positive semi-definite~\citep[e.g.][Lemma 4.1.3]{nourdin2012normal}! By H{\"o}lder's inequality for matrix inner products,
	\begin{align}\label{eq:theorem:GaussianComparison-3}
		\mathrm{E}\left[h(Y) - h(Z)\right] \leq  \left( \max_{j,k} |\Omega_{jk} - \Sigma_{jk}| \right) \left(\sum_{j,k} \mathrm{E}\left[\left|\frac{\partial^2 G_h}{\partial x_j \partial x_k }(Y) \right| \right] \right).
	\end{align}
	To complete the proof, we now bound the second derivative $D^2 G_h$. By Lemma~\ref{lemma:DifferentiatingUnderIntegral} (i), for arbitrary indices $1 \leq j, k \leq d$, 
	\begin{align*}
		\frac{\partial^2 G_h}{\partial y_j \partial y_k }(Y) = - \int_0^\infty e^{-2t} \mathrm{E}_Z\left[\frac{\partial^2 h}{\partial y_j \partial y_k} \left( e^{-t} y + \sqrt{1- e^{-2t}}Z \right) \Big|_{ y = Y} \right]dt,
	\end{align*}
	and, hence, by Lemma~\ref{lemma:DifferentiatingUnderIntegral} (ii),
	\begin{align*}
		& e^{-2t}\sum_{j,k}  \mathrm{E}_Z\left[\left|\frac{\partial^2 h}{\partial y_j \partial y_k} \left( e^{-t} y + \sqrt{1- e^{-2t}}Z \right) \Big|_{ y=Y} \right|\right]\\
		&\quad{} \leq C_2\lambda^{-2}e^{-2t} \mathrm{E}_Z\left[\sum_{j=1}^d \mathbf{1}_{[s, s+ 3\lambda]} \left(\|V^t\|_\infty \right) \mathbf{1}\left\{|V_j^t| \geq | V_m^t|, \: m \neq j \right\} \right],
	\end{align*}
	where
	\begin{align*}
		V^t :=	e^{-t}Y + \sqrt{1- e^{-2t}}Z.
	\end{align*}
	Since $\Omega$ is positive definite, no pair of entries in $Z$ can be perfectly (positively or negatively) correlated. Therefore, $\sum_{j=1}^d \mathbf{1}_{[s, s+ 3\lambda]} \left(\|V^t\|_\infty \right) \mathbf{1}\left\{|V_j^t| \geq | V_m^t|, \: m \neq j \right\} = 1$ almost surely. Hence,
	\begin{align}\label{eq:theorem:GaussianComparison-4}
		e^{-2t}\sum_{j,k}  \mathrm{E}\left[\left|\frac{\partial^2 h}{\partial y_j \partial y_k} \left( e^{-t} y + \sqrt{1- e^{-2t}}Z \right) \Big|_{ y=Y} \right|\right] \leq C_2  \lambda^{-2}e^{-2t}.
	\end{align}
	Combine eq.~\eqref{eq:theorem:GaussianComparison-3}--\eqref{eq:theorem:GaussianComparison-4}, integrate over $t \in (0, \infty)$, and conclude that
	\begin{align}\label{eq:theorem:GaussianComparison-6}
		\mathrm{E}\left[h(Y) - h(Z)\right]  \leq C_2\lambda^{-2} \left( \max_{j,k} |\Omega_{jk} - \Sigma_{jk}|  \right).
	\end{align}
	To conclude, combine eq.~\eqref{eq:theorem:GaussianComparison-1} and eq.~\eqref{eq:theorem:GaussianComparison-6} and optimize over $\lambda > 0$.
	
	\vspace{10pt}
	\noindent
	\textbf{The case of positive semi-definite $\Omega \neq \mathbf{0}$.} 
	\vspace{10pt}
	
	\noindent
	Suppose that both, $\Sigma$ and $\Omega$ are positive semi-definite. To avoid trivialities, we assume that $\Omega$ is not identical to zero. 
	
	Take $W \sim N(0, I_d)$ independent of $Z$ and, for $\eta > 0$ arbitrary, define $Z^\eta := Z + \eta W$. Now, consider eq.~\eqref{eq:theorem:GaussianComparison-1}, i.e.
	\begin{align*}
		\Delta \leq \sup_{s \in \mathbb{R}} \big|\mathrm{E}[h(Y) - h(Z)] \big| + \frac{2 \sqrt{3}\lambda }{ \sqrt{\mathrm{Var}(\|Y\|_\infty) \vee \mathrm{Var}(\|Z\|_\infty)  + \lambda^2/12}}.
	\end{align*} 
	Expand above inequality yields
	\begin{align}\label{eq:theorem:GaussianComparison-7}
		\Delta &\leq \sup_{s \in \mathbb{R}} \big|\mathrm{E}[h(Y) - h(Z^\eta)] \big| + \frac{2 \sqrt{3}\lambda }{ \sqrt{\mathrm{Var}(\|Y\|_\infty) \vee \mathrm{Var}(\|Z\|_\infty)  + \lambda^2/12}} + \sup_{s \in \mathbb{R}} \big|\mathrm{E}[h(Z^\eta) - h(Z)] \big| \nonumber\\
		&\leq\sup_{s \in \mathbb{R}} \big|\mathrm{E}[h(Y) - h(Z^\eta)] \big| + \frac{2 \sqrt{3}\lambda }{ \sqrt{\mathrm{Var}(\|Y\|_\infty) \vee \mathrm{Var}(\|Z\|_\infty)  + \lambda^2/12}} + \frac{\eta}{\lambda} \sqrt{2 \log 2d},
	\end{align}
	where the last inequality holds because $h$ is $\lambda^{-1}$-Lipschitz continuous w.r.t. the metric induced by the $\ell_\infty$-norm.
	
	Since $Z^\eta \sim N(0, \Omega + \eta^2 I_d)$ has positive definite covariance matrix, we can bound the first term on the far right hand side in above display using eq.~\eqref{eq:theorem:GaussianComparison-6}, i.e.
	\begin{align}\label{eq:theorem:GaussianComparison-8}
		\sup_{s \in \mathbb{R}} \big|\mathrm{E}[h(Y) - h(Z^\eta)] \big|  \leq  C_2\lambda^{-2} \left( \max_{j,k} |\Omega_{jk} - \Sigma_{jk} + \eta^2 \mathbf{1}\{j=k\}|  \right).
	\end{align}
	Combine  eq.~\eqref{eq:theorem:GaussianComparison-7} and eq.~\eqref{eq:theorem:GaussianComparison-8} to obtain
	\begin{align*}
		\Delta &\leq C_2\lambda^{-2} \left( \max_{j,k} |\Omega_{jk} - \Sigma_{jk} + \eta^2 \mathbf{1}\{j=k\}|  \right) + \frac{2 \sqrt{3}\lambda }{ \sqrt{\mathrm{Var}(\|Y\|_\infty) \vee \mathrm{Var}(\|Z\|_\infty)  + \lambda^2/12}}\\
		&\quad{} + \frac{\eta}{\lambda} \sqrt{2 \log 2d}.
	\end{align*}
	Letting $\eta \rightarrow 0$ and optimizing over $\lambda > 0$ gives the desired bound on $\Delta$.
\end{proof}

\subsection{Proofs of Propositions~\ref{lemma:Bootstrap-Max-Norm} and~\ref{lemma:Bootstrap-Max-Norm-Asymptotic-Size}}

\begin{proof}[\textbf{Proof of Proposition~\ref{lemma:Bootstrap-Max-Norm}}]
	By the triangle inequality,
	\begin{align*}
		&\sup_{s \geq 0} \Big|\mathbb{P}\left\{\|S_n\|_\infty \leq s \right\} -\mathbb{P}\left\{\|\widehat{\Sigma}_n^{1/2} Z\|_\infty \leq s \mid X_1, \ldots, X_n \right\} \Big|\\
		&\quad{}\leq \sup_{s \geq 0} \Big|\mathbb{P}\left\{\|S_n\|_\infty \leq s \right\} -\mathbb{P}\left\{\|\Sigma_n^{1/2} Z\|_\infty \leq s \right\} \Big|\\
		&\quad{}\quad{} + \sup_{s \geq 0} \Big|\mathbb{P}\left\{\|\Sigma_n^{1/2} Z\|_\infty \leq s \right\} -\mathbb{P}\left\{\|\widehat{\Sigma}_n^{1/2} Z\|_\infty \leq s \mid X_1, \ldots, X_n \right\} \Big|.
	\end{align*}
	Now, apply Proposition~\ref{lemma:CLT-Max-Norm} to the first summand and Proposition~\ref{lemma:GaussianComparison} to the second summand. This completes the proof.
\end{proof}

\begin{proof}[\textbf{Proof of Proposition~\ref{lemma:Bootstrap-Max-Norm-Asymptotic-Size}}]
	The proof is an adaptation of the proof of Theorem 3.1 in~\cite{chernozhukov2013GaussianApproxVec} to our setup. To simplify notation we write $c_n^*(\alpha) := c_n(\alpha; \widehat{\Sigma}_n)$ and $c_n(\alpha) := c_n(\alpha; \Sigma)$; see also eq.~\eqref{eq:subsec:QuantileComparison-0}. Note that
	\begin{align}\label{eq:theorem:Bootstrap-Max-Norm-Asymptotic-Size-1}
		&\sup_{\alpha \in (0,1)} \left|\mathbb{P}\left\{\|S_n\|_\infty  + \Theta \leq c_n^*(\alpha) \right\} - \alpha\right| \nonumber \\
		\begin{split}
			&\quad{}\leq\sup_{\alpha \in (0,1)}  \left|\mathbb{P}\left\{\|S_n\|_\infty  + \Theta \leq c_n^*(\alpha)\right\} - \mathbb{P}\left\{\|S_n\|_\infty  + \Theta \leq c_n(\alpha) \right\} \right|\\
			&\quad{}\quad{}+ \sup_{\alpha \in (0,1)} \left|\mathbb{P}\left\{\|S_n\|_\infty  + \Theta \leq c_n(\alpha) \right\} - \mathbb{P}\left\{\|\Sigma^{1/2}Z\|_\infty \leq c_n(\alpha)\right\} \right|.
		\end{split}
	\end{align}
	For $\delta > 0$ arbitrary, the first term can be upper bounded by Lemma~\ref{lemma:Bootstrap-Max-Norm-Quantil-Comparison} as
	\begin{align}\label{eq:theorem:Bootstrap-Max-Norm-Asymptotic-Size-2}
		&\sup_{\alpha \in (0,1)} \mathbb{P}\Big\{  c_n\big(\alpha - \pi_n(\delta) \big) < \|S_n\|_\infty + \Theta \leq  c_n \big(\alpha + \pi_n(\delta)\big)\Big\}  + 2 \mathbb{P}\left\{\max_{j,k} |\widehat{\Sigma}_{n,jk} - \Sigma_{jk}|  > \delta\right\}\nonumber\\
		&\quad \leq \sup_{\alpha \in (0,1)} \mathbb{P}\Big\{  c_n\big(\alpha - \pi_n(\delta) \big) - \eta < \|S_n\|_\infty \leq  c_n \big(\alpha + \pi_n(\delta)\big) + \eta \Big\} \\
		&\quad\quad  + \mathbb{P}\left\{|\Theta| > \eta \right\} + 2 \mathbb{P}\left\{\max_{j,k} |\widehat{\Sigma}_{n,jk} - \Sigma_{jk}|  > \delta\right\} \nonumber\\
		&\quad{}\leq \sup_{\alpha \in (0,1)} \mathbb{P} \Big\{c_n\big(\alpha - \pi_n(\delta)\big) - \eta <\|\Sigma^{1/2}Z\|_\infty \leq c_n\big(\alpha + \pi_n(\delta)\big) + \eta \Big\}  + \mathbb{P}\left\{|\Theta| > \eta \right\} \nonumber\\
		&\quad{}\quad{}+ 2\sup_{s \geq 0} \left|\mathbb{P}\big\{\|S_n\|_\infty \leq s\big\} - \mathbb{P}\big\{\|\Sigma^{1/2}Z\|_\infty \leq s\big\} \right| + 2\mathbb{P}\left\{\max_{j,k} |\widehat{\Sigma}_{n,jk} - \Sigma_{jk}|  > \delta\right\}\nonumber\\
		&\quad{} \leq \sup_{\alpha \in (0,1)} \mathbb{P} \Big\{c_n\big(\alpha - \pi_n(\delta)\big) <\|\Sigma^{1/2}Z\|_\infty \leq c_n\big(\alpha + \pi_n(\delta)\big) \Big\}  \nonumber\\
		&\quad{}\quad + \frac{\eta 8\sqrt{3}}{ \sqrt{\mathrm{Var}( \|\Sigma^{1/2}Z\|_\infty) + \eta^2/3} } + \mathbb{P}\left\{|\Theta| > \eta \right\} \nonumber\\
		&\quad{}\quad{}+ 2\sup_{s \geq 0} \left|\mathbb{P}\big\{\|S_n\|_\infty \leq s\big\} - \mathbb{P}\big\{\|\Sigma^{1/2}Z\|_\infty \leq s\big\} \right| + 2\mathbb{P}\left\{\max_{j,k} |\widehat{\Sigma}_{n,jk} - \Sigma_{jk}|  > \delta\right\} \nonumber\\
		\begin{split}
			&\quad{} \leq 2\pi_n(\delta) + \frac{\eta 8\sqrt{3}}{ \sqrt{\mathrm{Var}( \|\Sigma^{1/2}Z\|_\infty) + \eta^2/3} } +\mathbb{P}\left\{|\Theta| > \eta \right\}\\
			&\quad{}\quad{}+ 2\sup_{s \geq 0} \left|\mathbb{P}\big\{\|S_n\|_\infty \leq s\big\} - \mathbb{P}\big\{\|\Sigma^{1/2}Z\|_\infty \leq s\big\} \right| + 2\mathbb{P}\left\{\max_{j,k} |\widehat{\Sigma}_{n,jk} - \Sigma_{jk}|  > \delta\right\},\\
		\end{split}
	\end{align}
	where the second inequality follows from (several applications of) Lemma~\ref{lemma:AntiConcentration-SeparableProcess} and the third from the definition of quantiles and because $\|\Sigma^{1/2}Z\|_\infty$ has no point masses.
		
	Let $\eta > 0$ be arbitrary. We now bound the second term on the right hand side of eq.~\eqref{eq:theorem:Bootstrap-Max-Norm-Asymptotic-Size-1} by
	\begin{align}\label{eq:theorem:Bootstrap-Max-Norm-Asymptotic-Size-4}
		&\sup_{\alpha \in (0,1)} \left|\mathbb{P}\left\{\|S_n\|_\infty + \Theta \leq c_n(\alpha) \right\} - \mathbb{P}\left\{\|S_n\|_\infty\leq c_n(\alpha) \right\} \right|  \nonumber\\
		&\quad{}\quad{}\quad{}+ \sup_{s \geq 0} \left|\mathbb{P}\left\{\|S_n\|_\infty \leq s\right\} -\mathbb{P}\left\{\|\Sigma^{1/2}Z\|_\infty\leq s \right\} \right|\nonumber\\
		&\quad{}\quad{} \lesssim \mathbb{P}\left\{|\Theta| > \eta \right\} + \sup_{s \geq 0} \mathbb{P}\left\{s - \eta \leq \|S_n\|_\infty  \leq s+ \eta \right\} \nonumber\\
		&\quad{}\quad{}\quad{}\quad{} + \sup_{s \geq 0} \left| \mathbb{P}\left\{\|S_n\|_\infty \leq s\right\} - \mathbb{P}\left\{\|\Sigma^{1/2}Z\|_\infty \leq s \right\} \right|\nonumber\\
		&\quad{}\quad{} \lesssim \mathbb{P}\left\{|\Theta| > \eta \right\} + \sup_{s \geq 0} \mathbb{P}\left\{s- \eta \leq \|\Sigma^{1/2}Z\|_\infty \leq s + \eta \right\} \nonumber\\
		&\quad{}\quad{}\quad{}\quad{} + \sup_{s \geq 0} \left| \mathbb{P}\left\{\|S_n\|_\infty \leq s\right\} - \mathbb{P}\left\{\|\Sigma^{1/2}Z\|_\infty\leq s \right\} \right|\nonumber\\
		\begin{split}
		&\quad{}\quad{}\lesssim \mathbb{P} \left\{|\Theta| > \eta \right\} + \frac{\eta 4\sqrt{3}}{ \sqrt{\mathrm{Var}( \|\Sigma^{1/2}Z\|_\infty) + \eta^2/3} } \\
		&\quad{}\quad{}\quad{}\quad{} + \sup_{s\geq 0} \left|\mathbb{P}\left\{\|S_n\|_\infty \leq s\right\} -\mathbb{P}\left\{\|\Sigma^{1/2}Z\|_\infty\leq s\right\} \right|,
		\end{split}
	\end{align}	
	where the third inequality follows Lemma~\ref{lemma:AntiConcentration-SeparableProcess}. 
	
	Combine eq.~\eqref{eq:theorem:Bootstrap-Max-Norm-Asymptotic-Size-1}--\eqref{eq:theorem:Bootstrap-Max-Norm-Asymptotic-Size-4} to obtain
	\begin{align*}
		&\sup_{\alpha \in (0,1)} \left|\mathbb{P}\left\{\|S_n\|_\infty + \Theta \leq c_n^*(\alpha)\right\} - \alpha\right|\nonumber\\
		&\quad{} \lesssim  \sup_{s \geq 0} \left| \mathbb{P}\left\{\|S_n\|_\infty \leq s\right\} - \mathbb{P}\left\{\|\Sigma^{1/2}Z\|_\infty \leq s \right\} \right| + \inf_{\delta > 0}\left\{ \pi_n(\delta) + \mathbb{P}\left\{\max_{j,k} |\widehat{\Sigma}_{n,jk} - \Sigma_{jk}|  > \delta\right\}\right\} \nonumber\\
		&\quad{}\quad{}+ \inf_{\eta > 0} \left\{ \frac{\eta }{ \sqrt{\mathrm{Var}(\|\Sigma^{1/2}Z\|_\infty)} } + \mathbb{P}\left\{|\Theta| > \eta \right\}\right\}.
	\end{align*}
	To complete the proof bound the first term on the right hand side by Proposition~\ref{lemma:CLT-Max-Norm}.
\end{proof}

\newpage 

\section{Proofs of the results in Section~\ref{sec:ResultsSupremaEP}}\label{sec:Proofs-ResultsSupremaEP}

\subsection{Proofs of Theorem~\ref{theorem:CLT-Max-Norm-Simultaneous} and Corollary~\ref{corollary:theorem:CLT-Max-Norm-Simultaneous}} 
\begin{proof}[\textbf{Proof of Theorem~\ref{theorem:CLT-Max-Norm-Simultaneous}}]
	Let $\delta > 0$ be arbitrary and define $r_n(\delta) :=  \psi_n(\delta) \vee \phi_n(\delta)$. Let $\mathcal{H}_{n,\delta} \subset \mathcal{F}_n$ be a $\delta\|F_n\|_{P,2}$-net of $\mathcal{F}_n$ and set $\mathcal{F}_{n, \delta}' = \{ f - g: f, g \in \mathcal{F}_n, \: \rho(f, g) < \delta\|F_n\|_{P,2}\}$. Since $\mathcal{F}_n$ is totally bounded with respect to $\rho$, the $\mathcal{H}_{n,\delta}$'s are finite. Moreover, for each $\delta > 0$ there exists a map $\pi_\delta: \mathcal{F}_n \rightarrow \mathcal{H}_{n,\delta}$ such that $\rho(f, \pi_\delta f) < \delta\|F_n\|_{P,2}$ for all $f \in \mathcal{F}_n$ and, hence, 
	\begin{align*}
		\big|\|G_P\|_{\mathcal{H}_{n,\delta}} - \|G_P\|_{\mathcal{F}_n} \big| \leq  \|G_P \circ (\pi_\delta - id)\|_{\mathcal{F}_n} \leq \|G_P\|_{\mathcal{F}_{n,\delta}'},
	\end{align*}
	where the first inequality holds by the reverse triangle inequality and the prelinearity of the Gaussian $P$-bridge process~\citep[e.g.][p. 65, eq. 2.4]{dudley2014uniform}. (Here and in the following $id$ stands for the identity map.) Since the map $f \mapsto (P_n - P) f$ is obviously linear, the same inequality holds for the empirical process $\{ \mathbb{G}_n(f): f \in \mathcal{F}_n \}$. 
	
	By Lemma~\ref{lemma:Kolmogorov-Coupling-AntiConcentration} and Lemma~\ref{lemma:AntiConcentration-SeparableProcess}, 
	\begin{align*} 
		\begin{split}
			&\sup_{s \geq 0} \Big|\mathbb{P}\left\{\|\mathbb{G}_n\|_{\mathcal{F}_n} \leq s \right\} -\mathbb{P}\left\{\|G_P\|_{\mathcal{F}_n} \leq s \right\} \Big| \\
			&\quad{}\quad{}\leq \mathbb{P} \left\{ \big| \|\mathbb{G}_n\|_{\mathcal{F}_n} -  \|G_P\|_{\mathcal{F}_n} \big| > 3 \sqrt{r_n(\delta)\mathrm{Var}(\|G_P\|_{\mathcal{F}_n})} \right\} + 6 \sqrt{3}  \sqrt{r_n(\delta)}.
		\end{split}
	\end{align*}
	By the triangle inequality,
	\begin{align*} 
		\begin{split}
			&\mathbb{P} \left\{ \big| \|\mathbb{G}_n\|_{\mathcal{F}_n} -  \|G_P\|_{\mathcal{F}_n} \big| > 3  \sqrt{r_n(\delta)\mathrm{Var}(\|G_P\|_{\mathcal{F}_n})} \right\}\\
			&\quad{} \quad{}\leq \mathbb{P} \left\{ \big| \|\mathbb{G}_n\|_{\mathcal{H}_{n,\delta}} -  \|G_P\|_{\mathcal{H}_{n,\delta}} \big| > \sqrt{r_n(\delta)\mathrm{Var}(\|G_P\|_{\mathcal{F}_n})} \right\} \\
			&\quad{}\quad{}\quad{} + \mathbb{P} \left\{\|G_P\|_{\mathcal{F}_{n, \delta}'} >  \sqrt{r_n(\delta)\mathrm{Var}(\|G_P\|_{\mathcal{F}_n})} \right\} +  \mathbb{P} \left\{ \|\mathbb{G}_n\|_{\mathcal{F}_{n,\delta}'} >  \sqrt{r_n(\delta)\mathrm{Var}(\|G_P\|_{\mathcal{F}_n})} \right\},
		\end{split}
	\end{align*}
	Since $\mathcal{H}_\delta$ is finite, Proposition~\ref{lemma:CLT-Max-Norm} implies, for all $s \geq 0$,
	\begin{align*}
		\mathbb{P}\left\{ \|\mathbb{G}_n\|_{\mathcal{H}_{n,\delta}}  \leq s \right\} \leq \mathbb{P}\left\{ \| G_P\|_{\mathcal{H}_{n,\delta}}  \leq s \right\} + K B_n(\delta),
	\end{align*}	
	where $K > 0$ is an absolute constant and 
	\begin{align*}
		B_n(\delta) := \frac{\|F_n\|_{P, 3}}{n^{1/6}\sqrt{\mathrm{Var}(\|G_P\|_{\mathcal{H}_{n,\delta}})}} + \frac{\|F_n\mathbf{1}\{F_n > M_n\}\|_{P,3}^3}{\|F_n\|_{P,3}^3} + \frac{ \mathrm{E}\|G_P\|_{\mathcal{H}_{n,\delta}}+ M_n}{\sqrt{n\mathrm{Var}(\|G_P\|_{\mathcal{H}_{n,\delta}})}}.
	\end{align*}
	Therefore, by Strassen's theorem~\citep[e.g.][Theorem 11.6.2]{dudley2002real},
	\begin{align}\label{eq:theorem:CLT-Max-Norm-Simultaneous-1}
		\begin{split}
			&\sup_{s \geq 0} \Big|\mathbb{P}\left\{\|\mathbb{G}_n\|_{\mathcal{F}_n} \leq s \right\} -\mathbb{P}\left\{\|G_P\|_{\mathcal{F}_n} \leq s \right\} \Big| \\
			&\quad{}\quad{} \leq K B_n(\delta_n) + 6 \sqrt{3} \sqrt{r_n(\delta)} + \mathbb{P} \left\{\|G_P\|_{\mathcal{F}_{n,\delta}'} > \sqrt{r_n(\delta)\mathrm{Var}(\|G_P\|_{\mathcal{F}_n})} \right\}\\
			&\quad{}\quad{}\quad{}   +  \mathbb{P} \left\{ \|\mathbb{G}_n\|_{\mathcal{F}_{n,\delta}'} > \sqrt{r_n(\delta)\mathrm{Var}(\|G_P\|_{\mathcal{F}_n})} \right\},
		\end{split}
	\end{align}
	where $K > 0$ is an absolute constant and
	\begin{align*}
		B_n(\delta) = \frac{\|F_n\|_{P, 3}}{n^{1/6}\sqrt{\mathrm{Var}(\|G_P\|_{\mathcal{H}_{n,\delta}})}} + \frac{\|F_n\mathbf{1}\{F_n > M_n\}\|_{P,3}^3}{\|F_n\|_{P,3}^3} + \frac{\mathrm{E} \|G_P\|_{\mathcal{H}_{n, \delta}} + M_n}{\sqrt{n\mathrm{Var}(\|G_P\|_{\mathcal{H}_{n, \delta}})}}.
	\end{align*}
	By the hypothesis on the moduli of continuity of the Gaussian $P$-bridge and the empirical process, we can upper bound the right hand side in eq.~\eqref{eq:theorem:CLT-Max-Norm-Simultaneous-1} by
	\begin{align}\label{eq:theorem:CLT-Max-Norm-Simultaneous-2}
		K B_n(\delta) + 6 \sqrt{3} \sqrt{r_n(\delta)} + 2\frac{\psi_n(\delta) \vee \phi_n(\delta)}{r_n(\delta)} &< K B_n(\delta) + (2 + 6 \sqrt{3}) \sqrt{r_n(\delta) } \nonumber\\
		&< K B_n(\delta) + 13 \sqrt{r_n(\delta)}.
	\end{align}
	Since $r_n(\delta)$ is non-increasing in $\delta > 0$, we can assume without loss of generality that there exist $N  \geq 1$ such that for all $\delta > 0$ and all $n > N$, $\sqrt{r_n(\delta)} < 1/13$; otherwise the upper bound is trivial by eq.~\eqref{eq:theorem:CLT-Max-Norm-Simultaneous-2}. Thus, since $\|G_P\|_{\mathcal{F}_n} \geq \|G_P\|_{\mathcal{H}_{n,\delta}}$, Lemma~\ref{lemma:VarianceMaximum} with $X = \|G_P\|_{\mathcal{H}_{n,\delta}}$ and $Z = \|G_P\|_{\mathcal{F}_n}$ yields, for all $n > N$,
	\begin{align*}
		\sqrt{\mathrm{Var}(\|G_P\|_{\mathcal{F}_n})} &\leq \sqrt{\mathrm{Var}(\|G_P\|_{\mathcal{H}_{n,\delta}})}  + \sqrt{\mathrm{Var}(\|G_P\|_{\mathcal{F}_n} -\|G_P\|_{\mathcal{H}_{n,\delta}} ) }\\
		&\overset{(a)}{\leq} \sqrt{\mathrm{Var}(\|G_P\|_{\mathcal{H}_{n,\delta}})} + \sqrt{\mathrm{E} \|G_P\|_{\mathcal{F}_{n,\delta}'}^2}\\
		&\overset{(b)}{\leq} \sqrt{\mathrm{Var}(\|G_P\|_{\mathcal{H}_{n,\delta}})} + 12/13\sqrt{\mathrm{Var}(\|G_P\|_{\mathcal{F}_n})},
	\end{align*}
	where (a) holds because $0 \leq \|G_P\|_{\mathcal{F}_n} -\|G_P\|_{\mathcal{H}_{n,\delta}} \leq \|G_P\|_{\mathcal{F}_{n, \delta}'}$ and (b) holds because of $\sqrt{r_n(\delta)} < 1/13$ and for suprema of Gaussian processes we can reverse Liapunov's inequality (i.e. Lemma~\ref{lemma:ReverseLiapunov}, the unknown constant can be absorbed into $r_n(\delta)$). Conclude that for all $n > N$,
	\begin{align}\label{eq:theorem:CLT-Max-Norm-Simultaneous-3}
		(1-12/13) \sqrt{\mathrm{Var}(\|G_P\|_{\mathcal{F}_n})} \leq  \sqrt{\mathrm{Var}(\|G_P\|_{\mathcal{H}_{n,\delta}})}.
	\end{align}
	Combine eq.~\eqref{eq:theorem:CLT-Max-Norm-Simultaneous-1}--\eqref{eq:theorem:CLT-Max-Norm-Simultaneous-3} and obtain for all $n > N$,
	\begin{align*}
		\sup_{s \geq 0} \Big|\mathbb{P}\left\{\|\mathbb{G}_n\|_{\mathcal{F}_n} \leq s \right\} -\mathbb{P}\left\{\|G_P\|_{\mathcal{F}_n} \leq s \right\} \Big| \leq 13 \left(K B_n +  \sqrt{r_n(\delta)}\right),
	\end{align*}
	where
	\begin{align*}
		B_n = \frac{\|F_n\|_{P, 3}}{n^{1/6}\sqrt{\mathrm{Var}(\|G_P\|_{\mathcal{F}_n})}} + \frac{\|F_n\mathbf{1}\{F_n > M_n\}\|_{P,3}^3}{\|F_n\|_{P,3}^3} + \frac{\mathrm{E} \|G_P\|_{\mathcal{F}_n} + M_n}{\sqrt{n\mathrm{Var}(\|G_P\|_{\mathcal{F}_n})}}.
	\end{align*}
	Increase the constant $13$ on the right hand side until the bound holds also for all $ 1 \leq n \leq N$. Since $\delta > 0$ is arbitrary, set $\delta = \delta^*$ such that $r_n(\delta^*) = r_n := \inf\big\{ r_n(\delta) : \delta > 0\big\}$.
	
\end{proof}

\begin{proof}[\textbf{Proof of Corollary~\ref{corollary:theorem:CLT-Max-Norm-Simultaneous}}]
	We only need to verify that under the stated assumptions Theorem~\ref{theorem:CLT-Max-Norm-Simultaneous} applies with $r_n \equiv 0$.
	Let $n \geq 1$ be arbitrary and define  $\mathcal{F}_{n, \delta}' := \mathcal{F}_{n, \delta}'' \cap \mathcal{F}_{n, \delta}'''$, where $ \mathcal{F}_{n, \delta}'' = \{ f - g: f, g \in \mathcal{F}_n, \: d_P(f, g) < \delta\|F_n\|_{P,2}\}$ and $ \mathcal{F}_{n, \delta}''' = \{ f - g: f, g \in \mathcal{F}_n, \: d_{P_n}(f, g) < \delta\|F_n\|_{P,2}\}$. Note that $\mathcal{F}_{n, \delta}' =  \{ f - g: f, g \in \mathcal{F}_n, \: (d_P \vee d_{P_n})(f, g) < \delta\|F_n\|_{P,2}\}$.
	First, by Lemma~\ref{lemma:Sudakov}, $\mathcal{F}_n$ is totally bounded w.r.t. the pseudo-metric $d_P$. Second, by Lemma~\ref{lemma:ModulusContinuity}, $\mathrm{E}\|G_P\|_{\mathcal{F}_{n, \delta}'} \leq  \mathrm{E}\|G_P\|_{\mathcal{F}_{n, \delta}''}  \rightarrow 0$ as $\delta \rightarrow \infty$. Thus, for each $n \geq 1$, there exists some $\psi_n$ such that $\psi_n(\delta) \rightarrow 0$ as $\delta \rightarrow 0$. Third, $\big\|\|\mathbb{G}_n\|_{\mathcal{F}_{n,\delta}'} \big\|_{P,1}  \leq \big\|\|\mathbb{G}_n\|_{\mathcal{F}_{n,\delta}'''} \big\|_{P,1} \leq 2 \sqrt{n} \|F_n\|_{P,2} \delta$. Thus, for each $n \geq 1$, there also exists $\phi_n$ such that  $\psi_n(\delta) \rightarrow 0$ as $\delta \rightarrow 0$. This completes the proof.
	
	(Note that the last argument does not imply that the empirical process is asymptotically $d_P$-equicontinuous: Not only do we use a different norm, but we keep the sample size $n$ fixed and only consider the limit $\delta \rightarrow 0$. This shows once again that Gaussian approximability is a strictly weaker concept than weak convergence, see also Remark~\ref{remark:lemma:CLT-Max-Norm-3} and the discussion in Section~\ref{subsec:RelationPreviousWork-Math}.)
\end{proof}

\subsection{Proofs of Theorem~\ref{theorem:GaussianComparison-PQ-Simultaneous} and Corollaries~\ref{corollary:theorem:GaussianComparison-PQ-Simultaneous} and~\ref{corollary:theorem:GaussianComparison-PQ-Simultaneous-2}}

\begin{proof}[\textbf{Proof of Theorem~\ref{theorem:GaussianComparison-PQ-Simultaneous}}]
	
	The key idea is to construct ``entangled $\delta$-nets'' of the function classes $\mathcal{F}_n$ and $\mathcal{G}_n$. The precise meaning of this term and our reasoning for this idea will become clear in the course of the proof.
	
	Let $\delta > 0$ be arbitrary and define $r_n(\delta) :=  \psi_n(\delta) \vee \phi_n(\delta)$. Let $\mathcal{H}_{n,\delta} \subseteq \mathcal{F}_n$ and $\mathcal{I}_{n,\delta} \subseteq \mathcal{G}_n$ be $\delta \|F_n\|_{P,2}$- and $\delta \|G_n\|_{P,2}$-nets of $\mathcal{F}_n$ and $\mathcal{G}_n$ w.r.t. $\rho_1$ and $\rho_2$, respectively. Since $\mathcal{F}_n$ and $\mathcal{G}_n$ are totally bounded w.r.t. $\rho_1$ and $\rho_2$, respectively, the $\mathcal{H}_{n,\delta}$'s and $\mathcal{I}_{n,\delta}$'s are finite. 
	
	
	Next, let $\pi^1_{\mathcal{G}_n}$ be the projection from $\mathcal{F}_n \cup \mathcal{G}_n $ onto $\mathcal{G}_n$ defined by $\rho_1(h, \pi^1_{\mathcal{G}_n}h) = \inf_{g \in \mathcal{G}_n}\rho_1(h, g)$ for all $h \in \mathcal{F}_n \cup \mathcal{G}_n $. If this projection is not unique for some $h \in \mathcal{F}_n \cup \mathcal{G}_n$, choose any of the equivalent points. Define the projection $\pi^2_{\mathcal{F}_n}$ from $\mathcal{F}_n \cup \mathcal{G}_n $ onto $\mathcal{F}_n$ analogously. Finally, define the ``entangled $\delta$-nets'' of $\mathcal{F}_n$ and $\mathcal{G}_n$ as
	\begin{align*}
		\widetilde{\mathcal{H}}_{n,\delta} &: = \pi^2_{\mathcal{F}_n} \big(\mathcal{H}_{n,\delta} \cup \mathcal{I}_{n,\delta} \big) = \mathcal{H}_{n,\delta} \cup \pi^2_{\mathcal{F}_n} \big( \mathcal{I}_{n,\delta} \big) \subset \mathcal{F}_n, \quad \text{and} \\
		\widetilde{\mathcal{I}}_{n,\delta} &: = \pi^1_{\mathcal{G}_n} \big(\mathcal{H}_{n,\delta} \cup \mathcal{I}_{n,\delta} \big) = \pi^1_{\mathcal{G}_n} \big(\mathcal{H}_{n,\delta} \big) \cup \mathcal{I}_{n,\delta} \subset \mathcal{G}_n.
	\end{align*}
	Note that the entangled sets $\widetilde{\mathcal{H}}_{n,\delta}$ and $\widetilde{\mathcal{I}}_{n,\delta}$ are still $\delta \|F_n\|_{P,2}$- and $\delta \|G_n\|_{P,2}$-nets of $\mathcal{F}_n$ and $\mathcal{G}_n$ w.r.t. $\rho_1$ and $\rho_2$. Hence, for each $\delta > 0$, there exist maps $\tilde{\pi}_\delta: \mathcal{F}_n \rightarrow \widetilde{\mathcal{H}}_{n,\delta}$ and $\tilde{\theta}_\delta: \mathcal{G}_n \rightarrow \widetilde{\mathcal{I}}_{n,\delta}$ such that $\rho_1(f, \tilde{\pi}_\delta f) \leq \delta \|F_n\|_{P,2}$ for all $f \in \mathcal{F}_n$ and $\rho_2(g, \tilde{\theta}_\delta g) \leq \delta \|G_n\|_{Q,2}$ for all $g \in \mathcal{G}_n$. Therefore, by the reverse triangle inequality and the linearity of Gaussian $P$-bridge processes and $Q$-motions,
	\begin{align*}
		&\big|\|G_P\|_{\widetilde{\mathcal{H}}_{n, \delta}} - \|G_P\|_{\mathcal{F}_n} \big| \leq  \|G_P \circ (\tilde{\pi}_\delta - id)\|_{\mathcal{F}_n} \leq \|G_P\|_{\mathcal{F}_{n,\delta}'}, \quad{} \text{and}\\
		&\big|\|Z_Q\|_{\widetilde{\mathcal{I}}_{n, \delta}} - \|Z_Q\|_{\mathcal{G}_n} \big| \leq  \|Z_Q\circ (\tilde{\theta}_\delta - id)\|_{\mathcal{G}_n} \leq \|Z_Q\|_{\mathcal{G}_{n, \delta}'},
	\end{align*}
	where $\mathcal{F}_{n, \delta}' = \{ f_1 - f_2: f_1, f_2 \in \mathcal{F}_n, \: \rho_1(f_1, f_2) \leq \delta\|F_n\|_{P,2}\} $ and $\mathcal{G}_{n, \delta}' = \{ g_1 - g_2: g_1, g_2 \in \mathcal{G}_n, \: \rho_2(g_1, g_2) \leq \delta \|G_n\|_{P,2}\}$. 
	
	By Lemma~\ref{lemma:Kolmogorov-Coupling-AntiConcentration} and Lemma~\ref{lemma:AntiConcentration-SeparableProcess},
	\begin{align}\label{eq:theorem:GaussianComparison-Simultaneous-1}
		\begin{split}
			&\sup_{s \geq 0} \Big|\mathbb{P}\left\{\|G_P\|_{\mathcal{F}_n} \leq s \right\} -\mathbb{P}\left\{\|Z_Q\|_{\mathcal{G}_n} \leq s \right\} \Big| \\
			&\quad{}\quad{}\leq \mathbb{P} \left\{ \big| \|G_P\|_{\mathcal{F}_n} -  \|Z_Q\|_{\mathcal{G}_n} \big| > 3 \sqrt{r_n(\delta)}\sqrt{\mathrm{Var}(\|G_P\|_{\mathcal{F}_n})  \vee \mathrm{Var}(\|Z_Q\|_{\mathcal{G}_n})}   \right\} + 6 \sqrt{3} \sqrt{r_n(\delta)}.
		\end{split}
	\end{align}
	By the triangle inequality,
	\begin{align}\label{eq:theorem:GaussianComparison-Simultaneous-2}
		\begin{split}
			&\mathbb{P} \left\{ \big| \|G_P\|_{\mathcal{F}_n} -  \|Z_Q\|_{\mathcal{G}_n} \big| > 3  \sqrt{r_n(\delta)}\sqrt{\mathrm{Var}(\|G_P\|_{\mathcal{F}_n})  \vee \mathrm{Var}(\|Z_Q\|_{\mathcal{G}_n})} \right\}\\
			&\quad{} \quad{}\leq \mathbb{P} \left\{ \big| \|G_P\|_{\widetilde{\mathcal{H}}_{n,\delta}} -  \|Z_Q\|_{\widetilde{\mathcal{I}}_{n,\delta}} \big| > \sqrt{r_n(\delta)}\sqrt{\mathrm{Var}(\|G_P\|_{\mathcal{F}_n})  \vee \mathrm{Var}(\|Z_Q\|_{\mathcal{G}_n})}  \right\} \\
			&\quad{}\quad{}\quad{} + \mathbb{P} \left\{\|G_P\|_{\mathcal{F}_{\delta}'} >  \sqrt{r_n(\delta)}\sqrt{\mathrm{Var}(\|G_P\|_{\mathcal{F}_n})} \right\} +  \mathbb{P} \left\{ \|Z_Q\|_{\mathcal{G}_{\delta}'} >  \sqrt{r_n(\delta)}\sqrt{\mathrm{Var}(\|Z_Q\|_{\mathcal{G}_n})} \right\}.
		\end{split}
	\end{align}
	Since $	\|G_P\|_{\widetilde{\mathcal{H}}_{n,\delta}} \equiv \| G_P \circ \pi^2_{\mathcal{F}_n} \|_{\mathcal{H}_{n,\delta} \cup \mathcal{I}_{n,\delta}}$ and $\|Z_Q\|_{\widetilde{\mathcal{I}}_{n,\delta}} \equiv \| Z_Q \circ \pi^1_{\mathcal{G}_n} \|_{\mathcal{H}_{n,\delta} \cup \mathcal{I}_{n,\delta} }$ (note: two  finite dimensional mean zero Gaussian vectors with the same (!) index sets; to achieve this we needed the ``entangled $\delta$-nets'') we have by Proposition~\ref{lemma:GaussianComparison}, for all $s \geq 0$,
	\begin{align*}
		\mathbb{P}\left\{ \|G_P\|_{\widetilde{\mathcal{H}}_{n,\delta}}  \leq s \right\} \leq \mathbb{P}\left\{ \|Z_Q\|_{\widetilde{\mathcal{I}}_{n,\delta}} \leq s \right\} + K  \left( \frac{\Delta_n(\delta)}{ \mathrm{Var}(\|G_P\|_{\widetilde{\mathcal{H}}_{n,\delta}} ) \vee \mathrm{Var}(\|Z_Q\|_{\widetilde{\mathcal{I}}_{n,\delta}}  )} \right)^{1/3},
	\end{align*}
	where $K > 0$ is an absolute constant and, for all $\delta > 0$,
	\begin{align*}
		\Delta_n(\delta) &:= \sup_{h_1, h_2 \in \mathcal{H}_{n,\delta} \cup \mathcal{I}_{n,\delta}  } \big| \mathrm{E}[ (G_P \circ \pi^2_{\mathcal{F}_n}) (h_1) (G_P \circ \pi^2_{\mathcal{F}_n}) (h_2)] - \mathrm{E}[(Z_Q \circ \pi^1_{\mathcal{G}_n} )(h_1) (Z_Q \circ \pi^1_{\mathcal{G}_n} )(h_2)] \big|\\
		& \leq 	\sup_{f_1, f_2 \in \mathcal{F}_n} \big| \mathrm{E}[G_P(f_1)G_P(f_2)] - \mathrm{E}[Z_Q(\pi^1_{\mathcal{G}_n} f_1)Z_Q(\pi^1_{\mathcal{G}_n}f_2)] \big|\\
		&\quad \quad \quad\quad\quad \quad  \bigvee \sup_{g_1, g_2 \in \mathcal{G}_n }  \big| \mathrm{E}[Z_Q(g_1)Z_Q(g_2)] - \mathrm{E}[G_P(\pi^2_{\mathcal{F}_n} g_1)G_P(\pi^2_{\mathcal{F}_n} g_2)] \big|\\
		& = : \Delta_{P, Q}(\mathcal{F}_n, \mathcal{G}_n).
	\end{align*}
	Therefore, by eq.~\eqref{eq:theorem:GaussianComparison-Simultaneous-1} and eq.~\eqref{eq:theorem:GaussianComparison-Simultaneous-2}, Strassen's theorem~\citep[e.g.][Theorem 11.6.2; see also proof of Theorem~\ref{theorem:CLT-Max-Norm-Simultaneous}]{dudley2002real}, and Markov's inequality,
	\begin{align}\label{eq:theorem:GaussianComparison-Simultaneous-3}
		\sup_{s \geq 0} \Big|\mathbb{P}\left\{\|G_P\|_{\mathcal{F}_n} \leq s \right\} -\mathbb{P}\left\{\|Z_Q\|_{\mathcal{G}_n} \leq s \right\} \Big| \leq K  \left( \frac{\Delta_{P,Q}(\mathcal{F}_n, \mathcal{G}_n)}{  \mathrm{Var}(\|G_P\|_{\widetilde{\mathcal{H}}_{n,\delta}} ) \vee \mathrm{Var}(\|Z_Q\|_{\widetilde{\mathcal{I}}_{n,\delta}})} \right)^{1/3} + 13 \sqrt{r_n(\delta)}.
	\end{align}
	By the same arguments as in the proof of Theorem~\ref{theorem:CLT-Max-Norm-Simultaneous}, there exists $N > 0$ such that for all $\delta > 0$ and all $n \geq N$ we have $\sqrt{r_n(\delta)} < 1/13$ and, hence,
	\begin{align}\label{eq:theorem:GaussianComparison-Simultaneous-4}
		 \left( \frac{\Delta_{P,Q}(\mathcal{F}_n, \mathcal{G}_n)}{  \mathrm{Var}(\|G_P\|_{\widetilde{\mathcal{H}}_{n,\delta}} ) \vee \mathrm{Var}(\|Z_Q\|_{\widetilde{\mathcal{I}}_{n,\delta}})} \right)^{1/3} \leq 13 \left( \frac{\Delta_{P,Q}(\mathcal{F}_n, \mathcal{G}_n)}{ \mathrm{Var}(\|G_P\|_{\mathcal{F}_n} ) \vee \mathrm{Var}(\|Z_Q\|_{\mathcal{G}_n} )} \right)^{1/3}.
	\end{align}
	Combine eq.~\eqref{eq:theorem:GaussianComparison-Simultaneous-3} and eq.~\eqref{eq:theorem:GaussianComparison-Simultaneous-4} to conclude that the claim of the theorem holds for all $n \geq N$. Since for $\sqrt{r_n(\delta)} \geq 1/13$ the upper bound is trivial by eq.~\eqref{eq:theorem:GaussianComparison-Simultaneous-3}, the theorem holds in fact for all $n \geq 1$. Since $\delta > 0$ is arbitrary, set $\delta = \delta^*$ such that $r_n(\delta^*) = r_n := \inf\big\{ r_n(\delta) : \delta > 0\big\}$.
\end{proof}

\begin{proof}[\textbf{Proof of Corollary~\ref{corollary:theorem:GaussianComparison-PQ-Simultaneous}}]
	We only need to verify that under the stated assumptions Theorem~\ref{theorem:GaussianComparison-PQ-Simultaneous} applies with $r_n \equiv 0$. This follows by the same argument as in proof of Corollary~\ref{corollary:theorem:CLT-Max-Norm-Simultaneous}. We therefore omit a proof.
\end{proof}

\begin{proof}[\textbf{Proof of Corollary~\ref{corollary:theorem:GaussianComparison-PQ-Simultaneous-2}}]
	First, apply Theorem~\ref{theorem:GaussianComparison-PQ-Simultaneous} (and Remark~\ref{remark:theorem:GaussianComparison-PQ-Simultaneous}) with $\rho = d_P$, $Q = P$, and $\mathcal{G}_n \subseteq \mathcal{F}_n$ a $\delta \|F\|_{P,2}$-net of $\mathcal{F}_n$ with respect to $d_P$. Next, compute
	\begin{align*}
		\Delta_{P,P}(\mathcal{F}_n, \mathcal{G}_n) &= \sup_{f_1, f_2 \in \mathcal{F}_n} \big| \mathrm{E}[G_P(f_1)G_P(f_2)] - \mathrm{E}[G_P(\pi_{\mathcal{G}_n} f_1)G_P(\pi_{\mathcal{G}_n}f_2)] \big|\\
		&\leq \sup_{f_1, f_2 \in \mathcal{F}_n} \mathrm{E}\left[G_P^2(f_1)\big(G_P(f_2) - G_P(\pi_{\mathcal{G}_n} f_2)\big)^2  \right]^{1/2}\\
		&\quad +\sup_{f_1, f_2 \in \mathcal{F}_n} \mathrm{E}\left[G_P^2(\pi_{\mathcal{G}_n} f_2)\big(G_P(f_1) - G_P(\pi_{\mathcal{G}_n} f_1)\big)^2  \right]^{1/2}\\
		&\leq 2 \delta\|F\|_{P,2} \sup_{f \in \mathcal{F}_n} \sqrt{Pf^2}.
	\end{align*}
	Lastly, verify that under the stated assumptions Theorem~\ref{theorem:GaussianComparison-PQ-Simultaneous} applies with $r_n \equiv 0$. But this follows by the same argument as in proof of Corollary~\ref{corollary:theorem:CLT-Max-Norm-Simultaneous}. This completes the proof.
\end{proof}

\subsection{Proofs of Theorems~\ref{theorem:Abstract-Bootstrap-Sup-Empirical-Process} and~\ref{theorem:Abstract-Bootstrap-Sup-Empirical-Process-Quantiles}, Proposition~\ref{lemma:KarhunenLoeve-GP}, and Corollary~\ref{corollary:theorem:Abstract-Bootstrap-Sup-Empirical-Process-1}}

\begin{proof}[\textbf{Proof of Theorem~\ref{theorem:Abstract-Bootstrap-Sup-Empirical-Process}}]
	By the triangle inequality
	\begin{align*}
		&\sup_{s \geq 0} \Big|\mathbb{P}\left\{\|\mathbb{G}_n\|_{\mathcal{F}_n} \leq s \right\} -\mathbb{P}\left\{\|Z_Q\|_{\mathcal{G}_n} \leq s \mid X_1, \ldots, X_n \right\} \Big|\\
		&\quad{} \leq \sup_{s \geq 0} \Big|\mathbb{P}\left\{\|\mathbb{G}_n\|_{\mathcal{F}_n} \leq s \right\} -\mathbb{P}\left\{\|G_P\|_{\mathcal{F}_n} \leq s \right\} \Big|\\
		&\quad{} \quad{} + \sup_{s \geq 0} \Big|\mathbb{P}\left\{\|G_P\|_{\mathcal{F}_n} \leq s \right\} -\mathbb{P}\left\{\|Z_Q\|_{\mathcal{G}_n} \leq s \mid X_1, \ldots, X_n \right\} \Big|.
	\end{align*}
	To complete the proof, apply Theorem~\ref{theorem:CLT-Max-Norm-Simultaneous} to the first summand and Theorem~\ref{theorem:GaussianComparison-PQ-Simultaneous} to the second summand.
\end{proof}

\begin{proof}[\textbf{Proof of Proposition~\ref{lemma:KarhunenLoeve-GP}}]
	
	Our proof is modeled after the proof of Theorem 1 in~\cite{jain1970note}. However, unlike them we do not argue via the reproducing kernel Hilbert space associated to $\mathcal{E}_Q$. Instead, we leverage the fact that under the stated assumptions $Z_Q$ has a version that is both a mean-square continuous stochastic process and a random element on some Hilbert space. 
	
	Since $\mathcal{F}_n$ is compact it is totally bounded and, by Lemma~\ref{lemma:Sudakov}, separable w.r.t $e_Q$. Hence, the process $Z_Q = \{Z_Q(f) : f \in \mathcal{F}_n\}$ has a separable and jointly measurable version $\widetilde{Z}_Q$~\citep[e.g.][Proposition 2.1.12, note that the $e_Q$ is the intrinsic standard deviation metric of the Gaussian $Q$-motion $Z_Q$]{gine2016mathematical}. Since $\mathcal{F}_n$ is compact and $\mathcal{E}_Q$ is continuous, $T_{\mathcal{E}_Q}$ is a bounded linear operator. Let $\{\lambda_k,\varphi_k )\}_{k=1}^\infty$ be the eigenvalue and eigenfunction pairs of $T_{\mathcal{E}_Q}$ and define
	\begin{align}\label{eq:lemma:KarhunenLoeve-GP-4}
		\widetilde{Z}_Q^m(f) := \sum_{k=1}^m \langle \widetilde{Z}_Q, \varphi_k \rangle \varphi_k(f), \quad f \in \mathcal{F}_n.
	\end{align}
	By continuity of $\mathcal{E}_Q$, $\widetilde{Z}_Q$ is a mean-square continuous stochastic process in $L_2(\mathcal{F}_n, \mathcal{B}_n, \mu)$. Moreover, by joint measurability, $\widetilde{Z}_Q$ is also a random element on the Hilbert space $\big(L_2(\mathcal{F}_n, \mathcal{B}_n, \mu), \langle \cdot, \cdot \rangle\big)$. Hence, 
	Theorems 7.3.5 and 7.4.3 in~\cite{hsing2015theoretical} apply, and we conclude that the partial sums $\widetilde{Z}_Q^m(f)$ converge to $\widetilde{Z}_Q(f)$ in $L_2(\Omega, \mathcal{A}, \mathbb{P})$ as $m \rightarrow \infty$ pointwise in $f \in \mathcal{F}_n$.
	
	Now, observe that
	\begin{align}\label{eq:lemma:KarhunenLoeve-GP-3}
		\mathrm{E}[ \langle \widetilde{Z}_Q, \varphi_k \rangle ] = 0 \quad \quad \mathrm{and} \quad \quad \mathrm{E}[ \langle \widetilde{Z}_Q, \varphi_k \rangle \langle \widetilde{Z}_Q, \varphi_j \rangle] = \lambda_k \mathbf{1}\{j = k\}, \quad \forall j, k \in \mathbb{N}.
	\end{align}
	Thus, the $\langle \widetilde{Z}_Q, \varphi_k \rangle$'s are uncorrelated random variables. Since the $\langle \widetilde{Z}_Q, \varphi_k \rangle$'s are necessarily Gaussian (inner product of a Gaussian random element with a deterministic function $\varphi_k$!), they are in fact independent Gaussian random variables with mean zero and variance $\lambda_k$. Therefore, by L{\'e}vy's theorem, $\widetilde{Z}_Q^m(f)$ converges to $\widetilde{Z}_Q(f)$ almost surely pointwise in $f \in \mathcal{F}_n$. Thus, for $f \in \mathcal{F}$,
	\begin{align}\label{eq:lemma:KarhunenLoeve-GP-2}
		\widetilde{Z}_Q(f)= \widetilde{Z}_Q^\infty(f) \quad a.s.
	\end{align}
	
	Since $\mathrm{E}\|Z_Q\|_{\mathcal{F}_n} < \infty$, Lemma~\ref{lemma:AS-Bounded} implies that $\widetilde{Z}_Q$ is almost surely bounded. Moreover, since $\mathcal{E}_Q$ is continuous and $\mathcal{F}_n$ compact, Lemma 4.6.6 (4) in~\cite{hsing2015theoretical} implies that $\widetilde{Z}_Q$ is uniformly continuous. Consequently, by construction, $\widetilde{Z}_Q^m$ is almost surely bounded and uniformly continuous, too. Thus, $\widetilde{Z}_Q$ and $\widetilde{Z}_Q^m$ can be considered random elements on the Banach space $C(\mathcal{F}_n)$ of bounded continuous functions on $\mathcal{F}_n$
	equipped with the supremum norm. 
	
	Recall that the dual space $C^*(\mathcal{F}_n)$ is the space of finite signed Borel measures $\mu$ on $\mathcal{F}_n$ equipped with the total variation norm. The dual pairing between $X \in C(\mathcal{F}_n)$ and $\mu \in C^*(\mathcal{F}_n)$ is $\langle X, \mu\rangle := \int_{\mathcal{F}_n} X(f) d \mu(f)$. We compute
	\begin{align}\label{eq:lemma:KarhunenLoeve-GP-1}
		\mathrm{E}\left|\langle \widetilde{Z}_Q^m, \mu\rangle  - \langle \widetilde{Z}_Q, \mu\rangle\right| & = \mathrm{E} \left| \int_{\mathcal{F}_n} \big(  \widetilde{Z}_Q^m(f) -  \widetilde{Z}_Q(f)\big) d \mu(f)\right| \nonumber\\
		&\leq \sup_{f \in \mathcal{F}_n}  \mathrm{E} \left[ \big| \widetilde{Z}_Q^m(f) -  \widetilde{Z}_Q(f)\big| \right] \|\mu\|_{TV}\nonumber\\
		&\leq \sup_{f \in \mathcal{F}_n}  \mathrm{E} \left[ \big( \widetilde{Z}_Q^m(f) -  \widetilde{Z}_Q(f)\big)^2 \right]^{1/2} \|\mu\|_{TV}\nonumber\\
		&\overset{(a)}{=} \sup_{f \in \mathcal{F}_n}  \mathrm{E} \left[ \big( \widetilde{Z}_Q^m(f) -  \widetilde{Z}_Q^\infty(f)\big)^2 \right]^{1/2} \|\mu\|_{TV}\nonumber\\
		&= \sup_{f \in \mathcal{F}_n}   \left( \sum_{k=m+1}^\infty \lambda_k \varphi_k^2(f)\right)^{1/2} \|\mu\|_{TV}\nonumber\\
		&\overset{(b)}{\rightarrow} 0 \quad \text{as} \quad m \rightarrow \infty,
	\end{align}
	where (a) holds by~\eqref{eq:lemma:KarhunenLoeve-GP-2} and (b) holds by Lemma 4.6.6 (3) in~\cite{hsing2015theoretical}. Since~\eqref{eq:lemma:KarhunenLoeve-GP-1} implies convergence of the dual pairings in probability for every $\mu \in C^*(\mathcal{F}_n)$, we have, by It{\^o}-Nisio's theorem,
	\begin{align*}
		\big\| \widetilde{Z}_Q^m - \widetilde{Z}_Q \big\|_{\mathcal{F}_n} \rightarrow 0 \quad \text{as} \quad m \rightarrow \infty
	\end{align*}
	almost surely. Recall that $\widetilde{Z}_Q$ is a version of $Z_Q$. This completes the proof.
\end{proof}

\begin{proof}[\textbf{Proof of Corollary~\ref{corollary:theorem:Abstract-Bootstrap-Sup-Empirical-Process-1}}]
	By Theorem~\ref{theorem:Abstract-Bootstrap-Sup-Empirical-Process} and arguments as in the proof of Corollary~\ref{corollary:theorem:CLT-Max-Norm-Simultaneous} we have, for each $M \geq 0$,
	\begin{align*}
		&\sup_{s \geq 0} \Big|\mathbb{P}\left\{\|\mathbb{G}_n\|_{\mathcal{F}_n} \leq s \right\} -\mathbb{P}\left\{\|\widehat{Z}_n^m\|_{\mathcal{F}_n} \leq s \mid X_1, \ldots, X_n \right\}\Big| \\
		&\quad{}\quad{}\quad{} \lesssim \frac{\|F_n\|_{P, 3}}{\sqrt{n^{1/3}\mathrm{Var}(\|G_P\|_{\mathcal{F}_n})}} + \frac{\|F_n\mathbf{1}\{F_n > M\}\|_{P,3}^3}{\|F_n\|_{P,3}^3} + \frac{ \mathrm{E} \|G_P\|_{\mathcal{F}_n} + M}{\sqrt{n\mathrm{Var}(\|G_P\|_{\mathcal{F}_n})}}\\
		&\quad{}\quad{}\quad{}\quad{} \quad{} \quad{} + \left( \frac{\sup_{f, g \in \mathcal{F}_n} \big| \mathcal{C}_P(f,g) - \widehat{\mathcal{C}}_n^m(f, g) \big|}{ \mathrm{Var}(\|G_P\|_{\mathcal{F}_n} )} \right)^{1/3},
	\end{align*}
	where $\lesssim$ hides an absolute constant independent of $n, m, M, \mathcal{F}_n, F_n, P_n$, and $P$.

	Above approximation error can be further upper bounded by
	\begin{align*}
		\sup_{f, g \in \mathcal{F}_n} \big| \mathcal{C}_P(f,g) - \widehat{\mathcal{C}}_n^m(f, g) \big| \leq \sup_{f, g \in \mathcal{F}_n} \big| \mathcal{C}_P(f,g) - \widehat{\mathcal{C}}_n(f, g) \big| + \sup_{f, g \in \mathcal{F}_n} \big| \widehat{\mathcal{C}}_n(f,g) - \widehat{\mathcal{C}}_n^m(f, g) \big|.
	\end{align*}
	If $\mathcal{F}_n$ is compact w.r.t. $d_{\widehat{\mathcal{C}}_n}$, then by Mercer's theorem~\citep[e.g.][Lemma 4.6.6]{hsing2015theoretical}, 
	\begin{align*}
		\sup_{f, g \in \mathcal{F}_n} \big| \widehat{\mathcal{C}}_n(f,g) - \widehat{\mathcal{C}}_n^m(f, g) \big| = \sup_{f, g \in \mathcal{F}_n}\sum_{k=m+1}^\infty  \left|\widehat{\lambda}_k \widehat{\varphi}_k(f)\widehat{\varphi}_k(g)\right| \rightarrow 0 \quad \text{as} \quad m \rightarrow \infty.
	\end{align*}
	This completes the proof of the corollary. 
\end{proof}

\begin{proof}[\textbf{Proof of Theorem~\ref{theorem:Abstract-Bootstrap-Sup-Empirical-Process-Quantiles}}]
	The proof is identical to the proof of Proposition~\ref{lemma:Bootstrap-Max-Norm-Asymptotic-Size}. The only differences are the notation and that we use Lemma~\ref{lemma:Bootstrap-Sup-Empirical-Process-Quantil-Comparison} instead of Lemma~\ref{lemma:Bootstrap-Max-Norm-Quantil-Comparison} and Corollary~\ref{corollary:theorem:CLT-Max-Norm-Simultaneous}/ Theorem~\ref{theorem:Abstract-Bootstrap-Sup-Empirical-Process} instead of Proposition~\ref{lemma:CLT-Max-Norm}.
\end{proof}

\newpage
\section{Proofs of the results in Section~\ref{sec:Applications}}\label{sec:Proofs-Applications}

\begin{proof}[\textbf{Proof of Proposition~\ref{theorem:ConfidenceEllipsoid}}]
	Since the problem is finite dimensional we do not have to develop the Karhunen-Lo{\`e}ve expansion from Section~\ref{subsec:GaussianProcessBootstrap-Implementation-Consistency}. The covariance function $\mathcal{C}_\psi$ of the empirical process $ \frac{1}{\sqrt{n}} \sum_{i=1}^n \psi_i'u$ with $u \in \{v \in \mathbb{R}^d :v \in S^{d-1}\}$ has the explicit form $(u,v) \mapsto \mathcal{C}_\psi(u,v) = u' \mathrm{E}[\psi_1\psi_1'] v$. Hence, a natural estimate of $\mathcal{C}_\psi$ is $(u,v) \mapsto \widehat{\mathcal{C}}_\psi(u,v) = u' \widehat{\Omega}_\psi v$. A version of a centered Gaussian process defined on $\{v \in \mathbb{R}^d : \|v\|_2 =1\}$ and with covariance function $\widehat{\mathcal{C}}_\psi$ is $\{\widehat{Z}_\psi'u : u \in S^{d-1}\}$ where $\widehat{Z}_\psi \mid  X_1, \ldots, X_n \sim N(0, \widehat{\Omega}_\psi)$. It is standard to verify that, under the assumptions of the theorem, the metric entropy integrals associated to the Gaussian processes $G_P$ and $\widehat{Z}_\psi$ are finite for every fixed $n$ and $d$. Thus, by Lemma~\ref{lemma:MetricEntropCondition} there exist versions of $G_P$ and $\widehat{Z}_\psi$ which are almost surely bounded and almost surely uniformly continuous. 
	Hence, the modulus of continuity condition~\eqref{eq:theorem:Abstract-Bootstrap-Sup-Empirical-Process-0} holds for these versions and we can take $r_n = 0$. It follows by Theorem~\ref{theorem:Abstract-Bootstrap-Sup-Empirical-Process} that, for all $\alpha \in (0,1)$,
	\begin{align*}
		&\sup_{\alpha \in (0,1)} \Big|\mathbb{P} \left\{	\sqrt{n} \|\hat{\theta}_n - \theta_0\|_2 \leq c_n(\alpha; \widehat{\Omega}_\psi)  \right\}  - \alpha \Big| \\
		&\quad{}\quad{}\lesssim \frac{(\mathrm{E}[\|\psi_1\|_2^3])^{1/3}}{n^{1/6}\sqrt{\mathrm{Var}(\| \|Z_\psi\|_2\|_p)}} + \frac{\mathrm{E} \left[ \|\psi_1\|_2^3 \mathbf{1}\{\|\psi_1\|_2^3 > n\: \mathrm{E}[\|\psi_1\|_2^3]\}\right]}{\mathrm{E}\left[ \|\psi_1\|_2^3\right]}+  \frac{ \mathrm{E}[\|Z_\psi\|_2]}{\sqrt{n\mathrm{Var}(\|Z_\psi\|_2)}}\\
		& \quad{}\quad{}\quad{}\quad{}+ \inf_{\delta > 0}\left\{ \left(\frac{\delta }{ \mathrm{Var}( \|Z_\psi\|_2)} \right)^{1/3}  + \mathbb{P} \left\{ \|\widehat{\Omega}_\psi-  \mathrm{E}[\psi_1\psi_1']\|_{op} > \delta\right\}\right\} \\
		&\quad{}\quad{}\quad{} \quad{} \quad{}+ \inf_{\eta > 0} \left\{ \frac{\eta }{ \sqrt{\mathrm{Var}( \|Z_\psi\|_2)} } + \mathbb{P} \left\{|\Theta_n| > \eta \right\}\right\}.
	\end{align*}
	Under Assumption~\ref{assumption:SubGaussianInfluenceFunction} or~\ref{assumption:HeavyTailedInfluenceFunction} the right hand side in above inequality vanishes as $n$ diverges. For details we refer to Appendix A.1 in~\cite{giessing2023bootstrap}. This completes the proof.
\end{proof}

\begin{proof}[\textbf{Proof of Proposition~\ref{theorem:ConsistencyOperatorNorm}}]
	Since the problem is finite dimensional we do not have to develop the Karhunen-Lo{\`e}ve expansion from Section~\ref{subsec:GaussianProcessBootstrap-Implementation-Consistency}. Consider $T_n \equiv \sqrt{n}\|Q_n - Q\|_{{\mathcal{F}_n}}$ where $Q_n$  is the empirical measure of the collection of $Y_i = \mathrm{vech}'(X_iX_i') H_d'$, $1 \leq i \leq n$, $Q$ the pushforward of $P$ under the map $X \mapsto \mathrm{vech}'(XX') H_d'$, and $\mathcal{F}_n =\{ y \mapsto  y(v \otimes u) : u, v \in S^{d-1}\}$. Since each $f \in \mathcal{F}_n$ has a (not necessarily unique) representation in terms of $u, v \in S^{d-1}$, we will identify $f \in \mathcal{F}_n$ with pairs $(u,v) \in S^{d-1} \times S^{d-1}$ when there is no danger of confusion. Clearly, $F_n(y) = \|y\|_2^2$ is an envelope function for $\mathcal{F}_n$. 
	
	Define the Gaussian processes $\{\widehat{Z}_n^m(f) : f \in \mathcal{F}_n\} \equiv \{ \widehat{Z}_n' H_d' (v\otimes u): u, v \in S^{d-1}\}$, where $\widehat{Z}_n \mid X_1, \ldots, X_n \sim N(0, \widehat{\Omega}_n)$ with $\widehat{\Omega}_n = n^{-1} \sum_{i=1}^n \mathrm{vech}(X_iX_i'-\widehat{\Sigma}_n) \mathrm{vech}'(X_iX_i'-\widehat{\Sigma}_n) \in  \mathbb{R}^{d(d+1)/2 \times d(d+1)/2}$, and $\{G_Q(f) : f \in \mathcal{F}_n\} \equiv \{Z'H_d' (v\otimes u): u, v \in S^{d-1}\}$, where $Z \sim N(0, \Omega)$ with $\Omega = \mathrm{E}[ \mathrm{vech}(XX'-\Sigma) \mathrm{vech}'(XX'-\Sigma) ] \in  \mathbb{R}^{d(d+1)/2 \times d(d+1)/2}$. 
	
	Without loss of generality, we can assume that $\mathrm{tr}(\Omega)$ and $\|\Omega\|_{op}$ are finite since otherwise the statement is trivially true. Thus, under the assumptions of the theorem, the metric entropy integrals associated to the Gaussian processes $\widehat{Z}_n^m$ and $G_Q$ are finite for every fixed $n$ and $d$. Therefore by Lemma~\ref{lemma:MetricEntropCondition} there exist versions of $\widehat{Z}_n^m$ and $G_Q$ (which we also denote by $\widehat{Z}_n^m$ and $G_Q$) that are almost surely bounded and almost surely uniformly continuous. The modulus of continuity condition~\eqref{eq:theorem:Abstract-Bootstrap-Sup-Empirical-Process-0} holds for these versions and, in particular, we can take $r_n = 0$. By Theorem~\ref{theorem:Abstract-Bootstrap-Sup-Empirical-Process},
	\begin{align}\label{eq:theorem:QuantilesOperatorNorm-1}
		&\sup_{s \geq 0} \Big|\mathbb{P}\left\{T_n \leq s \right\} -\mathbb{P}\left\{ \|\widehat{S}_n\|_{op}\leq s \mid X_1, \ldots, X_n \right\}\Big| \nonumber\\
		\begin{split}
			&\quad{}\quad{}\quad{} \lesssim \frac{\|F_n\|_{P, 3}}{\sqrt{n^{1/3}\mathrm{Var}(\|G_P\|_{\mathcal{F}_n})}} + \frac{\|F_n\mathbf{1}\{F_n > n^{1/3} \|F_n\|_{P, 3}\}\|_{P,3}^3}{\|F_n\|_{P,3}^3} + \frac{ \mathrm{E} \|G_P\|_{\mathcal{F}_n}}{\sqrt{n\mathrm{Var}(\|G_P\|_{\mathcal{F}_n})}}\\		
			&\quad{}\quad{}\quad{}\quad{} \quad{} \quad{} + \left( \frac{ \sup_{u,v \in S^{d-1} } \big| (u \otimes v)'H_d(\Omega - \widehat{\Omega}_n)H_d'(u \otimes v) \big|}{ \mathrm{Var}(\|G_P\|_{\mathcal{F}_n} )} \right)^{1/3}.
		\end{split}
	\end{align}
	
	In the remainder of the proof we derive upper bounds on the quantities on the right-hand side in above display.
	
	Note that
	\begin{align}\label{eq:theorem:QuantilesOperatorNorm-2}
		\mathrm{tr}(\Sigma) = \mathrm{E}\big[ \|X\|_2^2\big] \overset{(a)}{\leq} \mathrm{E}\big[ \|X\|_2^6\big]^{1/3} = \|F_n\|_{P,3}  \lesssim \| F_n(X)\|_{\psi_1} \overset{(b)}{\lesssim} \mathrm{tr}(\Sigma),
	\end{align}
	where (a) follows from H{\"o}lder's inequality and (b) holds because $X = (X_1, \ldots, X_d)' \in \mathbb{R}^d$ is sub-Gaussian with mean zero and covariance $\Sigma$ and therefore
	\begin{align*}
		\|F_n(X)\|_{\psi_1} = \left\| X'X \right\|_{\psi_1} \leq \sum_{k=1}^d \left\| X_k^2  \right\|_{\psi_1} \leq \sum_{k=1}^d \left\| X_k \right\|_{\psi_2}^2 \lesssim \sum_{k=1}^d \mathrm{Var}(X_k) = \mathrm{tr}(\Sigma).
	\end{align*}
	
	Next, note that $\|u \otimes v\|_2^2  = (u' \otimes v' ) (u \otimes v) = u'u \otimes v'v = 1$ for all $u,v \in S^{d-1}$. Thus, we compute
	\begin{align*}
		\mathrm{E} \|G_P\|_{\mathcal{F}_n}  = \mathrm{E} \|S\|_{op} &= \mathrm{E} \left[\sup_{u , v \in S^{d-1}} | \mathrm{vec}(S)'(u \otimes v) | \right]\\
		&\leq \mathrm{E} \| \mathrm{vec}(S)\|_2 \leq \left(\mathrm{E} \| \mathrm{vec}(S)\|_2^2\right)^{1/2} = \sqrt{\mathrm{tr}(H_d\Omega H_d')}.
	\end{align*}
	Also, since $\mathrm{vec}(A) = H_d\mathrm{vech}(A)$ for all symmetric matrices $A$, we have
	\begin{align*}
		\mathrm{tr}(H_d\Omega H_d') = \mathrm{E}\left[ \mathrm{vec}(XX' - \Sigma)' \mathrm{vec}(XX' - \Sigma) \right] = \mathrm{E}\left[ \mathrm{tr}\left( (XX' - \Sigma)(XX' - \Sigma) \right)\right]\\
		= \mathrm{E}\left[ \mathrm{tr}\left( XX'XX'\right)\right] - \mathrm{tr}(\Sigma^2)= \mathrm{E}[\|X\|_2^4] - \mathrm{tr}(\Sigma^2) \lesssim \mathrm{tr}^2(\Sigma) - \mathrm{tr}(\Sigma^2),
	\end{align*}
	where the last inequality follows from similar arguments as used to derive the upper bound in~\eqref{eq:theorem:QuantilesOperatorNorm-2}. Thus,
	\begin{align}\label{eq:theorem:QuantilesOperatorNorm-4}
		\mathrm{E} \|G_P\|_{\mathcal{F}_n}  \lesssim \sqrt{ \mathrm{tr}^2(\Sigma) - \mathrm{tr}(\Sigma^2)}.
	\end{align}
	Moreover, since $(u \otimes v)' H_d \mathrm{vech}(A) = v'Au$ for all symmetric matrices $A$ and matching vectors $u,v$, we have, for arbitrary $f \in \mathcal{F}_n$, 
	\begin{align}\label{eq:theorem:QuantilesOperatorNorm-5}
		\mathrm{Var}(G_P(f)) &\geq \inf_{u, v \in S^{d-1}} \mathrm{Var}\left( \mathrm{vec}(S)'(u \otimes v) \right) \nonumber\\
		&= \inf_{u, v \in S^{d-1}} (u \otimes v)'H_d\Omega H_d' (u \otimes v) 	= \inf_{u, v \in S^{d-1}} \mathrm{E}\left[ \big(u'(XX'- \Sigma)v\big)^2  \right]\nonumber\\
		&= \inf_{u \in S^{d-1}}  \mathrm{E}[(X'u)^4] - \|\Sigma u\|_2^2 =  \inf_{u \in S^{d-1}} \mathrm{Var}\big( (X'u)^2 \big) \nonumber\\
		&\geq \kappa.
	\end{align}
	Combine eq.~\eqref{eq:theorem:QuantilesOperatorNorm-4}--\eqref{eq:theorem:QuantilesOperatorNorm-5} with Lemma~\ref{lemma:BoundsVariance-SeparableProcess} to obtain
	\begin{align}\label{eq:theorem:QuantilesOperatorNorm-6}
		\mathrm{Var}(\|G_P\|_{\mathcal{F}_n}) \gtrsim \left(\frac{ \kappa}{1 + \mathrm{E}[\|G_P\|_{\mathcal{F}_n}/\kappa]}\right)^2 \gtrsim \left(\frac{ \kappa^2}{\kappa + \sqrt{\mathrm{tr}^2(\Sigma) - \mathrm{tr}(\Sigma^2)}}\right)^2  \gtrsim \left(\frac{ \kappa^2}{\kappa + \mathrm{tr}(\Sigma)}\right)^2.
	\end{align}
	
	Next, since (again!) $(u \otimes v)' H_d \mathrm{vech}(A) = v'Au$ for all symmetric matrices $A$ and matching vectors $u,v$, we compute
	\begin{align}\label{eq:theorem:QuantilesOperatorNorm-7}
		&\sup_{u,v \in S^{d-1} } \left| (u \otimes v)'H_d(\Omega - \widehat{\Omega}_n)H_d(u \otimes v) \right| \nonumber\\
		&\quad \leq \sup_{u,v \in S^{d-1} } \left| \frac{1}{n}\sum_{i=1}^n (X_i'u)^2(X_i'v)^2 - \mathrm{E}[(X_i'u)^2(X_i'v)^2 ] \right| + \sup_{u,v \in S^{d-1} } u' ( \widehat{\Sigma} - \Sigma ) v \nonumber\\
		&\quad = O_p\left(  \mathrm{r}(\Sigma) \|\Sigma\|_{op}^2  \left( \sqrt{\frac{(\log en)^2  \mathrm{r}(\Sigma)}{n}} \vee \frac{(\log en)^2  \mathrm{r}(\Sigma)}{n} \right) \right),
	\end{align}
	where the last line holds since by Markov's inequality and Lemma~\ref{lemma:Tesseract},
	\begin{align*}
		&\sup_{u,v \in S^{d-1} } \left| \frac{1}{n}\sum_{i=1}^n (X_i'u)^2(X_i'v)^2 - \mathrm{E}[(X_i'u)^2(X_i'v)^2 ] \right|\\
		&\quad\quad = O_p \left(  \mathrm{r}(\Sigma) \|\Sigma\|_{op}^2  \left( \sqrt{\frac{(\log en)^2  \mathrm{r}(\Sigma)}{n}} \vee \frac{(\log en)^2  \mathrm{r}(\Sigma)}{n} \right) \right),
	\end{align*}
	and by Theorem 4 in~\cite{koltchinskiiConcentration2017}
	\begin{align*}
		\sup_{u,v \in S^{d-1} } u' ( \widehat{\Sigma} - \Sigma ) v = O_p \left( \|\Sigma\|_{op}  \left( \sqrt{\frac{\mathrm{r}(\Sigma)}{n}} \vee \frac{\mathrm{r}(\Sigma)}{n} \right) \right).
	\end{align*}
	
	Lastly, for $r > 0$, by the upper and lower bounds in~\eqref{eq:theorem:QuantilesOperatorNorm-2},
	\begin{align}\label{eq:theorem:QuantilesOperatorNorm-8}
		\frac{\|F_n\mathbf{1}\{F_n > n^{1/3} \|F_n\|_{P, 3}\}\|_{P,3}^3}{\|F_n\|_{P,3}^3} \leq \frac{\mathrm{E}[F_n^{3+r}]}{n^{r/3} \mathrm{E}[F_n^3]^{r/3+1}} =  \frac{\|F_n\|_{P, 3+r}^{3+r}}{n^{r/3} \|F_n\|_{P, 3}^{3+r}} \lesssim n^{-r/3}.
	\end{align}
	
	Now, combining~\eqref{eq:theorem:QuantilesOperatorNorm-1}--\eqref{eq:theorem:QuantilesOperatorNorm-8}, we obtain for $r =1$ and under Assumption~\ref{assumption:SubGaussianData},
	\begin{align*}
		&\sup_{s \geq 0} \Big|\mathbb{P}\left\{T_n \leq s \right\} -\mathbb{P}\left\{ \|\widehat{S}_n\|_{op}\leq s \mid X_1, \ldots, X_n \right\}\Big| \\
		&\quad{}\lesssim \frac{\mathrm{tr}(\Sigma)}{ n^{1/6} \kappa } +  \frac{\mathrm{tr}^2(\Sigma)}{ n^{1/6} \kappa^2 } + n^{-1/3} + \frac{ \sqrt{ \mathrm{tr}^2(\Sigma) - \mathrm{tr}(\Sigma^2)} }{n^{1/2} \kappa } + \frac{ \mathrm{tr}^2(\Sigma) - \mathrm{tr}(\Sigma^2) }{n^{1/2}\kappa^2 }\\	
		& \quad{}\quad{} + O_p\left(  \left( \mathrm{r}(\Sigma) \|\Sigma\|_{op}^2\left( \sqrt{\frac{(\log en)^2  \mathrm{r}(\Sigma)}{n}} \vee \frac{(\log en)^2  \mathrm{r}(\Sigma)}{n} \right) \right)^{1/3} \left(\frac{\kappa + \mathrm{tr}(\Sigma)}{\kappa^2}\right)^{2/3} \right)\\
		&\quad{}\lesssim  \frac{\|\Sigma\|_{op}}{\kappa } \frac{r(\Sigma)}{ n^{1/6}}  \vee \frac{\|\Sigma\|_{op}^2}{\kappa^2 } \frac{r^2(\Sigma)}{ n^{1/6} }\\
		&\quad{}\quad{} + O_p\left(\left( \frac{r^{1/3}(\Sigma) \|\Sigma\|_{op}^{2/3}}{\kappa^{2/3}} \vee \frac{r(\Sigma) \|\Sigma\|_{op}^{4/3}}{\kappa^{4/3}}  \right) \left( \sqrt{\frac{(\log en)^2  \mathrm{r}(\Sigma)}{n}} \vee \frac{(\log en)^2  \mathrm{r}(\Sigma)}{n} \right)^{1/3} \right).
	\end{align*}
	To complete the proof, adjust some constants.
\end{proof}

\begin{proof}[\textbf{Proof of Proposition~\ref{theorem:Uniform-CI-Bands-RKHS}}]	
	Let $\mathcal{F}_n = \{ v \mapsto \langle  v, u\rangle_\mathcal{H} : u \in S^* \}$. Further, let $Q_n$ be the empirical measure based on the $V_i = (T + \lambda)^{-2} T \big(Y_i - f_0(X_i)\big) k_{X_i}$, $1 \leq i \leq n$ and $Q$ the pushforward measure of $P$ under the map $(Y, X) \mapsto V =  (T + \lambda)^{-2} T \big(Y - f_0(X)\big) k_{X}$. Since $f_0$ is the best approximation in square loss, $\mathrm{E}[Y - f_0(X)]= 0$ and $V$ and the $V_i$'s have mean zero. Hence,
	\begin{align*}
		\sqrt{n} \|\widehat{f}^{\mathrm{bc}}_n - f_0\|_{\infty} = \sup_{u \in S^*} \left| \left \langle \frac{1}{\sqrt{n}} \sum_{i=1}^n V_i, u \right\rangle_\mathcal{H} \right|  + \Theta_n \equiv \sqrt{n}\|Q_n - Q\|_{\mathcal{F}_n} + \Theta_n,
	\end{align*}	
	where $\Theta_n$ is a reminder term which satisfies $|\Theta_n| \leq \sqrt{n}\|R_n\|_\infty$. We plan to apply Theorem~\ref{theorem:Abstract-Bootstrap-Sup-Empirical-Process-Quantiles} to the far right hand side in above display. To this end, we need to (i) find an envelope for the function class $\mathcal{F}_n$, (ii) establish (ii) compactness of $\mathcal{F}_n$, (iii) construct certain Gaussian processes (to be defined below) with continuous covariance functions, and (iv) establish almost sure uniform continuity of these processes w.r.t. their intrinsic standard deviation metrics.
	
	Recall that $\sqrt{k(x, y)} \leq \kappa$ for all $x, y \in S$ and that $\langle v, x^*\rangle_{\mathcal{H}} = \langle v, k_x\rangle_{\mathcal{H}} $ for all $v \in \mathcal{H}$ and $x \in S$. Hence, $|\langle v, u \rangle_{\mathcal{H}}| \leq \|v\|_{\mathcal{H}} \sup_{x \in \mathcal{H}} \| k_x\|_{\mathcal{H}} \leq \kappa \|v\|_{\mathcal{H}}$ for all $v \in \mathcal{H}$ and $u \in S^*$. Thus, $v\mapsto \kappa\|v\|_{\mathcal{H}}$ is an envelope of the function class $\mathcal{F}_n$. 	
	
	Next, let $Z$ and $\widehat{Z}_n^m$ be centered Gaussian random elements on $\mathcal{H}$ with covariance operators $\Omega$ and $\widehat{\Omega}_n^m$, respectively, i.e. for all $u, v \in \mathcal{H}$,
	\begin{align*}
		&\langle Z, u\rangle_{\mathcal{H}} \sim N\big(0, \mathcal{C}(u,v) \big), \quad \quad \text{where} \quad\quad \mathcal{C}(u,v) = \langle \Omega u, v\rangle_{\mathcal{H}}, \\
		&\langle \widehat{Z}_n, u\rangle_{\mathcal{H}} \sim N\big(0, \mathcal{\widehat{C}}_n^m(u,v)\big),\quad \quad \text{where} \quad \quad  \mathcal{\widehat{C}}_n^m(u,v) = \langle \widehat{\Omega}_n^m u, v\rangle_{\mathcal{H}}.
	\end{align*}
	By Cauchy-Schwarz these covariance functions are obviously continuous w.r.t. $\| \cdot\|_{\mathcal{H}}$. Denote the standard deviation metrics associated with above Gaussian random elements by $d_\mathcal{C}$ and $d_{\widehat{\mathcal{C}}_n^m}$. Then, for all $x, y \in S$,
	\begin{align*}
		d_\mathcal{C}^2(x^*, y^*) &= \mathrm{E}\big[(Z(x^*) - Z(y^*))^2\big] =  \mathrm{E}\big[\langle Z, x^* - y^*\rangle_{\mathcal{H}}^2\big] \leq \mathrm{E}\|Z\|_{\mathcal{H}}^2\|x^* - y^*\|_{\mathcal{H}}^2 =  \mathrm{tr}(\Omega)d_k^2(x,y),
	\end{align*}
	and, completely analogous,
	\begin{align*}
		d_{\widehat{\mathcal{C}}_n^m}^2(x^*, y^*)  \leq \mathrm{tr}(\widehat{\Omega}_n^m)d_k^2(x,y).
	\end{align*}		
	Since $(S, d_k)$ is compact and $\mathrm{tr}(\Omega) \vee \mathrm{tr}(\widehat{\Omega}_n) < \infty$ by assumption, these inequalities imply that both $S^*$ and (by continuity of the inner product) $\mathcal{F}_n$ are compact w.r.t. the standard deviation metrics $d_\mathcal{C}$ and $d_{\widehat{\mathcal{C}}_n^m}$. Moreover, since $\int_0^\infty \sqrt{N(S, d_k, \varepsilon)} d\varepsilon < \infty$, these inequalities also imply that $\int_0^\infty \sqrt{N(\mathcal{F}_n, d_\mathcal{C}, \varepsilon)} d\varepsilon < \infty$ and $\int_0^\infty \sqrt{N(\mathcal{F}_n, d_{\widehat{\mathcal{C}}_n}, \varepsilon)} d\varepsilon < \infty$. Hence, by Lemma~\ref{lemma:MetricEntropCondition} there exist versions of $Z$ and $\widehat{Z}_n$ that are almost surely bounded and has almost surely uniformly $d_\mathcal{C}$- and $d_{\widehat{\mathcal{C}}_n}$-continuous sample paths. In the following, we keep using $Z$ and $\widehat{Z}_n^m$ to denote these versions.
	
	Hence, Theorem~\ref{theorem:Abstract-Bootstrap-Sup-Empirical-Process-Quantiles} applies and we have
	\begin{align}\label{eq:theorem:Uniform-CI-Bands-RKHS-11}
		\begin{split}
		&\sup_{\alpha \in (0,1)} \Big|\mathbb{P}\left\{ \sqrt{n}\|\widehat{f}^{\mathrm{bc}}_n - f_0\|_\infty \leq 	c_n^m(\alpha) \right\}  - \alpha \Big| \\
		&\quad{}\quad{}\lesssim \frac{\kappa \mathrm{E}[\|V\|_{\mathcal{H}}^3]^{1/3}}{\sqrt{n^{1/3}\mathrm{Var}(\|Z\|_{\mathcal{F}_n})}} + \frac{\mathrm{E}[\|V\|_{\mathcal{H}}^3\mathbf{1}\{\|V\|_{\mathcal{H}}>  n^{1/3} \mathrm{E}[\|V\|_{\mathcal{H}}^3]^{1/3}\}]}{ \mathrm{E}[\|V\|_{\mathcal{H}}^3]} + \frac{ \mathrm{E} \|Z\|_{\mathcal{F}_n} }{\sqrt{n\mathrm{Var}(\|Z\|_{\mathcal{F}_n})}}\\
		& \quad{}\quad{}\quad{}  + \inf_{\delta > 0}\left\{ \left(\frac{\delta }{\mathrm{Var}(\|Z\|_{\mathcal{F}_n})} \right)^{1/3}  + \mathbb{P}\left\{ \big\| \Omega - \widehat{\Omega}_n^m\big\|_{op} > \delta\right\}\right\}\\
		&\quad{}\quad{}\quad{} + \inf_{\eta > 0} \left\{ \frac{\eta }{ \sqrt{\mathrm{Var}(\|Z\|_{\mathcal{F}_n})} } + \mathbb{P}\left\{ \sqrt{n} \|R_n\|_\infty > \eta \right\} \right\}.
		\end{split}
	\end{align}	
	In the remainder of the proof we derive upper bounds on the quantities on the right-hand side in above display.
	
	First, from the proof of Lemma~\ref{lemma:BC-KRR-Expansion} we know that
	\begin{align*}
		V = (T+\lambda)^{-1} \big( Y_i - f_0(X_i)\big)  k_X - \lambda (T+\lambda)^{-2} \big( Y_i - f_0(X_i)\big) k_X.
	\end{align*}
	Hence, by the proof of Lemma~\ref{lemma:BC-KRR-Consistency-Operator} 
	\begin{align*}
		\|V\|_{\mathcal{H}} &\leq  \left| Y_i - f_0(X_i)\right| \left\|(T + \lambda)^{-1}k_{X_i}\right\|_\mathcal{H} + \lambda \left| Y_i - f_0(X_i)\right| \left\|(T + \lambda)^{-2}k_{X_i}\right\|_\mathcal{H} \\
		&\leq \lambda^{-1} 2\kappa (B + \kappa\|f_0\|_{\mathcal{H}}) \quad a.s.
	\end{align*}
	and
	\begin{align*}
		\mathrm{E}\|V\|_{\mathcal{H}}^2 & \leq 2\sigma_0^2 \mathrm{E}\left[\left\|(T + \lambda)^{-1}k_{X_i}\right\|_{\mathcal{H}}^2 \right] + 2 \lambda \sigma_0^2 \mathrm{E}\left[\left\|(T + \lambda)^{-2}k_{X_i}\right\|_{\mathcal{H}}^2 \right]\\
		& \leq 2\sigma_0^2\mathrm{tr}\left((T+ \lambda)^{-2}T\right) + 2 \sigma_0^2 \mathrm{E}\left[\left\|(T + \lambda)^{-1}k_{X_i}\right\|_{\mathcal{H}}^2 \right]\\
		& \leq 4 \sigma_0^2\mathrm{tr}\left((T+ \lambda)^{-2}T\right).
	\end{align*}
	Therefore,
	\begin{align}\label{eq:theorem:Uniform-CI-Bands-RKHS-2}
		\mathrm{E}\|V\|_{\mathcal{H}}^3 & \leq  \lambda^{-1} 8 \sigma_0^2\mathrm{tr}\left((T+ \lambda)^{-2}T\right) \kappa (B + \kappa\|f_0\|_{\mathcal{H}}) \lesssim \lambda^{-1} \bar{\sigma}^3 \mathfrak{n}_2^2(\lambda),
	\end{align}
	where $\bar{\sigma}^2 \geq \sigma_0^2 \vee \kappa^2(B + \kappa\|f_0\|_{\mathcal{H}})^2 \vee 1$ and $\mathfrak{n}_\alpha^2(\lambda) = \mathrm{tr}\left((T+ \lambda)^{-2\alpha}T\right)$, $\alpha \in \mathbb{N}_0$.
	
	Second, by Cauchy-Schwarz,
	\begin{align}\label{eq:theorem:Uniform-CI-Bands-RKHS-1}
		\mathrm{E}\|Z\|_{\mathcal{F}_n} = \mathrm{E}\left[ \sup_{u \in S^*} \big|\langle Z, u \rangle_{\mathcal{H}} \big| \right] \leq \mathrm{E}\|Z\|_{\mathcal{H}} \sup_{u \in S^*} \|u\|_{\mathcal{H}} \leq \sqrt{\mathrm{tr}(\Omega)} \kappa < \infty.
	\end{align}
	
	Third, for $r > 0$,
	\begin{align}\label{eq:theorem:Uniform-CI-Bands-RKHS-3}
		\frac{\mathrm{E}[\|V\|_{\mathcal{H}}^3\mathbf{1}\{\|V\|_{\mathcal{H}}> n^{1/3}\mathrm{E}[\|V\|_{\mathcal{H}}^3]^{1/3}\}]}{ \mathrm{E}[\|V\|_{\mathcal{H}}^3]} &\leq \frac{\mathrm{E}[\|V\|_{\mathcal{H}}^{3+r}]}{n^{r/3} \mathrm{E}[\|V\|_{\mathcal{H}}^3]^{r/3+1}}  \nonumber\\
		&=  n^{-r/3}\left( \frac{\mathrm{E}[\|V\|_{\mathcal{H}}^{3+r}]^{1/(3+r)}}{ \mathrm{E}[\|V\|_{\mathcal{H}}^3]^{1/3}} \right)^{3+ r} \nonumber\\
		&\lesssim n^{-r/3} \left(\lambda^{-1} 2\kappa (B + \kappa\|f_0\|_{\mathcal{H}}) \right)^r.
	\end{align}
	
	Fourth, for $u \in S^*$ arbitrary,
	\begin{align}\label{eq:theorem:Uniform-CI-Bands-RKHS-4}
		\mathrm{Var}\big(Z(u)\big) \geq \inf_{u \in S^*} \mathrm{Var}\big(\langle Z, u \rangle_{\mathcal{H}}\big)  =  \inf_{u \in S^*} \langle \Omega u, u \rangle_{\mathcal{H}} \geq \omega_S > 0. 
	\end{align}
	Combined with~\eqref{eq:theorem:Uniform-CI-Bands-RKHS-1} and Lemma~\ref{lemma:BoundsVariance-SeparableProcess} this lower bound yields 
	\begin{align}\label{eq:theorem:Uniform-CI-Bands-RKHS-5}
		\mathrm{Var}(\|Z\|_{\mathcal{F}_n}) \gtrsim \left(\frac{ \omega_S }{1 + \mathrm{E}\|Z\|_{\mathcal{F}_n}/\omega_S}\right)^2 \gtrsim \left(\frac{ \omega_S^2}{\omega_S  + \sqrt{\mathrm{tr}(\Omega)} \kappa }\right)^2 \gtrsim \frac{ \omega_S^4}{\mathrm{tr}(\Omega) \kappa^2 }.
	\end{align}
	
	Fifth, by Lemma~\ref{lemma:BC-KRR-Consistency-Operator}, with probability at least $1-\delta$,
	\begin{align}\label{eq:theorem:Uniform-CI-Bands-RKHS-6}
		\big\| \Omega - \widehat{\Omega}_n^m\big\|_{op} &\leq \big\| \Omega - \widehat{\Omega}_n\big\|_{op} + \big\| \widehat{\Omega}_n - \widehat{\Omega}_n^m\big\|_{op} \nonumber\\
		& \lesssim \big\| \widehat{\Omega}_n - \widehat{\Omega}_n^m\big\|_{op} + \|T^3(T + \lambda)^{-4}\|_{op} \sqrt{\frac{ (\kappa^4 +\kappa^2  \mathfrak{n}_1^2(\lambda) + \bar{\sigma}^4)\log(2/\delta)}{n\lambda^2}}.
	\end{align}
	
	Finally, Lemma~\ref{lemma:BC-KRR-Remainder} combined with~\eqref{eq:theorem:Uniform-CI-Bands-RKHS-11}--\eqref{eq:theorem:Uniform-CI-Bands-RKHS-6}, $r =1$ and $\delta \in (0,1/n)$ implies that, with probability at least $1 - \delta$,
	\begin{align}
		\begin{split}\label{eq:theorem:Uniform-CI-Bands-RKHS-8}
		&\sup_{\alpha \in (0,1)} \Big|\mathbb{P}\left\{ \sqrt{n}\|\widehat{f}^{\mathrm{bc}}_n - f_0\|_\infty \leq 	c_n^m(\alpha) \right\}  - \alpha \Big| \\
		&\quad{}\quad{}\lesssim \kappa \left(\frac{\bar{\sigma}^3 \mathfrak{n}_2^2(\lambda)}{ \sqrt{n}\lambda}  \right)^{1/3} \frac{\sqrt{ \mathrm{tr}(\Omega)} \kappa }{\omega_S^2} + \delta +  \frac{\bar{\sigma}}{n^{1/3}\lambda}  + \frac{\mathrm{tr}(\Omega) \kappa^2 }{\sqrt{n}  \omega_S^2} \\
		& \quad{}\quad{}\quad{}  + \mathbb{P}\left\{ \big\| \widehat{\Omega}_n - \widehat{\Omega}_n^m\big\|_{op} > \|T^3(T + \lambda)^{-4}\|_{op} \sqrt{\frac{ (\kappa^4 +\kappa^2  \mathfrak{n}_1^2(\lambda) + \bar{\sigma}^4)\log(2/\delta)}{n\lambda^2}}\right\}\\
		&\quad \quad \quad + \left(\|T^3(T + \lambda)^{-4}\|_{op} \sqrt{\frac{ (\kappa^4 +\kappa^2  \mathfrak{n}_1^2(\lambda) + \bar{\sigma}^4)\log(2/\delta)}{n\lambda^2}} \frac{\mathrm{tr}(\Omega) \kappa^2}{\omega_S^4} \right)^{1/3}  + \delta + \\
		&\quad{}\quad{}\quad{} + \left(\sqrt{\frac{ \bar{\sigma}^2  \mathfrak{n}_1^2(\lambda)  }{n\lambda^2 } } \vee \frac{ \bar{\sigma}^2 }{n \lambda^2}\right)\left(\frac{\kappa^4\log^3(2/\delta)}{\sqrt{n} \lambda} \vee  \kappa^2 \log^2(2/\delta) \right)\frac{\sqrt{\mathrm{tr}(\Omega)} \kappa}{\omega_S^2} \\
		& \quad\quad\quad +  \left( \frac{\kappa^4\log^2(2/\delta)}{\sqrt{n}} \vee \sqrt{n} \lambda^2 \right) \kappa \| (T + \lambda)^{-2}f_0\|_{\mathcal{H}} \frac{\sqrt{\mathrm{tr}(\Omega)} \kappa}{\omega_S^2}
		\end{split}\\
		&\overset{(a)}{=} o(1) + \mathbb{P}\left\{ \big\| \widehat{\Omega}_n - \widehat{\Omega}_n^m\big\|_{op} > \|T^3(T + \lambda)^{-4}\|_{op} \sqrt{\frac{ (\kappa^4 +\kappa^2  \mathfrak{n}_1^2(\lambda) + \bar{\sigma}^4)\log n}{n\lambda^2}}\right\} \nonumber,
	\end{align}
	where (a) holds provided that sample size $n$, regularization parameter $\lambda$, kernel $k$, and operators $\Omega$ and $T$ are such that 
	\begin{align}\label{eq:theorem:Uniform-CI-Bands-RKHS-7}
		\begin{split}
			(i)&\quad \quad  \bar{\sigma}^2 \kappa^4 (\log n)^3 \mathfrak{n}_1(\lambda) \sqrt{\mathrm{tr}(\Omega)} \|T^3(T + \lambda)^{-4}\|_{op}= o( \sqrt{n} \lambda \omega_S^2)\\
			(ii)&\quad \quad \kappa^6 (\log n)^2 \sqrt{\mathrm{tr}(\Omega)} \|(T + \lambda)^{-2}f_0\|_{\mathcal{H}} = o(\sqrt{n}\omega_S^2)\\
			(iii)&\quad \quad  \lambda^2 \kappa^2\sqrt{\mathrm{tr}(\Omega)}\|(T + \lambda)^{-2}f_0\|_{\mathcal{H}} = o(\omega_S^2)\\
			(iv) &\quad \quad \bar{\sigma}^2 \kappa^2\sqrt{\log n} \mathrm{tr}(\Omega)\mathfrak{n}_1(\lambda)  \|T^3(T + \lambda)^{-4}\|_{op}  = o(\sqrt{n} \lambda \omega_S^4)\\
			(v) & \quad \quad \bar{\sigma}^3 \kappa^3 \big(\mathrm{tr}(\Omega) \big)^{3/2}\mathfrak{n}_2^2(\lambda) = o(\sqrt{n} \lambda \omega_S^6)\\
			(vi) & \quad \quad \bar{\sigma} = o(n^{1/3}\lambda),
		\end{split}
	\end{align}
	where $\bar{\sigma}^2 \geq \sigma_0^2 \vee \kappa^2(B + \kappa\|f_0\|_{\mathcal{H}})^2 \vee 1$ and $\mathfrak{n}_\beta^2(\lambda) = \mathrm{tr}\left((T+ \lambda)^{-2\beta}T\right)$, $\beta \in \mathbb{N}_0$.
	
	Simplifying these rates is beyond the scope of this illustrative example. The interested reader may consult Appendix H in~\cite{singh2023kernel} and Sections 6 and 7 in~\cite{lopes2022improved} for potentially useful results. If the Hilbert space $\mathcal{H}$ is finite dimensional and $T$ is invertible, then these rates are satisfied for $\lambda \rightarrow 0$ and $n^{1/3} \lambda \rightarrow \infty$ (the exact rates feature some $(\log n)^c$- factor for some $c \geq 1$).
\end{proof}

\newpage
\section{Proofs of the results in Section~\ref{sec:AuxiliaryResults}}\label{sec:Proofs-AuxiliaryResults}

\subsection{Proofs of Lemmas~\ref{lemma:Smooth-Lipschitz-Approx},~\ref{lemma:ChainRule-SecondOrder}, and~\ref{lemma:DerivativesMaxNorm}}

\begin{proof}[\textbf{Proof of Lemma~\ref{lemma:Smooth-Lipschitz-Approx}}]
	For $s,t \in \mathbb{R}$ and $\lambda > 0$ arbitrary, define
	\begin{align}\label{eq:lemma:Smooth-Lipschitz-Approx-1}
		g_{s, \lambda}(t) := \left( \left(1 + \frac{s-t}{\lambda}\right) \wedge 1\right) \vee 0 \quad{}\quad{}\mathrm{and} \quad{}\quad{} g_{s,0}(t): = \mathbf{1}_{(-\infty, s]}(t).
	\end{align}
	Since  $g_{s, 0} \leq g_{s, \lambda} \leq g_{s + \lambda, 0}$ for all $s \in \mathbb{R}$ and $\lambda > 0$, it follows that (draw a sketch!)
	\begin{align}\label{eq:lemma:Smooth-Lipschitz-Approx-2}
		g_{s, 0} \leq g_{s+\lambda, \lambda} \leq g_{s +\lambda, \lambda} \ast \varrho_\lambda \leq g_{s + 2 \lambda, \lambda} \leq g_{s + 3\lambda, 0},
	\end{align}
	where $\varrho_\lambda(\cdot) := \lambda^{-1} \varrho\left(\:\cdot\: \lambda^{-1}\right)$ and
	\begin{align*}
		\varrho(t) = C_0\exp\left(\frac{1}{t^2 - 1}\right) \mathbf{1}_{[-1,1]}(t),
	\end{align*}
	where $C_0 > 0$ is an absolute constant such that $\int \varrho(t) dt = 1$. Since $\varrho \in C_c^\infty(\mathbb{R})$ with support $[-1,1]$, we easily verify that the map $ t \mapsto (g_{s+ \lambda, \lambda} \ast \varrho_\lambda)(t)$ is continuously differentiable and its $k$th derivative satisfies
	\begin{align}\label{eq:lemma:Smooth-Lipschitz-Approx-3}
		\left| D^k \left(g_{s+ \lambda, \lambda} \ast \varrho_\lambda\right)(t) \right| =  \left|\left(g_{s+ \lambda, \lambda} \ast (D^k  \varrho_\lambda) \right)(t)\right| \leq C_k \lambda^{-k} \mathbf{1}_{[s, s+3\lambda]}(t),
	\end{align}
	where $C_k > 0$ is a constant depending only on $k \in \mathbb{N}_0$. This establishes the first claim of the lemma with $h_{s, \lambda} \equiv g_{s+ \lambda, \lambda}  \ast \varrho_\lambda$. 
	To prove the second claim of the lemma we proceed in two steps: First, we show that
	\begin{align*}
		\sup_{s \in \mathbb{R}} \left|\mathbb{P}\left\{X \leq s \right\} - \mathbb{P}\left\{Z \leq s\right\}\right| &\leq \sup_{s \in \mathbb{R}} \big|\mathrm{E}[(g_{s+ \lambda, \lambda}\ast \varrho_\lambda)(X) - (g_{s+ \lambda, \lambda}\ast \varrho_\lambda)(Z)] \big|  + \zeta_{3\lambda}(X) \wedge \zeta_{3\lambda}(Z).
	\end{align*} 
	By the chain of inequalities in eq.~\eqref{eq:lemma:Smooth-Lipschitz-Approx-2}, for $s \in \mathbb{R}$ arbitrary,
	\begin{align}\label{eq:lemma:Smooth-Lipschitz-Approx-4}
		&\mathbb{P}\left\{X \leq s\right\}-\mathbb{P}\left\{Z \leq s\right\} \nonumber\\
		&\quad{}= \mathrm{E} \left[g_{s, 0}(X) - g_{s, 0}(Z)\right] \nonumber\\
		&\quad{}\leq \left| \mathrm{E}\left[(g_{s+ \lambda, \lambda} \ast \varrho_\lambda)(X) - (g_{s+ \lambda, \lambda} \ast \varrho_\lambda)(Z)\right]\right| + \left| \mathrm{E}\left[(g_{s+ \lambda, \lambda} \ast \varrho_\lambda)(Z) - g_{s, 0}(Z)\right]\right| \nonumber\\
		&\quad{}\leq \left| \mathrm{E}\left[(g_{s+ \lambda, \lambda}\ast \varrho_\lambda)(X) - (g_{s+ \lambda, \lambda} \ast \varrho_\lambda)(Z)\right]\right| + \left| \mathrm{E}\left[g_{s + 3\lambda, 0}(Z) - g_{s, 0}(Z)\right]\right|\nonumber\\
		&\quad{}\leq \left| \mathrm{E}\left[(g_{s+ \lambda, \lambda}\ast \varrho_\lambda)(X) - (g_{s+ \lambda, \lambda} \ast \varrho_\lambda)(Z)\right]\right| + \mathbb{P}\left\{ s \leq Z \leq s + 3 \lambda \right\}.
	\end{align}
	Similarly,
	\begin{align}\label{eq:lemma:Smooth-Lipschitz-Approx-5}
		&\mathbb{P}\left\{Z \leq s + 3 \lambda \right\}-\mathbb{P}\left\{X \leq s  + 3 \lambda \right\}\nonumber\\
		&\quad{}= \mathrm{E} \left[g_{s + 3 \lambda, 0}(Z) - g_{s + 3 \lambda, 0}(X)\right] \nonumber\\
		&\quad{} \leq \left| \mathrm{E}\left[(g_{s+ \lambda, \lambda}  \ast \varrho_\lambda)(Z) - (g_{s+ \lambda, \lambda}  \ast \varrho_\lambda)(X)\right]\right| + \left| \mathrm{E}\left[(g_{s+ \lambda, \lambda}  \ast \varrho_\lambda)(Z) - g_{s+3\lambda, 0}(Z)\right]\right| \nonumber\\
		&\quad{} \leq \left| \mathrm{E}\left[(g_{s+ \lambda, \lambda}  \ast \varrho_\lambda)(Z) - (g_{s+ \lambda, \lambda} \ast \varrho_\lambda)(X)\right]\right| + \left| \mathrm{E}\left[g_{s+3\lambda} (Z) - g_{s, 0}(Z)\right]\right|\nonumber\\
		&\quad{}\leq \left| \mathrm{E}\left[(g_{s+ \lambda, \lambda}  \ast \varrho_\lambda)(Z) - (g_{s+ \lambda, \lambda}  \ast \varrho_\lambda)(X)\right]\right| + \mathbb{P}\left\{s \leq Z \leq s + 3 \lambda \right\}.
	\end{align}
	Now, take the supremum over $s \in \mathbb{R}$ in eq.~\eqref{eq:lemma:Smooth-Lipschitz-Approx-4} and~\eqref{eq:lemma:Smooth-Lipschitz-Approx-5} and switch the roles of $X$ and $Z$. Next, we show
	\begin{align*}
		\sup_{s \in \mathbb{R}} \big|\mathrm{E}[(g_{s+ \lambda, \lambda} \ast \varrho_\lambda)(X) - (g_{s+ \lambda, \lambda} \ast \varrho_\lambda)(Z)] \big| &\leq \sup_{s \in \mathbb{R}} \left|\mathbb{P}\left\{X \leq s \right\} - \mathbb{P}\left\{Z \leq s\right\}\right| + \zeta_{3\lambda}(X) \wedge \zeta_{3\lambda}(Z).
	\end{align*}
	As before, we compute
	\begin{align}\label{eq:lemma:Smooth-Lipschitz-Approx-6}
		&\mathrm{E}\left[(g_{s+ \lambda, \lambda} \ast \varrho_\lambda)(X) - (g_{s+ \lambda, \lambda} \ast \varrho_\lambda)(Z)\right]\nonumber\\
		&\quad{} \leq \mathrm{E}\left[g_{s + 3\lambda, 0}(X) - g_{s + 3 \lambda, 0}(Z)\right] + \mathrm{E}\left[g_{s + 3\lambda, 0}(Z) - g_{s, 0}(Z)\right]\nonumber\\
		&\quad{} \leq \left|\mathbb{P}\left\{X \leq s + 3\lambda \right\} - \mathbb{P}\left\{Z \leq s + 3 \lambda\right\}\right| + \mathbb{P}\left\{ s \leq Z \leq s + 3 \lambda \right\},
	\end{align}
	and 
	\begin{align}\label{eq:lemma:Smooth-Lipschitz-Approx-7}
		&\mathrm{E}\left[(g_{s+ \lambda, \lambda} \ast \varrho_\lambda)(Z) - (g_{s+ \lambda, \lambda} \ast \varrho_\lambda)(X)\right]\nonumber\\
		&\quad{} \leq \mathrm{E}\left[g_{s, 0}(Z) - g_{s, 0}(X)\right] + \mathrm{E}\left[g_{s + 3\lambda, 0}(Z) - g_{s, 0}(Z)\right]\nonumber\\
		&\quad{} \leq \left|\mathbb{P}\left\{X \leq s \right\} - \mathbb{P}\left\{Z \leq s \right\}\right| + \mathbb{P}\left\{ s \leq Z \leq s + 3 \lambda \right\}.
	\end{align}
	To conclude, take the supremum over $s \in \mathbb{R}$ in~\eqref{eq:lemma:Smooth-Lipschitz-Approx-6} and~\eqref{eq:lemma:Smooth-Lipschitz-Approx-7} and switch the roles of $X$ and $Z$.
\end{proof}

\begin{proof}[\textbf{Proof of Lemma~\ref{lemma:ChainRule-SecondOrder}}]
	The set on which $g \circ f$ is not differentiable is contained in the null set $N$ on which $f$ is not differentiable. Thus, $g \circ f$ is $k$-times differentiable almost everywhere. The derivatives now follow from~\cite{hardy2006combinatorics}.
\end{proof}

\begin{proof}[\textbf{Proof of Lemma~\ref{lemma:DerivativesMaxNorm}}]
	Since $f$ is a piecewise linear function, it has partial derivatives of any order at all its continuity points. Since the set of discontinuity points of $f$ forms a null set with respect to the Lebesgue measure on $\mathbb{R}^d$, $f$ differentiable almost everywhere on $\mathbb{R}^d$. (These partial derivatives need not to be continuous!) The expressions of these partial derivatives follow from direct calculation. 
	
	Even more is true: Since $f$ is Lipschitz continuous, Rademacher's theorem implies that $f$ is in fact totally differentiable almost everywhere on $\mathbb{R}^d$. Furthermore, since $f$ is convex, Alexandrov's theorem implies that $f$ satisfies a second order quadratic expansion almost everywhere on $\mathbb{R}^d$, i.e. for all $x, a \in \mathbb{R}^d$,
	\begin{align*}
		\left|f(x + a) - f(x) - D f(x) a - \frac{1}{2}a'Ha  \right| = o(\|a\|_2^2).
	\end{align*}
	According to~\cite{rockafellar1999second} the matrix $H \in \mathbb{R}^{d \times d}$ is symmetric, positive semi-definite, and equal to the Jacobian of $Df(x)$ for all $x \in \mathbb{R}^d$ at which $Df(x)$ exists. Thus, $H$ can be identified with the second derivative $D^2f(x)$ for almost all $x \in \mathbb{R}^d$.
\end{proof}

\subsection{Proofs of Lemmas~\ref{lemma:DifferentiatingUnderIntegral} and~\ref{lemma:L1ConvergencePartialRegularization}}

\begin{proof}[\textbf{Proof of Lemma~\ref{lemma:DifferentiatingUnderIntegral}}]
	Throughout the proof we write $h(x) = (h_{s, \lambda} \circ f) (x)$ with $h_{s, \lambda} \in C^\infty_c(\mathbb{R})$ as defined in Lemma~\ref{lemma:Smooth-Lipschitz-Approx} and $f(x) =\|x\|_\infty $. To simplify the notation, we denote partial derivatives w.r.t $x_j$ by $\partial_j$, e.g. we write $\partial_{i_1} h$ for $\frac{\partial h}{\partial x_{i_1} }$, $\partial_{i_2} f$ for $\frac{\partial f}{\partial x_{i_2} }$, etc. We use $\mathcal{L}^d$ to denote the Lebesgue measure on $\mathbb{R}^d$.
	
	\vspace{10pt}
	\noindent
	\textbf{Proof of part (i).} 
	\vspace{10pt}
	
	\noindent
	\textbf{Special case $k=1$.} 
	\vspace{10pt}
	
	\noindent
	Let $\varepsilon > 0$ and $e_{i_1} \in \mathbb{R}^d$ the $i_1$-th standard unit vector in $\mathbb{R}^d$. Recall that $h = h_{s, \lambda} \circ f$, where $h_{s, \lambda}$ is $\lambda^{-1}$-Lipschitz w.r.t. the norm induced by the absolute value and $f(x) = \|x\|_\infty$ is $1$-Lipschitz w.r.t. the metric induced by the $\ell_\infty$-norm. Thus, for $x_0, z \in \mathbb{R}^d$ and $t, \epsilon > 0$ arbitrary,
	\begin{align*}
		\Delta_{i_1}(x_0, z, t; \epsilon) &:= \epsilon^{-1}\left[ h\left(e^{-t}x_0 + e^{-t}\epsilon e_{i_1} + \sqrt{1 - e^{-2t}} z\right) - h\left(e^{-t}x_0 + \sqrt{1 - e^{-2t}} z\right)\right]\\ 
		&\leq  \epsilon^{-1} \lambda^{-1}\|e^{-t}\epsilon e_{i_1}\|_\infty \\
		&= \lambda^{-1} e^{-t}.
	\end{align*}
	Hence, the difference quotient $\Delta_{i_1}(x_0, z, t; \epsilon)$ is bounded uniformly in $x_0, z \in \mathbb{R}^d$ and $t, \epsilon > 0$. Furthermore, by Lemma~\ref{lemma:DerivativesMaxNorm} $h$ is differentiable $\mathcal{L}^d$-a.e. Therefore, the conditions of Corollary A.5 in~\cite{dudley2014uniform} are met and we can pass the derivative through both integrals (over $t > 0$ and $z \in \mathbb{R}^d$) to obtain
	\begin{align}\label{eq:lemma:DifferentiatingUnderIntegral-3}
		\partial_{i_1} \left( \int_0^\infty P_t h(x) dt\right) \Big|_{x=x_0} = \int_0^\infty e^{-t}P_t \partial_{i_1} h (x_0)dt.
	\end{align}
	
	\vspace{10pt}
	\noindent
	\textbf{The cases $k \geq 2$.} 
	\vspace{10pt}
	
	\noindent
	``Off-the-shelf" differentiating under the integral sign seems to be only possible in the case $k=1$. In all other cases, Corollary A.5 and related results in~\cite{dudley2002real} do not apply, because the higher-order difference quotients of $h$ are not uniformly integrable (considered as a collection of random variables indexed by a null sequence $(\epsilon_n)_{n \geq 1}$). Therefore, to prove the claim for $k \geq 2$ we develop a more specific inductive argument which is tailored explicitly to the map $x \mapsto \|x\|_\infty$ and the fact that we integrate w.r.t. a non-degenerate Gaussian measure.
	
	\vspace{10pt}
	\noindent
	\textbf{Base case $k=2$ and $i_1 \neq i_2$.} 
	\vspace{10pt}
	
	\noindent
	Let $Y \sim N(0, I_d)$, $A=[A_1, \ldots, A_d] = \Sigma^{-1/2}$, and $\varphi$ be the density of the law of $Y$.
	By a change of variable we can re-write eq.~\eqref{eq:lemma:DifferentiatingUnderIntegral-3} as 
	\begin{align}\label{eq:lemma:DifferentiatingUnderIntegral-4}
		\partial_{i_1} \left( \int_0^\infty P_t h(x) dt\right) \Big|_{x=x_0} &= \int_0^\infty e^{-t}P_t \partial_{i_1} h (x_0)dt \nonumber\\
		&= \int_0^\infty\int_{\mathbb{R}^d} e^{-t} \partial_{i_1}  h\left( e^{-t}x_0 + \sqrt{1- e^{-2t}}A^{-1} y \right) \varphi(y) dydt\nonumber\\
		& = \int_0^\infty\int_{\mathbb{R}^d} e^{-t}\partial_{i_1}  h\left(u \right)\psi (x_0, u, t) du dt,
	\end{align}
	where
	\begin{align*}
		\psi(x, u, t) := \left(1 - e^{-2t}\right)^{-d/2} \det(A) \varphi\left( \frac{A( u - e^{-t}x)}{\sqrt{1 - e^{-2t}}}\right).
	\end{align*}
	Since $A$ is full rank, $\psi$ is a smooth function of $\mathcal{L}^d$-a.e. $x \in \mathbb{R}^d$. Hence, the map $x \mapsto e^{-t}\partial_{i_1} h\left(u \right) \partial_{i_2}  \psi(x, u,t)$ exists and is continuous for $\mathcal{L}^d$-a.e. $x \in \mathbb{R}^d$. Furthermore, the map is uniformly bounded and integrable in $u \in \mathbb{R}^d$ and $t \geq 0$ for $\mathcal{L}^d$-a.e. $x \in \mathbb{R}^d$. Hence, by Corollary A.4 in~\cite{dudley2014uniform} we can pass a second partial derivative through the integral in~\eqref{eq:lemma:DifferentiatingUnderIntegral-4} and compute
	\begin{align}\label{eq:lemma:DifferentiatingUnderIntegral-5}
		\partial_{i_2}\partial_{i_1}\left( \int_0^\infty P_t h(x) dt\right) \Big|_{x=x_0}&=  \int_0^\infty\int_{\mathbb{R}^d} e^{-t} \partial_{i_1} h \left(u \right) \partial_{i_2} \psi(x_0, u, t)du dt.
	\end{align}
	In the following we show that one can ``pull back" the partial derivative $\partial_{i_2}$ from $\psi$ onto $\partial_{i_1} h$. Our argument involves a version of Stein's lemma~\citep[i.e.][Theorem 11.1]{chernozhukov2021NearlyOptimalCLT} and a regularization of $\partial_{i_1} h$. Several computations also crucially depend on the fact that $h(x) = h_{s, \lambda}(\|x\|_\infty)$.
	
	Let $\varrho \in C^\infty(\mathbb{R})$ be the standard mollifier
	\begin{align*}
		\varrho(r) = C_0\exp\left(\frac{1}{r^2 - 1}\right) \mathbf{1}\{-1 \leq r \leq 1\},
	\end{align*}
	where $C_0 > 0$ is an absolute constant such that $\int \varrho(r) dr = 1$. For $\eta > 0$, set $\varrho_\eta(\cdot) = \eta^{-1} \varrho\left(\:\cdot\: \eta^{-1}\right)$ and define the ``partial regularization" of a function $g$ on $\mathbb{R}^d$ in its $i$th coordinate by
	\begin{align}\label{eq:lemma:DifferentiatingUnderIntegral-5-2}
		(\varrho_\eta \ast_{(i)} g) := \int_{-1}^1 \varrho(r) \: g(u - r \eta e_i)dr,
	\end{align}
	where $e_i$ is the $i$th standard unit vector in $\mathbb{R}^d$. With this notation, define
	\begin{align*}
		h_{i_1}^{\eta, i_2}(u) &:= 	(D h_{s, \lambda} \circ f)(u) (\varrho_\eta \ast_{(i_2)} \partial_{i_1} f) (u)\\
		&\equiv (D h_{s, \lambda} \circ f)(u) \int_{-1}^1 \varrho(r) \: \partial_{i_1}f(u - r \eta e_{i_2})dr\\
		&\equiv (D h_{s, \lambda} \circ f)(u) \: \mathbf{1}\{|u_{i_1}| \geq |u_j|, \:\forall j\neq i_2\}  \big(\varrho_\eta \ast \mathbf{1}\{|u_{i_1}| \geq |\cdot|\}\big)(u_{i_2}).
	\end{align*}
	Obviously, the regularization $h_{i_1}^{\eta, i_2}(u)$ is just $\partial_{i_1}h$ with the discontinuity at $u_{i_2}$ smoothed out. (We comment on the rationale behind this partial regularization after eq.~\eqref{eq:lemma:DifferentiatingUnderIntegral-11}.) In particular, $u \mapsto \partial_{i_2} h_{i_1}^{\eta, i_2}(u)$ exists for all $u \in \mathbb{R}^d$, and, by Lemmas~\ref{lemma:ChainRule-SecondOrder} and~\ref{lemma:DerivativesMaxNorm}, 
	\begin{align}
		\partial_{i_2} h_{i_1}^{\eta, i_2} (u) &=\big(D^2 h_{s, \lambda} \circ f\big)(u)\: (\varrho_\eta \ast_{(i_2)} \partial_{i_1} f) (u) \: \partial_{i_2} f(u)\label{eq:lemma:DifferentiatingUnderIntegral-6}\\
		&\quad{} + \big(D^2 h_{s, \lambda} \circ f\big)(u)\: \partial_{i_2}(\varrho_\eta \ast_{(i_2)} \partial_{i_1}f) (u).\label{eq:lemma:DifferentiatingUnderIntegral-7}
	\end{align}
	Observe that by Leibniz's integral rule the partial derivative $\partial_{i_2}$ in line~\eqref{eq:lemma:DifferentiatingUnderIntegral-7} takes the form
	\begin{align*}
		\partial_{i_2}(\varrho_\eta \ast_{(i_2)} \partial_{i_1}f) (u) = \mathrm{sign}(u_{i_1}) \mathbf{1}\{|u_{i_1}| \geq |u_j|,\: j \neq i_2 \} \: \partial_{i_2} \big(\varrho_\eta \ast \mathbf{1}\{|u_{i_1}| \geq |\cdot|\}\big)(u_{i_2}),
	\end{align*}
	where
	\begin{align}\label{eq:lemma:DifferentiatingUnderIntegral-8}
		\partial_{i_2} \big(\varrho_\eta \ast \mathbf{1}\{|u_{i_1}| \geq |\cdot|\}\big)(u_{i_2}) 
		&= \partial_{i_2}\int_{ (u_{i_2} - |u_{i_1}|)/\eta}^{(u_{i_2} + |u_{i_1}|)/\eta } \varrho(r) dr \nonumber\\
		&= \varrho_\eta(u_{i_2} + |u_{i_1}|) -  \varrho_\eta(u_{i_2} - |u_{i_1}|)\nonumber\\
		&= \left(\varrho_\eta \ast \mathbf{1}_{-|u_{i_1}|} \right)(u_{i_2}) - \left(\varrho_\eta \ast \mathbf{1}_{|u_{i_1}|} \right)(u_{i_2}),
	\end{align}
	where $x \mapsto \mathbf{1}_a(x) = \mathbf{1}\{a = x\}$, $a, x \in \mathbb{R}$. Thus, by Lemma~\ref{lemma:Smooth-Lipschitz-Approx} (i) we easily verify that $h_{i_1}^{\eta, i_2}$ and $\partial_{i_2} h_{i_1}^{\eta, i_2}$ are bounded and integrable.
	
	We return to eq.~\eqref{eq:lemma:DifferentiatingUnderIntegral-5}. Adding and subtracting $h_{i_1}^{\eta, i_2}$ we expand its right hand side as 
	\begin{align}
		&\int_0^\infty\int_{\mathbb{R}^d} e^{-t} h_{i_1}^{\eta, i_2}(u) \partial_{i_2}\psi (x_0, u, t)du dt    \label{eq:lemma:DifferentiatingUnderIntegral-9}\\
		&\quad{} + \int_0^\infty\int_{\mathbb{R}^d}  e^{-t} \left(\partial_{i_1}h (u)- h_{i_1}^{\eta, i_2}(u) \right) \partial_{i_2}\psi  (x_0, u, t)du dt.\label{eq:lemma:DifferentiatingUnderIntegral-10}
	\end{align}
	Consider the integral in line~\eqref{eq:lemma:DifferentiatingUnderIntegral-9}. A change of variable yields
	\begin{align*}
		&\int_0^\infty\int_{\mathbb{R}^d} e^{-t} h_{i_1}^{\eta, i_2}(u) \partial_{i_2}\psi (x_0, u, t)du dt\\
		&\quad{}=  \int_0^\infty \int_{\mathbb{R}^d}  \frac{e^{-2t}}{\sqrt{1 - e^{-2t}}} h_{i_1}^{\eta, i_2}(u) \left(1 - e^{-2t}\right)^{-d/2} \det(A) \varphi\left( \frac{A( u - e^{-t}x_0)}{\sqrt{1 - e^{-2t}}}\right)\left( \frac{A( u - e^{-t}x_0)}{\sqrt{1 - e^{-2t}}}\right)'  A_{i_2} du dt\\
		&\quad{} = \int_0^\infty\int_{\mathbb{R}^d} \frac{e^{-2t}}{\sqrt{1 - e^{-2t}}}  h_{i_1}^{\eta, i_2} \left( e^{-t}x_0 + \sqrt{1- e^{-2t}}A^{-1}y\right)  \varphi(y)y'(A^{-1})' A' A_{i_2}dy dt\\
		&\quad{}= \int_0^\infty \frac{e^{-2t}}{\sqrt{1 - e^{-2t}}} \mathrm{E}\left[ h_{i_1}^{\eta, i_2} \left( e^{-t}x_0 + \sqrt{1- e^{-2t}}Z\right) Z'\right] dt  \: A' A_{i_2}.
	\end{align*}
	By Stein's lemma for non-differentiable (but bounded and integrable) functions~\citep[i.e.][Theorem 11.1]{chernozhukov2021NearlyOptimalCLT} the last line in above display is equal to
	\begin{align*}
		&\int_0^\infty e^{-2t} D \mathrm{E}\left[ h_{i_1}^{\eta, i_2} \left( x + \sqrt{1- e^{-2t}}Z\right)\right] \Big|_{x=e^{-t}x_0} dt  \:\Sigma A' A_{i_2} \nonumber\\
		&\quad{}= \int_0^\infty  e^{-2t} D \mathrm{E}\left[ h_{i_1}^{\eta, i_2} \left( x + \sqrt{1- e^{-2t}}Z\right)\right] \Big|_{x=e^{-t}x_0} dt  \:A^{-1} A_{i_2}\nonumber\\
		&\quad{}= \int_0^\infty  e^{-2t} D \mathrm{E}\left[ h_{i_1}^{\eta, i_2} \left( x + \sqrt{1- e^{-2t}}Z\right)\right] \Big|_{x=e^{-t}x_0} dt  \:e_{i_2}\nonumber\\
		&\quad{}= \int_0^\infty  e^{-2t} \partial_{i_2} \mathrm{E}\left[ h_{i_1}^{\eta, i_2} \left( x + \sqrt{1- e^{-2t}}Z\right)\right] \Big|_{x=e^{-t}x_0} dt,
	\end{align*}
	where $e_{i_2}$ denotes the $i_2$-th standard unit vector in $\mathbb{R}^d$.
	
	Consider the expression in the previous line. Since $ x \mapsto \partial_{i_2}h_{i_1}^{\eta, i_2}(x + u)$ exists and is uniformly bounded in $x, u \in \mathbb{R}^d$,  we can push the partial derivative $\partial_{i_2}$ through the expectation~\citep[e.g.][Theorem 2.27 (b); the integrability condition is satisified because we integrate w.r.t. the non-degenerate law of $N(0, \Sigma)$]{folland1999real}. Thus, we have shown that the integral in line~\eqref{eq:lemma:DifferentiatingUnderIntegral-9} equals
	\begin{align}\label{eq:lemma:DifferentiatingUnderIntegral-11}
		\int_0^\infty  e^{-2t}  \mathrm{E}\left[\partial_{i_2} h_{i_1}^{\eta, i_2} \left( e^{-t}x_0 + \sqrt{1- e^{-2t}}Z\right)\right] dt.
	\end{align}
	
	
	Obviously, we could have derived eq.~\eqref{eq:lemma:DifferentiatingUnderIntegral-11} under any kind of smoothing. However, under the ``partial regularization'' of $h_{i_1}$ in its $i_2$th coordinate the partial derivative $\partial_{i_2} h_{i_1}^{\eta, i_2}$ takes on a simple closed-form expression (eq.~\eqref{eq:lemma:DifferentiatingUnderIntegral-6}--\eqref{eq:lemma:DifferentiatingUnderIntegral-8}). The simple form of the expression in eq.~\eqref{eq:lemma:DifferentiatingUnderIntegral-8} is particularly important as it strongly suggests that $\partial_{i_2} h_{i_1}^{\eta, i_2}$ does not ``blow up'' as $\eta \downarrow 0$. Indeed, if we had ``fully'' regularized $h_{i_1}$ in all its coordinates, the expression in eq.~\eqref{eq:lemma:DifferentiatingUnderIntegral-8} would not be as simple and its asymptotic behavior (as $\eta \downarrow 0$) would be less clear.
	
	We now show that the integral in eq.~\eqref{eq:lemma:DifferentiatingUnderIntegral-11} converges to
	\begin{align}\label{eq:lemma:DifferentiatingUnderIntegral-12}
		\int_0^\infty  e^{-2t}  \mathrm{E}\left[\big(D^2 h_{s, \lambda} \circ f\big)(V_0^t)\: \partial_{i_1} f(V_0^t) \: \partial_{i_2} f(V_0^t)\right] dt,
	\end{align}
	as $\eta \downarrow 0$, where 
	\begin{align*}
		V_0^t := e^{-t}x_0 + \sqrt{1 - e^{-2t}}Z.
	\end{align*}
	We record the following useful bounds: First, by Lemma~\ref{lemma:Smooth-Lipschitz-Approx} (i), for all $u \in \mathbb{R}^d$,
	\begin{align}
		\Big| (D h_{s, \lambda} \circ f)(u)\Big| &\leq \lambda^{-1} \mathbf{1} \{s \leq \|u\|_\infty \leq s + 3\lambda\}, \label{eq:lemma:DifferentiatingUnderIntegral-13}\\
		\Big| (D^2 h_{s, \lambda} \circ f)(u) \Big| &\leq \lambda^{-2} \mathbf{1} \{s \leq \|u\|_\infty \leq s + 3\lambda\}.\label{eq:lemma:DifferentiatingUnderIntegral-14}
	\end{align}
	Second, there exists $C_{A,d}> 0$ (depending on $A$ and $d$) such that for all $u \in \mathbb{R}^d$ and $t \geq 0$,
	\begin{align}\label{eq:lemma:DifferentiatingUnderIntegral-15}
		\left|\psi \right| (x_0, u, t) \vee	\left|\partial_{i_2}\psi \right| (x_0, u, t) \leq C_{A,d} \quad{} \text{for\:\:\:}\mathcal{L}^d\text{-a.e.\:\:\:} x_0 \in \mathbb{R}^d.
	\end{align}
	By eq.~\eqref{eq:lemma:DifferentiatingUnderIntegral-14} and~\eqref{eq:lemma:DifferentiatingUnderIntegral-15},
	\begin{align}\label{eq:lemma:DifferentiatingUnderIntegral-16}
		&\left| \int_0^\infty  e^{-2t}  \mathrm{E}\left[\big(D^2 h_{s, \lambda} \circ f\big)(V_0^t)\: (\varrho_\eta \ast_{(i_2)} \partial_{i_1}f) (V_0^t) \: \partial_{i_2} f(V_0^t) \right] dt \right. \nonumber\\
		&\quad{}\quad{}\quad{}\quad{} \left. -  \int_0^\infty  e^{-2t}  \mathrm{E}\left[\big(D^2 h_{s, \lambda} \circ f\big)(V_0^t)\: \partial_{i_1} f(V_0^t) \: \partial_{i_2} f(V_0^t)\right] dt\right| \nonumber\\
		&\quad{} =\left| \int_0^\infty \int_{\mathbb{R}^d} e^{-2t}  \big(D^2 h_{s, \lambda} \circ f\big)(u) \: \Big( (\varrho_\eta \ast_{(i_2)} \partial_{i_1}f) (u) -\partial_{i_1} f(u) \Big)   \: \partial_{i_2} f(u) \psi(x_0, u, t)  du dt \right| \nonumber\\
		&\quad{}\leq  \frac{C_{A,d}}{2\lambda^2} \int_{\mathbb{R}^d} \Big| (\varrho_\eta \ast_{(i_2)} \partial_{i_1}f) (u) -\partial_{i_1} f(u) \Big| \mathbf{1}\{s \leq \|u\|_\infty \leq s + 3 \lambda\} du \nonumber\\
		&\quad{}\rightarrow 0 \quad{} \mathrm{as} \quad{} \eta \rightarrow 0,
	\end{align}
	where the limit in the last line follows from Lemma~\ref{lemma:L1ConvergencePartialRegularization}.
	
	Since the law of $N(0, \Sigma)$ is non-degenerate,
	\begin{align*}
		\big|Corr(V_{0,i_1}^t, V_{0,i_2}^t)\big| = \big|Corr(Z_{i_1}, Z_{i_2})\big| < 1,
	\end{align*}
	and, hence, the event $\{|V_{0,i_1}^t| =|V_{0,i_2}^t|\}$ is a $N(0, \Sigma)$-null set for all $t \geq 0$ and $\mathcal{L}^d$-a.e. $x_0 \in \mathbb{R}^d$ (for a more detailed argument, see also below proof of part (ii) of this lemma). Thus,
	\begin{align*}
		&\int_0^\infty  e^{-2t}  \mathrm{E}\Big[\big(D h_{s, \lambda} \circ f\big)(V_0^t)\: \mathrm{sign}(V_{0,i_1}^t) \mathbf{1}\{|V_{0,i_1}^t| \geq |V_{0,j}^t|,\: j \neq i_2 \} \left( - \mathbf{1}_{-|V_{0, i_1}^t|}(V_{0,i_2}^t) + \mathbf{1}_{|V_{0, i_1}^t|}(V_{0,i_2}^t)\right)\Big] dt\\
		&\quad{} = 0.
	\end{align*} 
	Therefore, by eq.~\eqref{eq:lemma:DifferentiatingUnderIntegral-8},~\eqref{eq:lemma:DifferentiatingUnderIntegral-13}, and~\eqref{eq:lemma:DifferentiatingUnderIntegral-15},
	\begin{align}\label{eq:lemma:DifferentiatingUnderIntegral-17}
		& \left|\int_0^\infty  e^{-2t}  \mathrm{E}\Big[\big(D h_{s, \lambda} \circ f\big)(V_0^t)\: \partial_{i_2}(\varrho_\eta \ast_{(i_2)} \partial_{i_1}f) (V_0^t) \Big] dt \right| \nonumber\\	
		&\quad{}=\left|\int_0^\infty \int_{\mathbb{R}^d} e^{-2t} \big(D h_{s, \lambda} \circ f\big)(u)\:  \mathrm{sign}(u_{i_1}) \mathbf{1}\{|u_{i_1}| \geq |u_j|,\: j \neq i_2 \} \right. \nonumber\\
		&\quad{}\quad{} \quad{}\quad{} \left. \phantom{\int_0^\infty} \times \left[\left(\varrho_\eta \ast \left[\mathbf{1}_{-|u_{i_1}|} - \mathbf{1}_{|u_{i_1}|} \right] \right)(u_{i_2}) - \left(\mathbf{1}_{-|u_{i_1}|}(u_{i_2}) - \mathbf{1}_{|u_{i_1}|}(u_{i_2}) \right) \right] \psi(x_0, u, t) dt \right| \nonumber\\
		&\quad{}\leq \frac{C_{A,d}}{2\lambda}  \int_{\mathbb{R}^d} \left|\left(\varrho_\eta \ast \mathbf{1}_{-|u_{i_1}|}\right)(u_{i_2}) - \mathbf{1}_{-|u_{i_1}|} (u_{i_2}) \right| \mathbf{1}\{s \leq \|u\|_\infty \leq s + 3 \lambda\} du  \nonumber\\
		&\quad{}\quad{} +\frac{C_{A,d}}{2\lambda}  \int_{\mathbb{R}^d} \left|\left(\varrho_\eta \ast \mathbf{1}_{|u_{i_1}|}\right)(u_{i_2}) - \mathbf{1}_{|u_{i_1}|} (u_{i_2}) \right|  \mathbf{1}\{s \leq \|u\|_\infty \leq s + 3 \lambda\}du  \nonumber\\
		&\quad{}\rightarrow 0 \quad{} \mathrm{as} \quad{} \eta \rightarrow 0,
	\end{align}
	where the limit follows from Lemma~\ref{lemma:L1ConvergencePartialRegularization}. The limit in line~\eqref{eq:lemma:DifferentiatingUnderIntegral-12} now follows by combining eq.~\eqref{eq:lemma:DifferentiatingUnderIntegral-16} and~\eqref{eq:lemma:DifferentiatingUnderIntegral-17} with eq.~\eqref{eq:lemma:DifferentiatingUnderIntegral-6} and~\eqref{eq:lemma:DifferentiatingUnderIntegral-7}.
	
	Lastly, we return to the integral in line~\eqref{eq:lemma:DifferentiatingUnderIntegral-10}. By eq.~\eqref{eq:lemma:DifferentiatingUnderIntegral-14} and~\eqref{eq:lemma:DifferentiatingUnderIntegral-15},
	\begin{align}\label{eq:lemma:DifferentiatingUnderIntegral-18}
		&\left|\int_0^\infty\int_{\mathbb{R}^d}  e^{-t} \left(\partial_{i_1}h (u)- h_{i_1}^{\eta, i_2}(u) \right) \partial_{i_2}\psi  (x_0, u, t)du dt\right|\nonumber\\
		&\quad{}= \left| \int_0^\infty\int_{\mathbb{R}^d} e^{-t} (D h_{s, \lambda} \circ f)(u) \Big(  \partial_{i_1} f(u) - (\varrho_\eta \ast_{(i_2)} \partial_{i_1} f) (u) \Big)\partial_{i_2} \psi (x_0, u, t)du dt \right|\nonumber\\
		&\quad{}\leq  \frac{C_{A,d}}{\lambda} \int_{\mathbb{R}^d} \Big|  \partial_{i_1} f(u) - (\varrho_\eta \ast_{(i_2)} \partial_{i_1} f) (u) \Big|  \mathbf{1}\{s \leq \|u\|_\infty \leq s + 3 \lambda\}du \nonumber\\
		&\quad{}\rightarrow 0 \quad{} \mathrm{as} \quad{} \eta \rightarrow 0,
	\end{align}
	where the limit follows from Lemma~\ref{lemma:L1ConvergencePartialRegularization}.
	
	Combine eq.~\eqref{eq:lemma:DifferentiatingUnderIntegral-12} and~\eqref{eq:lemma:DifferentiatingUnderIntegral-18} with eq.~\eqref{eq:lemma:DifferentiatingUnderIntegral-5},~\eqref{eq:lemma:DifferentiatingUnderIntegral-9}, and~\eqref{eq:lemma:DifferentiatingUnderIntegral-10} to conclude via Lemma~\eqref{lemma:DerivativesMaxNorm} that for $\mathcal{L}^d$-a.e. $x_0 \in \mathbb{R}^d$,
	\begin{align}\label{eq:lemma:DifferentiatingUnderIntegral-19}
		\begin{split}
			\partial_{i_2}\partial_{i_1}\left( \int_0^\infty P_t h(x) dt\right) \Big|_{x=x_0} &= \int_0^\infty  e^{-2t}  \mathrm{E}\left[\big(D^2 h_{s, \lambda} \circ f\big)(V_0^t)\: \partial_{i_1} f(V_0^t) \: \partial_{i_2} f(V_0^t)\right] dt\\
			&\equiv \int_0^\infty\int_{\mathbb{R}^d} e^{-2t} \partial_{i_2} \partial_{i_1} h(u) \psi(x_0, u, t)du dt\\	
			&\equiv \int_0^\infty e^{-2t}P_t \left(\partial_{i_2}\partial_{i_1} h\right) (x_0)dt.
		\end{split}
	\end{align}
	Notice that the order in which we take the partial derivatives $\partial_{i_1}$ and $\partial_{i_2}$ does not matter.
	
	\vspace{10pt}
	\noindent
	\textbf{Base case $k=2$ and $i_1 = i_2$.} 
	\vspace{10pt}
	
	\noindent
	The strategy is identical to the one of the preceding case $i_1 \neq i_2$. The only difference is the regularization of $\partial_{i_1} h$. Recall the notion of ``partial regularization'' from eq.~\eqref{eq:lemma:DifferentiatingUnderIntegral-5-2} and define, for $\eta > 0$,
	\begin{align*}
		h_{i_1}^{\eta, i_1}(u) &:= 	(D h_{s, \lambda} \circ f)(u) (\varrho_\eta \ast_{(i_1)} \partial_{i_1} f) (u)\\
		&\equiv (D h_{s, \lambda} \circ f)(u) \int_{-1}^1 \varrho(r) \: \partial_{i_1}f(u - r \eta e_{i_1})dr\\
		&\equiv (D h_{s, \lambda} \circ f)(u) \big(\varrho_\eta \ast \mathrm{sign}(\cdot) \mathbf{1}\{|\cdot| \geq |u_j|, \:\forall j\}\big)(u_{i_1}).
	\end{align*}
	The map $u \mapsto \partial_{i_1} h_{i_1}^{\eta, i_1}(u)$ exists for all $u \in \mathbb{R}^d$, and, by Lemmas~\ref{lemma:ChainRule-SecondOrder} and~\ref{lemma:DerivativesMaxNorm}, 
	\begin{align}
		\partial_{i_1} h_{i_1}^{\eta, i_2} (u) &=\big(D^2 h_{s, \lambda} \circ f\big)(u)\: (\varrho_\eta \ast_{(i_1)} \partial_{i_1}f) (u) \: \partial_{i_1} f(u)\label{eq:lemma:DifferentiatingUnderIntegral-20}\\
		&\quad{} + \big(D^2 h_{s, \lambda} \circ f\big)(u)\: \partial_{i_1}(\varrho_\eta \ast_{(i_1)} \partial_{i_1}f) (u).\label{eq:lemma:DifferentiatingUnderIntegral-21}
	\end{align}
	By Leibniz's integral rule we find that the partial derivative $\partial_{i_1}$ in line~\eqref{eq:lemma:DifferentiatingUnderIntegral-21} equals
	\begin{align}\label{eq:lemma:DifferentiatingUnderIntegral-22}
		&\partial_{i_1}(\varrho_\eta \ast_{(i_1)} \partial_{i_1}f) (u) \nonumber\\
		&\quad{}=\partial_{i_1}  \big(\varrho_\eta \ast \mathrm{sign}(\cdot) \mathbf{1}\{|\cdot| \geq |u_j|, \:\forall j\}\big)(u_{i_1})\nonumber\\
		&\quad{}= \partial_{i_1} \left( \int_{(u_{i_1} + \max_{j \neq i_1} |u_j| )/\eta}^1 \varrho(r) \mathrm{sign}(u_{i_1}-r\eta)  dr + \int_{-1}^{(u_{i_1} - \max_{j \neq i_1} |u_j| )/\eta} \varrho(r) \mathrm{sign}(u_{i_1}-r\eta)  dr \right)  \nonumber\\
		&\quad{}= \partial_{i_1} \left( -\int_{(u_{i_1} + \max_{j \neq i_1} |u_j| )/\eta}^1 \varrho(r) dr + \int_{-1}^{(u_{i_1} - \max_{j \neq i_1} |u_j| )/\eta} \varrho(r) dr \right)  \nonumber\\
		&\quad{}= \varrho_\eta(u_{i_1} + \max_{j \neq i_1} |u_j|) -  \varrho_\eta(u_{i_1} - \max_{j \neq i_1} |u_j|)\nonumber\\
		&\quad{}= \left(\varrho_\eta \ast \mathbf{1}_{-\max_{j \neq i_1} |u_j|} \right)(u_{i_1}) - \left(\varrho_\eta \ast \mathbf{1}_{\max_{j \neq i_1} |u_j|} \right)(u_{i_1}),
	\end{align}
	where $x \mapsto \mathbf{1}_a(x) = \mathbf{1}\{a = x\}$, $a, x \in \mathbb{R}$. Thus, from Lemma~\ref{lemma:Smooth-Lipschitz-Approx} (i) we infer that $h_{i_1}^{\eta, i_1}$ and $\partial_{i_1} h_{i_1}^{\eta, i_1}$ are both bounded and integrable.
	
	The same arguments that led to eq.~\eqref{eq:lemma:DifferentiatingUnderIntegral-5},~\eqref{eq:lemma:DifferentiatingUnderIntegral-9}, and~\eqref{eq:lemma:DifferentiatingUnderIntegral-10} also yield
	\begin{align}
		\frac{\partial^2}{\partial x_{i_1}^2}\left( \int_0^\infty P_t h(x) dt\right) \Big|_{x=x_0}
		&=\int_0^\infty\int_{\mathbb{R}^d} e^{-t} h_{i_1}^{\eta, i_1}(u) \partial_{i_1}\psi (x_0, u, t)du dt    \label{eq:lemma:DifferentiatingUnderIntegral-23}\\
		&\quad{} + \int_0^\infty\int_{\mathbb{R}^d}  e^{-t} \left(\partial_{i_1}h (u)- h_{i_1}^{\eta, i_1}(u) \right) \partial_{i_1}\psi  (x_0, u, t)du dt.\label{eq:lemma:DifferentiatingUnderIntegral-24}
	\end{align}
	Repeating the arguments that gave eq.~\eqref{eq:lemma:DifferentiatingUnderIntegral-11} we find that the integral in line~\eqref{eq:lemma:DifferentiatingUnderIntegral-23} equals
	\begin{align}\label{eq:lemma:DifferentiatingUnderIntegral-25}
		\int_0^\infty  e^{-2t}  \mathrm{E}\left[\partial_{i_1} h_{i_1}^{\eta, i_1} \left( e^{-t}x_0 + \sqrt{1- e^{-2t}}Z\right)\right] dt.
	\end{align}
	We now study the behavior of the integrals in lines~\eqref{eq:lemma:DifferentiatingUnderIntegral-24} and \eqref{eq:lemma:DifferentiatingUnderIntegral-25} as $\eta \downarrow 0$. The arguments are similar to those used in the case $i_1 \neq i_2$; we provide them for completeness only. As before, let $V_0^t = e^{-t}x_0 + \sqrt{1 - e^{-2t}}Z$. By eq.~\eqref{eq:lemma:DifferentiatingUnderIntegral-14} and~\eqref{eq:lemma:DifferentiatingUnderIntegral-15} and Lemma~\ref{lemma:L1ConvergencePartialRegularization},
	\begin{align}\label{eq:lemma:DifferentiatingUnderIntegral-27}
		&\left| \int_0^\infty  e^{-2t}  \mathrm{E}\left[\big(D^2 h_{s, \lambda} \circ f\big)(V_0^t)\: (\varrho_\eta \ast_{(i_1)} \partial_{i_1}f) (V_0^t) \: \partial_{i_1} f(V_0^t) \right] dt \right. \nonumber\\
		&\quad{}\quad{}\quad{}\quad{} \left. -  \int_0^\infty  e^{-2t}  \mathrm{E}\left[\big(D^2 h_{s, \lambda} \circ f\big)(V_0^t)\: \partial_{i_1} f(V_0^t) \: \partial_{i_1} f(V_0^t)\right] dt\right| \nonumber\\
		&\quad{} =\left| \int_0^\infty \int_{\mathbb{R}^d} e^{-2t}  \big(D^2 h_{s, \lambda} \circ f\big)(u) \: \Big( (\varrho_\eta \ast_{(i_1)} \partial_{i_1}f) (u) -\partial_{i_1} f(u) \Big)   \: \partial_{i_1} f(u) \psi(x_0, u, t)  du dt \right| \nonumber\\
		&\quad{}\leq  \frac{C_{A,d}}{2\lambda^2} \int_{\mathbb{R}^d} \Big| (\varrho_\eta \ast_{(i_1)} \partial_{i_1}f) (u) -\partial_{i_1} f(u) \Big| \mathbf{1}\{s \leq \|u\|_\infty \leq s + 3 \lambda\} du \nonumber\\
		&\quad{}\rightarrow 0 \quad{} \mathrm{as} \quad{} \eta \rightarrow 0,
	\end{align}
	Since the law of $N(0, \Sigma)$ is non-degenerate, by eq.~\eqref{eq:lemma:DifferentiatingUnderIntegral-22},~\eqref{eq:lemma:DifferentiatingUnderIntegral-13}, and~\eqref{eq:lemma:DifferentiatingUnderIntegral-15} and Lemma~\ref{lemma:L1ConvergencePartialRegularization},
	\begin{align}\label{eq:lemma:DifferentiatingUnderIntegral-28}
		&\left|\int_0^\infty  e^{-2t}  \mathrm{E}\Big[\big(D h_{s, \lambda} \circ f\big)(V_0^t)\: \partial_{i_1}(\varrho_\eta \ast_{(i_1)} \partial_{i_1}) (V_0^t) \Big] dt \right| \nonumber\\	
		&\quad{}=\left|\int_0^\infty \int_{\mathbb{R}^d} e^{-2t} \big(D h_{s, \lambda} \circ f\big)(u) \right. \nonumber\\
		&\left. \phantom{\int_0^\infty} \times \left[\left(\varrho_\eta \ast \left[\mathbf{1}_{-\max_{j \neq i_1}|u_j|}-\mathbf{1}_{\max_{j \neq i_1}|u_j|}\right] \right)(u_{i_2}) - \left(\mathbf{1}_{-\max_{j \neq i_1}|u_j|}- \mathbf{1}_{\max_{j \neq i_1}|u_j|} \right)(u_{i_1}) \right] \psi(x_0, u, t) dt \right| \nonumber\\
		&\quad{}\leq \frac{C_{A,d}}{2\lambda}  \int_{\mathbb{R}^d} \mathbf{1}\{s \leq \|u\|_\infty \leq s + 3 \lambda\} \left|\left(\varrho_\eta \ast \mathbf{1}_{-\max_{j \neq i_1}|u_j|}\right)(u_{i_1}) - \mathbf{1}_{-\max_{j \neq i_1}|u_j|} (u_{i_1}) \right|  du  \nonumber\\
		&\quad{}\quad{} +\frac{C_{A,d}}{2\lambda}  \int_{\mathbb{R}^d} \mathbf{1}\{s \leq \|u\|_\infty \leq s + 3 \lambda\} \left|\left(\varrho_\eta \ast \mathbf{1}_{\max_{j \neq i_1}|u_j|}\right)(u_{i_1}) - \mathbf{1}_{\max_{j \neq i_1}|u_j|} (u_{i_1}) \right| du  \nonumber\\
		&\quad{}\rightarrow 0 \quad{} \mathrm{as} \quad{} \eta \rightarrow 0.
	\end{align}
	Combine eq.~\eqref{eq:lemma:DifferentiatingUnderIntegral-27} and~\eqref{eq:lemma:DifferentiatingUnderIntegral-28} to conclude that, as $\eta \downarrow 0$, the integral in~\eqref{eq:lemma:DifferentiatingUnderIntegral-25} converges to
	\begin{align}\label{eq:lemma:DifferentiatingUnderIntegral-26}
		\int_0^\infty  e^{-2t}  \mathrm{E}\left[\big(D^2 h_{s, \lambda} \circ f\big)(V_0^t)\: \partial_{i_1} f(V_0^t) \: \partial_{i_1} f(V_0^t)\right] dt.
	\end{align}
	Lastly, we turn to the integral in line~\eqref{eq:lemma:DifferentiatingUnderIntegral-24}. By eq.~\eqref{eq:lemma:DifferentiatingUnderIntegral-14} and~\eqref{eq:lemma:DifferentiatingUnderIntegral-15} and Lemma~\ref{lemma:L1ConvergencePartialRegularization},
	\begin{align}\label{eq:lemma:DifferentiatingUnderIntegral-29}
		&\left|\int_0^\infty\int_{\mathbb{R}^d}  e^{-t} \left(\partial_{i_1}h (u)- h_{i_1}^{\eta, i_1}(u) \right) \partial_{i_1}\psi  (x_0, u, t)du dt\right|\nonumber\\
		&\quad{}= \left| \int_0^\infty\int_{\mathbb{R}^d} e^{-t} (D h_{s, \lambda} \circ f)(u) \Big(  \partial_{i_1} f(u) - (\varrho_\eta \ast_{(i_1)} \partial_{i_1} f) (u) \Big)\partial_{i_1} \psi (x_0, u, t)du dt \right|\nonumber\\
		&\quad{}\leq  \frac{C_{A,d}}{\lambda} \int_{\mathbb{R}^d} \Big|  \partial_{i_1} f(u) - (\varrho_\eta \ast_{(i_1)} \partial_{i_1} f) (u) \Big|  \mathbf{1}\{s \leq \|u\|_\infty \leq s + 3 \lambda\}du \nonumber\\
		&\quad{}\rightarrow 0 \quad{} \mathrm{as} \quad{} \eta \rightarrow 0,
	\end{align}
	Thus, combining eq.~\eqref{eq:lemma:DifferentiatingUnderIntegral-26} and~\eqref{eq:lemma:DifferentiatingUnderIntegral-29} with eq.~\eqref{eq:lemma:DifferentiatingUnderIntegral-23} and~\eqref{eq:lemma:DifferentiatingUnderIntegral-24} and invoking Lemma~\eqref{lemma:DerivativesMaxNorm} we conclude for $\mathcal{L}^d$-a.e. $x_0 \in \mathbb{R}^d$,
	\begin{align}\label{eq:lemma:DifferentiatingUnderIntegral-30}
		\begin{split}
			\partial_{i_1}\partial_{i_1}\left( \int_0^\infty P_t h(x) dt\right) \Big|_{x=x_0}&= \int_0^\infty  e^{-2t}  \mathrm{E}\left[\big(D^2 h_{s, \lambda} \circ f\big)(V_0^t)\: \partial_{i_1} f(V_0^t) \: \partial_{i_1} f(V_0^t)\right] dt\\
			&\equiv \int_0^\infty\int_{\mathbb{R}^d} e^{-2t} \partial_{i_1} \partial_{i_1} h(u) \psi(x_0, u, t)du dt\\	
			&\equiv \int_0^\infty e^{-2t}P_t \left(\partial_{i_1}\partial_{i_1} h\right) (x_0)dt.
		\end{split}
	\end{align}

	\vspace{10pt}
	\noindent
	\textbf{Inductive step from $k$ to $k+1$.} 
	\vspace{10pt}
	
	\noindent
	Suppose that for arbitrary indices $1 \leq i_1, \ldots, i_k \leq d$, $k \geq 2$,
	\begin{align}\label{eq:lemma:DifferentiatingUnderIntegral-31}
		\begin{split}
			\partial_{i_k} \cdots \partial_{i_1}\left( \int_0^\infty P_t h(x) dt\right) \Big|_{x=x_0} 
			&= \int_0^\infty\int_{\mathbb{R}^d} e^{-kt} \partial_{i_k} \cdots \partial_{i_1} h(u) \psi(x_0, u, t)du dt\\		
			&\equiv \int_0^\infty e^{-kt}P_t \left(\partial_{i_k} \cdots \partial_{i_1} h \right) (x_0)dt.
		\end{split}
	\end{align}
	Under the induction hypothesis~\eqref{eq:lemma:DifferentiatingUnderIntegral-31}  the argument that gave identity~\eqref{eq:lemma:DifferentiatingUnderIntegral-5} also gives
	\begin{align}\label{eq:lemma:DifferentiatingUnderIntegral-31-2}
		\partial_{i_{k+1}} \cdots \partial_{i_1}\left( \int_0^\infty P_t h(x) dt\right) \Big|_{x=x_0} =  \int_0^\infty\int_{\mathbb{R}^d} e^{-kt} \partial_{i_k} \cdots \partial_{i_1} h(u) \partial_{i_{k+1}} \psi(x_0, u, t)du dt.
	\end{align}
	As in above case with $k=2$, we now show that one can ``pull back'' the partial derivative $\partial_{i_{k+1}}$ from $\psi$ onto $\partial_{i_k} \cdots \partial_{i_1} h$. Let $1 \leq i_{k+1} \leq d$ be an arbitrary index and, for $\eta > 0$, define
	\begin{align*}
		h_{i_1, \ldots, i_k}^{\eta, i_{k+1}}(u) &:= (D^k h_{s, \lambda} \circ f)(u) \prod_{j=1}^k (\varrho_\eta \ast_{(i_{k+1})} \partial_{i_j} f) (u),
	\end{align*}
	where $(\varrho_\eta \ast_{(i_{k+1})} \partial_{i_j} f) (u)$ denotes the ``partial regularization'' in the $i_{k+1}$th coordinate as defined in eq.~\eqref{eq:lemma:DifferentiatingUnderIntegral-5-2}. The map $u \mapsto \partial_{i_{k+1}} h_{i_1, \ldots, i_k}^{\eta, i_{k+1}}(u)$ exists for all $u \in \mathbb{R}^d$, and satisfies, by the chain rule,
	\begin{align}
		&\partial_{i_{k+1}} h_{i_1, \ldots, i_k}^{\eta, i_{k+1}}(u) \nonumber\\
		&\quad{}= \big(D^{k+1} h_{s, \lambda} \circ f\big)(u) \left(\prod_{j=1}^k (\varrho_\eta \ast_{(i_{k+1})} \partial_{i_j} f) (u)\right) \partial_{i_{k+1}} f(u)\label{eq:lemma:DifferentiatingUnderIntegral-32}\\
		&\quad{}\quad{} + \big(D^{k+1} h_{s, \lambda} \circ f\big)(u) \sum_{j : i_j \neq i_{k+1}} \left(\prod_{\ell \neq j} (\varrho_\eta \ast_{(i_{k+1})} \partial_{i_\ell} f) (u)\right) \partial_{i_{k+1}}(\varrho_\eta \ast_{(i_{k+1})} \partial_{i_j}f) (u)\label{eq:lemma:DifferentiatingUnderIntegral-33}\\
		&\quad{}\quad{} + \big(D^{k+1} h_{s, \lambda} \circ f\big)(u) \sum_{j : i_j = i_{k+1}} \left(\prod_{\ell \neq j} (\varrho_\eta \ast_{(i_{k+1})} \partial_{i_\ell} f) (u)\right) \partial_{i_{k+1}}(\varrho_\eta \ast_{(i_{k+1})} \partial_{i_{k+1}}f) (u).\label{eq:lemma:DifferentiatingUnderIntegral-34}
	\end{align}
	Notice that the partial derivative $\partial_{i_{k+1}}(\varrho_\eta \ast_{(i_{k+1})} \partial_{i_j}f) (u)$ in eq.~\eqref{eq:lemma:DifferentiatingUnderIntegral-33} and $\partial_{i_{k+1}}(\varrho_\eta \ast_{(i_{k+1})} \partial_{i_{k+1}}f) (u)$ in eq.~\eqref{eq:lemma:DifferentiatingUnderIntegral-34} follow the patterns derived in eq.~\eqref{eq:lemma:DifferentiatingUnderIntegral-8} and~\eqref{eq:lemma:DifferentiatingUnderIntegral-22}, respectively.
	
	Next, expand the right hand side of~\eqref{eq:lemma:DifferentiatingUnderIntegral-31-2} as
	\begin{align}
		&\int_0^\infty\int_{\mathbb{R}^d} e^{-kt} h_{i_1, \ldots, i_k}^{\eta, i_{k+1}}(u) \partial_{i_{k+1}}\psi (x_0, u, t)du dt    \label{eq:lemma:DifferentiatingUnderIntegral-35}\\
		&\quad{} + \int_0^\infty\int_{\mathbb{R}^d}  e^{-kt} \left(\partial_{i_1, \ldots, i_k}h (u)- h_{i_1, \ldots, i_k}^{\eta, i_{k+1}}(u)  \right) \partial_{i_{k+1}}\psi  (x_0, u, t)du dt.\label{eq:lemma:DifferentiatingUnderIntegral-36}
	\end{align}
	Recall the argument developed to establish the limit~\eqref{eq:lemma:DifferentiatingUnderIntegral-12} for $\eta \downarrow 0$. The same argument, now combined with~\eqref{eq:lemma:DifferentiatingUnderIntegral-32}--\eqref{eq:lemma:DifferentiatingUnderIntegral-34},  also yields that the integral in~\eqref{eq:lemma:DifferentiatingUnderIntegral-35} converges to 
	\begin{align}\label{eq:lemma:DifferentiatingUnderIntegral-37}
		&\int_0^\infty  e^{-(k+1)t}  \mathrm{E}\left[\big(D^{k+1} h_{s, \lambda} \circ f\big)(V_0^t)\: \prod_{j=1}^{k+1}\partial_{i_j} f(V_0^t)\right] dt \nonumber\\
		&\quad{}\equiv \int_0^\infty\int_{\mathbb{R}^d} e^{-(k+1)t} \partial_{i_{k+1}} \cdots \partial_{i_1} h(u) \psi(x_0, u, t)du dt\nonumber\\
		&\quad{}\equiv \int_0^\infty e^{-(k+1)t}P_t \left(\partial_{i_{k+1}} \cdots \partial_{i_1} h \right) (x_0)dt.
	\end{align}	
	Similarly, the same argument used to show that the integral in~\eqref{eq:lemma:DifferentiatingUnderIntegral-10} vanishes for $\eta \downarrow 0$ also guarantees that the integral in~\eqref{eq:lemma:DifferentiatingUnderIntegral-36} vanishes. Combine eq.~\eqref{eq:lemma:DifferentiatingUnderIntegral-31-2} and~\eqref{eq:lemma:DifferentiatingUnderIntegral-35}--\eqref{eq:lemma:DifferentiatingUnderIntegral-37} to conclude the inductive step from $k$ to $k+1$ for $k \geq 2$.
	
	\vspace{10pt}
	\noindent
	\textbf{Proof of part (ii).}
	\vspace{10pt}
	
	\noindent	
	Combine Lemma~\ref{lemma:Smooth-Lipschitz-Approx} (i),~\ref{lemma:ChainRule-SecondOrder}, and~\ref{lemma:DerivativesMaxNorm} and conclude that, for all $x \in \mathbb{R}^d \setminus \mathcal{N}$,
	\begin{align*}
		&\left| \partial_{i_k} \ldots \partial_{i_1} h \right| (x) \nonumber\\
		&\quad{} =\left| \left(D^k h_{s, \lambda} \circ f\right) \right|(x)    \left| \partial_{i_1} f \right|(x) \cdots  \left|  \partial_{i_k}\right| (x)\nonumber\\
		&\quad{}\leq C_k \lambda^{-k} \mathbf{1}\left\{ s \leq \|x\|_\infty \leq s + 3 \lambda \right\} \mathbf{1}\left\{|x_{i_1}| \geq |x_\ell|, \: \forall \ell \right\} \cdots \mathbf{1}\left\{|x_{i_k}| \geq |x_\ell|, \: \forall \ell \right\},
	\end{align*}
	where $C_k  >0$ is the absolute constant from Lemma~\ref{lemma:Smooth-Lipschitz-Approx}.
	Notice that
	\begin{align*}
		\mathbf{1}\left\{|x_{i_1}| \geq |x_\ell|, \: \forall \ell \right\} \cdots \mathbf{1}\left\{|x_{i_k}| \geq |x_\ell|, \: \forall \ell \right\} = 1 \quad{} \Longleftrightarrow \quad{} |x_{i_1}| = \ldots = |x_{i_k}|,
	\end{align*}
	and, hence,
	\begin{align*}
		\mathcal{L}^d\left( \left\{x \in \mathbb{R}^d :\mathbf{1}\left\{|x_{i_1}| \geq |x_\ell|, \: \forall \ell \right\} \cdots \mathbf{1}\left\{|x_{i_k}| \geq |x_\ell|, \: \forall \ell \right\} = 1 \right\}\right) > 0\:\: \Leftrightarrow \:\: i_1= \ldots = i_k.
	\end{align*}
	Thus, there exists a $\mathcal{L}^d$-null set $\mathcal{N}' \supseteq \mathcal{N}$ such that for all $x \in \mathbb{R}^d \setminus \mathcal{N}'$,
	\begin{align}\label{eq:lemma:DifferentiatingUnderIntegral-1}
		\left| \partial_{i_k} \ldots \partial_{i_1} h \right|(x) \leq C_k \lambda^{-k} \mathbf{1}\left\{ s \leq \|x\|_\infty \leq s + 3 \lambda \right\} \mathbf{1}\left\{|x_{i_1}| \geq |x_{\ell}|, \: \ell \neq {i_1} \right\} \mathbf{1}\{i_1 = \ldots = i_k\}.
	\end{align}
	Since $\Sigma$ is positive definite, $N(0, \Sigma)$ is absolutely continuous with respect to $\mathcal{L}^d$. Thus, the $\mathcal{L}^d$-a.e. upper bound~\eqref{eq:lemma:DifferentiatingUnderIntegral-1} continues to hold when evaluated at $V_0^t = e^{-t}x_0 + \sqrt{1- e^{-2t}}Z$ and integrated over $Z \sim N(0, \Sigma)$ and $t \sim Exp(k)$. We conclude that for all $x_0 \in \mathbb{R}^d$,
	\begin{align}\label{eq:lemma:DifferentiatingUnderIntegral-2}
		\begin{split}
			&\left| \int_0^\infty e^{-kt} P_t \left( \partial_{i_k} \ldots \partial_{i_1} h \right) (x_0)dt \right|\\
			&\:\:\: \leq C_k\lambda^{-k} \int_0^\infty e^{-kt}\mathrm{E}\left[\mathbf{1}\left\{ s \leq \|V_0^t\|_\infty \leq s + 3 \lambda \right\} \mathbf{1}\left\{|V_{0{i_1}}^t| \geq |V_{0\ell}^t|, \: \forall \ell \right\} \right] dt \:  \mathbf{1}\{i_1 = \ldots = i_k\}.
		\end{split}
	\end{align}
\end{proof}

\begin{proof}[\textbf{Proof of Lemma~\ref{lemma:L1ConvergencePartialRegularization}}]
	For any function $h$ on $\mathbb{R}^d$ and vector $y \in \mathbb{R}^d$, we define the translation operator $\tau$ by $\tau_yh(x) = h(x-y)$. With this notation,
	\begin{align*}
		(\varrho_\eta \ast_{(i)} h)  = \int \varrho(r) h(x - r\eta e_i) dr = \int \varrho(r)\tau_{ r\eta e_i} h(x) dr
	\end{align*}
	Without loss of generality we can assume that the $f_j$'s are bounded by one. Also, since $\varrho$ integrates to one, we have 
	\begin{align*}
		(\varrho_\eta \ast_{(i_j)} f_j) g(x) - f_j(x)g(x) = \int \varrho(r) \left(\tau_{ r\eta e_{i_j}} f_j(x) - f_j(x)  \right)g(x) dr.
	\end{align*}
	Hence, by the product comparison inequality,
	\begin{align}\label{eq:lemma:L1ConvergencePartialRegularization}
		\left\| \prod_{j=1}^k(\varrho_\eta \ast_{(i_j)} f_j) g - \prod_{j=1}^kf_jg \right\|_1 &= \int \left| \int \varrho(r) \prod_{j=1}^k\left(\tau_{ r\eta e_{i_j}} f_j(x)g(x)  - f_j(x) g(x)  \right)dr \right| dx \nonumber\\
		&\leq \int \int |\varrho(r) | \left|  \prod_{j=1}^k\left(\tau_{ r\eta e_{i_j}} f_j(x)g(x) - f_j(x)g(x)  \right)\right| dr dx\nonumber\\
		&\leq \sum_{k=1}^k  \int \int |\varrho(r)|  \left| \tau_{ r\eta e_{i_j}} f_j(x)g(x) - f_j(x)g(x)\right| dr dx\nonumber\\
		&\leq \sum_{k=1}^k \int |\varrho(r)| \| (\tau_{ r\eta e_{i_j}} f_j)g - f_jg \|_1 dr.
	\end{align}
	Next, compute
	\begin{align*}
		\int |\varrho(r)| \|\tau_{ r\eta e_{i_j}} f_j g - f_j g \|_1 dr	&\leq \int  |\varrho(r)|  \| (\tau_{r\eta e_{i_j}} f_j) g - \tau_{r\eta e_{i_j}}(f_jg) \|_1 dr + \int \varrho(r) \| \tau_{r\eta e_{i_j}}(f_jg) - f_jg \|_1 dr\\
		&=\int  |\varrho(r)|  \| (\tau_{r\eta e_{i_j}} f_j) ( g - \tau_{r\eta e_{i_j}}g) \|_1 dr + \int  |\varrho(r)|  \| \tau_{r\eta e_{i_j}}(f_jg) - f_jg \|_1 dr\\
		&\leq\int  |\varrho(r)|  \|g - \tau_{r\eta e_{i_j}}g \|_1 dr + \int  |\varrho(r)|  \| \tau_{r\eta e_{i_j}}(f_jg) - f_jg \|_1 dr.
	\end{align*}	
	Since $\|g - \tau_{z}g \|_1$ and $\| \tau_{z}(f_jg) - f_jg \|_1$ are both bounded by $2 \|g\|_1 < \infty$ for all $z \in \mathbb{R}$, Proposition 8.5 in~\cite{folland1999real} implies that, for all $r \in \mathbb{R}$,
	\begin{align*}
		\|g - \tau_{r\eta e_i}g \|_1  \vee  \| \tau_{r\eta e_i}(f_jg) - f_jg \|_1 \rightarrow 0 \quad{} \mathrm{as} \quad{}\eta \rightarrow 0.
	\end{align*}
	Hence, by the dominated convergence theorem
	\begin{align*}
		\int |\varrho(r)| \|\tau_{ r\eta e_{i_j}} f_j g - f_j g \|_1 dr \rightarrow 0 \quad{} \mathrm{as} \quad{}\eta \rightarrow 0.
	\end{align*}
	Conclude that each summand in eq.~\eqref{eq:lemma:L1ConvergencePartialRegularization} vanishes as $\eta \downarrow 0$. This completes the proof.
\end{proof}

\subsection{Proofs of Lemmas~\ref{lemma:Kolmogorov-Coupling-AntiConcentration},~\ref{lemma:VarianceMaximum},~\ref{lemma:Bootstrap-Max-Norm-Quantil-Comparison}, and~\ref{lemma:Bootstrap-Max-Norm-Quantil-Comparison-Non-Gaussian}}

\begin{proof}[\textbf{Proof of Lemma~\ref{lemma:Kolmogorov-Coupling-AntiConcentration}}]
	We provide a full proof for completeness only. Incidentally, this result has already been established by~\cite{lecam1986asymptotic}, p. 402 (Lemma 2). We compute,
	\begin{align*}
		\mathbb{P}\left\{X \leq s \right\} &= \mathbb{P}\left\{X  \leq s, \: Z  \leq  s+ \varepsilon \right\} +\mathbb{P}\left\{X \leq s, \: Z  > s+ \varepsilon \right\} \\
		&\leq\mathbb{P}\left\{Z  \leq  s+ \varepsilon \right\} + \mathbb{P}\left\{ | Z - X|  > \varepsilon \right\}\\
		& \leq  \mathbb{P}\left\{ Z \leq s \right\} +  \mathbb{P}\left\{s \leq Z \leq s + \varepsilon \right\}  + \mathbb{P}\left\{ | Z  - X |  > \varepsilon \right\}.
	\end{align*}	
	For the reverse inequality,
	\begin{align*}
		\mathbb{P}\left\{Z \leq s - \varepsilon \right\}  &\leq  \mathbb{P}\left\{ Z \leq  s-\varepsilon, \: X  \leq s \right\} + 	
		\mathbb{P}\left\{ Z \leq  s- \varepsilon, \:X  > s \right\} \\
		&\leq \mathbb{P}\left\{ X \leq s \right\} + \mathbb{P}\left\{| Z -X |  > \varepsilon \right\},
	\end{align*}
	and, hence,
	\begin{align*}
		\mathbb{P}\left\{Z \leq s\right\}  \leq \mathbb{P}\left\{X \leq s \right\} + \mathbb{P}\left\{|Z - X |  > \varepsilon \right\} + \mathbb{P}\left\{ s -\varepsilon \leq Z\leq s\right\}.
	\end{align*}
	Take the supremum over $s \geq 0$ and combine both inequalities. Then switch the roles of $X$ and $Z$ to conclude the proof.
\end{proof}

\begin{proof}[\textbf{Proof of Lemma~\ref{lemma:VarianceMaximum}}]
	Observe that $X \vee Z = X + (Z - X)_+$. Now, by Cauchy-Schwarz,
	\begin{align*}
		\mathrm{Var}(X \vee Z) &= \mathrm{Var}(X) + \mathrm{Var}\big( (Z-X)_+\big) + 2 \mathrm{Cov}\big(X, (Z-X)_+\big)\\
		&\leq \mathrm{Var}(X) + \mathrm{Var}\big( (Z-X)_+\big) + 2\sqrt{ \mathrm{Var}(X) } \sqrt{\mathrm{Var}\big( (Z-X)_+\big)}\\
		&= \left(\sqrt{\mathrm{Var}(X)} + \sqrt{\mathrm{Var}\big( (Z-X)_+\big)}\right)^2.
	\end{align*}
	Similarly, use the estimate $2 \mathrm{Cov}\big(X, (Z-X)_+\big) \geq -2\sqrt{ \mathrm{Var}(X) } \sqrt{\mathrm{Var}\big( (Z-X)_+\big)}$ in the first line of above display, to obtain
	\begin{align*}
		\mathrm{Var}(X \vee Z) \geq \left(\sqrt{\mathrm{Var}(X)} - \sqrt{\mathrm{Var}\big( (Z-X)_+\big)}\right)^2.
	\end{align*}
	Combine both inequalities to obtain the desired two-sided inequality. Further, if $\mathrm{E}[Z-X] \geq 0$, then by convexity of the map $a \mapsto (a)_+$, $a \in \mathbb{R}$, and Jensen's inequality,
	\begin{align*}
		\mathrm{Var}\big((Z-X)_+\big) = \mathrm{E}\big[(Z-X)_+^2\big] - \mathrm{E}\big[(Z-X)_+\big]^2 \leq \mathrm{E}\big[(Z-X)^2\big] - \big(\mathrm{E}[Z-X])_+^2 = \mathrm{Var}(Z-X).
	\end{align*}
\end{proof}

\begin{proof}[\textbf{Proof of Lemma~\ref{lemma:Bootstrap-Max-Norm-Quantil-Comparison}}]
	The proof is identical to the one of Lemma 3.2 in~\cite{chernozhukov2013GaussianApproxVec}. Let $\pi_n(\delta) = \delta^{1/3} \left(\mathrm{Var}(\|\Sigma^{1/2} Z\|_\infty)\right)^{-1/3}$. By Lemma~\ref{lemma:GaussianComparison} there exists an absolute constant $K > 0$ such that on the event $\left\{\max_{j,k} |\widehat{\Sigma}_{n,jk} - \Sigma_{jk}| \leq \delta \right\}$, for all $s \geq 0$,
	\begin{align*}
		\Big|\mathbb{P}\left\{\|\widehat{\Sigma}_n^{1/2} Z\|_\infty \leq s \mid X_1, \ldots, X_n \right\} -\mathbb{P}\left\{\|\Sigma_n^{1/2} Z\|_\infty \leq s \right\} \Big| \leq K \pi_n(\delta).
	\end{align*}
	In particular, for $s = c_n(K\pi_n(\delta; \Sigma) + \alpha)$, we obtain
	\begin{align*}
		\mathbb{P}\left\{\|\widehat{\Sigma}_n^{1/2} Z\|_\infty \leq  c_n(K\pi_n(\delta) + \alpha; \Sigma) \mid X_1, \ldots, X_n \right\}  &\geq \mathbb{P}\left\{\|\Sigma_n^{1/2} Z\|_\infty \leq  c_n(K\pi_n(\delta; \Sigma) + \alpha) \right\} - K\pi_n(\delta) \\
		& \geq K\pi_n(\delta) + \alpha - K\pi_n(\delta) = \alpha.
	\end{align*}
	To conclude the proof of the first statement apply the definition of quantiles. The second claim follows in the same way.
\end{proof}

\begin{proof}[\textbf{Proof of Lemma~\ref{lemma:Bootstrap-Max-Norm-Quantil-Comparison-Non-Gaussian}}]
	Let $\delta, \eta > 0$ be arbitrary and set $\gamma_n = \sup_{s \geq 0} | \mathbb{P}\{S_n \leq s \mid X_1, \ldots, X_n\} - \mathbb{P}\{\|\Sigma_n^{1/2}Z\|_\infty \leq s \}|$, $\kappa_n(\delta) =  \delta (\mathrm{Var}(\|\Sigma_n^{1/2} Z\|_\infty))^{-1/2}$, $\rho_n(\delta)= \mathbb{P}\{|R_n| > \delta \mid X_1, \ldots, X_n\}$. By Lemma~\ref{lemma:AntiConcentration-SeparableProcess} there exists an absolute constant $K > 0$ such that on the event $\{\gamma_n + \rho(\delta) \leq \eta\}$,
	\begin{align*}
		&\left| \mathbb{P}\left\{T_n \leq s \mid X_1, \ldots, X_n\right\}  - \mathbb{P}\left\{\|\Sigma_n^{1/2}Z\|_\infty\leq s\right\}  \right|\\
		&\quad\leq \left| \mathbb{P}\left\{T_n \leq s, |R_n| \leq \delta \mid X_1, \ldots, X_n \right\}  - \mathbb{P}\left\{ \|\Sigma_n^{1/2}Z\|_\infty\leq s\right\}  \right| + \mathbb{P}\left\{|R_n| > \delta \mid X_1, \ldots, X_n \right\}\\
		&\quad \leq \sup_{s \geq 0} \left|  \mathbb{P}\left\{ S_n \leq s \mid X_1, \ldots, X_n \right\}  - \mathbb{P}\left\{ \|\Sigma_n^{1/2}Z\|_\infty\leq s\right\} \right| \\
		&\quad \quad + \sup_{s \geq 0} \mathbb{P}\left\{ s \leq\|\Sigma_n^{1/2} Z\|_\infty \leq s + \delta \right\} + \mathbb{P}\left\{|R_n| > \delta \mid X_1, \ldots, X_n \right\}\\
		&\quad \leq \gamma_n + \kappa_n(\delta) + \rho_n(\delta)\\
		&\quad \leq K \kappa_n(\delta) + \eta.
	\end{align*}
	Hence, for $s = c_n(K\kappa_n(\delta) + \eta + \alpha; \Sigma)$, we obtain
	\begin{align*}
		&\mathbb{P}\left\{T_n \leq  c_n(K\kappa_n(\delta) + \eta + \alpha; \Sigma) \right\} \\
		&\quad \geq \mathbb{P}\left\{\|\Sigma_n^{1/2} Z\|_\infty \leq  c_n(K\kappa_n(\delta) + \eta + \alpha; \Sigma) \right\} - K\kappa_n(\delta) - \eta \\
		&\quad \geq K\kappa_n(\delta) + \eta + \alpha - K\kappa_n(\delta) - \eta\\
		&\quad = \alpha.
	\end{align*}
	The first statement now follows from the definition of quantiles. The second claim follows in the same way.
\end{proof}

\subsection{Proof of Lemma~\ref{lemma:ModulusContinuity}}\label{subsec:ProofModulusContinuity}

\begin{proof}[\textbf{Proof of Lemma~\ref{lemma:ModulusContinuity}}]
	First, we show necessity: Suppose that the sample paths of $X$ on $U$ are almost surely uniformly continuous w.r.t. $d_X$, i.e. for almost every $\omega \in \Omega$,
	\begin{align*}
		\lim_{\delta \rightarrow 0} \sup_{d_X(u,v) < \delta} |X_u(\omega) - X_v(\omega)| = 0.
	\end{align*}
	Since $X$ is Gaussian, by Lemma~\ref{lemma:AS-Bounded}, for $\delta > 0$ arbitrary,
	\begin{align*}
		\mathrm{E}\left[ \sup_{d_X(u,v) < \delta} |X_u(\omega) - X_v(\omega)| \right] \leq 2\mathrm{E}[Z] < \infty.
	\end{align*}
	Hence, the claim follows from the dominated convergence theorem.
	
	Next, we show sufficiency: Given the premise, we can find a null sequence $\{\delta_n\}_{n \geq 1}$ such that 
	\begin{align}\label{eq:proof:lemma:ModulusContinuity-1}
		\mathrm{E} \left[ \sup_{d_X(u, v) < \delta_n} |X_u - X_v| \right] \leq 2^{-n}
	\end{align}
	Define events
	\begin{align*}
		A_n = \left\{  \sup_{d_X(u,v) < \delta_n \wedge 2^{-n} } |X_u - X_v| > 2^{-n/2}\right\}.
	\end{align*}
	By Markov's inequality and~\eqref{eq:proof:lemma:ModulusContinuity-1},
	\begin{align*}
		\sum_{i=1}^\infty \mathbb{P}\{A_n\} \leq \sum_{i=1}^\infty 2^{-n/2} = 1 + \sqrt{2} < \infty.
	\end{align*}
	Thus, by Borel-Cantelli, $\mathbb{P}\{\limsup_{n \rightarrow \infty} A_n\} = 0$, i.e. $X$ is almost surely uniformly continuous on $U$. (Note that we have not used the fact that $X$ is Gaussian and $d_X$ the intrinsic standard deviation. Hence, sufficiency holds for arbitrary stochastic processes on general metric spaces.)
\end{proof}

\subsection{Proofs of Lemmas~\ref{lemma:Tesseract},~\ref{lemma:BC-KRR-Expansion},~\ref{lemma:BC-KRR-Remainder},~\ref{lemma:BC-KRR-Consistency-Operator}, and~\ref{lemma:Consistency-Gram-Matrix}}

\begin{proof}[\textbf{Proof of Lemma~\ref{lemma:Tesseract}}]
	Throughout the proof $\otimes$ denotes the tensor product between linear maps.
	
	Let $M_n = \max_{1 \leq i \leq n} \|X_i\|_2^2$. Define the map $x \mapsto \varphi(x) := x^2\mathbf{1}\{|x| \leq M_n\} + M_n^2 \mathbf{1}\{|x| > M_n\}$. This map is Lipschitz continuous with Lipschitz constant $2M_n$ and $\varphi(0) = 0$. Let $\varepsilon_1, \ldots, \varepsilon_n$ be i.i.d. Rademacher random variables independent of $X_1, \ldots, X_n$. Then, symmetrization and contraction principle applied to $\varphi(\cdot)/(2M_n)$ yield
	\begin{align*}
		\mathrm{E}  \left\| \frac{1}{n}\sum_{i=1}^n X_i \otimes X_i \otimes X_i \otimes X_i - T_4 \right\|_{op} 
		&\leq 4 \mathrm{E} \left[ M_n \sup_{\substack{\|u\|_2 =1\\ \|v\|_2 =1} } \left| \frac{1}{n} \sum_{i=1}^n \varepsilon_i (X_i'u)(X_i'v)  \right| \right].
	\end{align*}
	We upper bound the right-hand side in above display using Cauchy-Schwarz, Hoffmann-J\o{}rgensen, and de-symmetrization inequalities by
	\begin{align*}
		& \left(\mathrm{E} \left[ M_n^2\right]\right)^{1/2} \left(\mathrm{E} \left[\sup_{\substack{\|u\|_2 =1\\ \|v\|_2 =1} } \left( \frac{1}{n} \sum_{i=1}^n \varepsilon_i (X_i'u)(X_i'v)  \right)^2 \right]\right)^{1/2}\\
		&\quad  \lesssim \left(\mathrm{E} \left[ M_n^2\right]\right)^{1/2} \left\{  \mathrm{E} \left\| \frac{1}{n} \sum_{i=1}^n X_i \otimes X_i - \Sigma \right\|_{op} + n^{-1} \left(\mathrm{E}\left[ M_n^2 \right]\right)^{1/2}  \right\}.
	\end{align*}
	Since $X = (X_1, \ldots, X_d)' \in \mathbb{R}^d$ is sub-Gaussian with mean zero and covariance $\Sigma$, we compute
	\begin{align*}
		\left\| X'X \right\|_{\psi_1} \leq \sum_{k=1}^d \left\| X_k^2  \right\|_{\psi_1} \leq \sum_{k=1}^d \left\| X_k \right\|_{\psi_2}^2 \lesssim \sum_{k=1}^d \mathrm{Var}(X_k) = \mathrm{tr}(\Sigma).
	\end{align*}
	Hence, by Lemma 35 in~\cite{giessing2021debiased} and Lemma 2.2.2 in~\cite{vandervaart1996weak} (see Exercise 2.14.1 (ibid.) for how to handle the non-convexity of the map inducing the $\psi_{1/2}$-norm)
	\begin{align*}
		\mathrm{E} \left[ M_n^2\right] = \mathrm{E} \left[ \max_{1 \leq i \leq n}  (X_i'X_i)^2\right] \lesssim ( \log en)^2 \max_{1 \leq i \leq n} \left\|  (X_i'X_i)^2 \right\|_{\psi_{1/2}} \lesssim (\log en)^2  \mathrm{tr}^2(\Sigma).
	\end{align*}
	Moreover, by Theorem 4 in~\cite{koltchinskiiConcentration2017},
	\begin{align*}
		\mathrm{E} \left\| \frac{1}{n} \sum_{i=1}^n X_i \otimes X_i - \Sigma \right\|_{op}\lesssim \|\Sigma\|_{op} \left( \sqrt{\frac{\mathrm{r}(\Sigma)}{n} } \vee \frac{\mathrm{r}(\Sigma)}{n} \right).
	\end{align*}
	Thus, we conclude that
	\begin{align*}
		&\mathrm{E}  \left\| \frac{1}{n}\sum_{i=1}^n X_i \otimes X_i \otimes X_i \otimes X_i - T_4 \right\|_{op} \\
		&\quad\lesssim (\log en) \mathrm{r}(\Sigma) \|\Sigma\|_{op}^2  \left( \sqrt{\frac{\mathrm{r}(\Sigma)}{n}} \vee \frac{ \mathrm{r}(\Sigma)}{n} \right) +\frac{ (\log en)^2 \mathrm{r}^2(\Sigma)  \|\Sigma\|_{op}^2 }{n}.
	\end{align*}
	Adjust some absolute constants to complete the proof.
\end{proof}

\begin{proof}[\textbf{Proof of Lemma~\ref{lemma:BC-KRR-Expansion}}]
	To simplify notation we write $\varepsilon_i = Y_i - f_0(X_i)$ for $1 \leq i \leq n$.
	
	We begin with the following fundamental identity: Let $I$ be the identity operator and $P$ be such that $I + P$ is invertible. Then, 
	\begin{align*}
		(I + P)^{-1} P - I = -(I + P)^{-1}.
	\end{align*}
	Indeed, we have $I = (I + P)^{-1}(I + P) = (I +P)^{-1} + (I + P)^{-1} P$. Now, re-arrange the terms on the far left and right hand side of this identity to conclude. Applied to $P = \lambda^{-1}\widehat{T}_n$, we obtain
	\begin{align}\label{eq:lemma:BC-KRR-Expansion-0}
		(\widehat{T}_n + \lambda)^{-1} \widehat{T}_n - I = (\lambda^{-1} \widehat{T}_n+ I)^{-1} \lambda^{-1}\widehat{T}_n- I = - (\lambda^{-1} \widehat{T}_n+ I)^{-1} 
		=  - \lambda (\widehat{T}_n + \lambda)^{-1}.
	\end{align}
	Hence, we compute
	\begin{align}\label{eq:lemma:BC-KRR-Expansion-1}
		\widehat{f}^{\mathrm{bc}}_n - f_0 &= (\widehat{T}_n + \lambda)^{-1}\frac{1}{n}\sum_{i=1}^n k_{X_i} Y_i + \lambda(\widehat{T}_n + \lambda)^{-1} \widehat{f}_n  - f_0\nonumber\\		
		&= (\widehat{T}_n + \lambda)^{-1} \frac{1}{n}\sum_{i=1}^n\varepsilon_i k_{X_i}  + (\widehat{T}_n + \lambda)^{-1} \widehat{T}_nf_0 + \lambda(\widehat{T}_n + \lambda)^{-1} \widehat{f}_n  - f_0\nonumber\\
		&\overset{(a)}{=} (\widehat{T}_n + \lambda)^{-1} \frac{1}{n}\sum_{i=1}^n\varepsilon_i k_{X_i} + \lambda(\widehat{T}_n + \lambda)^{-1} (\widehat{f}_n - f_0)\nonumber\\
		&= (\widehat{T}_n + \lambda)^{-1} \frac{1}{n}\sum_{i=1}^n \varepsilon_i k_{X_i}  + \lambda  (\widehat{T}_n + \lambda)^{-2} \frac{1}{n}\sum_{i=1}^n \varepsilon_i k_{X_i} \nonumber\\
		&\quad  + \lambda (\widehat{T}_n + \lambda)^{-1} \left( (\widehat{T}_n + \lambda)^{-1}  \widehat{T}_n f_0 - f_0 \right)  \nonumber\\
		&\overset{(b)}{=}  (\widehat{T}_n + \lambda)^{-1} \frac{1}{n}\sum_{i=1}^n \varepsilon_i k_{X_i}  - \lambda  (\widehat{T}_n + \lambda)^{-2} \frac{1}{n}\sum_{i=1}^n \varepsilon_i k_{X_i}  - \lambda^2 (\widehat{T}_n + \lambda)^{-2}f_0\nonumber\\
		& =  (\widehat{T}_n + \lambda)^{-1} \left( I - \lambda  (\widehat{T}_n + \lambda)^{-1} \right)\frac{1}{n}\sum_{i=1}^n \varepsilon_i k_{X_i} - \lambda^2 (\widehat{T}_n + \lambda)^{-2}f_0\nonumber\\
		&\overset{(c)}{=} (\widehat{T}_n + \lambda)^{-2}\widehat{T}_n \frac{1}{n}\sum_{i=1}^n \varepsilon_i k_{X_i} - \lambda^2 (\widehat{T}_n + \lambda)^{-2}f_0,
	\end{align}
	where (a), (b), and (c) follow from identity~\eqref{eq:lemma:BC-KRR-Expansion-0}.
	
	We now further upper bound the terms in~\eqref{eq:lemma:BC-KRR-Expansion-1}. Re-walking the steps from (c) to (b) we obtain
	\begin{align}\label{eq:lemma:BC-KRR-Expansion-2}
		&\left\| \left( (\widehat{T}_n + \lambda)^{-2}\widehat{T}_n - (T+ \lambda)^{-2}T\right) \frac{1}{n}\sum_{i=1}^n \varepsilon_i k_{X_i} \right\|_{\infty} \nonumber\\
		\begin{split}
			&\quad \quad  \leq  \left\| \left((\widehat{T}_n + \lambda)^{-1} - (T + \lambda)^{-1} \right) \frac{1}{n}\sum_{i=1}^n \varepsilon_i k_{X_i} \right\|_{\infty} \\
			&\quad \quad \quad + \lambda \left\| \left((\widehat{T}_n + \lambda)^{-2} - (T + \lambda)^{-2} \right) \frac{1}{n}\sum_{i=1}^n \varepsilon_i k_{X_i} \right\|_{\infty}.
		\end{split}
	\end{align}
	Note that
	\begin{align}\label{eq:lemma:BC-KRR-Expansion-3}
		\| (\widehat{T}_n + \lambda)^{-1}(T+ \lambda) - I \|_{op} &\leq \| (\widehat{T}_n + \lambda)^{-1}\|_{op}\|T - \widehat{T}_n\|_{op}  \leq \lambda^{-1} \|T - \widehat{T}_n\|_{op}.
	\end{align}
	Recall that $|a^2 - 1| \leq 3 \max\{|a-1|, |a-1|^2\}$ for all $a  \geq 0$. Let $U_{\mathcal{H}}$ be the unit ball of $\mathcal{H}$. We have
	\begin{align}\label{eq:lemma:BC-KRR-Expansion-4}
		&\| (\widehat{T}_n + \lambda)^{-2} (T+ \lambda)^2 - I \|_{op} \nonumber\\
		&\quad \overset{(a)}{=} \sup_{u \in U_{\mathcal{H}}} \left|\left\langle \left((\widehat{T}_n + \lambda)^{-2} (T+ \lambda)^2 - I \right) u, u \right\rangle_{\mathcal{H}} \right| \nonumber\\
		&\quad \overset{(b)}{=}  \sup_{u \in U_{\mathcal{H}}} \left| \left\langle (\widehat{T}_n + \lambda)^{-1} (T+ \lambda)u , (\widehat{T}_n + \lambda)^{-1} (T+ \lambda)u \right\rangle_{\mathcal{H}} - 1\right| \nonumber\\
		&\quad=  \sup_{u \in U_{\mathcal{H}}} \left| \left\| (\widehat{T}_n + \lambda)^{-1} (T+ \lambda) u \right\|_{\mathcal{H}}^2 - 1\right| \nonumber\\
		&\quad\leq 3\sup_{u \in U_{\mathcal{H}}} \left| \left\| (\widehat{T}_n + \lambda)^{-1} (T+ \lambda) u \right\|_{\mathcal{H}} - 1\right| \vee 3 \sup_{u \in U_{\mathcal{H}}} \left| \left\| (\widehat{T}_n + \lambda)^{-1} (T+ \lambda) u \right\|_{\mathcal{H}} - 1\right|^2 \nonumber\\
		&\quad \overset{(c)}{\leq}  3 \left\| (\widehat{T}_n + \lambda)^{-1} (T+ \lambda) - I\right\|_{op} \vee 3 \left\| (\widehat{T}_n + \lambda)^{-1} (T+ \lambda) - I \right\|_{op}^2\nonumber\\
		&\quad  \overset{(d)}{\leq}  3\lambda^{-1} \|T- \widehat{T}_n\|_{op} \vee 3\lambda^{-2} \|T- \widehat{T}_n\|_{op}^2,
	\end{align}
	where (a) and (b) hold because $(\widehat{T}_n + \lambda)^{-1}$ and $(T+ \lambda)$ are self-adjoint, (c) follows from the reverse triangle inequality applied $\|\cdot\|_{\mathcal{H}}$ and the definition of the operator norm, and (d) follows from~\eqref{eq:lemma:BC-KRR-Expansion-3}.
	
	By the reproducing property of the kernel $k$ (see Remark~\ref{remark:Assumptions-RE}) and~\eqref{eq:lemma:BC-KRR-Expansion-3},
	\begin{align}\label{eq:lemma:BC-KRR-Expansion-6}
		\begin{split}
		&\left\| \left( (\widehat{T}_n + \lambda)^{-1} - (T+ \lambda)^{-1}\right) \frac{1}{n}\sum_{i=1}^n \varepsilon_i k_{X_i} \right\|_{\infty} \leq \kappa \lambda^{-1} \|\widehat{T}_n - T\|_{op}  \left\| \frac{1}{n}\sum_{i=1}^n (T + \lambda)^{-1}\varepsilon_i k_{X_i}\right\|_{\mathcal{H}},
		\end{split}
	\end{align}
	and similarly, by~\eqref{eq:lemma:BC-KRR-Expansion-4},
	\begin{align}\label{eq:lemma:BC-KRR-Expansion-7}
		\begin{split}
		&\lambda \left\| \left((\widehat{T}_n + \lambda)^{-2} - (T + \lambda)^{-2} \right) \frac{1}{n}\sum_{i=1}^n \varepsilon_i k_{X_i} \right\|_{\infty} \\ 
		&\quad\leq \kappa\lambda^{-1} \|\widehat{T}_n - T\|_{op} \left( 3 + 3\lambda^{-1}\|\widehat{T}_n - T\|_{op} \right)  \left\| \frac{1}{n}\sum_{i=1}^n (T + \lambda)^{-1} \varepsilon_i k_{X_i}\right\|_{\mathcal{H}},
		\end{split}
	\end{align}
	and~\eqref{eq:lemma:BC-KRR-Expansion-3} again,
	\begin{align}\label{eq:lemma:BC-KRR-Expansion-8}
		\lambda^2 \|(\widehat{T}_n + \lambda)^{-2}f_0\|_\infty 	&\leq \lambda^2 \kappa \|(\widehat{T}_n + \lambda)^{-1}(T + \lambda)\|_{op}^2\| (T + \lambda)^{-2}f_0\|_{\mathcal{H}} \nonumber\\
		&\leq \lambda^2 \kappa \|(\widehat{T}_n + \lambda)^{-1}(T + \lambda) - I\|_{op}^2  \|(T + \lambda)^2f_0\|_{\mathcal{H}}  + \lambda^2 \kappa \| (T + \lambda)^{-2}f_0\|_{\mathcal{H}} \nonumber\\
		&\leq \kappa \left( \|\widehat{T}_n - T\|_{op}^2  + \lambda^2 \right) \| (T + \lambda)^{-2}f_0\|_{\mathcal{H}}.
	\end{align}
	Combine~\eqref{eq:lemma:BC-KRR-Expansion-1},~\eqref{eq:lemma:BC-KRR-Expansion-2}, and~\eqref{eq:lemma:BC-KRR-Expansion-6}--\eqref{eq:lemma:BC-KRR-Expansion-8} to conclude that
	\begin{align*}
		\sqrt{n} (\widehat{f}^{\mathrm{bc}}_n - f_0) = (T + \lambda)^{-2}T \frac{1}{\sqrt{n}}\sum_{i=1}^n \varepsilon_i k_{X_i} + \sqrt{n} R_n,
	\end{align*}
	where 
	\begin{align*}
		\sqrt{n}\|R_n\|_\infty &\lesssim  \kappa\lambda^{-1} \|\widehat{T}_n - T\|_{op} \left( 1 + \lambda^{-1}\|\widehat{T}_n - T\|_{op} \right)  \left\| \frac{1}{n}\sum_{i=1}^n (T + \lambda)^{-1} \varepsilon_i k_{X_i}\right\|_{\mathcal{H}}\\
		& + \kappa \left( \sqrt{n}\|\widehat{T}_n - T\|_{op}^2  + \sqrt{n}\lambda^2 \right) \| (T + \lambda)^{-2}f_0\|_{\mathcal{H}}.
	\end{align*}	
\end{proof}

\begin{proof}[\textbf{Proof of Lemma~\ref{lemma:BC-KRR-Remainder}}]
	Observe that $\mathrm{tr}(T) \leq \kappa^2$. Therefore, by Lemma~\ref{lemma:Consistency-Gram-Matrix}, with probability at least $1-\delta$,
	\begin{align*}
		\lambda^{-1} \|\widehat{T}_n - T\|_{op}  \left\| \frac{1}{\sqrt{n}}\sum_{i=1}^n (T + \lambda)^{-1}\varepsilon_i k_{X_i}\right\|_{\mathcal{H}} \lesssim \kappa^2 \log^2(2/\delta) \left(\sqrt{\frac{ \bar{\sigma}^2  \mathfrak{n}_1^2(\lambda)  }{n\lambda^2 } } \vee \frac{ \bar{\sigma}^2 }{n \lambda^2}\right),
	\end{align*}
	and
	\begin{align*}
		\lambda^{-2} \|\widehat{T}_n - T\|_{op}^2  \left\| \frac{1}{\sqrt{n}}\sum_{i=1}^n (T + \lambda)^{-1}\varepsilon_i k_{X_i}\right\|_{\mathcal{H}} \lesssim \frac{\kappa^4\log^3(2/\delta)}{\sqrt{n} \lambda} \left(\sqrt{\frac{ \bar{\sigma}^2  \mathfrak{n}_1^2(\lambda)  }{n\lambda^2 } } \vee \frac{ \bar{\sigma}^2 }{n \lambda^2}\right),
	\end{align*}
	and
	\begin{align*}
		\kappa \left( \sqrt{n}\|\widehat{T}_n - T\|_{op}^2  + \sqrt{n}\lambda^2 \right) \| (T + \lambda)^2f_0\|_{\mathcal{H}} \lesssim \kappa \left( \frac{\kappa^4\log^2(2/\delta)}{\sqrt{n}} + \sqrt{n} \lambda^2 \right) \| (T + \lambda)^{-2}f_0\|_{\mathcal{H}}.
	\end{align*}	
	We combine above three inequalities and conclude that, with probability at least $1 - \delta$,
	\begin{align*}
		\sqrt{n}\|R_n\|_\infty& \lesssim \left(\sqrt{\frac{ \bar{\sigma}^2  \mathfrak{n}_1^2(\lambda)  }{n\lambda^2 } } \vee \frac{ \bar{\sigma}^2 }{n \lambda^2}\right)\left(\kappa^2 \log^2(2/\delta) + \frac{\kappa^4\log^3(2/\delta)}{\sqrt{n} \lambda}\right) \\
		&+ \kappa \| (T + \lambda)^{-2}f_0\|_{\mathcal{H}} \left( \frac{\kappa^4\log^2(2/\delta)}{\sqrt{n}} + \sqrt{n} \lambda^2 \right).
	\end{align*}
\end{proof}

\begin{proof}[\textbf{Proof of Lemma~\ref{lemma:BC-KRR-Consistency-Operator}}]
	Note that $A$ and $(A + \lambda)^{-1}$ commute for any operator $A$. Hence, we compute
	\begin{align*}
		\|\widehat{\Omega}_n - \Omega\|_{op} \leq \left|\widehat{\sigma}_n^2 - \sigma_0^2 \right| \| (T + \lambda)^{-4} T^3  \|_{op} + \sigma_0^2 \| (\widehat{T}_n + \lambda)^{-4}\widehat{T}_n^3 -  (T + \lambda)^{-4} T^3 \|_{op} = \mathbf{I} + \mathbf{II}.
	\end{align*}
	
	Bound on term $\mathbf{I}$. Since $\mathrm{E}[(Y_i - f_0(X_i))^2 - \sigma_0^2] = 0$, $|(Y_i - f_0(X_i))^2 - \sigma_0^2| \leq (B + \kappa\|f_0\|_{\mathcal{H}})^2 + \sigma_0^2$ almost surely, and $\mathrm{Var}(Y_i - f_0(X_i))^2) \leq 2(B + \kappa\|f_0\|_{\mathcal{H}})^4 + 2\sigma_0^4$ Bernstein's inequality for real-valued random variables implies that, with probability at least $1 - \delta$,
	\begin{align*}
		\left|\widehat{\sigma}_n^2 - \sigma_0^2 \right|  \lesssim ( (B + \kappa\|f_0\|_{\mathcal{H}})^2 + \sigma_0^2) \left(\sqrt{\frac{\log(2/\delta)}{n}}  \vee \frac{\log(2/\delta)}{n} \right).
	\end{align*}
	Also,
	\begin{align*}
		\|(T + \lambda)^{-4}T^3\|_{op} \leq \lambda^{-1}\|(T + \lambda)^{-1}T\|_{op}^3.
	\end{align*}
	Hence, with probability at least $1- \delta$,
	\begin{align*}
		\mathbf{I} \lesssim ( (B + \kappa\|f_0\|_{\mathcal{H}})^2 + \sigma_0^2) \|(T + \lambda)^{-1}T\|_{op}^3\left(\sqrt{\frac{\log(2/\delta)}{n\lambda^2}}  \vee \frac{\log(2/\delta)}{n\lambda} \right).
	\end{align*}
	
	Bound on term $\mathbf{II}$. Compute
	\begin{align*}
		&\| (\widehat{T}_n + \lambda)^{-4}\widehat{T}_n^3 -  (T + \lambda)^{-4} T^3 \|_{op}\\
		&\quad \leq \|( \widehat{T}_n^3 - T^3 )(T + \lambda)^{-4} \|_{op} + \|  (\widehat{T}_n + \lambda)^{-4}(T + \lambda)^4 - I \|_{op}\|(T + \lambda)^{-4}T^3\|_{op}.
	\end{align*}
	Since $(a^3 - b^3)c^4 = (a-b)c^2\big( (a-b)^2c^2 + 3(a-b)cbc + 3b^2c^2  \big)$ for $a, b,c \in \mathbb{R}$, it follows that
	\begin{align*}
		&\|( \widehat{T}_n^3 - T^3 )(T + \lambda)^{-4} \|_{op} \\
		&\quad \leq \lambda^{-1} \| (\widehat{T}_n - T)(T + \lambda)^{-1} \|_{op} \left(  \| (\widehat{T}_n - T)(T + \lambda)^{-1} \|_{op}^2 \right.\\
		&\quad\quad \left. + 3\|(\widehat{T}_n - T)(T + \lambda)^{-1} \|_{op}\|T(T + \lambda)^{-1}\|_{op} + 3 \|T(T + \lambda)^{-1}\|_{op}^2\right)\\
		&\quad\leq 3\lambda^{-1} \| (\widehat{T}_n - T)(T + \lambda)^{-1} \|_{op} \left(\| (\widehat{T}_n - T)(T + \lambda)^{-1} \|_{op} +  \|T(T + \lambda)^{-1}\|_{op}  \right)^2.
	\end{align*}
	Hence, by Lemma~\ref{lemma:Consistency-Gram-Matrix}, with probability at least $1- \delta$,
	\begin{align}\label{eq:lemma:BC-KRR-Consistency-Operator-1}
		\begin{split}
			\|( \widehat{T}_n^3 - T^3 )(T + \lambda)^{-4} \|_{op} &\lesssim  \left(\sqrt{ \frac{\kappa^2  \mathfrak{n}_1^2(\lambda) \log(2/\delta)}{n\lambda^2}}  \vee \frac{\kappa^2 \log(2/\delta) }{n\lambda^2} \right) \\
			& \quad\times \left( \sqrt{ \frac{\kappa^2  \mathfrak{n}_1^2(\lambda) \log(2/\delta)}{n\lambda^2}}  \vee \frac{\kappa^2 \log(2/\delta) }{n\lambda^2} \vee \|T(T + \lambda)^{-1}\|_{op}\right)^2.
		\end{split}
	\end{align}
	Next, recall the approach that led to the bound in~\eqref{eq:lemma:BC-KRR-Expansion-4} in the proof of Lemma~\ref{lemma:BC-KRR-Expansion}. We iterate this approach to obtain
	\begin{align*}
		&\|  (\widehat{T}_n + \lambda)^{-4}(T + \lambda)^4 - I \|_{op}\\
		&\quad\leq  3 \| (\widehat{T}_n + \lambda)^{-2} (T+ \lambda)^2 - I\|_{op} \vee 3 \| (\widehat{T}_n + \lambda)^{-2} (T+ \lambda)^2 - I \|_{op}^2\\
		&\quad\leq 9 \| (\widehat{T}_n + \lambda)^{-1} (T+ \lambda) - I\|_{op} \vee 9 \| (\widehat{T}_n + \lambda)^{-1} (T+ \lambda) - I \|_{op}^4\\
		&\quad\leq 9\lambda^{-1} \|T- \widehat{T}_n\|_{op} \vee 9\lambda^{-4} \|T- \widehat{T}_n\|_{op}^4,
	\end{align*}
	Thus, by Lemma~\ref{lemma:Consistency-Gram-Matrix} and since $\mathrm{tr}(T) \leq \kappa^2$, with probability at least $1- \delta$,
	\begin{align}\label{eq:lemma:BC-KRR-Consistency-Operator-2}
		\begin{split}
			& \|  (\widehat{T}_n + \lambda)^{-4}(T + \lambda)^4 - I \|_{op}\|(T + \lambda)^{-1}T\|_{op}^3\\
			&\quad \lesssim  \left( \sqrt{\frac{\kappa^4 \log(2/\delta)}{ n \lambda^2} } \vee \left(\frac{\kappa^2 \log(2/\delta)}{n \lambda^2}\right)^4 \right) \|(T + \lambda)^{-4}T^3\|_{op}.
		\end{split}
	\end{align}
	Combine~\eqref{eq:lemma:BC-KRR-Consistency-Operator-1} and~\eqref{eq:lemma:BC-KRR-Consistency-Operator-2} and conclude that, with probability at least $1- \delta$,
	\begin{align*}
		\mathbf{II} &\lesssim \left(\sqrt{ \frac{\kappa^2  \mathfrak{n}_1^2(\lambda) \log(2/\delta)}{n\lambda^2}}  \vee \frac{\kappa^2 \log(2/\delta) }{n\lambda^2} \right) \|(T + \lambda)^{-1}T\|_{op}^2 \\
		& \quad + \left( \sqrt{ \frac{\kappa^2  \mathfrak{n}_1^2(\lambda) \log(2/\delta)}{n\lambda^2}}  \vee \frac{\kappa^2 \log(2/\delta) }{n\lambda^2}\right)^3\\
		&\quad +  \left( \sqrt{\frac{\kappa^4 \log(2/\delta)}{ n \lambda^2} } \vee \left(\frac{\kappa^2 \log(2/\delta)}{n \lambda^2}\right)^4 \right) \|T^3(T + \lambda)^{-4}\|_{op}.
	\end{align*}
	Since $(T+\lambda)^{-1}T \lesssim I$ it follows that $\|(T + \lambda)^{-1}T\|_{op} \leq 1$. Hence, under the rate conditions the upper bounds on terms $\mathbf{I}$ and $\mathbf{II}$ simplify as stated in the lemma.
\end{proof}

\begin{proof}[\textbf{Proof of Lemma~\ref{lemma:Consistency-Gram-Matrix}}]
	
	Proof of claim (i). Let $1 \leq i \leq n$ and $\alpha \in \mathbb{N}$ be arbitrary. Since $f_0$ minimizes the expected square loss, we have
	\begin{align*}
		\mathrm{E}\left[ (T + \lambda)^{-\alpha} \big(Y_i - f_0(X_i)\big) k_{X_i}\right] = 0.
	\end{align*}
	Therefore, as in the proof of Lemma G.4 in~\cite{singh2023kernel}, we compute,
	\begin{align*}
		\left\|(T + \lambda)^{-1} (T + \lambda)^{-\alpha} \big(Y_i - f_0(X_i)\big) k_{X_i}\right\|_{\mathcal{H}} &\leq  \left| Y_i - f_0(X_i)\right| \left\|(T + \lambda)^{-\alpha}k_{X_i}\right\|_\mathcal{H} \\
		&\leq \lambda^{-\alpha} \kappa (B + \kappa\|f_0\|_{\mathcal{H}}) \quad a.s.
	\end{align*}
	and
	\begin{align*}
		\mathrm{E}\left[\left\|(T + \lambda)^{-\alpha}\big(Y_i - f_0(X_i)\big) k_{X_i}\right\|_{\mathcal{H}}^2 \right] \leq \sigma_0^2 \mathrm{E}\left[\left\|(T + \lambda)^{-\alpha}k_{X_i}\right\|_{\mathcal{H}}^2 \right] \leq \sigma_0^2\mathrm{tr}\left((T+ \lambda)^{-2\alpha}T\right).
	\end{align*}
	Moreover, since $S$ is separable and $k$ is continuous, $\mathcal{H}$ is separable as well. Hence, the conditions of Lemma~\ref{lemma:Bernstein-RE} are satisfied with $\nu = \lambda^{-\alpha}\kappa B + \lambda^{-\alpha}\kappa^2\|f_0\|_{\mathcal{H}}$ and $\sigma =  \sigma_0 \mathfrak{n}_\alpha(\lambda) =\sigma_0 \sqrt{\mathrm{tr}\left((T+ \lambda)^{-2\alpha}T\right)}$. This completes the proof of the first claim.
	
	To proof of claim (ii). Let $1 \leq i \leq n$ be arbitrary. Then $\mathrm{E}[(T + \lambda)^{-\alpha} (( k_{X_i} \otimes k_{X_i}^*) - T)] = 0$. Moreover,
	\begin{align*}
		\left\| (T + \lambda)^{-\alpha} \left(( k_{X_i} \otimes k_{X_i}^*) - T \right)\right\|_{HS} &\leq  \left\|(T + \lambda)^{-\alpha}  k_{X_i} \otimes k_{X_i}^* \right\|_{HS} + \left\| (T + \lambda)^{-\alpha} T \right\|_{HS}\\
		& \leq \lambda^{-\alpha}\kappa^2 + \sqrt{\lambda^{-2\alpha}\mathrm{tr}(T^2)}\\
		& \leq 2 \lambda^{-\alpha}\kappa^2 \quad a.s.
	\end{align*}
	and
	\begin{align*}
		\mathrm{E}\left[\left\| (T + \lambda)^{-\alpha} \left(( k_{X_i} \otimes k_{X_i}^*) - T \right)\right\|_{HS}^2 \right] &= \mathrm{E}\left[ \mathrm{tr}\left(  (T + \lambda)^{-2\alpha}(k_{X_i} \otimes k_{X_i}^* - T)^2 \right)\right]\\
		& = \mathrm{E}\left[ \mathrm{tr}\left( (T + \lambda)^{-2\alpha} (k_{X_i} \otimes k_{X_i}^*)^2 \right)\right] - \mathrm{tr}((T + \lambda)^{-2\alpha}T^2)\\
		&\leq \kappa^2 \mathrm{tr}\left((T+ \lambda)^{-2\alpha}T\right).
	\end{align*}	
	The space of linear operators $A : \mathcal{H} \rightarrow \mathcal{H}$ equipped with the Hilbert-Schmidt norm $\|A\|_{HS} = \sqrt{\mathrm{tr}(A'A)}$ is a Hilbert space. This space can be identified as the tensor product $\mathcal{H} \otimes \mathcal{H}$. Since this tensor product space is separable whenever $\mathcal{H}$ is separable, the conditions of Lemma~\ref{lemma:Bernstein-RE} are satisfied $\nu = 2\lambda^{-\alpha}\kappa^2$ and $\sigma = \kappa \sqrt{\mathrm{tr}((T+ \lambda)^{-2\alpha}T)}$. This completes the proof of the second claim. 
\end{proof}

\end{document}